\documentclass[fleqn,twoside,12pt,a4paper]{article}

\usepackage{ddwmacs,latexsym,amssymb}
\usepackage{epsfig}
\usepackage{cite}

\begin{document}

\title{%
  Explicit Construction of the \\
  Representation of the Braid Generator $\sigma$ \\
  Associated with the One-Parameter Family of \\
  Minimal Typical Highest Weight $\(0,0\,|\,\alpha\)$ \\
  Representations of $U_q[gl(2|1)]$ \\
  and its use in the Evaluation of the \\
  Links--Gould \\
  Two-Variable Laurent Polynomial Invariant \\
  of Oriented $(1,1)$ Tangles
}

\author{%
  David~~De Wit
  \\
  BSc, BScAppGeophHons
  \\
  BSc, PGDipSci, MScSt, BScHons
  ~
  \\
  \vspace{10mm}
  \\
  A thesis submitted in the \\
  Department of Mathematics \\
  at \\
  The University of Queensland \\
  for the requirements of the degree of \\
  Doctor of Philosophy
}

\date{25 November 1998}

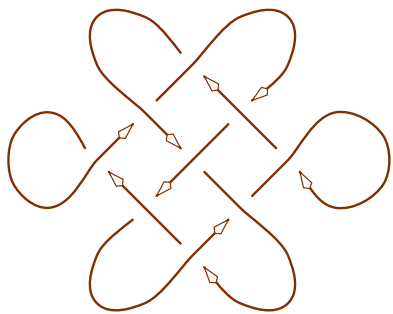
\pagenumbering{roman}


\maketitle

\begin{center}
  \begin{picture}(0,0)%
\epsfig{file=graphics/BuddhistKnot.pstex}%
\end{picture}%
\setlength{\unitlength}{0.00083300in}%
\begingroup\makeatletter\ifx\SetFigFont\undefined
\def\x#1#2#3#4#5#6#7\relax{\def\x{#1#2#3#4#5#6}}%
\expandafter\x\fmtname xxxxxx\relax \def\y{splain}%
\ifx\x\y   
\gdef\SetFigFont#1#2#3{%
  \ifnum #1<17\tiny\else \ifnum #1<20\small\else
  \ifnum #1<24\normalsize\else \ifnum #1<29\large\else
  \ifnum #1<34\Large\else \ifnum #1<41\LARGE\else
     \huge\fi\fi\fi\fi\fi\fi
  \csname #3\endcsname}%
\else
\gdef\SetFigFont#1#2#3{\begingroup
  \count@#1\relax \ifnum 25<\count@\count@25\fi
  \def\x{\endgroup\@setsize\SetFigFont{#2pt}}%
  \expandafter\x
    \csname \romannumeral\the\count@ pt\expandafter\endcsname
    \csname @\romannumeral\the\count@ pt\endcsname
  \csname #3\endcsname}%
\fi
\fi\endgroup
\begin{picture}(2423,1928)(2130,-1625)
\end{picture}

\end{center}

\vfill

\thispagestyle{empty}

\pagebreak


\begin{center}
  {\large \textbf{Statement of Originality}}
\end{center}

The work presented in this thesis is, to the best of my knowledge and
belief, original and my own work, except as acknowledged in the text,
and the material has not been submitted, either in whole or in part,
for a degree at this or any other university.

The material presented herein will be published in:

\begin{description}
\item[\cite{DeWitKauffmanLinks:99}]
  David De Wit, Louis H Kauffman and Jon R Links.
  An Invariant of Knots and Links Brought to You by $U_q[gl(2|1)]$.
  To appear in the Proceedings of the
  XXII International Colloquium on Group Theoretical Methods in
  Physics (Group 22), held at the University of Tasmania, Hobart from
  13--17 July 1998.
  International Press.

\item[\cite{DeWitKauffmanLinks:98}]
  David De Wit, Louis H Kauffman and Jon R Links.
  On the Links--Gould Invariant of Links.
  Accepted for publication by the
  \emph{Journal of Knot Theory and its Ramifications}
  in September 1998.
\end{description}

\vspace{30mm}

\begin{center}
  {\large \textbf{Acknowledgements}}
\end{center}

My PhD candidature at The University of Queensland over the period
1995--1998 was jointly supervised by Mark D Gould and Jon R Links. In
1992 these two published the original description of the Links--Gould
invariant, however its evaluation and hence the experimental deduction
of its properties, remained infeasible.

In 1997, Louis H Kauffman of The University of Illinois at Chicago
visited our group, and introduced us to state models for link
invariants. This theory, developed by him (amongst others) as an
elegant method of evaluation for the Jones polynomial, was also
appropriate for the Links--Gould invariant.

The confluence of these two streams led to the opportunity for me
to explore the properties of the Links--Gould invariant.

\pagebreak


\begin{abstract}
  \noindent
  The ``Links--Gould invariant'' is a two-variable Laurent polynomial
  invariant of oriented $(1,1)$ tangles, which is derived from the
  representation of the braid generator $\sigma$ associated with the
  one-parameter family of representations with highest weights
  $(0,0\,|\,\alpha)$ of the quantum superalgebra $U_q[gl(2|1)]$.  We
  use an abstract tensor state model to evaluate the invariant, as per
  the construction of the bracket polynomial state model used by Louis
  Kauffman to derive the Jones polynomial.  This model facilitates both
  computation and theoretical exploration.

  \vspace{3mm}\noindent
  Evaluation of the invariant for any particular link involves first
  drawing a presentation of the link as a $(1,1)$ tangle (open diagram)
  in a quasi-Morse function form. To this diagram is associated an
  abstract tensor expression that is a contraction over all the free
  indices of a product of rank $2$ and $4$ tensors describing the
  structure of the link.  This expression represents a summation over
  all states in a state space, and the primary (rank $4$) tensor
  involved is a representation of $\sigma$.

  \vspace{3mm}\noindent
  In principle, the quantum R matrix associated with any
  finite-dimensional irreducible representation of any quantum
  (super)algebra will suffice to yield an appropriate representation of
  $\sigma$. After obtaining a particular instance of such a tensor, the
  model also requires a set of four rank $2$ tensors (i.e. genuine
  matrices), and these are chosen using considerations of graphical
  consistency and the properties of quantum (super)algebras.

  \vspace{3mm}\noindent
  The family of four dimensional highest weight $U_q[gl(2|1)]$
  representations labelled $\(0,0\,|\,\alpha\)$ has a two-variable
  quantum R matrix (unique up to orthogonal transformations).  Choosing
  this R matrix to yield the representation of $\sigma$ ensures that
  our polynomial invariant will also have two variables, and hence
  might be expected to have a greater discriminatory power than an
  invariant with only one.  We construct this R matrix from first
  principles, invoking \textsc{Mathematica} to assist with the
  algebra.

  \vspace{3mm}\noindent
  We have evaluated the invariant for several critical link examples
  and numerous other links of special forms.  Performing the
  contraction over the abstract tensor involves quite lengthy symbolic
  computation, and is infeasible for a human to do.  Again, the
  assistance of \textsc{Mathematica} has been invoked.

  \vspace{3mm}\noindent
  From the experimental data, we observe that the Links--Gould
  invariant is distinct from the two-variable HOMFLY polynomial in that
  it detects the chirality of some links where the HOMFLY fails.  As we
  have no method of automatically computing the invariant for large
  numbers of links, we cannot make empirical observations about how
  often our invariant will be able to distinguish chirality.  Notably,
  it does \emph{not} distinguish inverses, which is not surprising as
  we are able to demonstrate that \emph{no} invariant of this type
  should be able to distinguish between inverses. It also does not
  distinguish between mutants.
\end{abstract}

\cleardoublepage

\tableofcontents

\pagebreak

\listoffigures

\listoftables

\cleardoublepage

\pagenumbering{arabic}


\section{Overview}
\addtocontents{toc}{\protect\vspace{-2.5ex}}

The announcement in 1985 of the discovery of the Jones polynomial
\cite{Jones:85} prompted vigorous research which has since yielded
many other link invariants. In particular, the quantum algebras as
defined by Drinfel'd \cite{Drinfeld:87} and Jimbo \cite{Jimbo:85},
being examples of quasi-triangular Hopf algebras, provide a systematic
means of solving the Yang--Baxter equation, and in turn may be employed
to construct representations of the braid group. From each of these
representations a prescription exists to compute invariants of oriented
knots and links \cite{Reshetikhin:87,Turaev:88,ZhangGouldBracken:91b}.
The simplest example of this process is that the Jones polynomial is
recoverable using the simplest quantum algebra $U_q[sl(2)]$ in its
minimal ($2$-dimensional) representation.

From such a large class of available invariants it is natural to ask
if generalisations exist, with the view to gaining a classification.
One possibility is to look to multiparametric extensions in order
to see which invariants occur as special cases. A notable example is
the HOMFLY%
\footnote{
  The HOMFLY polynomial is named by the conjunction of the initials of
  six of its discoverers
  \cite{FreydYetterHosteLickorishMilletOcneanu:85}, omitting those
  (``P'' and ``T'') of two independent discoverers
  \cite{PrzytyckiTraczyk:87}.  Przytycki, the omitted ``P'', has
  furthered the entymological spirit with the suggestion ``FLYPMOTH''
  \cite[p256]{Przytycki:89}, which includes all the discoverers and
  has a muted reference to the ``flyping'' operation of the Tait,
  Kirkwood and Little -- the original compilers of knot tables.  A
  similar possibility is ``HOMFLYPT''.  Bar-Natan (Prasolov and
  Sossinsky \cite[p36]{PrasolovSossinsky:96} cite Bar-Natan
  \cite{BarNatan:95}, who cites ``L Rudulph'') goes further, adding a
  ``U'' for good measure, to account for any unknown discoverers,
  yielding the unpalatable ``LYMPHTOFU''!
}
invariant \cite{FreydYetterHosteLickorishMilletOcneanu:85}
which includes both the Jones and Alexander--Conway invariants
\cite{Alexander:23,Conway:70} as particular cases as well as the
invariants arising from minimal representations of $U_q[sl(n)]$
\cite{Turaev:88}. Another is the Kauffman polynomial which includes the
Jones invariant as well as those obtained from the quantum algebras
$U_q[o(n)]$ and $U_q[sp(2n)]$ in the $q$-deformations of the
defining representations \cite{Turaev:88}.

The work of Turaev and Reshetikhin \cite{ReshetikhinTuraev:90} shows
that the algebraic properties of quantum algebras are such that an
extension of this method to produce invariants of oriented tangles is
permissible. (A tangle diagram is a link diagram with free ends.)  An
associated invariant takes the form of a tensor operator acting on a
product of vector spaces. Zhang \cite{Zhang:95} has extended this
formalism to the case of quantum superalgebras, which are
$\mathbb{Z}_2$ graded generalisations of quantum algebras.

Since quantum superalgebras give rise to nontrivial one-parameter
families of irreducible representations, it is possible to utilise them
for the construction of two-variable invariants. This was first shown
by Links and Gould \cite{LinksGould:92b} for the simplest case using the
family of four dimensional representations of $U_q[gl(2|1)]$. It was
also made known that a one-variable reduction of this invariant
coincides with a one-variable reduction of the Kauffman polynomial by
the use of the Birman--Wenzl--Murakami algebra. Extensions to more
general representations of quantum superalgebras are discussed in
\cite{GouldLinksZhang:96b}.

\clearpage

The original paper describing the Links--Gould invariant
\cite{LinksGould:92b} does not contain evaluations of it, for want of
an efficient method of computation. As such, none of its properties
were known. Here, we describe a `state space' method of evaluation and
investigations into the properties of the invariant.
We provide many examples and a complete description of the
state model for the invariant in abstract tensor form. This description
of the invariant directly facilitates the construction of a computer
program in \textsc{Mathematica} for evaluation of the invariant.  We
emphasise that evaluation of the invariant is computationally
expensive, and as such it would be infeasible to do manually for any
but the simplest knots.  The code is fairly straightforward, and to
some extent can be understood directly, without prior experience with
\textsc{Mathematica}. The readability is somewhat reduced in places
where computational efficiency has been improved by more sophisticated
constructs.

The method of evaluating the invariant involves a prior construction of
the quantum R matrix associated with a family of four dimensional
representations, and this is done in explicit detail in
\secref{Uqgl21}. We have used \textsc{Mathematica} to check the details.
(As an aside, we further illustrate the power of
\textsc{Mathematica} by also explicitly constructing the associated
trigonometric R matrix.) Having obtained this R matrix, the
construction of the invariant follows from properties of ribbon Hopf
(super)algebras and their representations.  Here we consider the
invariants of $(1,1)$ tangles for the following reason: for invariants
derived from representations of quantum superalgebras with zero
$q$-superdimension, the corresponding invariant is also zero.  If the
representation is irreducible, the quantum superalgebra symmetry of the
procedure ensures that the invariant of $(1,1)$ tangles takes the form
of some scalar multiple of the identity matrix.  (See
\cite{ReshetikhinTuraev:90} for a discussion of this symmetry.) We take
this scalar to be the invariant.

We prove that the Links--Gould invariant is not able to distinguish a
knot from its inverse (\propref{CantDetectInvertibility}), nor can it
distinguish between mutant links (\propref{CantDistinguishMutants}).
However it does distinguish many links from their
reflections (see Propositions \ref{prp:AmphichiralMeansPalindromic} and
\ref{prp:NotPalindromicImpliesChiral}), and it is distinguished from
the Kauffman and HOMFLY (and hence the Jones) polynomials by this
behaviour (see \secref{NineFortyTwoandTenFortyEight} and
\secref{NoninvertiblePretzelsarenotDistinguished} for specific
examples). Whether it is a complete invariant for chirality is an open
question.

\pagebreak


\section{Link Invariants from Combinatorial Topology}
\addtocontents{toc}{\protect\vspace{-2.5ex}}

\subsection{Basic Terminology of Knot Theory}
\addtocontents{toc}{\protect\vspace{-2.5ex}}

A good introduction to the basic concepts of knot theory is found in a
book by Adams \cite{Adams:94}, and aspects of the theory that are
particularly pertinent to this work may be found in books by Kauffman
\cite{Kauffman:87a,Kauffman:93}.  Between them, these sources describe
the relationship between physical and mathematical knots, and the
concepts of regular and ambient isotopy.  Here, we assume knowledge of
these basics.

To be certain, multicomponent knots are called \emph{links}, and almost
all of what we shall have to say under the heading of `knot theory'
actually applies to links.  We are in general interested in the
classification and properties of \emph{oriented} links, except that we
only investigate the invertibility properties of single-component links
(i.e. true knots).  It goes without apology that we are only
considering \emph{tame} (i.e.  non-pathological) links.

We shall write ``='' to denote the \emph{ambient isotopy} of link
diagrams, meaning that they are equivalent under \emph{all} of the
Reidemeister moves \cite{Kauffman:93,Kauffman:96}; and we shall use the
Alexander--Briggs notation and ordering for \emph{prime knots} (i.e.
knots not able to be decomposed into a `product' of simpler knots,
which are called \emph{composite}).

Substantial lists of data (i.e. diagrams, values of invariants) for all
prime knots of up to $10$ crossings are found in a book by Kawauchi
\cite[\S F]{Kawauchi:96}; this information is complemented by data for
braid presentations in a paper by Jones \cite[pp381-388]{Jones:87}.  A
smaller amount of similar data for some (multicomponent) `prime links'
is contained in a book by Adams \cite{Adams:94} and a paper by Doll and
Hoste \cite{DollHoste:91}.


\subsection{The Link Classification Problem}
\addtocontents{toc}{\protect\vspace{-2.5ex}}

Consider an arbitrary presentation (i.e. a link diagram) of an
arbitrary link. Under ambient isotopy, the diagram will in general
correspond to a composition of standard link presentations, from which
we can name the link.  Even in the simplest case, where the link is
prime, how do we discover \emph{which} standard presentation?  This is
a difficult problem when tackled with only the machinery of the
Reidemeister moves, and is only made worse when non-prime links are
permitted.  Thus, whilst we \emph{do} have a classification system for
prime links (and a way of decomposing composite into prime links), we
do not have an \emph{efficient method} for determining, from a given
link diagram, to which of the standard presentations our link
corresponds.

Link (polynomial) invariants attempt to address this deficiency.  To
each link is assigned an algebraic expression (typically a Laurent
polynomial), and the method of evaluation of this expression ensures
that it is invariant under the Reidemeister moves, i.e. that it is
independent of the presentation.

However, invariants may in general assign the \emph{same} polynomial to
\emph{distinct} links, which reduces their utility. A major goal of
knot theory is then to develop a \emph{complete} invariant that assigns
distinct polynomials to distinct links.  Currently, this goal is
unsatisfied.  The best-known and perhaps the most powerful polynomial
invariants currently known are the two-variable HOMFLY and Kauffman
polynomials, but they are far from complete.

For a particular link invariant to be complete, it is minimally
necessary that it must distinguish links under the transformations of
reflection, inversion (for oriented links) and mutation.


\subsection{Reflection, Inversion and Mutation of Links}
\addtocontents{toc}{\protect\vspace{-2.5ex}}

Denote by $K^*$ the \emph{reflection} (or mirror image) of a link $K$.
A link is \emph{chiral} if it is distinct from its reflection:
$K^*\neq K$; i.e.  there are actually two distinct links with the same
name.  Alternatively, a link is \emph{amphichiral} if it is ambient
isotopic to its reflection, i.e. $K^*=K$.  Note that this definition
doesn't require an orientation.  To illustrate, $3_1$ (the Trefoil
Knot, see \figref{HopfLinkTrefoil}) is chiral:  ${(3_1)}^* \neq 3_1$.
Perusal of the table of all prime knots of up to ten crossings shows
that most of these knots are chiral, and it is probably safe to claim
that most links are chiral.

The HOMFLY and Kauffman polynomials can distinguish
many (but not all) chiral knots from their reflections.  The first
chiral knot that neither the HOMFLY nor the Kauffman polynomial can
distinguish is $9_{42}$, i.e.  ${(9_{42})}^* \neq 9_{42}$, but the
polynomials are equal.  Similarly, $10_{48}$ is chiral, but the HOMFLY
polynomial fails to detect this, although the Kauffman does detect it
\cite[p218]{Kauffman:93}. (Note that the diagram for $10_{48}$ is
wrongly labeled $10_{79}$ in that volume. $10_{79}$ is in fact not
chiral.)

Denote by $K^{-1}$ the \emph{inverse} of an oriented (true) knot $K$,
obtained by reversing its orientation.  Whilst this is a simple concept
for a true knot, there are of course many possibilities for the
reversal of only some components of oriented, multi-component links; we
shall not go into these here.  Commonly, $K = K^{-1}$, and we say that
$K$ is \emph{invertible}.  For example, the Trefoil Knot is invertible
${(3_1)}^{-1}=3_1$.  Less commonly, $K \neq K^{-1}$, and we say that
$K$ is \emph{noninvertible}.  The first example of a noninvertible
prime knot is $8_{17}$.  To date, no known invariant of \emph{Vassiliev
type} can detect invertibility.  Our invariant is of this type.

Both the reflection and the inverse are automorphisms of order two,
i.e.  ${(K^*)}^*=K$ and ${(K^{-1})}^{-1}=K$.  The notions
combined yield:  ${(K^*)}^{-1}={(K^{-1})}^*$.
For example, for $3_1$, we have two equivalence classes:
$3_1={(3_1)}^{-1}$ and
${(3_1)}^*={({(3_1)}^{-1})}^*={({(3_1)}^*)}^{-1}$.

A third transformation of a link (diagram) is by \emph{mutation}, in
which a component is isolated and replaced with a $180^\circ$ rotation
of itself. In general, this process yields a new link.  The HOMFLY and
Kauffman polynomials \emph{always} fail to distinguish between mutant
links.

\pagebreak


\subsection{The Links--Gould Invariant}
\addtocontents{toc}{\protect\vspace{-2.5ex}}

The paper by Links and Gould \cite{LinksGould:92b} describes the
construction of link invariants from a representation of the braid
generator $\sigma$, corresponding to any particular finite-dimensional
irreducible representation of any particular quantum (super)algebra.
That paper demonstrates the existence of the Links--Gould invariant,
but does not supply an efficient method for its evaluation. Here, we
set up a state space model for this purpose. This model uses an
explicit representation of $\sigma$, and also requires `cap' and `cup'
matrices $\Omega^{\pm}$ and $\mho^{\pm}$. Below, we develop the state
space model for the general case of oriented links, before returning to
the specific choices for the Links--Gould invariant. It is emphasised
that we do not reiterate the material of \cite{LinksGould:92b} here, as
it is not germane to the state space model.


\subsection{State Space Models}
\addtocontents{toc}{\protect\vspace{-2.5ex}}

In order to construct a state space model for a given link, we must
first construct an \emph{abstract tensor} expression corresponding to
that link, which consists of a contraction (i.e. a sum) over the
indices of a product of tensors that reflects the link's structure.
This may be regarded as a summation of \emph{states}.

In order that the contraction yields a topological invariant of ambient
isotopy, the tensors included within the abstract tensor must satisfy a
set of minimal graphical consistency considerations; these constraints
are well-documented in a paper by Hennings \cite[pp59-60]{Hennings:91}.
That paper describes the situation for unoriented links, but the
material carries over naturally to the oriented case (see a book by
Kauffman \cite[pp~235-237]{Kauffman:93}).  It is well-known that the
use of a representation of the braid generator will suffice for the
most important of these tensors.

Further introductions to state space models include the expositions by
Kauffman, e.g.  \cite[pp395-399]{Kauffman:87b} and
\cite[pp204-209]{Kauffman:88}, which evaluate invariants by graphical
considerations, and \cite[\S\S4.1-4.2,~pp163-172]{Kauffman:97a}, which
evaluates the Jones polynomial using the minimal
$U_q[sl(2)]$ representation, using only abstract tensors.


\subsection{Morse and Quasi-Morse Functions}
\addtocontents{toc}{\protect\vspace{-2.5ex}}

For a link diagram, a \emph{Morse function} applied to that diagram is
a drawing of it in terms of vertically-aligned crossings (i.e. with all
arrows pointing upwards), vertical lines and horizontal arcs
\cite{Hennings:91}. This assignment is not unique. We shall here define
a `\emph{quasi-Morse' function} as one that does essentially the same
thing, but allows non-vertically-aligned crossings and non-straight
vertical lines. Our link diagrams (see Figures
\ref{fig:HopfLinkTrefoil} to \ref{fig:TenFortyEight}) are drawn
according to such functions; state space invariants are quickly seen to
be independent of the distinction between Morse and quasi-Morse
functions.


\subsection{Construction of Abstract Tensors}
\addtocontents{toc}{\protect\vspace{-2.5ex}}

To define the abstract tensor corresponding to a particular link, we
begin with a presentation of the link's diagram according to a
(quasi-)Morse function, and assign to the arms of each crossing in the
diagram a set of labels $a,b,c,d$, etc.

By a `\emph{positive oriented}' or `right-handed' crossing within the
link diagram, we intend a crossing such that if the right hand thumb
points in the direction of one of the arrows, the (curled!)
fingers of the right hand will point in the direction of the other
arrow. The opposite situation is naturally called a `\emph{negative
oriented}' or `left-handed' crossing.  If the arms of a positive
oriented crossing with upward-pointing arrows have labels $a,b,c,d$, we
shall associate a rank $4$ tensor $\sigma$ with that crossing, where
the positioning of the indices in the tensor corresponds to those on
the diagram. We shall do the same using $\sigma^{-1}$ for negative
oriented crossings.  Similarly, to each horizontally-oriented arc, we
shall assign a rank $2$ tensor (i.e. a true matrix) called a `cap'
$\Omega^\pm$ or `cup' $\mho^\pm$ as appropriate. This information is
summarised in \figref{CrossingsCapsCups}.

\begin{figure}[ht]
  \begin{center}
    \input{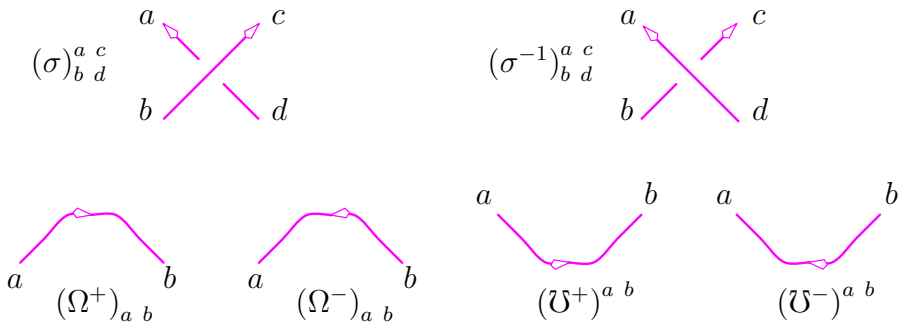}
    \caption[%
      The tensors $\sigma$ and $\sigma^{-1}$ and
      cap and cup matrices $\Omega^\pm$ and $\mho^\pm$.
    ]{%
      The tensors $\sigma$ and $\sigma^{-1}$ and
      cap and cup matrices $\Omega^\pm$ and $\mho^\pm$, and
      associated oriented crossings and horizontal arcs.
    }
    \addtocontents{lof}{\protect\vspace{-2.5ex}}
    \figlabel{CrossingsCapsCups}
  \end{center}
\end{figure}

The abstract tensor for our link is then formed as a contraction over
the indices of the product of the tensors representing the components
of the diagram.  This amounts to a summation over all possible
\emph{states}, where each state amounts to a particular assignment of
values to the dummy indices in the terms of the abstract tensor
expression prior to contraction.  To be a little more precise, the
construction of the bracket polynomial (and hence the Jones polynomial)
involves a contraction over \emph{all} the indices of the tensor,
whilst our invariant will require a contraction over all but two of the
indices. This corresponds to drawing the link as a $(1,1)$ tangle, with
free ends at the top (index $y$) and bottom (index $x$), to which we
will assign equal, constant values (see
\secref{FormationoftheInvariant}).

Each crossing in the link diagram will then have a copy of a rank $4$
tensor $\sigma$ or $\sigma^{-1}$ associated to it, in a manner
describing the entanglement of the link.  The choice of a suitable
$\sigma$ is not unique, but a representation of the braid generator
associated with a quantum R matrix (i.e. a solution of the quantum
Yang--Baxter equation) will suffice.  Evaluation of the abstract tensor
will require that $\sigma$ explicitly, and we will also require the set
of matrices $\Omega^\pm$ and $\mho^\pm$. The latter will be determined
from other properties of the representation used to represent
$\sigma$.


\subsection{Forming the Invariant from the Abstract Tensor}
\addtocontents{toc}{\protect\vspace{-2.5ex}}
\seclabel{FormationoftheInvariant}

Having formed the abstract tensor ${(T_K)}^y_x$ for any particular link
(where the indices $x$ and $y$ correspond to the lower and upper loose
ends of the tangle), we form an invariant by setting $x$ and $y$ to be
the \emph{same}:
\bse
  {LG}_K
  \defeq
  {(T_K)}^i_i
  \qquad
  \textrm{(no sum on $i$)},
  \eqlabel{Nerida}
\ese
for \emph{any} allowable index $i$; we shall typically choose $i=1$.
This situation is in contrast to the situation for the bracket
polynomial, where the invariant is formed as a trace (i.e. by summation
on $i$ in \eqref{Nerida}).  In our case, that trace would have to be
replaced with a \emph{super}trace as we are dealing with a quantum
\emph{super}algebra, which would lead to an uninteresting invariant
that would always be zero.  \figref{ClosingTangle} demonstrates how the
closure of a $(1,1)$ tangle (open diagram) to form a proper link
relates to the abstract tensors.

\begin{figure}[htbp]
  \begin{center}
    \input{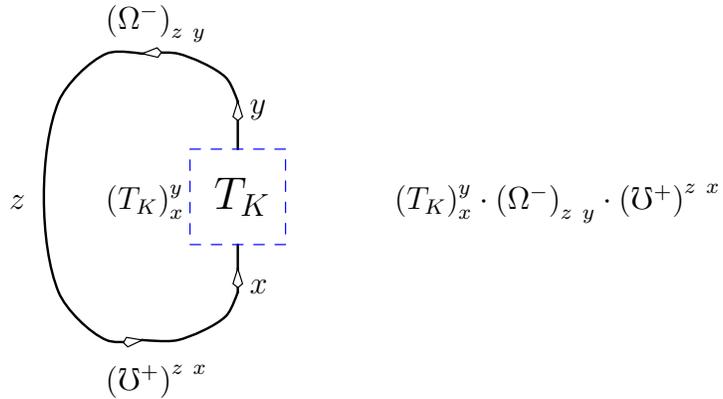}
    \caption{Closure of a $(1,1)$ tangle to form a closed link diagram.}
    \addtocontents{lof}{\protect\vspace{-2.5ex}}
    \figlabel{ClosingTangle}
  \end{center}
\end{figure}

In \secref{RepresentationofSigma}, we define the Links--Gould invariant
by \eqref{Nerida}, where we have made specific choices for the tensors
composing $T_K$. Such an invariant is typically a Laurent polynomial in
one or two variables; indeed our invariant is, although we supply no
proof.

Recall the choice of a representation of the braid generator for
$\sigma$. Deduction of appropriate cap and cup matrices follows from
graphical consistency considerations which ensure that our invariant is
automatically an invariant of ambient isotopy, i.e. it is
\emph{writhe-normalised} cf. the case of the bracket polynomial
\cite[p102]{Kawauchi:96}.  \figref{Loop} depicts removal of a single
twist from a diagram, showing that this requires that we must
choose $\sigma$, $\Omega^\pm$ and $\mho^\pm$ such that:
\bse
  {(\sigma)}^{y~a}_{x~b}
  \cdot
  {(\Omega^+)}_{a~c}
  \cdot
  {(\mho^-)}^{b~c}
  =
  \delta^y_x.
  \eqlabel{AmbientIsotopy}
\ese

\begin{figure}[ht]
  \begin{center}
    \input{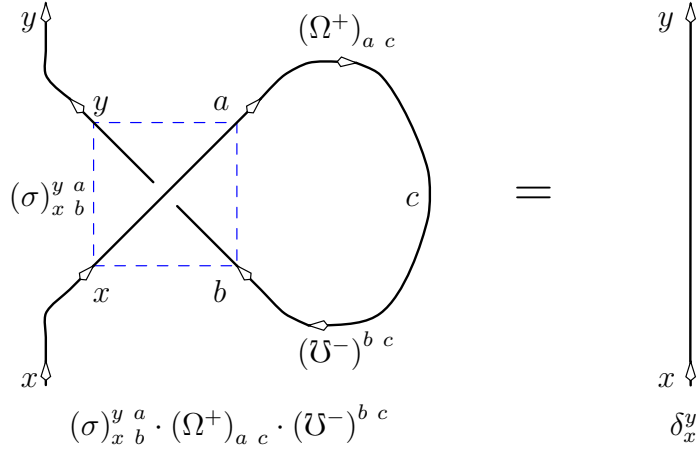}
    \caption{Equivalence of a single loop with an upward line.}
    \addtocontents{lof}{\protect\vspace{-2.5ex}}
    \figlabel{Loop}
  \end{center}
\end{figure}
We will return to this result in \secref{CapsandCups}.


\subsection{Auxiliary Tensors}
\addtocontents{toc}{\protect\vspace{-2.5ex}}

For any particular link, we may write down its corresponding abstract
tensor expression using only the tensors $\sigma$ and $\sigma^{-1}$
(as well as $\Omega^\pm$ and $\mho^\pm$).  However, subsequent
evaluation of the abstract tensor typically involves the repetition of
computations reflecting repeated patterns of crossings in the link. To
improve computational efficiency, we define auxiliary tensors.  In the
following, we use the notation $X$ to represent either $\sigma$ or
$\sigma^{-1}$, corresponding to a diagram component of a crossing with
upward pointing arrows.  The Einstein summation convention is
used throughout.

\begin{itemize}
\item
  The first auxiliary tensors represent crossings that have been
  `twisted' relative to $X$. The left, right, and upside-down-twisted
  versions of $X$ are called $X_l$, $X_r$ and $X_d$ respectively (see
  Figures \ref{fig:XLandXR} and \ref{fig:XD}).
  \bne
    {(X_l)}^{a~c}_{b~d}
    & \defeq &
    {(X)}^{e~a}_{d~h}
    \cdot
    {(\Omega^-)}_{b~e}
    \cdot
    {(\mho^-)}^{h~c}
    \nonumber
    \\
    {(X_r)}^{a~c}_{b~d}
    & \defeq &
    {(X)}^{c~g}_{f~b}
    \cdot
    {(\mho^+)}^{a~f}
    \cdot
    {(\Omega^+)}_{g~d}
    \eqlabel{XRdef}
    \\
    \nonumber
    {(X_d)}^{a~c}_{b~d}
    & \defeq &
    {(X)}^{e~g}_{f~h}
    \cdot
    {(\mho^+)}^{a~h}
    \cdot
    {(\Omega^+)}_{g~b}
    \cdot
    {(\mho^+)}^{c~f}
    \cdot
    {(\Omega^+)}_{e~d}.
    \nonumber
  \ene
  \begin{figure}[htbp]
    \begin{center}
      \input{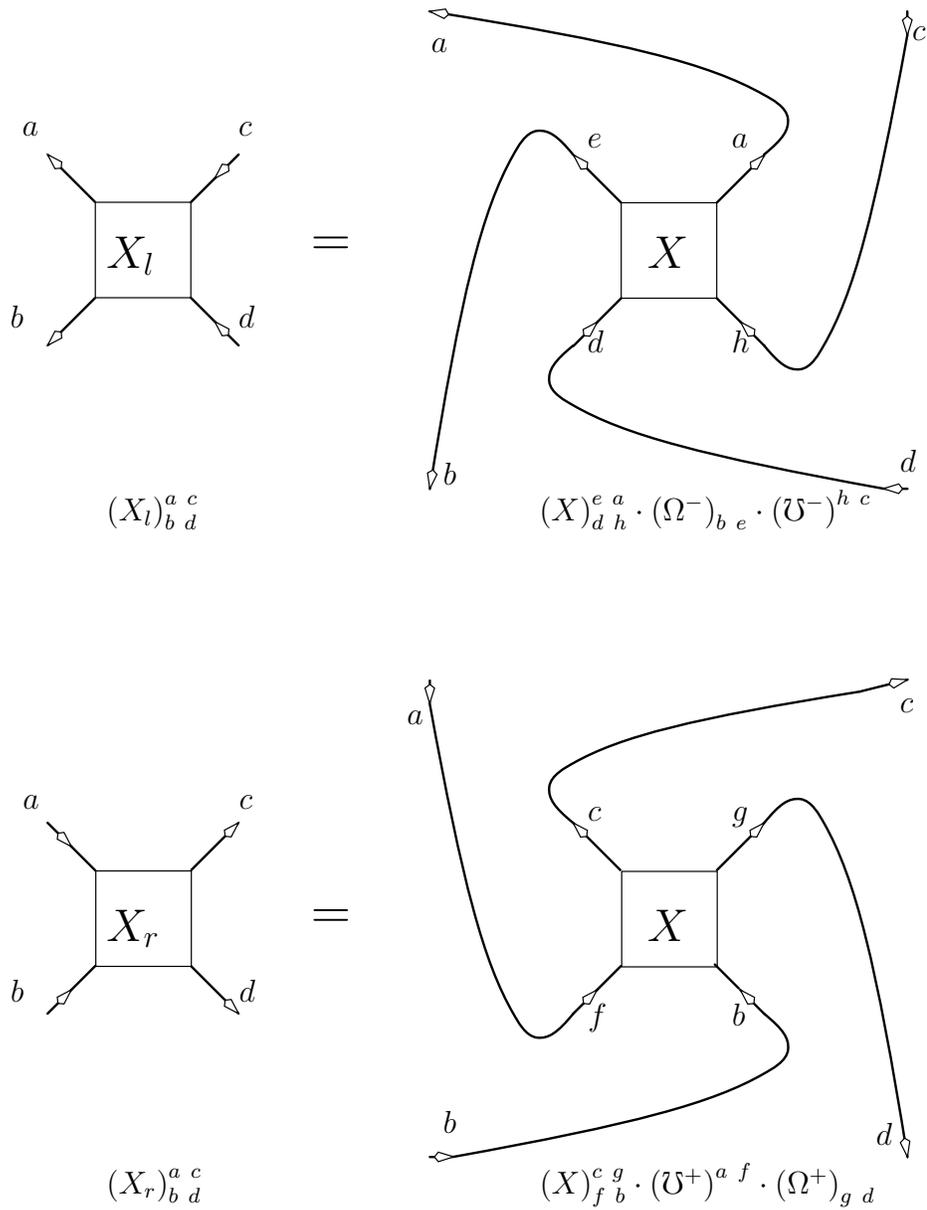}
      \caption[%
        The auxiliary tensors $X_l$ and $X_r$.
      ]{%
        The auxiliary tensors
        $X_l$ and $X_r$;
        $X$ is either $\sigma$ or $\sigma^{-1}$.
      }
      \addtocontents{lof}{\protect\vspace{-2.5ex}}
      \figlabel{XLandXR}
    \end{center}
  \end{figure}
  \begin{figure}[htbp]
    \begin{center}
      \input{graphics/XD.pstex_t}
      \caption[%
        The auxiliary tensors $X_d$.
      ]{%
	The auxiliary tensors $X_d$;
        $X$ is either $\sigma$ or $\sigma^{-1}$.
      }
      \addtocontents{lof}{\protect\vspace{-2.5ex}}
      \figlabel{XD}
    \end{center}
  \end{figure}
\item
  The second set of auxiliary tensors represent $N$ copies of the
  \emph{same} crossing $X$ atop one another (see \figref{XN}).
  We use the shorthand $\sigma^{-N}\equiv {(\sigma^{-1})}^N$.
  \be
    {(X^{N})}^{a~c}_{b~d}
    \defeq
    {(X)}^{a~c}_{e~f}
    \cdot
    {(X^{N-1})}^{e~f}_{b~d},
    \qquad
    N \geqslant 2.
  \ee
  \begin{figure}[htbp]
    \begin{center}
      \input{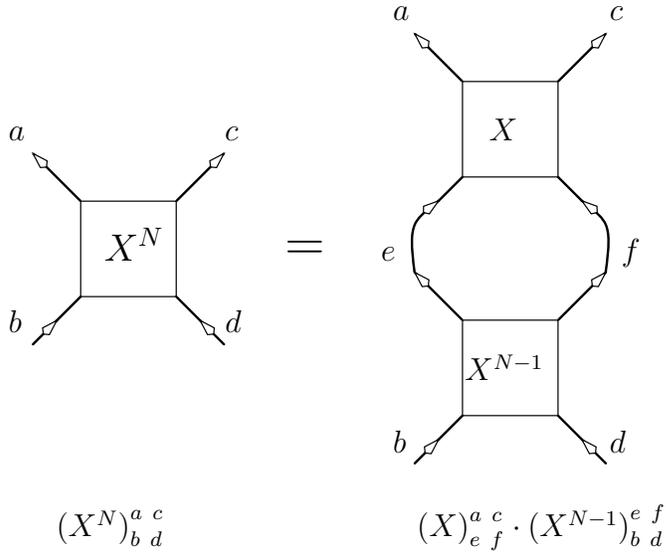}
      \caption[%
        The auxiliary tensors $X^{N}$.
      ]{%
        The auxiliary tensors
        $X^{N}$, expressed in terms of $X$ and $X^{N-1}$;
	$X$ is either $\sigma$ or $\sigma^{-1}$ and $N\geqslant2$.
	If all arrows are reversed, the definition also holds for $X$
	being $\sigma_d$ or $\sigma^{-1}_d$.
      }
      \addtocontents{lof}{\protect\vspace{-2.5ex}}
      \figlabel{XN}
    \end{center}
  \end{figure}

\pagebreak

\item
  The third set of auxiliary tensors correspond to the situation where
  a crossing $X$ (upright, so we write $X_u\equiv X$ for consistency)
  sits to the left or right of its own `upside-downness' $X_d$ (see
  \figref{XUXDandXDXU}).
  \be
    {(X_u X_d)}^{a~c}_{b~d}
    & \defeq &
    {(X)}^{a~e}_{b~f}
    \cdot
    {(X_d)}^{g~c}_{h~d}
    \cdot
    {(\Omega^+)}_{e~g}
    \cdot
    {(\mho^-)}^{f~h}
    \\
    {(X_d X_u)}^{a~c}_{b~d}
    & \defeq &
    {(X_d)}^{a~e}_{b~f}
    \cdot
    {(X)}^{g~c}_{h~d}
    \cdot
    {(\Omega^-)}_{e~g}
    \cdot
    {(\mho^+)}^{f~h}.
  \ee
  \begin{figure}[htbp]
    \begin{center}
      \input{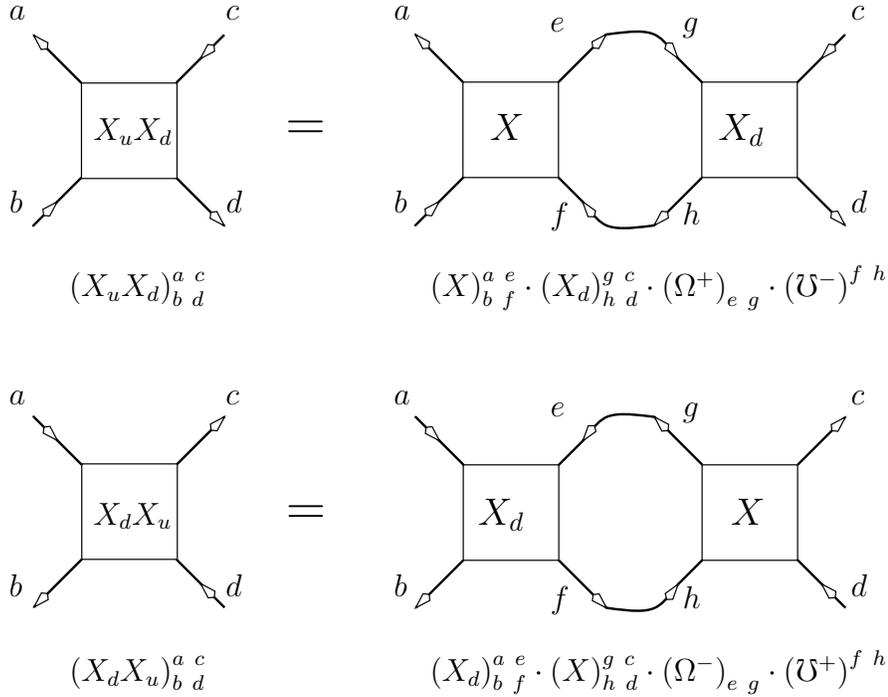}
      \caption[%
        The auxiliary tensors
        $X_d X_u$ and $X_u X_d$.
      ]{%
        The auxiliary tensors
        $X_d X_u$ and $X_u X_d$;
        $X$ is either $\sigma$ or $\sigma^{-1}$.
      }
      \addtocontents{lof}{\protect\vspace{-2.5ex}}
      \figlabel{XUXDandXDXU}
    \end{center}
  \end{figure}
\item
  The fourth set of auxiliary tensors describe the situation
  where a crossing $X_l$ is placed atop or below a crossing $X_r$
  (see \figref{XLXRandXRXL}).
  \be
    {(X_l X_r)}^{a~c}_{b~d}
    & \defeq &
    {(X_l)}^{a~c}_{e~f}
    \cdot
    {(X_r)}^{e~f}_{b~d}
    \\
    {(X_r X_l)}^{a~c}_{b~d}
    & \defeq &
    {(X_r)}^{a~c}_{e~f}
    \cdot
    {(X_l)}^{e~f}_{b~d}.
  \ee
  The diagram for $X_l X_r$ is of course a right rotation of the
  diagram for $X_d X_u$.
  \begin{figure}[htbp]
    \begin{center}
      \input{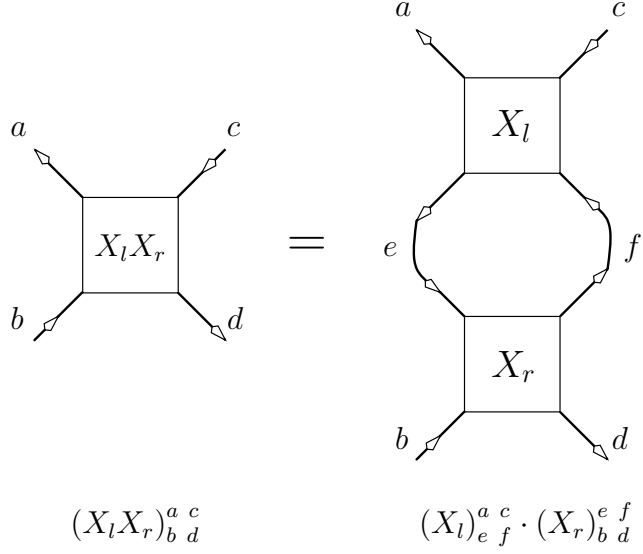}
      \caption[%
        The auxiliary tensors $X_l X_r$ and $X_r X_l$.
      ]{%
	The auxiliary tensors $X_l X_r$ and $X_r X_l$; $X$ is either
	$\sigma$ or $\sigma^{-1}$. The diagram for $X_r X_l$ is
	obtained by interchanging $r$ and $l$.
      }
      \addtocontents{lof}{\protect\vspace{-2.5ex}}
      \figlabel{XLXRandXRXL}
    \end{center}
  \end{figure}
\item
  The fifth set of auxiliary tensors describe the situation where a
  chain of $N$ crossings is formed by the recursive placing of $X$
  alongside $X_d$, with $X$ as the leftmost and rightmost crossings;
  where $X$ is either $\sigma$ or $\sigma^{-1}$. The recursive
  definition of such chains $X_{udu}^N$ is provided in \figref{XUDU}.

  \begin{figure}[htbp]
    \begin{center}
      \input{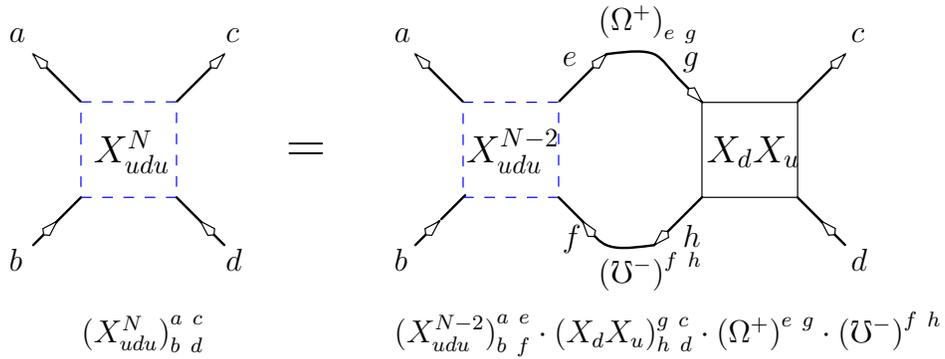}
      \caption[%
        Recursive definition of the chains $X_{udu}^N$.
      ]{%
        Recursive definition of the chains $X_{udu}^N$ used in the
        presentation of the
        pretzel knots; $X$ is either $\sigma$ or $\sigma^{-1}$.
        The minimum is the case $N=1$, which
        corresponds to $X$, i.e.  $X_{udu}^1 \defeq X$.
      }
      \addtocontents{lof}{\protect\vspace{-2.5ex}}
      \figlabel{XUDU}
    \end{center}
  \end{figure}
\end{itemize}

\cleardoublepage


\section{A Collection of Especially Interesting Links}
\addtocontents{toc}{\protect\vspace{-2.5ex}}

\subsection{Several Interesting Ordinary Links}
\addtocontents{toc}{\protect\vspace{-2.5ex}}

We have selected several links as test examples for evaluation of our
invariant. The first few examples illustrate the principles involved in
the construction of the abstract tensors, and the need for the use of
the auxiliary tensors. The latter three are important cases that we
must investigate to determine how our invariant compares with other
invariants.  In \tabref{linknamesandwrithes}, we list the links to be
studied, in order of their Alexander--Briggs notation
\cite{AlexanderBriggs:26}.  The data for proper (i.e. single component)
knots have been taken from \cite[\S F.2]{Kawauchi:96} and for
two-component links from \cite{Adams:94} (itself citing
\cite{Rolfsen:76} and \cite{DollHoste:91}).

\begin{table}[ht]
  \centering
  \begin{tabular}{||l||c|c||}
    \hline\hline
    & & \\[-3mm]
      \multicolumn{1}{||c||}{$K$} & Amphichiral? & Invertible? \\
    & & \\[-3mm]
    \hline\hline
    & & \\[-3mm]
    $ 0_1 $ (Unknot) &
      Yes & Yes \\[1mm]
    \hline
    & & \\[-3mm]
    $ 2^2_1 $ (Hopf Link) &
      No & Yes \\[1mm]
    \hline
    & & \\[-3mm]
    $ 3_1 $ (Trefoil Knot) &
      No & Yes
      \\[1mm]
    \hline
    & & \\[-3mm]
      $ 4_1 $ (Figure Eight Knot) & Yes & Yes \\[1mm]
    \hline
    & & \\[-3mm]
      $ 5_1 $ (Cinquefoil Knot) & No & Yes \\[1mm]
    \hline
    & & \\[-3mm]
      $ 5_2 $ & No & Yes \\[1mm]
    \hline
    & & \\[-3mm]
    $ 5^2_1 $ (Whitehead Link) & No \cite[pp49-50]{Kauffman:87a} & Yes
      \\[1mm]
    \hline
    & & \\[-3mm]
      $ 6_1 $ & No & Yes \\[1mm]
    \hline
    & & \\[-3mm]
      $ 6_2 $ & No & Yes \\[1mm]
    \hline
    & & \\[-3mm]
      $ 6_3 $ & Yes & Yes \\[1mm]
    \hline
    & & \\[-3mm]
      $ 7_1 $ (Septfoil Knot) & No & Yes \\[1mm]
    \hline
    & & \\[-3mm]
      $ 7_2 $ & No & Yes \\[1mm]
    \hline\hline
    & & \\[-3mm]
    $ 8_{17} $ & Yes & No \\[1mm]
    \hline\hline
    & & \\[-3mm]
    $ 9_{42} $ & No & Yes \\[1mm]
    \hline
    & & \\[-3mm]
    $ 10_{48} $ & No & Yes \\[1mm]
    \hline\hline
  \end{tabular}
  \caption[%
    Data for the links to be investigated.
  ]{%
    Data for the links to be investigated, see Figures
    \ref{fig:HopfLinkTrefoil} to \ref{fig:TenFortyEight}.
  }
  \addtocontents{lot}{\protect\vspace{-2.5ex}}
  \tablabel{linknamesandwrithes}
\end{table}

\pagebreak

Below, we describe each of our example links in terms of
$(1,1)$-tangle (open diagram) forms, where the
indices $x$ and $y$ correspond to the lower and upper loose ends of the
tangle. From these diagrams, their abstract tensors ${(T_K)}^y_x $ are
immediately obvious. Some of the diagrams have been constructed using a
braid presentation, as taken from \cite[pp109-110]{Jones:85} and
\cite[pp381-386]{Jones:87}.  For the Unknot ($0_1$), a braid
presentation is the trivial $e\in B_1$ as the tangle is a simple
vertical line, and we have ${(T_{0_1})}^y_x \defeq \delta^y_x$.

\begin{description}
\item[$\mathbf{2^2_1}$ (Hopf Link) and $\mathbf{3_1}$ (Trefoil Knot):]
  Diagrams, constructed from the braid presentations
  ${\sigma_1}^2$ and ${\sigma_1}^3 \in B_2$, are found in
  \figref{HopfLinkTrefoil}.  The abstract tensors follow immediately:
  \be
    {( T_{2^2_1} )}^y_x
    & \defeq &
    {\( \sigma^2 \)}^{y~a}_{x~b}
    \cdot
    {\( \Omega^+ \)}_{a~c}
    \cdot
    {\( \mho^- \)}^{b~c}
    \\
    {\( T_{3_1} \)}^y_x
    & \defeq &
    {\( \sigma^3 \)}^{y~a}_{x~b}
    \cdot
    {\( \Omega^+ \)}_{a~c}
    \cdot
    {\( \mho^- \)}^{b~c}.
  \ee

  \begin{figure}[htbp]
    \begin{center}
      \input{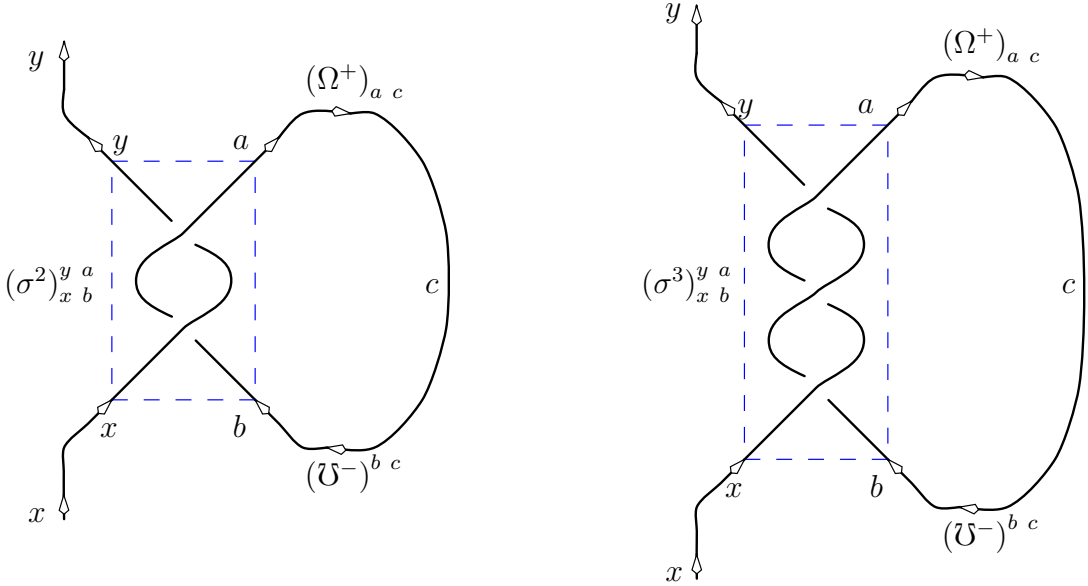}
      \caption{Tangle forms of $2^2_1$ (the Hopf Link) and
          $3_1$ (the Trefoil Knot).
      }
      \addtocontents{lof}{\protect\vspace{-2.5ex}}
      \figlabel{HopfLinkTrefoil}
    \end{center}
  \end{figure}

\item[$\mathbf{4_1}$ (Figure Eight Knot):]
  From the diagram found in \figref{FigureEight}, we have:
  \be
    {\( T_{4_1} \)}^y_x
    \defeq
    {\( \sigma_l \sigma_r \)}^{y~b}_{a~c}
    \cdot
    {\( \sigma^{-1}_u \sigma^{-1}_d \)}^{a~c}_{x~d}
    \cdot
    {\( \Omega^- \)}_{b~e}
    \cdot
    {\( \mho^+ \)}^{d~e}.
  \ee
  \begin{figure}[ht]
    \begin{center}
      \input{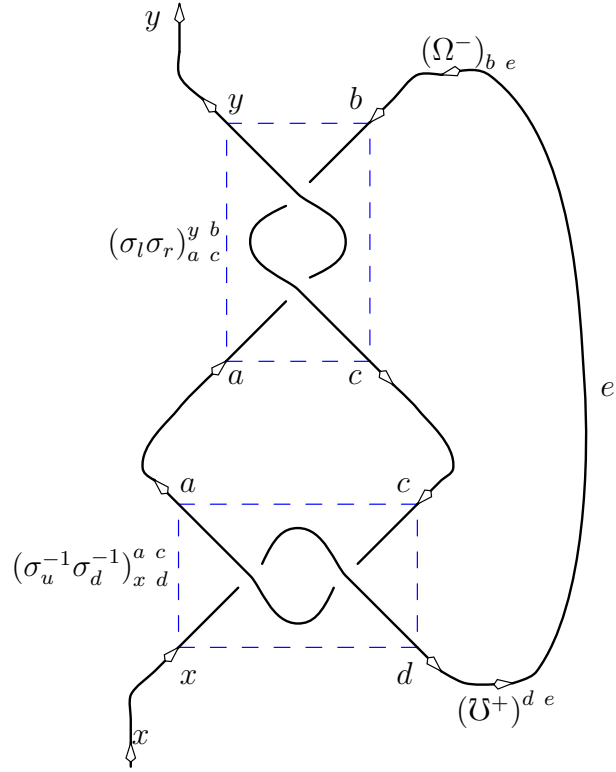}
      \caption{Tangle form of $4_1$ (the Figure Eight Knot).}
      \addtocontents{lof}{\protect\vspace{-2.5ex}}
      \figlabel{FigureEight}
    \end{center}
  \end{figure}

\item[$\mathbf{5_1}$ and $\mathbf{7_1}$
     (Cinquefoil Knot and Septfoil Knot):]
  Diagrams, constructed from the braid presentations
  $\sigma_1^5 \in B_2$ and $\sigma_1^7 \in B_2$,
  are found in \figref{CinqueSeptfoil}.  We have:
  \be
    {\( T_{5_1} \)}^y_x
    & \defeq &
    {\( \sigma^5 \)}^{y~a}_{x~b}
    \cdot
    {\( \Omega^+ \)}_{a~c}
    \cdot
    {\( \mho^- \)}^{b~c}
    \\
    {\( T_{7_1} \)}^y_x
    & \defeq &
    {\( \sigma^7 \)}^{y~a}_{x~b}
    \cdot
    {\( \Omega^+ \)}_{a~c}
    \cdot
    {\( \mho^- \)}^{b~c}.
  \ee

  \begin{figure}[htbp]
    \begin{center}
      \input{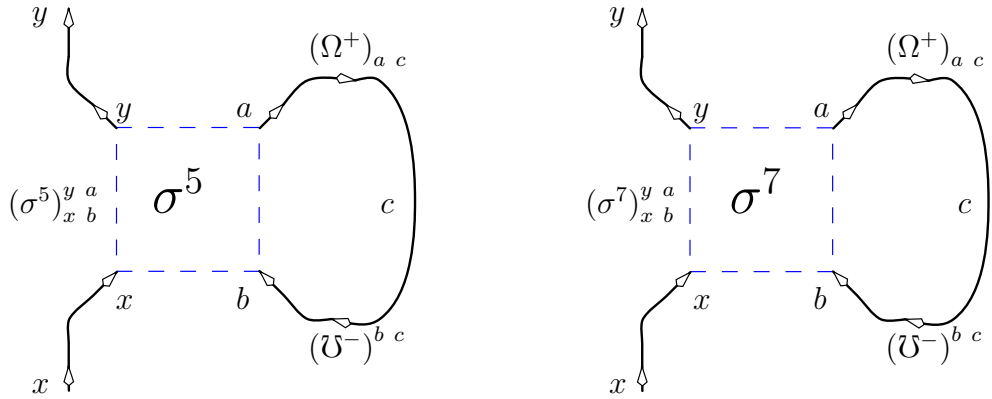}
      \caption{Tangle forms of $5_1$ and $7_1$
               (the Cinquefoil Knot and the Septfoil Knot).}
      \addtocontents{lof}{\protect\vspace{-2.5ex}}
      \figlabel{CinqueSeptfoil}
    \end{center}
  \end{figure}

\item[$\mathbf{5_2}$:]
  A diagram is found in \figref{FiveTwo}.  We have:
  \be
    {( T_{5_2} )}^y_x
    \defeq
    {\( \sigma_{udu}^3 \)}^{b~d}_{c~x}
    \cdot
    {\( \sigma^2 \)}^{a~y}_{b~d}
    \cdot
    {\( \Omega^- \)}_{e~a}
    \cdot
    {\( \mho^+ \)}^{e~c}.
  \ee

  \begin{figure}[htbp]
    \begin{center}
      \input{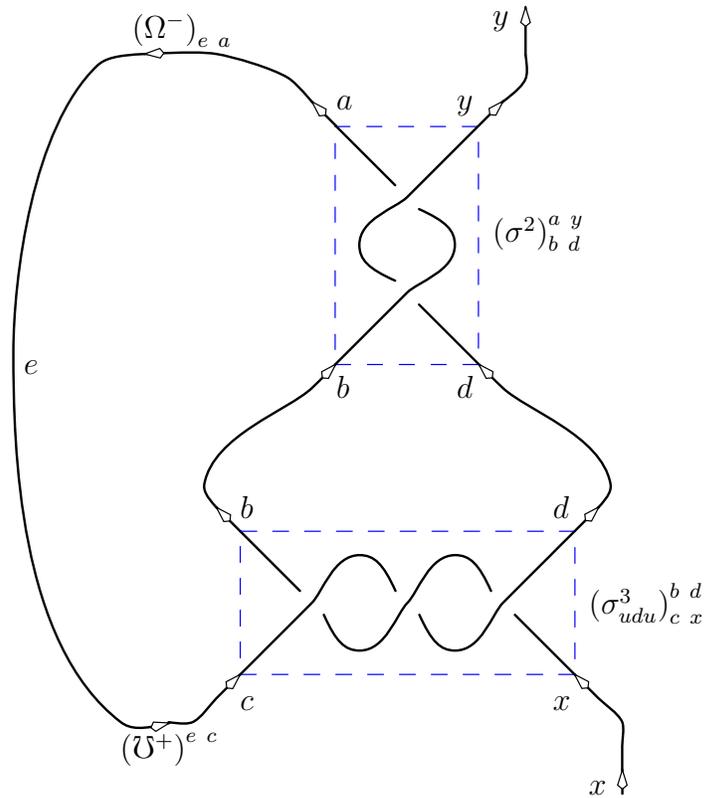}
      \caption{Tangle form of $5_2$.}
      \addtocontents{lof}{\protect\vspace{-2.5ex}}
      \figlabel{FiveTwo}
    \end{center}
  \end{figure}

\item[$\mathbf{5^2_1}$ (Whitehead Link):]
  This link is named after the topologist J H C Whitehead
  \cite[p200]{Kauffman:88}.  A diagram is found in
  \figref{WhiteheadLink}. To deduce the abstract tensor, we initially
  define a temporary tensor to reduce computation:
  \be
    {\( W \)}^{c~i}_{x~d}
    \defeq
    {\( \sigma^{-2} \)}^{c~e}_{x~f}
    \cdot
    {( \sigma_d^2 )}^{g~i}_{h~d}
    \cdot
    {\( \Omega^+ \)}_{e~g}
    \cdot
    {\( \mho^- \)}^{f~h}.
  \ee
  With this, we have:
  \be
    {( T_{5^2_1} )}^y_x
    \defeq
    {\( W \)}^{c~i}_{x~d}
    \cdot
    {\( \sigma_r \sigma_l \)}^{a~y}_{i~b}
    \cdot
    {\( \Omega^+ \)}_{c~a}
    \cdot
    {\( \mho^+ \)}^{d~b}.
  \ee

  \begin{figure}[htbp]
    \begin{center}
      \input{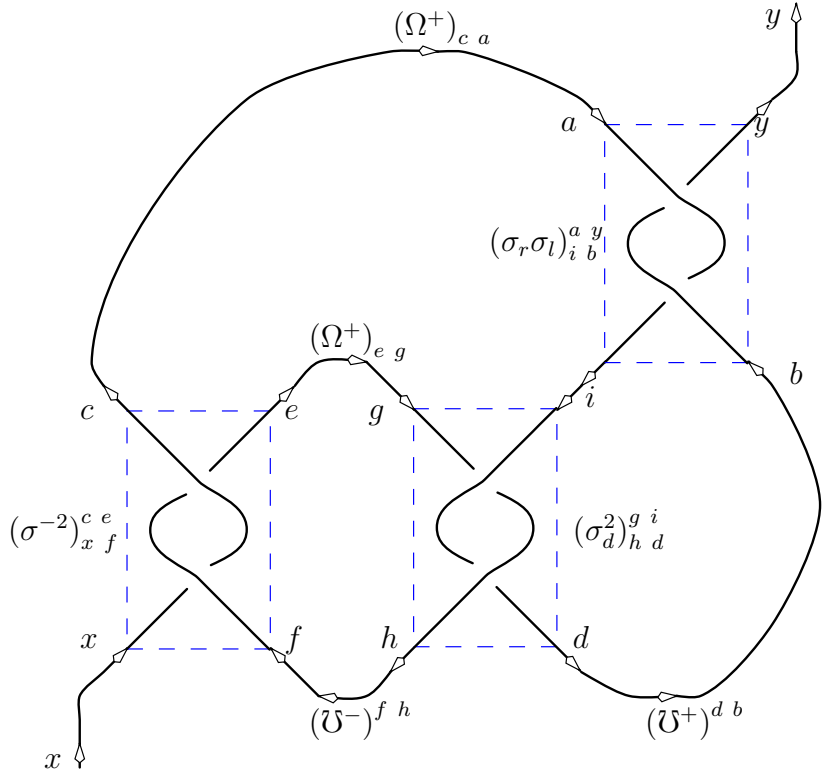}
      \caption{Tangle form of $5^2_1$ (the Whitehead Link).}
      \addtocontents{lof}{\protect\vspace{-2.5ex}}
      \figlabel{WhiteheadLink}
    \end{center}
  \end{figure}

\item[$\mathbf{6_1}$:]
  A diagram is found in \figref{SixOneTwo}.  We define an auxiliary
  tensor:
  \be
    {\( SOA \)}^{b~d}_{c~x}
    \defeq
    {\( \sigma^{-1}_d \)}^{b~f}_{c~h}
    \cdot
    {\( \sigma^{-3}_{udu} \)}^{g~d}_{i~x}
    \cdot
    {\( \Omega^- \)}_{f~g}
    \cdot
    {\( \mho^+ \)}^{h~i}.
  \ee
  With this, we have:
  \be
    {( T_{6_1} )}^y_x
    \defeq
    {\( SOA \)}^{b~d}_{c~x}
    \cdot
    {\( \sigma_r \sigma_l \)}^{a~y}_{b~d}
    \cdot
    {\( \Omega^+ \)}_{e~a}
    \cdot
    {\( \mho^- \)}^{e~c}.
  \ee

  \begin{figure}[htbp]
    \begin{center}
      \input{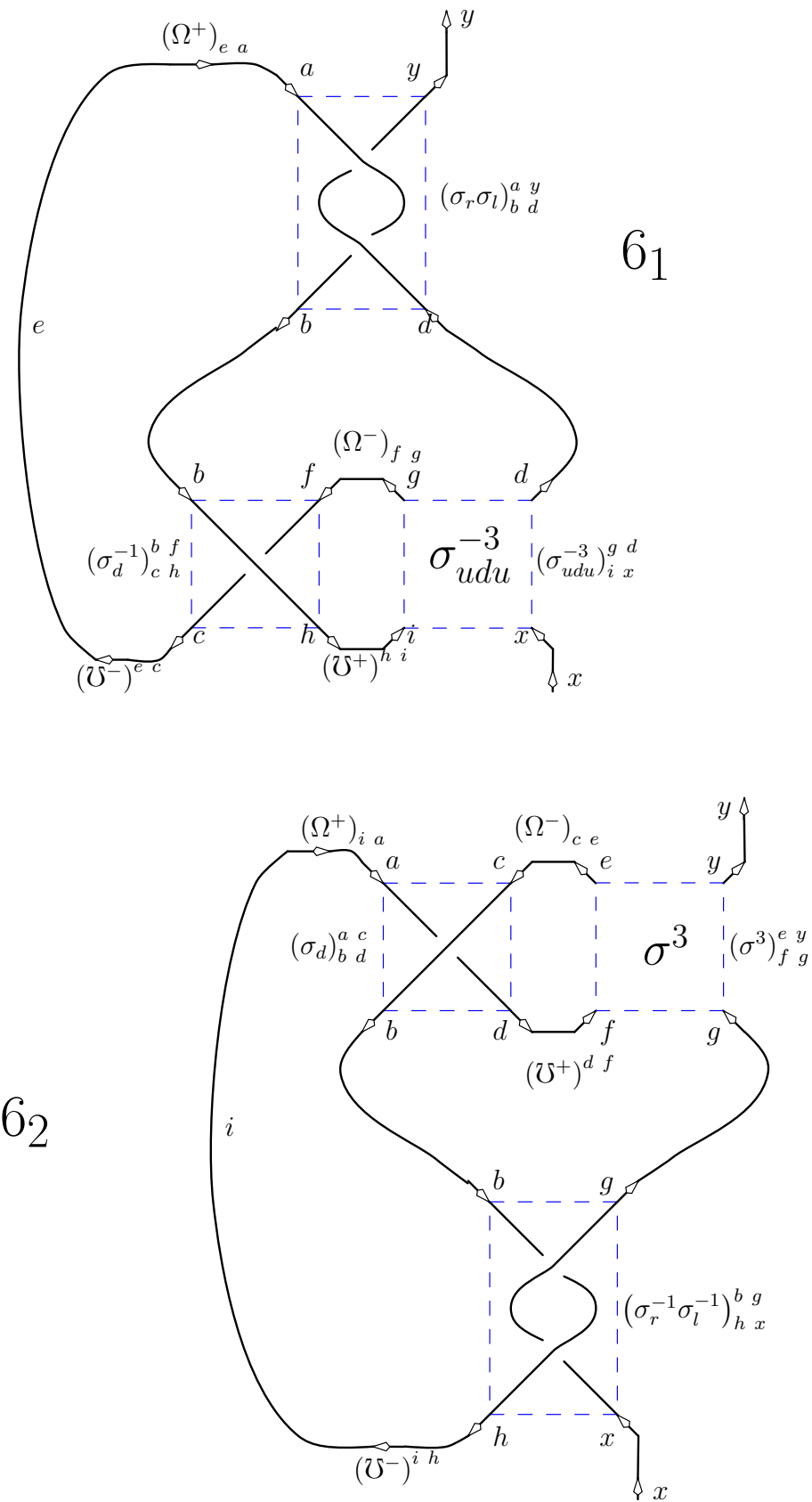}
      \caption{Tangle forms of $6_1$ and $6_2$.}
      \addtocontents{lof}{\protect\vspace{-2.5ex}}
      \figlabel{SixOneTwo}
    \end{center}
  \end{figure}

\item[$\mathbf{6_2}$:]
  A diagram is found in \figref{SixOneTwo}.  We define an auxiliary
  tensor:
  \be
    {\( STA \)}^{a~y}_{b~g}
    \defeq
    {\( \sigma_d \)}^{a~c}_{b~d}
    \cdot
    {\( \sigma^3 \)}^{e~y}_{f~g}
    \cdot
    {\( \Omega^- \)}_{c~e}
    \cdot
    {\( \mho^+ \)}^{d~f}.
  \ee
  With this, we have:
  \be
    {( T_{6_2} )}^y_x
    \defeq
    {\( STA \)}^{a~y}_{b~g}
    \cdot
    {\( \sigma^{-1}_r \sigma^{-1}_l \)}^{b~g}_{h~x}
    \cdot
    {\( \Omega^+ \)}_{i~a}
    \cdot
    {\( \mho^- \)}^{i~h}.
  \ee

\item[$\mathbf{6_3}$:]
  A diagram, drawn from the braid presentation
  $\sigma_1^{-1} \sigma_2^2 \sigma_1^{-2} \sigma_2 \in B_3$,
  is found in \figref{SixThree}.  We define some
  auxiliary tensors:
  \be
    {\( STA \)}^{a~d~y}_{b~f~i}
    & \defeq &
    {\( \sigma^{-2} \)}^{a~e}_{b~f}
    \cdot
    {\( \sigma \)}^{d~y}_{e~i}
    \\
    {\( STB \)}^{b~f~i}_{c~h~x}
    & \defeq &
    {{\( \sigma^{-1} \)}}^{b~g}_{c~h}
    \cdot
    {\( \sigma^2 \)}^{f~i}_{g~x}
    \\
    {\( ST \)}^{a~d~y}_{c~h~x}
    & \defeq &
    {\( STA \)}^{a~d~y}_{b~f~i}
    \cdot
    {\( STB \)}^{b~f~i}_{c~h~x}.
  \ee
  With these, we have:
  \be
    {( T_{6_3} )}^y_x
    \defeq
    {\( ST \)}^{a~d~y}_{c~h~x}
    \cdot
    {\( \Omega^- \)}_{j~d}
    \cdot
    {\( \mho^+ \)}^{j~h}
    \cdot
    {\( \Omega^- \)}_{k~a}
    \cdot
    {\( \mho^+ \)}^{k~c}.
  \ee

  \begin{figure}[htbp]
    \begin{center}
      \input{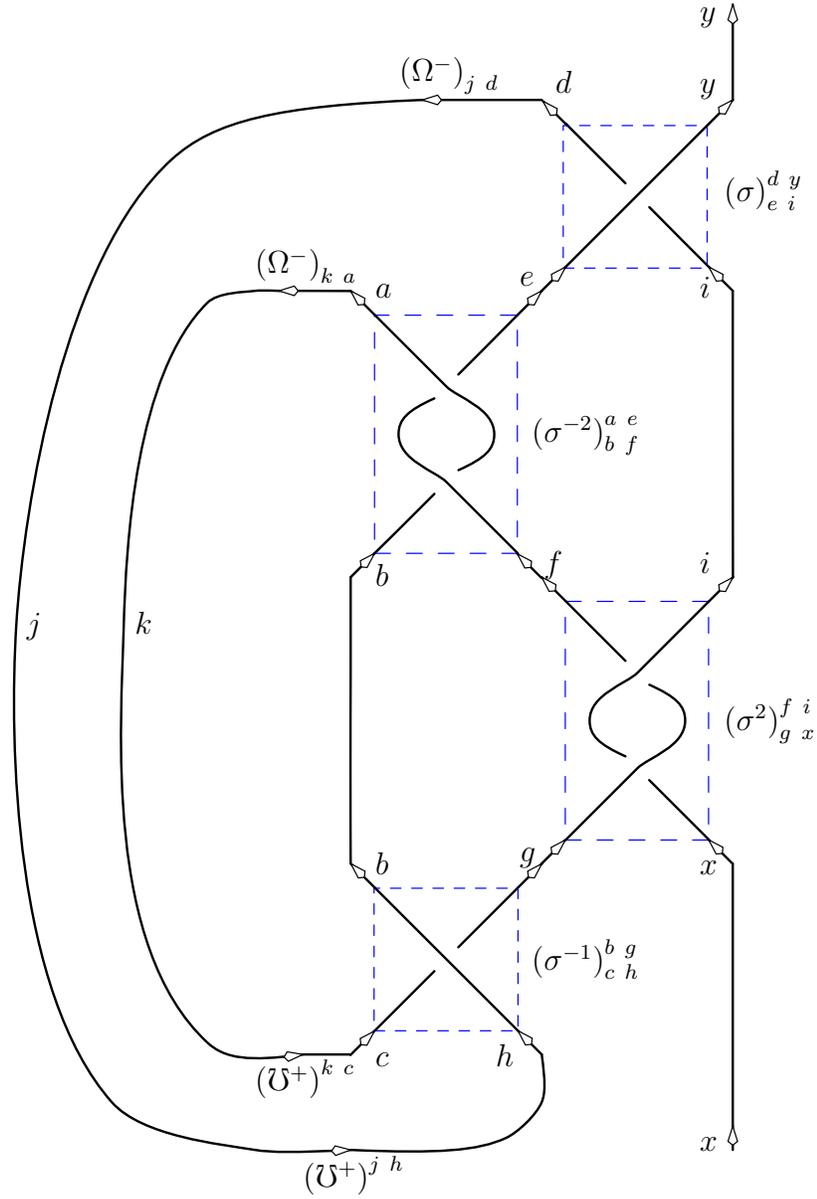}
      \caption{Tangle form of $6_3$.}
      \addtocontents{lof}{\protect\vspace{-2.5ex}}
      \figlabel{SixThree}
    \end{center}
  \end{figure}

\item[$\mathbf{7_2}$:]
  A diagram is found in \figref{SevenTwo}.  We have:
  \be
    {( T_{7_2} )}^y_x
    \defeq
    {\( \sigma_{udu}^5 \)}^{b~d}_{c~x}
    \cdot
    {\( \sigma^2 \)}^{a~y}_{b~d}
    \cdot
    {\( \Omega^- \)}_{e~a}
    \cdot
    {\( \mho^+ \)}^{e~c}.
  \ee

  \begin{figure}[htbp]
    \begin{center}
      \input{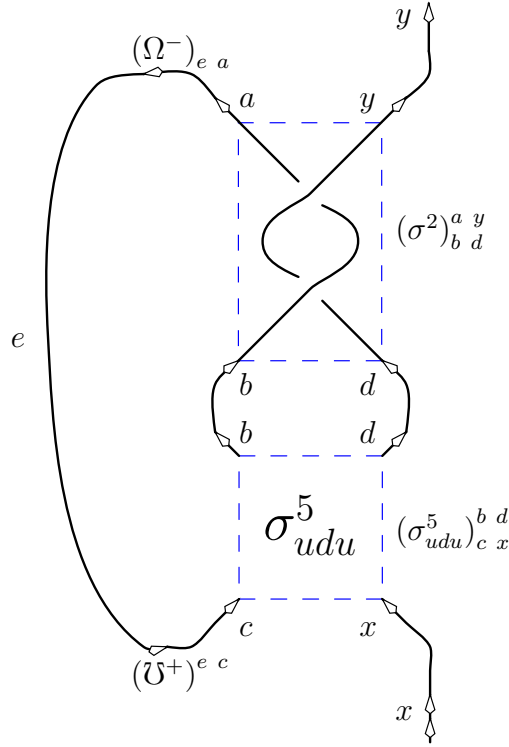}
      \caption{Tangle form of $7_2$.}
      \addtocontents{lof}{\protect\vspace{-2.5ex}}
      \figlabel{SevenTwo}
    \end{center}
  \end{figure}

\item[$\mathbf{8_{17}}$:]
  A braid presentation is
  $
    {\( \sigma_1^{-1} \sigma_2 \)}^2 \sigma_2^2 \sigma_1^{-2} \sigma_2
    \in B_3
  $,
  from which we draw the diagram found in \figref{EightSeventeen}.
  Again, we define some temporary (rank $6$) tensors to reduce
  computation:
  \be
    {\( EA \)}^{y~c~e}_{b~d~f}
    & \defeq &
    {\( \sigma^{-2} \)}^{c~e}_{g~f}
    \cdot
    {\( \sigma^2 \)}^{y~g}_{b~d}
    \\
    {\( EB \)}^{b~d~f}_{x~i~j}
    & \defeq &
    {{\( \sigma^{-1} \)}}^{d~f}_{k~l}
    \cdot
    {\( \sigma \)}^{b~k}_{m~n}
    \cdot
    {{\( \sigma^{-1} \)}}^{n~l}_{o~j}
    \cdot
    {\( \sigma \)}^{m~o}_{x~i}.
  \ee
  With these, we have:
  \be
    \hspace{-9pt}
    {\( T_{8_{17}} \)}^y_x
    \defeq
    {\( EA \)}^{y~c~e}_{b~d~f}
    \cdot
    {\( EB \)}^{b~d~f}_{x~i~j}
    \cdot
    {\( \Omega^+ \)}_{c~r}
    \cdot
    {\( \mho^- \)}^{i~r}
    \cdot
    {\( \Omega^+ \)}_{e~q}
    \cdot
    {\( \mho^- \)}^{j~q}.
  \ee
  To reduce computation, we may define even more auxiliary tensors:
  \be
    {\( EB \)}^{b~d~f}_{x~i~j}
    =
    {\( EC \)}^{b~d~f}_{m~n~l}
    \cdot
    {\( ED \)}^{m~n~l}_{x~i~j},
  \ee
  where:
  \be
    {\( EC \)}^{b~d~f}_{m~n~l}
    & \defeq &
    {{\( \sigma^{-1} \)}}^{d~f}_{k~l}
    \cdot
    {\( \sigma \)}^{b~k}_{m~n}
    \\
    {\( ED \)}^{m~n~l}_{x~i~j}
    & \defeq &
    {{\( \sigma^{-1} \)}}^{n~l}_{o~j}
    \cdot
    {\( \sigma \)}^{m~o}_{x~i}.
  \ee

  \begin{figure}[htbp]
    \begin{center}
      \input{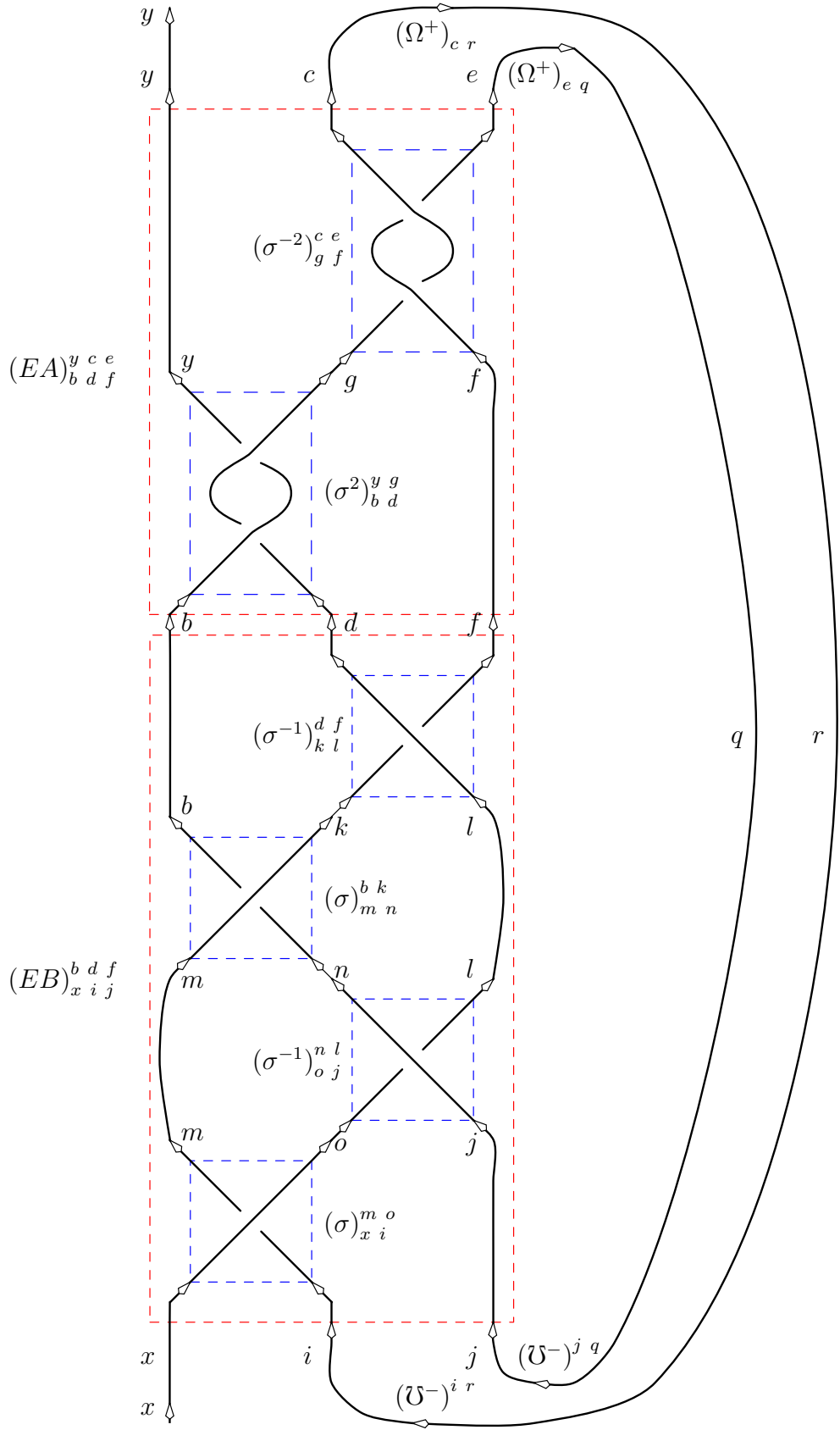}
      \caption{Tangle form of $ 8_{17} $.}
      \addtocontents{lof}{\protect\vspace{-2.5ex}}
      \figlabel{EightSeventeen}
    \end{center}
  \end{figure}

\item[$\mathbf{9_{42}}$:]
  A diagram is found in \figref{NineFortyTwo}.  Again, we define a
  temporary tensor to reduce computation:
  \be
    {\( N \)}^{a~y}_{b~h}
    \defeq
    {( \sigma_d^2 )}^{a~c}_{b~d}
    \cdot
    {\( \sigma^{-3} \)}^{e~y}_{f~h}
    \cdot
    {\( \Omega^- \)}_{c~e}
    \cdot
    {\( \mho^+ \)}^{d~f}.
  \ee
  With this, we have:
  \be
    {\( T_{9_{42}} \)}^y_x
    \defeq
    {\( N \)}^{a~y}_{b~h}
    \cdot
    {\( \sigma^{-1}_d \sigma^{-1}_u \)}^{b~h}_{i~j}
    \cdot
    {\( \sigma_u \sigma_d \)}^{k~i}_{x~m}
    \cdot
    {\( \mho^+ \)}^{m~j}
    \cdot
    {\( \Omega^+ \)}_{k~a}.
  \ee
  \begin{figure}[htbp]
    \begin{center}
      \input{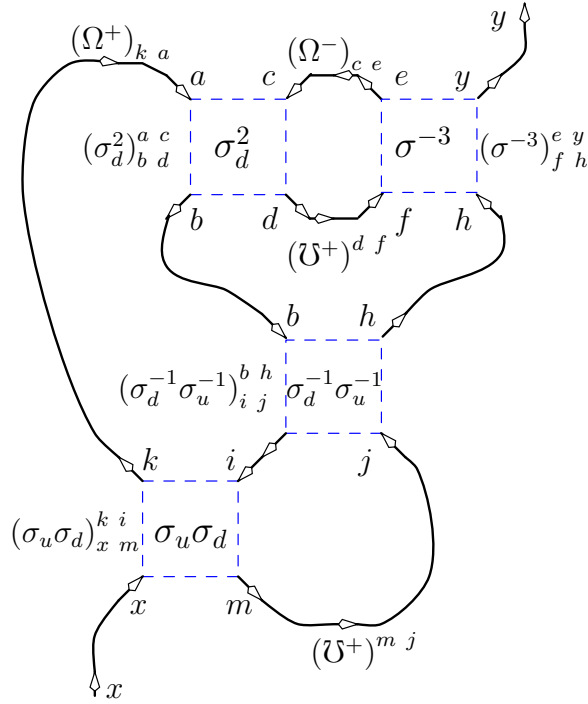}
      \caption{Tangle form of $ 9_{42} $.}
      \addtocontents{lof}{\protect\vspace{-2.5ex}}
      \figlabel{NineFortyTwo}
    \end{center}
  \end{figure}
\item[$\mathbf{10_{48}}$:]
  A braid presentation is
  $
    \sigma_1^{-2} \sigma_2^{4} \sigma_1^{-3} \sigma_2
    \in B_3
  $,
  from which we may draw the diagram found in \figref{TenFortyEight}.
  Again, we improve efficiency by the use of some auxiliary tensors:
  \be
    {\( TA \)}^{a~y~g}_{b~d~h}
    & \defeq &
    {{\( \sigma^{-2} \)}}^{a~y}_{b~f}
    \cdot
    {\( \sigma^{4} \)}^{f~g}_{d~h}
    \\
    {\( TB \)}^{b~d~h}_{c~x~i}
    & \defeq &
    {\( \sigma^{-3} \)}^{b~d}_{c~e}
    \cdot
    {\( \sigma \)}^{e~h}_{x~i}.
  \ee
  Then:
  \be
    \hspace{-13pt}
    {\( T_{10_{48}} \)}^y_x
    \defeq
    {\( TA \)}^{a~y~g}_{b~d~h}
    \cdot
    {\( TB \)}^{b~d~h}_{c~x~i}
    \cdot
    {\( \Omega^- \)}_{j~a}
    \cdot
    {\( \mho^+ \)}^{j~c}
    \cdot
    {\( \Omega^+ \)}_{g~k}
    \cdot
    {\( \mho^- \)}^{i~k}.
  \ee

  \begin{figure}[htbp]
    \begin{center}
      \input{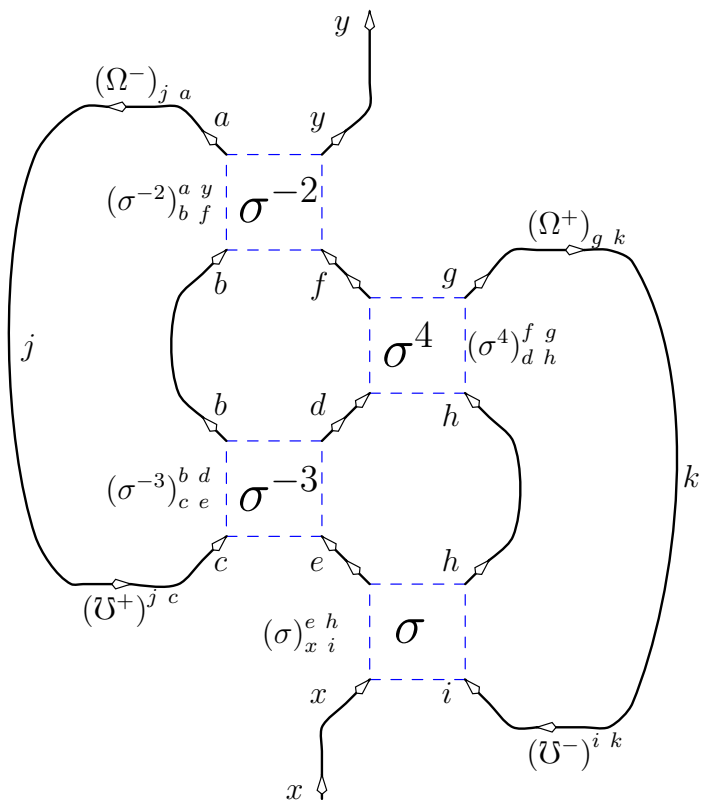}
      \caption{Tangle form of $10_{48}$.}
      \addtocontents{lof}{\protect\vspace{-2.5ex}}
      \figlabel{TenFortyEight}
    \end{center}
  \end{figure}
\end{description}

\clearpage


\subsection{A Class of Noninvertible Pretzel Knots}
\addtocontents{toc}{\protect\vspace{-2.5ex}}

The first demonstration of the existence of noninvertible knots
\cite[p25]{Livingston:93} was by Trotter in 1964 \cite{Trotter:64}, who
showed that the well-known class of \emph{pretzel knots} labelled
$(p,q,r)$ in the notation of Conway \cite{Conway:70}, were
noninvertible when $p,q$ and $r$ are distinct, odd, and greater than
$1$ (their ordering is irrelevant).  Trotter is of the opinion that the
knots are chiral (for some evidence, see
\secref{NoninvertiblePretzelsarenotDistinguished}).  This class of
pretzels provides an easily-programmable set of examples for us to
evaluate our invariant to see if it detects their noninvertibility.

The structure of the knots $(p,q,r)$ in this family is depicted in
\figref{Pretzel}, where the auxiliary tensor $X_{udu}^N$ is used.

\begin{figure}[htbp]
  \begin{center}
    \input{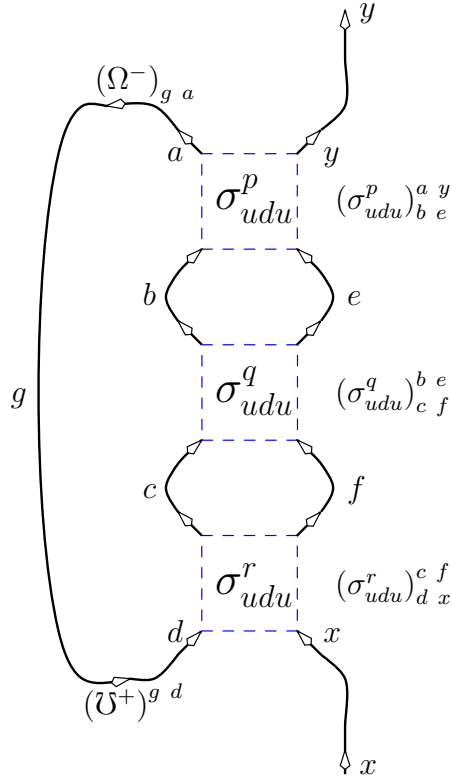}
    \caption{%
      Tangle form of the (noninvertible) pretzel knots $(p,q,r)$.
    }
    \addtocontents{lof}{\protect\vspace{-2.5ex}}
    \figlabel{Pretzel}
  \end{center}
\end{figure}

The abstract tensor associated with the pretzel is:
\be
  {(T_{(p,q,r)})}^y_x
  \defeq
  {(\sigma^{p}_{udu})}^{a~y}_{b~e}
  \cdot
  {(\sigma^{q}_{udu})}^{b~e}_{c~f}
  \cdot
  {(\sigma^{r}_{udu})}^{c~f}_{d~x}
  \cdot
  {(\Omega^-)}_{g~a}
  \cdot
  {(\mho^+)}^{g~d}.
\ee

\clearpage


\subsection{The Kinoshita--Terasaka Pair of Mutant Knots}
\addtocontents{toc}{\protect\vspace{-2.5ex}}
\seclabel{KTMutants}

The `Kinoshita--Terasaka Pair' is the best known example of a pair of
mutant knots that are known to be distinct (\cite[p106,174]{Adams:94}
and \cite{MortonCromwell:96}).  To be precise, the first of the pair is
usually known as the ``Kinoshita--Terasaka Knot'', and the second has
been called the ``Conway Knot'', as Conway used it in coining the
term ``mutant''.  In the original source by Kinoshita and Terasaka
\cite[p151]{KinoshitaTerasaka:57}, the knot involved is the one
labelled $ \kappa (2,2) $ (reproduced in \cite[p53]{Livingston:93}).
They had constructed this knot as an example of a nontrivial $11$
crossing knot with Alexander polynomial equal to $1$.  The source used
to draw our example is from \cite[p174]{Adams:94}; note that these
diagrams have $12$ crossings, so they are not minimal.

It is known that neither the HOMFLY nor the Kauffman polynomial can
distinguish \emph{any} pair of mutants \cite[p174]{Adams:94}, but
stronger statements can be made when these invariants are regarded as
being examples of Vassiliev invariants.  A Vassiliev invariant is
defined \cite[p229]{MortonCromwell:96} to be of \emph{type} $d$ if it
is zero on any link diagram of $d+1$ nodes, and of \emph{degree} $d$ if
it is of type $d$ but not of type $d-1$.  In 1994, Chmutov, Duzhin and
Lando \cite{ChmutovDuzhinLando:94} proved that \emph{all} Vassiliev
invariants of type at most $8$ will agree on \emph{any} pair of
mutants. The Links--Gould invariant is certainly a Vassiliev invariant,
although we do not know its type, but it seems unlikely that it is of
type greater than $8$. We shall investigate whether the Links--Gould
invariant distinguishes these mutants.

We illustrate $KT$, the first of Kinoshita--Terasaka pair, in
\figref{KT}, where the components $KTA$ and $KTB$ are defined below, in
Figures \ref{fig:KTA} and \ref{fig:KTC}. From $KT$, we may build the
mutant $KT'$ by replacing the component $KTA$ with $KTA'$ (depicted in
\figref{KTAprime}), which is formed by reflection of $KTA$ about a
horizontal line.

\begin{figure}[htbp]
  \begin{center}
    \input{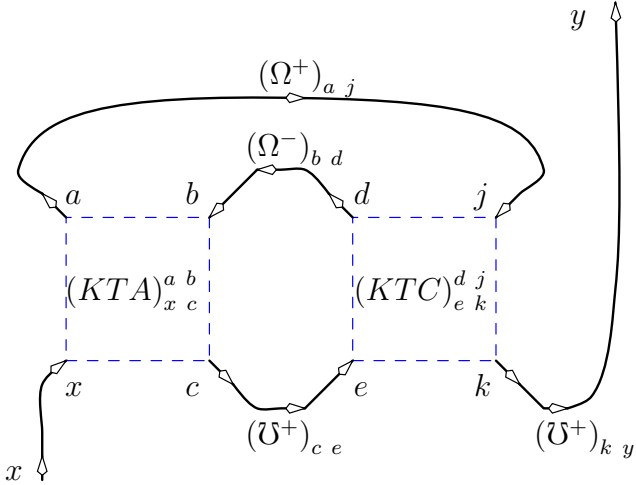}
    \caption[%
      $KT$, the first of the Kinoshita--Terasaka pair of mutant knots.
    ]{%
      $KT$, the first of the Kinoshita--Terasaka pair of mutant knots,
      where the subdiagrams $KTA$ and $KTC$ are found in
      Figures \ref{fig:KTA} and \ref{fig:KTC}
      respectively.  The mutant $KT'$ of $KT$ is obtained by
      exchanging $KTA$ with $KTA'$ (see \ref{fig:KTAprime}).
    }
    \addtocontents{lof}{\protect\vspace{-2.5ex}}
    \figlabel{KT}
  \end{center}
\end{figure}

The tensors associated with $KT$ and $KT'$ are:
\be
  {\( T_{KT} \)}^y_x
  & \defeq &
  {\( KTA \)}^{a~b}_{x~c}
  \cdot
  {\( KTC \)}^{d~j}_{e~k}
  \cdot
  {\( \Omega^- \)}_{b~d}
  \cdot
  {\( \mho^+ \)}^{c~e}
  \cdot
  {\( \Omega^+ \)}_{a~j}
  \cdot
  {\( \mho^+ \)}^{k~y}
  \\
  {\( T_{KT'} \)}^y_x
  & \defeq &
  {\( KTA' \)}^{a~b}_{x~c}
  \cdot
  {\( KTC \)}^{d~j}_{e~k}
  \cdot
  {\( \Omega^- \)}_{b~d}
  \cdot
  {\( \mho^+ \)}^{c~e}
  \cdot
  {\( \Omega^+ \)}_{a~j}
  \cdot
  {\( \mho^+ \)}^{k~y},
\ee
where:
\be
  {\( KTA \)}^{a~b}_{q~c}
  & \defeq &
  {\( \sigma_u \sigma_d \)}^{a~b}_{d~e}
  \cdot
  {\( \sigma^{-2} \)}^{d~f}_{q~g}
  \cdot
  {\( \sigma^{-1}_d \)}^{h~e}_{i~c}
  \cdot
  {\( \Omega^+ \)}_{f~h}
  \cdot
  {\( \mho^- \)}^{g~i}
  \\
  {\( KTA' \)}^{a~b}_{q~c}
  & \defeq &
  {\( \sigma^{-2} \)}^{a~f}_{d~g}
  \cdot
  {\( \sigma^{-1}_d \)}^{h~b}_{i~e}
  \cdot
  {\( \sigma_u \sigma_d \)}^{d~e}_{q~c}
  \cdot
  {\( \Omega^+ \)}_{f~h}
  \cdot
  {\( \mho^- \)}^{g~i}
  \\
  {\( KTC \)}^{d~j}_{e~k}
  & \defeq &
  {\( KTB \)}^{d~f}_{e~g}
  \cdot
  {\( \sigma^{-1}_l \sigma^{-1}_r \)}^{h~j}_{i~k}
  \cdot
  {\( \Omega^- \)}_{f~h}
  \cdot
  {\( \mho^+ \)}^{g~i}
  \\
  {\( KTB \)}^{d~f}_{e~g}
  & \defeq &
  {\( \sigma \)}^{d~b}_{a~c}
  \cdot
  {( \sigma_d^2 )}^{l~f}_{m~n}
  \cdot
  {\( \sigma^{-1}_u \sigma^{-1}_d \)}^{a~n}_{e~g}
  \cdot
  {\( \Omega^+ \)}_{b~l}
  \cdot
  {\( \mho^- \)}^{c~m}.
\ee

\begin{figure}[htbp]
  \begin{center}
    \input{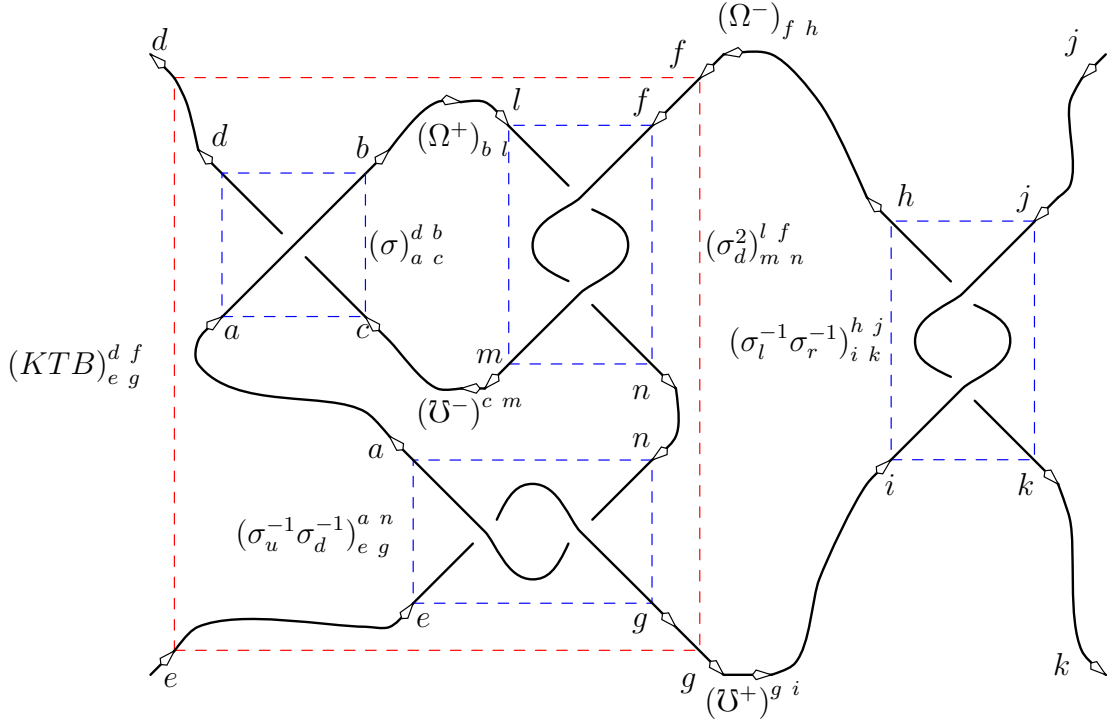}
    \caption[%
      The component $ KTC $ of the K--T pair.
    ]{%
      The component $ KTC $ of the Kinoshita--Terasaka pair.
    }
    \addtocontents{lof}{\protect\vspace{-2.5ex}}
    \figlabel{KTC}
  \end{center}
\end{figure}

\begin{figure}[htbp]
  \begin{center}
    \input{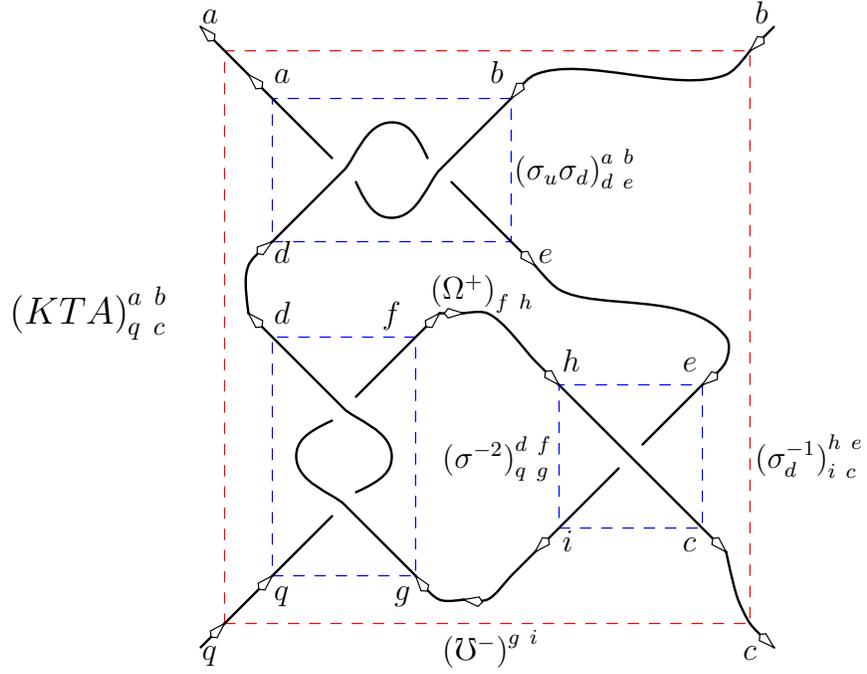}
    \caption[%
      The component $ KTA $ of $ KT $.
    ]{%
      The component $ KTA $ of $ KT $, the first of the
      Kinoshita--Terasaka pair.
    }
    \addtocontents{lof}{\protect\vspace{-2.5ex}}
    \figlabel{KTA}
  \end{center}
\end{figure}

\begin{figure}[htbp]
  \begin{center}
    \input{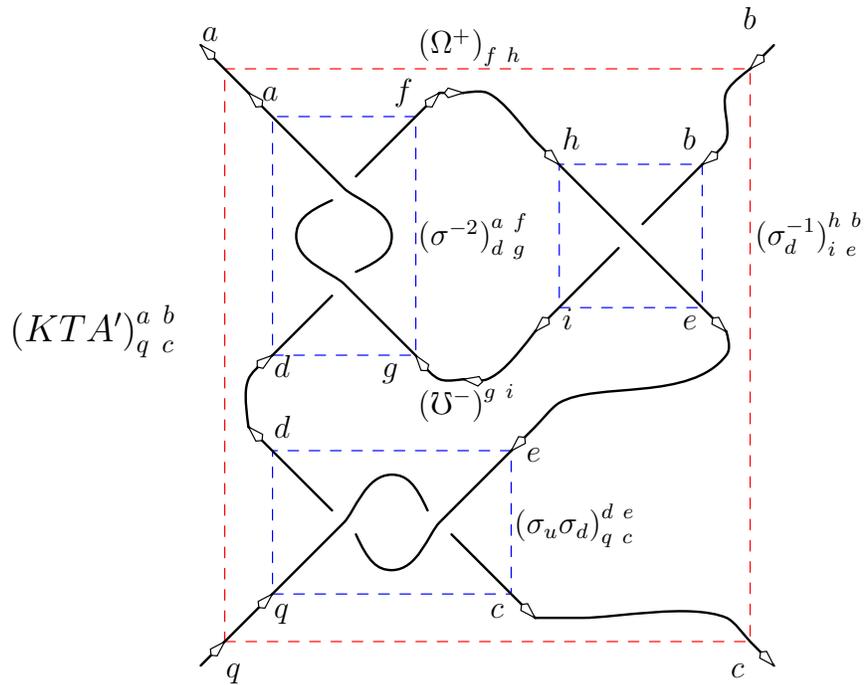}
    \caption[%
      The component $ KTA' $ of $ KT' $.
    ]{%
      The component $ KTA' $ of $ KT' $, the second of the
      Kinoshita--Terasaka pair.
    }
    \addtocontents{lof}{\protect\vspace{-2.5ex}}
    \figlabel{KTAprime}
  \end{center}
\end{figure}

\cleardoublepage


%
%

\section{The $U_q[gl(2|1)]$ Representation $\(0,0\,|\,\alpha\)$}
\addtocontents{toc}{\protect\vspace{-2.5ex}}

We consider the family of four dimensional representations of highest
weight $\(0,0\,|\,\alpha\)$ of the quantum superalgebra $U_q[gl(2|1)]$
(for nonzero complex $q$ `generic', i.e. not a root of unity). These
representations depend on a free complex parameter $\alpha$, and are
irreducible for $\alpha\neq 0,-1$, which we shall assume throughout.
Below, we use the notation of \cite{GouldHibberdLinksZhang:96}, which
itself is based on the nonquantised case $gl\(2|1\)$ found in
\cite{BrackenGouldLinksZhang:95}.  Material in this section has been
obtained by specialising that in \secref{Uqglmn}.

The material in this section is presented in great detail, with many
routine derivations presented in full. This is to ensure correctness,
and to illustrate the principles used to explicitly construct these
quantum superalgebra representations. In \secref{Conclusions}, we
discuss the idea of the automation of the construction of quantum
superalgebra representations.


\subsection{Some Standard Notation}
\addtocontents{toc}{\protect\vspace{-2.5ex}}

The following operations are actually defined for any quantum
superalgebra $U_q[\mathfrak{g}]$ (where $\mathfrak{g}$ is a Lie
superalgebra), in particular, for $U_q[gl(m|n)]$ (see
\secref{Uqglmn}).

\subsubsection{The Graded Commutator}

We shall employ the \emph{graded commutator}
$
  \[\cdot,\cdot\]
  :
  U_q[gl(2|1)] \times U_q[gl(2|1)]
  \to
  U_q[gl(2|1)]
$,
defined for homogeneous
$ x, y \in U_q[gl(2|1)] $ by:
\bse
  \[ x, y \]
  \defeq
  x y - {\(-\)}^{\[x\]\[y\]} y x,
  \eqlabel{GradedCommutator}
\ese
and extended by linearity.


\subsubsection{The $ q $ bracket}

We introduce the \emph{$ q $ bracket}, a shorthand notation defined for
$x\in U_q[gl(2|1)]$ by:
\bse
  {[x]}_q
  \defeq
  \frac{q^x - q^{-x}}{q - q^{-1}}.
  \eqlabel{qBracketDefinition}
\ese
The $q$ bracket should not be confused with the graded commutator
bracket!  Observe:
\be
  \lim_{q \to 1}{[x]}_q=x.
\ee

\pagebreak


\subsection{The Quantum Superalgebra $U_q[gl(2|1)]$}
\addtocontents{toc}{\protect\vspace{-2.5ex}}
\seclabel{Uqgl21}

\subsubsection{$ \protect \mathbb{Z}_2 $ Grading on $U_q[gl(2|1)]$
  Indices}

Firstly, we define a $\mathbb{Z}_2$ grading
$\[\cdot\] : \{ 1, 2, 3 \} \to \mathbb{Z}_2$
on the three $gl\(2|1\)$ \emph{indices}:
\be
  \begin{array}{ll}
    \[1\] = \[2\] = 0 \qquad & \textrm{(even)}
    \\
    \[3\] = 1                  & \textrm{(odd)}.
  \end{array}
\ee
Throughout, we shall use dummy indices $ a, b = 1, 2, 3 $ where
meaningful.


\subsubsection{$U_q[gl(2|1)]$ Generators}

We shall follow the notation used by Zhang
\cite[pp1237-1238]{Zhang:93} to define $U_q[gl(m|n)]$, except that
we shall substitute%
\footnote{%
  After this substitution, our notation for the Cartan generators
  approaches that of another of Zhang's papers
  \cite[pp1970-1971]{Zhang:92}, which lists ${E^a}_a$ as the Cartan
  generators. Strictly speaking, this causes problems with the
  meaning of $q^{{E^a}_a}$ as a formal power series.%
}
the notation $q^{{\(-\)}^{\[a\]}{E^a}_a}$ for $K_a$.
Specialising to the case $U_q[gl(2|1)]$, our superalgebra
has $ 7 $ \emph{simple} generators:
\be
  \left\{
    \begin{array}{ll}
      q^{{E^1}_1}, q^{{E^2}_2}, q^{{E^3}_3} \qquad
                                      & \textrm{(Cartan)} \\
      {E^1}_2, {E^2}_1                & \textrm{(raising)} \\
      {E^2}_3, {E^3}_2.               & \textrm{(lowering)}
    \end{array}
  \right\}
\ee
Note that we are \emph{not} writing the Cartan generators as
${E^a}_a$.  The fact that there are $2$ simple raising generators
indicates that $U_q[gl(2|1)]$ has rank $2$.

In addition to the above $7$ simple generators we have $2$ others
defined in terms of them \cite[p1971,~(3)]{Zhang:92} and
\cite[p1238,~(2)]{Zhang:93}:
\bse
  \left.
  \begin{array}{lrcll}
    (a) &
    {E^1}_3
    & \defeq &
    {E^1}_2
    {E^2}_3
    -
    q^{-1}
    {E^2}_3
    {E^1}_2
    \qquad
    & \textrm{(raising)}
    \\
    (b) &
    {E^3}_1
    & \defeq &
    {E^3}_2
    {E^2}_1
    -
    q
    {E^2}_1
    {E^3}_2
    &
    \textrm{(lowering)}.
  \end{array}
  \qquad \qquad \qquad
  \right\}
  \eqlabel{NonSimpleGenerators}
\ese
On \emph{all} of the $9$ \emph{generators} we define a natural
$\mathbb{Z}_2$ grading in terms of the grading on the indices:
\be
  [{E^a}_b]
  \defeq
  [a] + [b]
  \qquad
  \( \mod~2 \),
\ee
where the definition of the exponential shows that we have
$
  [q^{{E^a}_a}] = [{E^a}_a] = 0
$.
The product of homogeneous $x,y\in U_q[gl(2|1)]$ has degree:
\be
  \[x y\] = \[x\] + \[y\]
  \qquad
  \( \mod~2 \).
\ee

\vfill


\subsubsection{$U_q[gl(2|1)]$ Relations}
\seclabel{Uqgl21Relations}

For the \emph{simple} $U_q[gl(2|1)]$ generators, we have the
following $U_q[gl(2|1)]$ \emph{relations}:

\begin{itemize}
\item
  The Cartan generators all commute:
  \bse
    q^{\pm{E^a}_a}
    q^{\pm{E^b}_b}
    =
    q^{\pm{E^b}_b}
    q^{\pm{E^a}_a},
    \qquad \qquad
    q^{{E^a}_a}
    q^{-{E^a}_a}
    =
    1.
    \eqlabel{CartanGeneratorsCommute}
  \ese

\item
  The Cartan generators commute with the simple raising and lowering
  generators in the following manner:
  \bse
    q^{{\(-\)}^{\[a\]}{E^a}_a} {E^b}_{b\pm1} q^{-{\(-\)}^{\[a\]}{E^a}_a}
    =
    q^{{\(-\)}^{\[a\]} \( \delta^a_b - \delta^a_{b\pm1} \)}
    {E^b}_{b\pm1}.
    \eqlabel{CartanRaisingCommutation}
  \ese
  From \eqref{CartanRaisingCommutation}, we have the following useful
  interchange:
  \bse
    q^{{\(-\)}^{\[a\]}{E^a}_a}
    {E^b}_{b\pm1}
    =
    q^{{\(-\)}^{\[a\]} \( \delta^a_b - \delta^a_{b\pm1} \)}
    {E^b}_{b\pm1}
    q^{{\(-\)}^{\[a\]}{E^a}_a}.
    \eqlabel{UsefulInterchangePositive}
  \ese
  Replacing $q$ with $q^{-1}$ in \eqref{UsefulInterchangePositive}
  yields the equivalent:
  \bse
    q^{-{\(-\)}^{\[a\]}{E^a}_a}
    {E^b}_{b\pm1}
    =
    q^{- {\(-\)}^{\[a\]} \( \delta^a_b - \delta^a_{b\pm1} \)}
    {E^b}_{b\pm1}
    q^{-{\(-\)}^{\[a\]}{E^a}_a}.
    \eqlabel{UsefulInterchangeNegative}
  \ese
\item
  The squares of all odd generators are zero (this will also apply to
  ${E^3}_1$ and $ {E^1}_3 $, although we will never actually use this
  in practice):
  \bse
    {\( {E^2}_3 \)}^2
    =
    {\( {E^3}_2 \)}^2
    =
    0.
    \eqlabel{SquaresOfOddGeneratorsAreZero}
  \ese
\item
  The non-Cartan generators satisfy the following commutation relations
  (this is the really interesting part!):
  \bse
    \left.
    \begin{array}{rcl}
      \[ {E^1}_2, {E^2}_1 \]
      \eq
      {\[ {E^1}_1 - {E^2}_2 \]}_q
      \\
      \[ {E^2}_3, {E^3}_2 \]
      \eq
      {\[ {E^2}_2 + {E^3}_3 \]}_q
      \\
      \[ {E^1}_2, {E^3}_2 \]
      \eq
      0
      \\
      \[ {E^2}_3, {E^2}_1 \]
      \eq
      0.
    \end{array}
    \qquad \qquad \qquad \qquad \qquad \qquad \qquad \quad
    \right\}
    \eqlabel{nonCartanCommutation}
  \ese

\item
  There are two additional relations, known as the \emph{Serre
  relations}. Their inclusion ensures that the algebra is reduced
  enough to be \emph{simple}. We include them for completeness, but
  they will not concern us further.
  \be
    {\( {E^1}_2 \)}^2 {E^2}_3
    -
    \( q + q^{-1} \) {E^1}_2 {E^2}_3 {E^1}_2
    +
    {E^2}_3 {\( {E^1}_2 \)}^2
    \eq
    0
    \\
    {\( {E^2}_1 \)}^2 {E^3}_2
    -
    \( q + q^{-1} \) {E^2}_1 {E^3}_2 {E^2}_1
    +
    {E^3}_2 {\( {E^2}_1 \)}^2
    \eq
    0.
  \ee
  These relations may be more succinctly expressed using the two
  nonsimple generators, viz:
  $
    {E^1}_2 {E^1}_3 - q^{1} {E^1}_3 {E^1}_2
    =
    0
  $
  and
  $
    {E^3}_1 {E^2}_1 - q^{-1} {E^2}_1 {E^3}_1
    =
    0
  $.
\end{itemize}
More generally, expressions involving the two nonsimple generators
${E^3}_1$ and ${E^1}_3$ may be manipulated by prior expansion using
\eqref{NonSimpleGenerators}.


\subsubsection{$U_q[gl(2|1)]$ Simple Generator Interchange Rules}

Application of the definition of the graded commutator
\eqref{GradedCommutator} to the LHS of
\eqref{nonCartanCommutation} yields:
\bse
  [ {E^a}_{a+1}, {E^{b+1}}_b ]
  =
  {E^a}_{a+1} {E^{b+1}}_b
  -
  {\( - \)}^{\(\[a\]+\[a+1\]\)\(\[b\]+\[b+1\]\)}
  {E^{b+1}}_b {E^a}_{a+1}.
  \eqlabel{GradedCommutatoronGenerators}
\ese
Substitution of all possible choices for $a$ and $b$ in
\eqref{GradedCommutatoronGenerators} yields the following explicit
commutation expansions:
\bse
  \left.
  \begin{array}{rcl}
    \[ {E^1}_2, {E^2}_1 \]
    \eq
    {E^1}_2 {E^2}_1 - {E^2}_1 {E^1}_2
    \\
    \[ {E^2}_3, {E^3}_2 \]
    \eq
    {E^2}_3 {E^3}_2 + {E^3}_2 {E^2}_3
    \\
    \[ {E^1}_2, {E^3}_2 \]
    \eq
    {E^1}_2 {E^3}_2 - {E^3}_2 {E^1}_2
    \\
    \[ {E^2}_3, {E^2}_1 \]
    \eq
    {E^2}_3 {E^2}_1 - {E^2}_1 {E^2}_3.
  \end{array}
  \qquad \qquad \qquad
  \qquad \qquad \qquad \qquad
  \;
  \right\}
  \eqlabel{ExplicitCommutationRelations}
\ese
Combining the information of the commutators of
\eqref{nonCartanCommutation} with their expansions in
\eqref{ExplicitCommutationRelations} yields the following
interchanges (`commutations'):
\bse
  \left.
  \begin{array}{lrcl}
    (a) &
    {E^1}_2 {E^2}_1
    \eq
    {\[ {E^1}_1 - {E^2}_2 \]}_q
    +
    {E^2}_1 {E^1}_2
    \\
    (b) &
    {E^2}_3 {E^3}_2
    \eq
    {\[ {E^2}_2 + {E^3}_3 \]}_q
    -
    {E^3}_2 {E^2}_3
    \\
    (c) &
    {E^1}_2 {E^3}_2
    \eq
    {E^3}_2 {E^1}_2
    \\
    (d) &
    {E^2}_3 {E^2}_1
    \eq
    {E^2}_1 {E^2}_3.
  \end{array}
  \qquad \qquad \qquad \qquad
  \qquad \quad \;
  \right\}
  \eqlabel{ExplicitInterchanges1}
\ese
Lastly, expansion of \eqref{UsefulInterchangePositive} and
\eqref{UsefulInterchangeNegative} to all possible cases yields:
\bse
  \left.
  \begin{array}{lrclrcl}
    (a) &
    q^{{E^1}_1}
    {E^1}_2
    \eq
    q^{1}
    {E^1}_2
    q^{{E^1}_1}
    \qquad \quad
    &
    q^{- {E^1}_1}
    {E^1}_2
    \eq
    q^{-1}
    {E^1}_2
    q^{- {E^1}_1}
    \\
    (b) &
    q^{{E^1}_1}
    {E^2}_3
    \eq
    {E^2}_3
    q^{{E^1}_1}
    &
    q^{-{E^1}_1}
    {E^2}_3
    \eq
    {E^2}_3
    q^{-{E^1}_1}
    \\
    (c) &
    q^{{E^1}_1}
    {E^2}_1
    \eq
    q^{-1}
    {E^2}_1
    q^{{E^1}_1}
    &
    q^{-{E^1}_1}
    {E^2}_1
    \eq
    q^{1}
    {E^2}_1
    q^{-{E^1}_1}
    \\
    (d) &
    q^{{E^1}_1}
    {E^3}_2
    \eq
    {E^3}_2
    q^{{E^1}_1}
    &
    q^{-{E^1}_1}
    {E^3}_2
    \eq
    {E^3}_2
    q^{-{E^1}_1}
    \\
    (e) &
    q^{{E^2}_2}
    {E^1}_2
    \eq
    q^{-1}
    {E^1}_2
    q^{{E^2}_2}
    &
    q^{-{E^2}_2}
    {E^1}_2
    \eq
    q^{1}
    {E^1}_2
    q^{-{E^2}_2}
    \\
    (f) &
    q^{{E^2}_2}
    {E^2}_3
    \eq
    q^{1}
    {E^2}_3
    q^{{E^2}_2}
    &
    q^{-{E^2}_2}
    {E^2}_3
    \eq
    q^{-1}
    {E^2}_3
    q^{-{E^2}_2}
    \\
    (g) &
    q^{{E^2}_2}
    {E^2}_1
    \eq
    q^{1}
    {E^2}_1
    q^{{E^2}_2}
    &
    q^{-{E^2}_2}
    {E^2}_1
    \eq
    q^{-1}
    {E^2}_1
    q^{-{E^2}_2}
    \\
    (h) &
    q^{{E^2}_2}
    {E^3}_2
    \eq
    q^{-1}
    {E^3}_2
    q^{{E^2}_2}
    &
    q^{-{E^2}_2}
    {E^3}_2
    \eq
    q^{1}
    {E^3}_2
    q^{-{E^2}_2}
    \\
    (i) &
    q^{{E^3}_3}
    {E^1}_2
    \eq
    {E^1}_2
    q^{{E^3}_3}
    &
    q^{-{E^3}_3}
    {E^1}_2
    \eq
    {E^1}_2
    q^{-{E^3}_3}
    \\
    (j) &
    q^{{E^3}_3}
    {E^2}_3
    \eq
    q^{-1}
    {E^2}_3
    q^{{E^3}_3}
    &
    q^{-{E^3}_3}
    {E^2}_3
    \eq
    q^{1}
    {E^2}_3
    q^{-{E^3}_3}
    \\
    (k) &
    q^{{E^3}_3}
    {E^2}_1
    \eq
    {E^2}_1
    q^{{E^3}_3}
    &
    q^{-{E^3}_3}
    {E^2}_1
    \eq
    {E^2}_1
    q^{-{E^3}_3}
    \\
    (l) &
    q^{{E^3}_3}
    {E^3}_2
    \eq
    q^{1}
    {E^3}_2
    q^{{E^3}_3}
    &
    q^{-{E^3}_3}
    {E^3}_2
    \eq
    q^{-1}
    {E^3}_2
    q^{-{E^3}_3}.
  \end{array}
  \quad
  \right\}
  \eqlabel{ExplicitInterchanges2}
\ese
The definition of the nonsimple generators \eqref{NonSimpleGenerators},
together with the Cartan commutations
\eqref{CartanGeneratorsCommute}, the nilpotency in
\eqref{SquaresOfOddGeneratorsAreZero} and the exchanges in
\eqref{ExplicitInterchanges1} and \eqref{ExplicitInterchanges2},
provide all the necessary information for manipulation of vectors
encountered in the explicit construction of representations for
$U_q[gl(2|1)]$.

\pagebreak


\subsection{The $U_q[gl(2|1)]$ Module $V$}
\addtocontents{toc}{\protect\vspace{-2.5ex}}
\seclabel{TheUqgl21Module}

We shall be concerned with the $U_q[gl(2|1)]$ module
$V\equiv V_{(0,0\,|\,\alpha)}$, which is irreducible and typical
(recall that we insist $\alpha \neq 0,-1$). We may obtain a concrete
realisation of it using the Kac induced module construction (from which
we immediately discover that $V$ is $4$ dimensional).  This section
follows \cite{GouldHibberdLinksZhang:96}; with the caveat that there is
a small error in that work, and that we shall choose a slightly
different definition for our basis. Where differences occur, we shall
mention these.


\subsubsection{Dual Bases $\{ \ket{i} \}_{i=1}^4 $
  and $\{ \bra{i} \}_{i=1}^4 $ for $V$}

Let $ \{ \ket{i} \}_{i = 1}^4 $ denote a basis for the four dimensional
$U_q[gl(2|1)] $ module $ V $.  In accordance with the Kac induced
module construction, we shall define these basis vectors as created by
the action of all possible nonzero combinations of all (not just the
simple) the \emph{odd} lowering generators on $\ket{1}$, defined to be
the highest weight state:
\bse
  \left.
  \begin{array}{lrcl}
    (a) &
    \ket{2}
    & \defeq &
    \beta_2
    {E^3}_2
    \cdot
    \ket{1}
    \\
    (b) &
    \ket{3}
    & \defeq &
    \beta_3
    {E^3}_1
    \cdot
    \ket{1}
    \\
    (c) &
    \ket{4}
    & \defeq &
    \beta_4
    {E^3}_2
    {E^3}_1
    \cdot
    \ket{1},
  \end{array}
  \qquad \qquad \qquad \qquad
  \qquad \qquad \qquad \qquad \; \;
  \right\}
  \eqlabel{InitialKetDefinition}
\ese
where the $ \beta_i, i = 2, 3, 4 $ are suitable normalisation
constants, to be deduced below (implicitly, $\beta_1 \equiv 1$).
We shall use the symbol ``$\cdot$'' to represent algebra-module
multiplication.

The above definition \eqref{InitialKetDefinition} differs from that
(implicitly) used in \cite{GouldHibberdLinksZhang:96} in that the
definitions of $\ket{2}$ and $\ket{3}$ are interchanged. The choice
here is more natural in that the vectors are ordered in terms of
decreasing weight.  To be certain:
\bse
  \left.
  \begin{array}{
      ccl@{}r@{\hspace{1mm}}r@{\hspace{1mm}}c@{\hspace{1mm}}l@{}l@{}c
               }
    \ket{1} & \textrm{has~weight} & ( &  0, &  0 & | & \alpha   & ) & \\
    \ket{2} &                     & ( &  0, & -1 & | & \alpha+1 & ) & \\
    \ket{3} &                     & ( & -1, &  0 & | & \alpha+1 & ) & \\
    \ket{4} &                     & ( & -1, & -1 & | & \alpha+2 & ) & .
  \end{array}
  \qquad \qquad \qquad \qquad \qquad \qquad \; \;
  \right\}
  \eqlabel{WeightsofKets}
\ese

Consistent with the $\mathbb{Z}_2$ grading on $U_q[gl(2|1)]$,
following \cite[p158]{GouldHibberdLinksZhang:96}, we grade the basis
states by:
\be
  \[ \ket{1} \] = \[ \ket{4} \] \defeq 0,
  \qquad \qquad
  \[ \ket{2} \] = \[ \ket{3} \] \defeq 1.
\ee
By ``consistent'', we mean that the action of a homogeneous
$U_q[gl(2|1)]$ generator $x$ on a (homogeneous)
basis vector $\ket{i}$, yields a homogeneous
vector of degree $[x]+[\ket{i}]$:
\be
  [x \cdot \ket{i}]
  =
  [x]
  +
  [\ket{i}].
\ee
Thus, declaring $\ket{1}$ as even
means that $\ket{2}$ and $\ket{3}$, which are created from the action
of odd lowering generators on $\ket{1}$, must necessarily be odd.
Lastly, $\ket{4}$ may be regarded as being defined in terms of the
action of an odd lowering generator on the odd $\ket{2}$, hence is
necessarily even.

We define a basis $\{\bra{i}\}_{i=1}^4$ \emph{dual} to
$\{\ket{i}\}_{i=1}^4$, where $[\bra{i}]\defeq[\ket{i}]$.  Written in
component form, these dual vectors are represented by the transpose
complex conjugates of the original basis:
$\bra{i}\defeq\overline{\ket{i}}^T\equiv\ket{i}^\dagger$. Then we
have:
\be
  \bra{i} \ket{j} \equiv \langle i | j \rangle
  \defeq
  \delta_{i j}.
\ee


\subsubsection{Action of the Generators on the Basis}

The module $V$ has highest weight vector $\ket{1}$, of weight
$\(0,0\,|\,\alpha\)$.  Using the shorthand
$\lambda\defeq\(0,0\,|\,\alpha\)$, we have $\lambda_1=\lambda_2=0$ and
$\lambda_3=\alpha$.

\begin{itemize}
\item
  As $\ket{1}$ is a \emph{weight vector}, the action of the Cartan
  generators on $ \ket{1} $ is:%
  \footnote{%
    The equations of \eqref{ExplicitCartanActiononket1} are of course
    the relations that would define a $gl\(2|1\)$ representation as
    having weight $\(0,0\,|\,\alpha\)$.  As we are regarding
    $q^{{\(-\)}^{\[a\]} {E^a}_a}$ as the generators of
    $U_q[gl(2|1)]$, perhaps it would be more honest to label this
    the $\(1,1\,|\,q^{-\alpha}\)$ representation, as per
    \eqref{ExplicitExpCartanActiononket1}.
  }
  \bse
    {E^i}_i
    \cdot
    \ket{1}
    \defeq
    \lambda_i
    \ket{1},
    \eqlabel{CartanActiononket1}
  \ese
  i.e.:
  \bse
    {E^1}_1
    \cdot
    \ket{1}
    =
    0,
    \qquad \qquad
    {E^2}_2
    \cdot
    \ket{1}
    =
    0,
    \qquad \qquad
    {E^3}_3
    \cdot
    \ket{1}
    =
    \alpha
    \ket{1}.
    \eqlabel{ExplicitCartanActiononket1}
  \ese
  Expansion of the exponential, and the use of
  \eqref{CartanActiononket1} yields:
  \be
    q^{{E^i}_i}
    \cdot
    \ket{1}
    =
    q^{\lambda_i}
    \ket{1}.
  \ee
  We take a moment to justify this, using the observation
  that, recursively for all integers $ j \geqslant 2 $ we have
  $
    {\( {E^i}_i \)}^j
    \cdot
    \ket{1}
    =
    {\( {E^i}_i \)}^{j-1}
    {E^i}_i
    \cdot
    \ket{1}
    =
    \lambda_i
    {\( {E^i}_i \)}^{j-1}
    \cdot
    \ket{1}
  $,
  hence
  $
    {({E^i}_i)}^j \cdot \ket{1}
    =
    {(\lambda_i)}^j \ket{1}
  $.
  \be
    ~\hspace{-10mm}
    q^{{E^i}_i}
    \cdot
    \ket{1}
    \eq
    e^{\ln\(q\) {E^i}_i}
    \cdot
    \ket{1}
    =
    \sum_{j=0}^\infty
      \frac{{\( \ln\(q\) {E^i}_i \)}^j}{j!}
    \cdot
    \ket{1}
    =
    \sum_{j=0}^\infty
      \frac{{\ln\(q\)}^j}{j!}
      {\({E^i}_i\)}^j
      \cdot
      \ket{1}
    \\
    \eq
    \sum_{j=0}^\infty
      \frac{{\ln\(q\)}^j}{j!}
      {\(\lambda_i\)}^j
      \ket{1}
    =
    \sum_{j=0}^\infty
      \frac{{\( \ln\(q\) \lambda_i\)}^j}{j!}
    \ket{1}
    =
    e^{\ln\(q\)\lambda_i}
    \ket{1}
    =
    q^{\lambda_i}
    \ket{1}.
  \ee

  \vfill

  Thus, \eqref{ExplicitCartanActiononket1} is equivalent to:
  \bse
    q^{{E^1}_1}
    \cdot
    \ket{1}
    =
    q^{{E^2}_2}
    \cdot
    \ket{1}
    =
    q^0
    \ket{1}
    =
    \ket{1},
    \qquad \qquad
    q^{{E^3}_3}
    \cdot
    \ket{1}
    =
    q^{\alpha}
    \ket{1}.
    \eqlabel{ExplicitExpCartanActiononket1}
  \ese
  Temporarily writing ${E^0}_0\equiv1$ and $\lambda_0\equiv1$, for any
  complex scalars $\gamma_i$:
  \be
    q^{\gamma_i {E^i}_i}
    \cdot
    \ket{1}
    =
    q^{\lambda_i \gamma_i}
    \ket{1},
    \qquad
    i = 0, \dots, 3,
  \ee
  which may be extended to:
  \bse
    q^{\sum_i \gamma_i {E^i}_i}
    \cdot
    \ket{1}
    =
    q^{\sum_i \lambda_i \gamma_i}
    \ket{1}.
    \eqlabel{Botia111}
  \ese
  Using \eqref{Botia111} , we have:
  \be
    {\[
      {\textstyle \sum_{i=0}^3}
        \gamma_i {E^i}_i
    \]}_q
    \cdot
    \ket{1}
    \eq
    {\( q - q^{-1} \)}^{-1}
    \(
      {\textstyle \prod_{i=0}^3} q^{\gamma_i {E^i}_i}
      -
      {\textstyle \prod_{i=0}^3} q^{-\gamma_i {E^i}_i}
    \)
    \cdot
    \ket{1}
    \\
    \eq
    {\( q - q^{-1} \)}^{-1}
    \(
      {\textstyle \prod_{i=0}^3} q^{\gamma_i \lambda_i}
      -
      {\textstyle \prod_{i=0}^3} q^{-\gamma_i \lambda_i}
    \)
    \ket{1}
    \\
    \eq
    {\[ {\textstyle \sum}_{i=0}^3 \gamma_i \lambda_i \]}_q
    \ket{1}.
  \ee
  In our case, where $\lambda=(0,0\,|\,\alpha)$, this simplifies
  to:
  \bse
    {\[
      {\textstyle \sum_{i=0}^3}
        \gamma_i {E^i}_i
    \]}_q
    \cdot
    \ket{1}
    =
    {\[ \gamma_0 + \gamma_3 \alpha \]}_q
    \ket{1}.
    \eqlabel{SumExpCartanActiononket1}
  \ese
\item
  As $ \ket{1} $ is a \emph{highest} weight vector, it is annihilated
  by the action of the $U_q[gl(2|1)]$ \emph{raising}
  generators:
  \bse
    {E^1}_2
    \cdot
    \ket{1}
    =
    {E^2}_3
    \cdot
    \ket{1}
    =
    {E^1}_3
    \cdot
    \ket{1}
    =
    0.
    \eqlabel{Annihilationofket1byRaisingGens}
  \ese
\item
  The even subalgebra of $gl(2|1)$ is
  $ {gl(2|1)}_{\overline{0}} = gl(2) \oplus gl(1)$,
  and the Kac induced module construction builds a $gl(2|1)$
  module by acting on a ${gl \(2| 1 \)}_{\overline{0}}$ module with
  all the odd lowering generators. This naturally extends to the
  $U_q[gl(2|1)]$ case.

  As $\ket{1}$ spans a one dimensional
  ${U_q[gl(2|1)]}_0=U_q\[gl\(2\)\]\oplus U_q\[gl\(1\)\]$
  module, it is annihilated by the action of the even $U_q\[gl\(2\)\]$
  \emph{lowering} generator $ {E^2}_1 $:
  \bse
    {E^2}_1
    \cdot
    \ket{1}
    =
    0.
    \eqlabel{E21Annihilatesket1}
  \ese
\end{itemize}

The definitions of the duals imply that \eqref{CartanActiononket1},
\eqref{Annihilationofket1byRaisingGens}  and \eqref{E21Annihilatesket1}
have natural analogues:
\be
  \bra{1} \cdot {E^i}_i
  =
  \lambda_i \bra{1},
  \qquad
  \bra{1} \cdot {E^2}_1
  =
  \bra{1} \cdot {E^3}_2
  =
  \bra{1} \cdot {E^3}_1
  =
  0,
  \quad
  \bra{1}
  \cdot
  {E^1}_2
  =
  0.
\ee
We shall use these implicitly, without further comment.

\pagebreak


\subsubsection{The Normalisation Constants $ \beta_i $}
\seclabel{NormalisationConstants}

We demand that the basis $\{\ket{i}\}_{i=1}^4$ is normalised, viz:
\bse
  \braket{i}{i}
  =
  1,
  \eqlabel{Braketii=1}
\ese
and this defines $ \beta_i $, for $ i = 2, 3, 4 $.  Let us assume for
now that $ \alpha $ and $ q $ are real and positive. (Our results for
the R matrix will hold later by analytic continuation.) Initially,
we declare that $ \ket{1} $ is already normalised, i.e.
$\braket{1}{1}=1$, equivalently $ \beta_1 = 1 $. Before proceeding to
deduce the normalisation constants, we deduce some auxiliary
results, which will also be of use when we come to discover the
representation itself.
\bne
  \lefteqn{
    \hspace{-5mm}
    {\[ {E^2}_2 + {E^3}_3 \]}_q
    {E^2}_1
    {E^3}_2
    \cdot
    \ket{1}
  }
  \nonumber
  \\
  & \stackrel{(\ref{eq:qBracketDefinition},
               \ref{eq:CartanGeneratorsCommute})}{=} &
  {\( q - q^{-1} \)}^{-1}
  \(
    q^{{E^2}_2}
    q^{{E^3}_3}
    {E^2}_1
    -
    q^{-{E^2}_2}
    q^{-{E^3}_3}
    {E^2}_1
  \)
  {E^3}_2
  \cdot
  \ket{1}
  \nonumber
  \\
  & \stackrel{(\ref{eq:ExplicitInterchanges2}k)}{=} &
  {\( q - q^{-1} \)}^{-1}
  \(
    q^{{E^2}_2}
    {E^2}_1
    q^{{E^3}_3}
    -
    q^{-{E^2}_2}
    {E^2}_1
    q^{-{E^3}_3}
  \)
  {E^3}_2
  \cdot
  \ket{1}
  \nonumber
  \\
  & \stackrel{(\ref{eq:ExplicitInterchanges2}g)}{=} &
  {E^2}_1
  {\( q - q^{-1} \)}^{-1}
  \(
    q^{1}
    q^{{E^2}_2}
    q^{{E^3}_3}
    -
    q^{-1}
    q^{-{E^2}_2}
    q^{-{E^3}_3}
  \)
  {E^3}_2
  \cdot
  \ket{1}
  \nonumber
  \\
  & \stackrel{(\ref{eq:ExplicitInterchanges2}l)}{=} &
  {E^2}_1
  {\( q - q^{-1} \)}^{-1}
  \(
    q^{2}
    q^{{E^2}_2}
    {E^3}_2
    q^{{E^3}_3}
    -
    q^{-2}
    q^{-{E^2}_2}
    {E^3}_2
    q^{-{E^3}_3}
  \)
  \cdot
  \ket{1}
  \nonumber
  \\
  & \stackrel{(\ref{eq:ExplicitInterchanges2}h)}{=} &
  {E^2}_1
  {E^3}_2
  {\( q - q^{-1} \)}^{-1}
  \(
    q^{1}
    q^{{E^2}_2}
    q^{{E^3}_3}
    -
    q^{-1}
    q^{-{E^2}_2}
    q^{-{E^3}_3}
  \)
  \cdot
  \ket{1}
  \nonumber
  \\
  \eq
  {E^2}_1
  {E^3}_2
  {\[ 1 + {E^2}_2 + {E^3}_3 \]}_q
  \cdot
  \ket{1}
  \stackrel{\eqref{SumExpCartanActiononket1}}{=}
  {\[\alpha+1\]}_q
  {E^2}_1
  {E^3}_2
  \cdot
  \ket{1}.
  \eqlabel{Botia38}
\ene
\bne
  \lefteqn{
    \hspace{-5mm}
    {\[ {E^1}_1 - {E^2}_2 \]}_q
    {E^3}_2
    \cdot
    \ket{1}
  }
  \nonumber
  \\
  & \stackrel{(\ref{eq:qBracketDefinition},
               \ref{eq:CartanGeneratorsCommute})}{=} &
  {\( q - q^{-1} \)}^{-1}
  \(
    q^{{E^1}_1} q^{- {E^2}_2}
    {E^3}_2
    -
    q^{- {E^1}_1} q^{{E^2}_2}
    {E^3}_2
  \)
  \cdot
  \ket{1}
  \nonumber
  \\
  & \stackrel{(\ref{eq:ExplicitInterchanges2}h)}{=} &
  {\( q - q^{-1} \)}^{-1}
  \(
    q^{1}
    q^{{E^1}_1}
    {E^3}_2
    q^{- {E^2}_2}
    -
    q^{-1}
    q^{- {E^1}_1}
    {E^3}_2
    q^{{E^2}_2}
  \)
  \cdot
  \ket{1}
  \nonumber
  \\
  & \stackrel{(\ref{eq:ExplicitInterchanges2}d)}{=} &
  {\( q - q^{-1} \)}^{-1}
  \(
    q^{1}
    {E^3}_2
    q^{{E^1}_1}
    q^{- {E^2}_2}
    -
    q^{-1}
    {E^3}_2
    q^{- {E^1}_1}
    q^{{E^2}_2}
  \)
  \cdot
  \ket{1}
  \nonumber
  \\
  \eq
  {E^3}_2
  {\[1 + {E^1}_1 - {E^2}_2 \]}_q
  \cdot
  \ket{1}
  \stackrel{\eqref{SumExpCartanActiononket1}}{=}
  {\[1\]}_q
  {E^3}_2
  \cdot
  \ket{1}
  =
  {E^3}_2
  \cdot
  \ket{1}.
  \eqlabel{Botia31}
\ene
\bse
  ~\hspace{-6mm}
  {E^3}_1
  \cdot
  \ket{1}
  \stackrel{(\ref{eq:NonSimpleGenerators}b)}{=}
  \( {E^3}_2 {E^2}_1 - q {E^2}_1 {E^3}_2 \)
  \cdot
  \ket{1}
  \stackrel{\eqref{E21Annihilatesket1}}{=}
  - q {E^2}_1 {E^3}_2
  \cdot
  \ket{1}.
  \eqlabel{Botia23}
\ese
\bne
  ~\hspace{-6mm}
  {E^1}_2
  {E^3}_1
  \cdot
  \ket{1}
  & \stackrel{\eqref{Botia23}}{=} &
  - q
  {E^1}_2
  {E^2}_1
  {E^3}_2
  \cdot
  \ket{1}
  \nonumber
  \\
  & \stackrel{(\ref{eq:ExplicitInterchanges1}a)}{=} &
  - q
  \(
    {\[ {E^1}_1 - {E^2}_2 \]}_q
    +
    {E^2}_1 {E^1}_2
  \)
  {E^3}_2
  \cdot
  \ket{1}
  \nonumber
  \\
  & \stackrel{(\ref{eq:ExplicitInterchanges1}c)}{=} &
  - q
  \(
    {\[ {E^1}_1 - {E^2}_2 \]}_q
    {E^3}_2
    +
    {E^2}_1
    {E^3}_2
    {E^1}_2
  \)
  \cdot
  \ket{1}
  \nonumber
  \\
  & \stackrel{\eqref{Annihilationofket1byRaisingGens}}{=} &
  - q
  {\[ {E^1}_1 - {E^2}_2 \]}_q
  {E^3}_2
  \cdot
  \ket{1}
  \stackrel{\eqref{Botia31}}{=}
  - q
  {E^3}_2
  \cdot
  \ket{1}.
  \eqlabel{Botia46}
\ene
\bne
  {E^2}_3
  {E^3}_1
  \cdot
  \ket{1}
  & \stackrel{(\ref{eq:NonSimpleGenerators}b)}{=} &
  {E^2}_3
  \(
    {E^3}_2
    {E^2}_1
    -
    q
    {E^2}_1
    {E^3}_2
  \)
  \cdot
  \ket{1}
  \nonumber
  \\
  & \stackrel{\eqref{E21Annihilatesket1}}{=} &
  -
  q
  {E^2}_3
  {E^2}_1
  {E^3}_2
  \cdot
  \ket{1}
  \stackrel{(\ref{eq:ExplicitInterchanges1}d)}{=}
  -
  q
  {E^2}_1
  {E^2}_3
  {E^3}_2
  \cdot
  \ket{1}
  \nonumber
  \\
  & \stackrel{(\ref{eq:ExplicitInterchanges1}b)}{=} &
  -
  q
  {E^2}_1
  \(
    {\[ {E^2}_2 + {E^3}_3 \]}_q
    -
    {E^3}_2 {E^2}_3
  \)
  \cdot
  \ket{1}
  \nonumber
  \\
  & \stackrel{\eqref{Annihilationofket1byRaisingGens}}{=} &
  -
  q
  {E^2}_1
  {\[ {E^2}_2 + {E^3}_3 \]}_q
  \cdot
  \ket{1}
  \stackrel{\eqref{SumExpCartanActiononket1}}{=}
  -
  q
  {\[\alpha\]}_q
  {E^2}_1
  \cdot
  \ket{1}
  \stackrel{\eqref{E21Annihilatesket1}}{=}
  0.
  \eqlabel{Botia42}
\ene
\bne
  {E^2}_3 {E^3}_2
  \cdot
  \ket{1}
  & \stackrel{(\ref{eq:ExplicitInterchanges1}b)}{=} &
  \(
    {\[{E^2}_2 + {E^3}_3 \]}_q
    -
    {E^3}_2 {E^2}_3
  \)
  \cdot
  \ket{1}
  \stackrel{\eqref{Annihilationofket1byRaisingGens}}{=}
  {\[{E^2}_2 + {E^3}_3 \]}_q
  \cdot
  \ket{1}
  \nonumber
  \\
  & \stackrel{\eqref{SumExpCartanActiononket1}}{=} &
  {\[\alpha\]}_q
  \ket{1}.
  \eqlabel{Botia32}
\ene
\bne
  \! \! \! \! \!
  {E^2}_3
  {E^1}_2
  {E^2}_1
  {E^3}_2
  \cdot
  \ket{1}
  & \stackrel{(\ref{eq:ExplicitInterchanges1}a)}{=} &
  {E^2}_3
  \(
    {\[ {E^1}_1 - {E^2}_2 \]}_q
    +
    {E^2}_1 {E^1}_2
  \)
  {E^3}_2
  \cdot
  \ket{1}
  \nonumber
  \\
  & \stackrel{(\ref{eq:ExplicitInterchanges1}c,d)}{=} &
  \(
    {E^2}_3
    {\[ {E^1}_1 - {E^2}_2 \]}_q
    {E^3}_2
    +
    {E^2}_1
    {E^2}_3
    {E^3}_2
    {E^1}_2
  \)
  \cdot
  \ket{1}
  \nonumber
  \\
  & \stackrel{\eqref{Annihilationofket1byRaisingGens}}{=} &
  {E^2}_3
  {\[ {E^1}_1 - {E^2}_2 \]}_q
  {E^3}_2
  \cdot
  \ket{1}
  \nonumber
  \\
  & \stackrel{\eqref{Botia31}}{=} &
  {E^2}_3
  {E^3}_2
  \cdot
  \ket{1}
  \stackrel{\eqref{Botia32}}{=}
  {\[\alpha\]}_q
  \ket{1}.
  \eqlabel{Botia12}
\ene

We now compute the normalisation coefficients:
\be
  \braket{2}{2}
  =
  \bra{2} \cdot \ket{2}
  \stackrel{(\ref{eq:InitialKetDefinition}b)}{=}
  \beta_2 \beta_2^*
  \bra{1}
  \cdot
  {E^2}_3 {E^3}_2
  \cdot
  \ket{1}
  \stackrel{\eqref{Botia32}}{=}
  {\[\alpha\]}_q
  \beta_2 \beta_2^*.
\ee
Demanding that $ \braket{2}{2} = 1 $ then amounts to selecting
$\beta_2$ as real and equal to ${\[\alpha\]}_q^{-1/2}$ (up to
a phase factor, which we shall set as $1$).
\be
  \braket{3}{3}
  \eq
  \bra{3} \cdot \ket{3}
  \stackrel{(\ref{eq:InitialKetDefinition}a)}{=}
  \beta_3
  \beta_3^*
  \bra{1}
  \cdot
  {E^1}_3
  {E^3}_1
  \cdot
  \ket{1}
  \\
  & \stackrel{\eqref{Botia23}}{=} &
  q^2
  \beta_3
  \beta_3^*
  \bra{1}
  \cdot
  {E^2}_3
  {E^1}_2
  {E^2}_1
  {E^3}_2
  \cdot
  \ket{1}
  \stackrel{\eqref{Botia12}}{=}
  q^2
  {\[\alpha\]}_q
  \beta_3
  \beta_3^*.
\ee
Demanding that $ \braket{3}{3} = 1 $ then amounts to selecting
$\beta_3$ as real and equal to $q^{-1}{\[\alpha \]}_q^{-1/2}$
(again, up to a phase factor, which we shall set as $1$).

\be
  \hspace{-28pt}
  \begin{array}{cl}
    \hspace{-5pt}
    \braket{4}{4}
    \\
    =
    &
    \bra{4} \cdot \ket{4}
    \stackrel{(\ref{eq:InitialKetDefinition}c)}{=}
    \beta_4 \beta_4^*
    \bra{1}
    \cdot
    {E^1}_3
    {E^2}_3
    {E^3}_2
    {E^3}_1
    \cdot
    \ket{1}
    \\
    \stackrel{(\ref{eq:ExplicitInterchanges1}b)}{=} &
    \beta_4 \beta_4^*
    \bra{1}
    \cdot
    {E^1}_3
    \(
      {\[{E^2}_2 + {E^3}_3 \]}_q
      -
      {E^3}_2 {E^2}_3
    \)
    {E^3}_1
    \cdot
    \ket{1}
    \\
    \stackrel{\eqref{Botia23}}{=} &
    q^2
    \beta_4 \beta_4^*
    \bra{1}
    \cdot
    {E^2}_3
    {E^1}_2
    \(
      {\[{E^2}_2 + {E^3}_3 \]}_q
      -
      {E^3}_2 {E^2}_3
    \)
    {E^2}_1 {E^3}_2
    \cdot
    \ket{1}
    \\
    \stackrel{(\ref{eq:ExplicitInterchanges1}d)}{=} &
    q^2
    \beta_4 \beta_4^*
    \bra{1}
    \cdot
    \(
      {E^2}_3
      {E^1}_2
      {\[{E^2}_2 + {E^3}_3 \]}_q
      {E^2}_1
      {E^3}_2
       -
      {E^2}_3
      {E^3}_2
      {E^1}_2
      {E^2}_1
      {E^2}_3
      {E^3}_2
    \)
    \cdot
    \ket{1}
    \\
    \stackrel{\eqref{Botia32}}{=} &
    q^2
    \beta_4 \beta_4^*
    \(
      \bra{1}
      \cdot
      {E^2}_3
      {E^1}_2
      {\[{E^2}_2 + {E^3}_3 \]}_q
      {E^2}_1
      {E^3}_2
      \cdot
      \ket{1}
       -
      {\[\alpha\]}_q^2
      \bra{1}
      \cdot
      {E^1}_2
      {E^2}_1
      \cdot
      \ket{1}
    \)
    \\
    \stackrel{\eqref{E21Annihilatesket1}}{=} &
    q^2
    \beta_4 \beta_4^*
    \bra{1}
    \cdot
    {E^2}_3
    {E^1}_2
    {\[{E^2}_2 + {E^3}_3 \]}_q
    {E^2}_1
    {E^3}_2
    \cdot
    \ket{1}
    \\
    \stackrel{\eqref{Botia38}}{=} &
    q^2
    {\[\alpha+1\]}_q
    \beta_4 \beta_4^*
    \bra{1}
    \cdot
    {E^2}_3
    {E^1}_2
    {E^2}_1
    {E^3}_2
    \cdot
    \ket{1}
    \\
    \stackrel{\eqref{Botia12}}{=} &
    q^2
    {\[\alpha\]}_q
    {\[\alpha+1\]}_q
    \beta_4 \beta_4^*.
  \end{array}
\ee
Thus, up to a phase factor, we have:
$
  \beta_4
  =
  q^{-1}
  {\[\alpha\]}_q^{-1/2}
  {\[\alpha+1\]}_q^{-1/2}
$.
In conclusion, we have the normalisation constants:
\bse
  \left.
  \begin{array}{lrcl}
    (a) &
    \beta_2
    \eq
    {\[\alpha\]}_q^{-1/2}
    \\
    (b) &
    \beta_3
    \eq
    q^{-1}
    {\[\alpha\]}_q^{-1/2}
    \\
    (c) &
    \beta_4
    \eq
    q^{-1}
    {\[\alpha\]}_q^{-1/2}
    {\[\alpha+1\]}_q^{-1/2},
  \end{array}
  \qquad \qquad \qquad \qquad \qquad \qquad \qquad \; \;
  \right\}
  \eqlabel{NormalisationConstants}
\ese
that is, the definitions of \eqref{InitialKetDefinition} may be
replaced by:
\bne
  \left.
  \begin{array}{lrcl}
    (a) &
    \ket{2}
    & \defeq &
    {\[\alpha\]}_q^{-1/2}
    {E^3}_2
    \cdot
    \ket{1}
    \\
    (b) &
    \ket{3}
    & \defeq &
    q^{-1}
    {\[\alpha\]}_q^{-1/2}
    {E^3}_1
    \cdot
    \ket{1}
    \\
    (c) &
    \ket{4}
    & \defeq &
    q^{-1}
    {\[\alpha\]}_q^{-1/2}
    {\[\alpha+1\]}_q^{-1/2}
    {E^3}_2
    {E^3}_1
    \cdot
    \ket{1}.
  \end{array}
  \qquad \qquad \qquad \qquad \; \;
  \right\}
  \eqlabel{FinalKetDefinition}
\ene

\pagebreak


\subsection{The Representation $\pi \equiv \pi_\lambda$}
\addtocontents{toc}{\protect\vspace{-2.5ex}}

In terms of these dual bases ${\{\ket{i}\}}_{i=1}^4$ and
${\{\bra{i}\}}_{i=1}^4$, we define a representation
$\pi\equiv\pi_\lambda$ on $U_q[gl(2|1)]$. For the generators, we
have:
\bse
  \lefteqn{
    \hspace{-5mm}
    \left.
    \begin{array}{lrcl}
      (a) &
      \pi \( {E^1}_1 \)
      \eq
      - \ket{3} \bra{3} - \ket{4} \bra{4}
      \\
      (b) &
      \pi \( {E^2}_2 \)
      \eq
      - \ket{2} \bra{2} - \ket{4} \bra{4}
      \\
      (c) &
      \pi \( {E^3}_3 \)
      \eq
      \alpha \ket{1} \bra{1}
      +
      \(\alpha+1\) \ketbra{2}{2}
      +
      \(\alpha+1\) \ketbra{3}{3}
      +
      \(\alpha+2\) \ketbra{4}{4}
      \qquad
      \\
      (d) &
      \pi \( {E^1}_2 \)
      \eq
      -
      \ket{2} \bra{3}
      \\
      (e) &
      \pi \( {E^2}_1 \)
      \eq
      -
      \ket{3} \bra{2}
      \\
      (f) &
      \pi \( {E^2}_3 \)
      \eq
      {\[\alpha\]}_q^{1/2} \ket{1} \bra{2}
      +
      {\[\alpha+1\]}_q^{1/2} \ket{3} \bra{4}
      \\
      (g) &
      \pi \( {E^3}_2 \)
      \eq
      {\[\alpha\]}_q^{1/2} \ket{2} \bra{1}
      +
      {\[\alpha+1\]}_q^{1/2} \ket{4} \bra{3}.
    \end{array}
    \! \! \! \! \! \! \! \! \! \!
    \! \! \!
    \right\}
  }
  \eqlabel{TheRepresentation}
\ese
This comes from \cite[p157,~(3)]{GouldHibberdLinksZhang:96}, although
that paper does not describe any of the conventions concerning the
definition of the basis vectors, or the phase factors in their
normalisation constants. The interchange of the definitions of
$\ket{2}$ and $\ket{3}$ between that paper and this work generates the
same changes in \eqref{TheRepresentation}. (Our paper
\cite[p3]{DeWitKauffmanLinks:98} uses the conventions of
\cite{GouldHibberdLinksZhang:96} rather than the above.)

Observe that we have used the expression ${E^a}_a$ here, rather than
the more technically correct $q^{{E^a}_a}$ (or even
$q^{{(-)}^{[a]}{E^a}_a}$).  The representation is of course extended to
the rest of $U_q[gl(2|1)]$ by linearity.  In particular, for future
reference, we have:
\bse
  \lefteqn{
    \hspace{-5mm}
    \left.
    \begin{array}{lrcl}
      (a) &
      \pi ( q^{{E^1}_1} )
      \eq
      \ketbra{1}{1}
      +
      \ketbra{2}{2}
      +
      q^{-1}
      \ketbra{3}{3}
      +
      q^{-1}
      \ketbra{4}{4}
      \\
      (b) &
      \pi ( q^{{E^2}_2} )
      \eq
      \ketbra{1}{1}
      +
      q^{-1}
      \ketbra{2}{2}
      +
      \ketbra{3}{3}
      +
      q^{-1}
      \ketbra{4}{4}
      \\
      (c) &
      \pi ( q^{{E^3}_3} )
      \eq
      q^{\alpha}
      \ketbra{1}{1}
      +
      q^{\alpha+1}
      \ketbra{2}{2}
      +
      q^{\alpha+1}
      \ketbra{3}{3}
      +
      q^{\alpha+2}
      \ketbra{4}{4}
      \\
      (d) &
      \pi \( {E^1}_3 \)
      \eq
      -
      {\[\alpha+1\]}_q^{1/2}
      \ketbra{2}{4}
      +
      q^{-1}
      {\[\alpha\]}_q^{1/2} \ketbra{1}{3}
      \\
      (e) &
      \pi \( {E^3}_1 \)
      \eq
      -
      {\[\alpha+1\]}_q^{1/2} \ketbra{4}{2}
      +
      q
      {\[\alpha\]}_q^{1/2}
      \ketbra{3}{1}.
    \end{array}
    \qquad
    \right\}
  }
  \eqlabel{TheRepresentationExtended}
\ese
The expressions presented in \eqref{TheRepresentation} may be
constructed using the following definition, where $X \in U_q[gl(2|1)]$:
\bne
  \pi (X)
  & \defeq &
  X \cdot \pi (I)
  =
  X
  \cdot
  (
    \ketbra{1}{1}
    +
    \ketbra{2}{2}
    +
    \ketbra{3}{3}
    +
    \ketbra{4}{4}
  )
  \nonumber
  \\
  \eq
  (X \cdot \ket{1}) \bra{1}
  +
  (X \cdot \ket{2}) \bra{2}
  +
  (X \cdot \ket{3}) \bra{3}
  +
  (X \cdot \ket{4}) \bra{4}.
  \eqlabel{ConstructionRule}
\ene
In \secref{RepresentationoftheCartanGenerators} and
\secref{RepresentationoftheNon-CartanGenerators}, we
evaluate the action of the $U_q[gl(2|1)]$ simple generators on the
basis vectors, then use \eqref{ConstructionRule} to construct the
explicit representations presented in \eqref{TheRepresentation}.

\pagebreak


%

\subsubsection{Representation of the Cartan Generators}
\seclabel{RepresentationoftheCartanGenerators}

Firstly we consider the Cartan generators.
Initially, we deduce the form of $ \pi ( q^{{E^1}_1} ) $:
\bne
  q^{{E^1}_1}
  \cdot
  \ket{1}
  & \stackrel{\eqref{ExplicitExpCartanActiononket1}}{=} &
  \ket{1}.
  \eqlabel{E11onket1}
  \\
  q^{{E^1}_1}
  \cdot
  \ket{2}
  & \stackrel{(\ref{eq:InitialKetDefinition}a)}{=} &
  \beta_2
  q^{{E^1}_1}
  {E^3}_2
  \cdot
  \ket{1}
  \stackrel{(\ref{eq:ExplicitInterchanges2}d)}{=}
  \beta_2
  {E^3}_2
  q^{{E^1}_1}
  \cdot
  \ket{1}
  \stackrel{\eqref{ExplicitExpCartanActiononket1}}{=}
  \beta_2
  {E^3}_2
  \cdot
  \ket{1}
  \nonumber
  \\
  & \stackrel{(\ref{eq:InitialKetDefinition}a)}{=} &
  \ket{2}.
  \eqlabel{E11onket2}
  \\
  q^{{E^1}_1}
  \cdot
  \ket{3}
  & \stackrel{(\ref{eq:InitialKetDefinition}b)}{=} &
  \beta_3
  q^{{E^1}_1}
  {E^3}_1
  \cdot
  \ket{1}
  \stackrel{\eqref{Botia23}}{=}
  -
  q^{1}
  \beta_3
  q^{{E^1}_1}
  {E^2}_1
  {E^3}_2
  \cdot
  \ket{1}
  \nonumber
  \\
  & \stackrel{(\ref{eq:ExplicitInterchanges2}c)}{=} &
  -
  \beta_3
  {E^2}_1
  q^{{E^1}_1}
  {E^3}_2
  \cdot
  \ket{1}
  \stackrel{(\ref{eq:ExplicitInterchanges2}d)}{=}
  -
  \beta_3
  {E^2}_1
  {E^3}_2
  q^{{E^1}_1}
  \cdot
  \ket{1}
  \nonumber
  \\
  & \stackrel{\eqref{ExplicitExpCartanActiononket1}}{=} &
  -
  \beta_3
  {E^2}_1
  {E^3}_2
  \cdot
  \ket{1}
  \stackrel{\eqref{Botia23}}{=}
  q^{-1}
  \beta_3
  {E^3}_1
  \cdot
  \ket{1}
  \stackrel{(\ref{eq:InitialKetDefinition}b)}{=}
  q^{-1}
  \ket{3}.
  \eqlabel{E11onket3}
  \\
  q^{{E^1}_1}
  \cdot
  \ket{4}
  & \stackrel{(\ref{eq:InitialKetDefinition}c)}{=} &
  \beta_4
  q^{{E^1}_1}
  {E^3}_2
  {E^3}_1
  \cdot
  \ket{1}
  \stackrel{(\ref{eq:ExplicitInterchanges2}d)}{=}
  \beta_4
  {E^3}_2
  q^{{E^1}_1}
  {E^3}_1
  \cdot
  \ket{1}
  \nonumber
  \\
  & \stackrel{\eqref{E11onket3}}{=} &
  q^{-1}
  \beta_4
  \beta_3^{-1}
  {E^3}_2
  \cdot
  \ket{3}
  \stackrel{(\ref{eq:InitialKetDefinition}b)}{=}
  q^{-1}
  \beta_4
  {E^3}_2
  {E^3}_1
  \cdot
  \ket{1}
  \stackrel{(\ref{eq:InitialKetDefinition}c)}{=}
  q^{-1}
  \ket{4}.
  \eqlabel{E11onket4}
\ene
Combining the results of \eqref{E11onket1} to \eqref{E11onket4},
we have, using \eqref{ConstructionRule}, in accordance with
(\ref{eq:TheRepresentationExtended}$a$):
\be
  \pi ( q^{{E^1}_1} )
  =
  \ketbra{1}{1}
  +
  \ketbra{2}{2}
  +
  q^{-1}
  \ketbra{3}{3}
  +
  q^{-1}
  \ketbra{4}{4},
\ee
equivalently, in accordance with (\ref{eq:TheRepresentation}$a$):
\be
  \pi \( {E^1}_1 \)
  =
  -
  \ketbra{3}{3}
  -
  \ketbra{4}{4}.
\ee
Let us pause to justify this statement. Starting with
(\ref{eq:TheRepresentationExtended}$a$), firstly, we have, as ${E^1}_1$
is a diagonal operator/matrix:
$
  \pi ( q^{{E^1}_1} )
  =
  q^{\pi \( {E^1}_1 \)}
$,
thus:
\be
  \pi ( q^{{E^1}_1} )
  =
  e^{\ln \( q \) \pi \( {E^1}_1 \)}
  =
  \sum_{i=0}^\infty
  \frac{
    {\[ \ln \( q \) \pi \( {E^1}_1 \) \]}^i
  }{
    i!
  }
  =
  \sum_{i=0}^\infty
  \frac{
    {\[ \ln \( q \) \]}^i {\[ \pi \( {E^1}_1 \) \]}^i
  }{
    i!
  },
\ee
by the definition of the exponential as a formal power series. Next,
we have, by inspection:
\be
  {\[ \pi ( {E^1}_1 ) \]}^i
  =
  \left\{
  \begin{array}{cl}
    - {\(-\)}^i \pi \( {E^1}_1 \) \qquad & \mathrm{if~} i > 0
    \\
    I & \mathrm{if~} i = 0.
  \end{array}
  \right.
\ee

\pagebreak

Thus, we have:
\be
  \! \! \! \! \!
  \pi ( q^{{E^1}_1} )
  \eq
  \sum_{i=1}^\infty
    \frac{
      {\[ \ln \( q \) \]}^i {\[ \pi \( {E^1}_1 \) \]}^i
    }{
      i!
    }
  +
  I
  =
  -
  \pi \( {E^1}_1 \)
  \sum_{i=1}^\infty
    \frac{
      {\[ \ln \( q \) \]}^i {\(-\)}^i
    }{
      i!
    }
  +
  I
  \\
  \eq
  -
  \pi \( {E^1}_1 \)
  \sum_{i=1}^\infty
    \frac{
      {\[ - \ln \( q \) \]}^i
    }{
      i!
    }
  +
  I
  =
  -
  \pi \( {E^1}_1 \)
  \[
    \sum_{i=0}^\infty
      \frac{
        {\[ {\ln \( q \)}^{-1} \]}^i
      }{
        i!
      }
    -
    1
  \]
  +
  I
  \\
  \eq
  -
  \pi \( {E^1}_1 \)
  \[ e^{{\ln \( q \)}^{-1}} - 1 \]
  +
  I
  =
  \( 1 - q^{-1} \)
  \pi \( {E^1}_1 \)
  +
  I
  \\
  \eq
  \( 1 - q^{-1} \)
  \(
    -
    \ketbra{3}{3}
    -
    \ketbra{4}{4}
  \)
  +
  \ketbra{1}{1}
  +
  \ketbra{2}{2}
  +
  \ketbra{3}{3}
  +
  \ketbra{4}{4}
  \\
  \eq
  \ketbra{1}{1}
  +
  \ketbra{2}{2}
  +
  q^{-1}
  \ketbra{3}{3}
  +
  q^{-1}
  \ketbra{4}{4}.
\ee
A similar expansion may be made in the opposite direction. Below,
we shall simply report results for the other Cartan generators without
further ado.

We repeat the process for the other Cartan generators.
For $ \pi ( q^{{E^2}_2} ) $:
\bne
  q^{{E^2}_2}
  \cdot
  \ket{1}
  & \stackrel{\eqref{ExplicitExpCartanActiononket1}}{=} &
  \ket{1}.
  \eqlabel{E22onket1}
  \\
  q^{{E^2}_2}
  \cdot
  \ket{2}
  & \stackrel{(\ref{eq:InitialKetDefinition}a)}{=} &
  \beta_2
  q^{{E^2}_2}
  {E^3}_2
  \cdot
  \ket{1}
  \stackrel{(\ref{eq:ExplicitInterchanges2}h)}{=}
  q^{-1}
  \beta_2
  {E^3}_2
  q^{{E^2}_2}
  \cdot
  \ket{1}
  \stackrel{\eqref{ExplicitExpCartanActiononket1}}{=}
  q^{-1}
  \beta_2
  {E^3}_2
  \cdot
  \ket{1}
  \nonumber
  \\
  & \stackrel{(\ref{eq:InitialKetDefinition}a)}{=} &
  q^{-1}
  \ket{2}.
  \eqlabel{E22onket2}
  \\
  q^{{E^2}_2}
  \cdot
  \ket{3}
  & \stackrel{(\ref{eq:InitialKetDefinition}b)}{=} &
  \beta_3
  q^{{E^2}_2}
  {E^3}_1
  \cdot
  \ket{1}
  \stackrel{\eqref{Botia23}}{=}
  -
  q^{1}
  \beta_3
  q^{{E^2}_2}
  {E^2}_1
  {E^3}_2
  \cdot
  \ket{1}
  \nonumber
  \\
  & \stackrel{(\ref{eq:ExplicitInterchanges2}g)}{=} &
  -
  q^2
  \beta_3
  {E^2}_1
  q^{{E^2}_2}
  {E^3}_2
  \cdot
  \ket{1}
  \stackrel{(\ref{eq:ExplicitInterchanges2}h)}{=}
  -
  q^{1}
  \beta_3
  {E^2}_1
  {E^3}_2
  q^{{E^2}_2}
  \cdot
  \ket{1}
  \nonumber
  \\
  & \stackrel{\eqref{ExplicitExpCartanActiononket1}}{=} &
  -
  q^{1}
  \beta_2
  {E^2}_1
  {E^3}_2
  \cdot
  \ket{1}
  \stackrel{\eqref{Botia23}}{=}
  \beta_3
  {E^3}_1
  \cdot
  \ket{1}
  \stackrel{(\ref{eq:InitialKetDefinition}b)}{=}
  \ket{3}.
  \eqlabel{E22onket3}
  \\
  q^{{E^2}_2}
  \cdot
  \ket{4}
  & \stackrel{(\ref{eq:InitialKetDefinition}c)}{=} &
  \beta_4
  q^{{E^2}_2}
  {E^3}_2
  {E^3}_1
  \cdot
  \ket{1}
  \stackrel{(\ref{eq:ExplicitInterchanges2}h)}{=}
  q^{-1}
  \beta_4
  {E^3}_2
  q^{{E^2}_2}
  {E^3}_1
  \cdot
  \ket{1}
  \nonumber
  \\
  & \stackrel{\eqref{E22onket3}}{=} &
  q^{-1}
  \beta_4
  \beta_3^{-1}
  {E^3}_2
  \cdot
  \ket{3}
  \stackrel{(\ref{eq:InitialKetDefinition}b)}{=}
  q^{-1}
  \beta_4
  {E^3}_2
  {E^3}_1
  \cdot
  \ket{1}
  \stackrel{(\ref{eq:InitialKetDefinition}c)}{=}
  q^{-1}
  \ket{4}.
  \eqlabel{E22onket4}
\ene
Combining the results of \eqref{E22onket1} to \eqref{E22onket4}, we
have, in accordance with (\ref{eq:TheRepresentationExtended}$b$):
\be
  \pi ( q^{{E^2}_2} )
  =
  \ketbra{1}{1}
  +
  q^{-1}
  \ketbra{2}{2}
  +
  \ketbra{3}{3}
  +
  q^{-1}
  \ketbra{4}{4},
\ee
thus, we have, in accordance with (\ref{eq:TheRepresentation}$b$):
\be
  \pi \( {E^2}_2 \)
  =
  -
  \ketbra{2}{2}
  -
  \ketbra{4}{4}.
\ee

\pagebreak

Lastly, for $ \pi ( q^{{E^3}_3} ) $:
\bne
  q^{{E^3}_3}
  \cdot
  \ket{1}
  & \stackrel{\eqref{ExplicitExpCartanActiononket1}}{=} &
  q^{\alpha}
  \ket{1}.
  \eqlabel{E33onket1}
  \\
  q^{{E^3}_3}
  \cdot
  \ket{2}
  & \stackrel{(\ref{eq:InitialKetDefinition}a)}{=} &
  \beta_2
  q^{{E^3}_3}
  {E^3}_2
  \cdot
  \ket{1}
  \stackrel{(\ref{eq:ExplicitInterchanges2}l)}{=}
  q^{1}
  \beta_2
  {E^3}_2
  q^{{E^3}_3}
  \cdot
  \ket{1}
  \stackrel{\eqref{ExplicitExpCartanActiononket1}}{=}
  q^{1}
  q^{\alpha}
  \beta_2
  {E^3}_2
  \cdot
  \ket{1}
  \nonumber
  \\
  & \stackrel{(\ref{eq:InitialKetDefinition}a)}{=} &
  q^{\alpha+1}
  \ket{2}.
  \eqlabel{E33onket2}
  \\
  q^{{E^3}_3}
  \cdot
  \ket{3}
  & \stackrel{(\ref{eq:InitialKetDefinition}b)}{=} &
  \beta_3
  q^{{E^3}_3}
  {E^3}_1
  \cdot
  \ket{1}
  \stackrel{\eqref{Botia23}}{=}
  -
  \beta_3
  q^{1}
  q^{{E^3}_3}
  {E^2}_1
  {E^3}_2
  \cdot
  \ket{1}
  \nonumber
  \\
  & \stackrel{(\ref{eq:ExplicitInterchanges2}k)}{=} &
  -
  q^{1}
  \beta_3
  {E^2}_1
  q^{{E^3}_3}
  {E^3}_2
  \cdot
  \ket{1}
  \stackrel{(\ref{eq:ExplicitInterchanges2}l)}{=}
  -
  q^{2}
  \beta_3
  {E^2}_1
  {E^3}_2
  q^{{E^3}_3}
  \cdot
  \ket{1}
  \nonumber
  \\
  & \stackrel{\eqref{ExplicitExpCartanActiononket1}}{=} &
  -
  q^{2}
  q^{\alpha}
  \beta_3
  {E^2}_1
  {E^3}_2
  \cdot
  \ket{1}
  \stackrel{\eqref{Botia23}}{=}
  q^{\alpha+2}
  q^{-1}
  \beta_3
  {E^3}_1
  \cdot
  \ket{1}
  \stackrel{(\ref{eq:InitialKetDefinition}b)}{=}
  q^{\alpha+1}
  \ket{3}.
  \eqlabel{E33onket3}
  \\
  q^{{E^3}_3}
  \cdot
  \ket{4}
  & \stackrel{(\ref{eq:InitialKetDefinition}c)}{=} &
  \beta_4
  q^{{E^3}_3}
  {E^3}_2
  {E^3}_1
  \cdot
  \ket{1}
  \stackrel{(\ref{eq:ExplicitInterchanges2}l)}{=}
  q^{1}
  \beta_4
  {E^3}_2
  q^{{E^3}_3}
  {E^3}_1
  \cdot
  \ket{1}
  \nonumber
  \\
  & \stackrel{\eqref{E33onket3}}{=} &
  q^{1}
  q^{\alpha+1}
  \beta_4
  \beta_3^{-1}
  {E^3}_2
  \cdot
  \ket{3}
  \stackrel{(\ref{eq:InitialKetDefinition}b)}{=}
  q^{\alpha+2}
  \beta_4
  {E^3}_2
  {E^3}_1
  \cdot
  \ket{1}
  \stackrel{(\ref{eq:InitialKetDefinition}c)}{=}
  q^{\alpha+2}
  \ket{4}.
  \eqlabel{E33onket4}
\ene
Combining the results of \eqref{E33onket1} to \eqref{E33onket4}, we
have, in accordance with (\ref{eq:TheRepresentationExtended}$c$):
\be
  \pi ( q^{{E^3}_3} )
  =
  q^{\alpha}
  \ketbra{1}{1}
  +
  q^{\alpha+1}
  \ketbra{2}{2}
  +
  q^{\alpha+1}
  \ketbra{3}{3}
  +
  q^{\alpha+2}
  \ketbra{4}{4},
\ee
thus, in accordance with (\ref{eq:TheRepresentation}$c$):
\be
  \pi \( {E^3}_3 \)
  =
  \alpha
  \ketbra{1}{1}
  +
  \( \alpha + 1 \)
  \ketbra{2}{2}
  +
  \( \alpha + 1 \)
  \ketbra{3}{3}
  +
  \( \alpha + 2 \)
  \ketbra{4}{4}.
\ee


\subsubsection{Representation of the Non-Cartan Generators}
\seclabel{RepresentationoftheNon-CartanGenerators}

Now, we deduce the forms of the simple raising and lowering
generators.  Firstly, let us deduce the form of $\pi\({E^1}_2\)$:
\bne
  {E^1}_2 \cdot \ket{1}
  & \stackrel{\eqref{Annihilationofket1byRaisingGens}}{=} &
  0.
  \eqlabel{E12onket1}
  \\
  {E^1}_2 \cdot \ket{2}
  & \stackrel{(\ref{eq:InitialKetDefinition}a)}{=} &
  \beta_2
  {E^1}_2
  {E^3}_2
  \cdot
  \ket{1}
  \stackrel{(\ref{eq:ExplicitInterchanges1}c)}{=}
  \beta_2
  {E^3}_2
  {E^1}_2
  \cdot
  \ket{1}
  \stackrel{\eqref{Annihilationofket1byRaisingGens}}{=}
  0.
  \eqlabel{E12onket2}
  \\
  {E^1}_2 \cdot \ket{3}
  & \stackrel{(\ref{eq:InitialKetDefinition}b)}{=} &
  \beta_3
  {E^1}_2
  {E^3}_1
  \cdot
  \ket{1}
  \stackrel{\eqref{Botia46}}{=}
  -
  q
  \beta_3
  {E^3}_2
  \cdot
  \ket{1}
  \stackrel{(\ref{eq:InitialKetDefinition}a)}{=}
  -
  q
  \beta_3
  \beta_2^{-1}
  \ket{2}
  \nonumber
  \\
  & \stackrel{(\ref{eq:NormalisationConstants}a,b)}{=} &
  -
  \ket{2}.
  \eqlabel{E12onket3}
  \\
  {E^1}_2
  \cdot
  \ket{4}
  & \stackrel{(\ref{eq:InitialKetDefinition}c)}{=} &
  \beta_4
  {E^1}_2
  {E^3}_2
  {E^3}_1
  \cdot
  \ket{1}
  \stackrel{(\ref{eq:ExplicitInterchanges1}c)}{=}
  \beta_4
  {E^3}_2
  {E^1}_2
  {E^3}_1
  \cdot
  \ket{1}
  \nonumber
  \\
  & \stackrel{\eqref{Botia46}}{=} &
  -
  q
  \beta_4
  {E^3}_2
  {E^3}_2
  \cdot
  \ket{1}
  \stackrel{\eqref{SquaresOfOddGeneratorsAreZero}}{=}
  0.
  \eqlabel{E12onket4}
\ene
Combining the results of \eqref{E12onket1} to \eqref{E12onket4}, we
have, in accordance with (\ref{eq:TheRepresentation}$d$):
\be
  \pi \( {E^1}_2 \)
  =
  -
  \ketbra{2}{3}.
\ee
As we know that the representation is unitary
\cite{GouldHibberdLinksZhang:96,GouldScheunert:95}, we have
\be
  \pi \( {E^2}_1 \)
  =
  {\pi \( {E^1}_2 \)}^\dagger,
\ee
where, in matrix notation, $ X^\dagger $ is the conjugate transpose of
$X$.  Thus, immediately, in accordance with
(\ref{eq:TheRepresentation}$e$):
\be
  \pi \( {E^2}_1 \)
  =
  -
  \ketbra{3}{2}.
\ee

Secondly, $\pi \( {E^2}_3 \)$ is more work:
\bne
  {E^2}_3 \cdot \ket{1}
  & \stackrel{\eqref{Annihilationofket1byRaisingGens}}{=} &
  0.
  \eqlabel{E23onket1}
  \\
  {E^2}_3 \cdot \ket{2}
  & \stackrel{(\ref{eq:InitialKetDefinition}a)}{=} &
  \beta_2
  {E^2}_3 {E^3}_2
  \cdot
  \ket{1}
  \stackrel{(\ref{eq:Botia32},\ref{eq:NormalisationConstants}a)}{=}
  {\[ \alpha \]}_q^{1/2}
  \ket{1}.
  \eqlabel{E23onket2}
  \\
  {E^2}_3 \cdot \ket{3}
  & \stackrel{(\ref{eq:InitialKetDefinition}b)}{=} &
  \beta_3
  {E^2}_3 {E^3}_1
  \cdot
  \ket{1}
  \stackrel{\eqref{Botia42}}{=}
  0.
  \eqlabel{E23onket3}
  \\
  {E^2}_3
  \cdot
  \ket{4}
  & \stackrel{(\ref{eq:InitialKetDefinition}c)}{=} &
  \beta_4
  {E^2}_3
  {E^3}_2
  {E^3}_1
  \cdot
  \ket{1}
  \nonumber
  \\
  & \stackrel{(\ref{eq:ExplicitInterchanges1}b,\ref{eq:Botia23})}{=} &
  -
  q
  \beta_4
  \(
    {\[ {E^2}_2 + {E^3}_3 \]}_q
    -
    {E^3}_2
    {E^2}_3
  \)
  {E^2}_1
  {E^3}_2
  \cdot
  \ket{1}
  \nonumber
  \\
  & \stackrel{(\ref{eq:Botia38},\ref{eq:ExplicitInterchanges1}d)}{=} &
  -
  q
  \beta_4
  \(
    {\[\alpha+1\]}_q
    {E^2}_1
    {E^3}_2
    \cdot
    \ket{1}
    -
    {E^3}_2
    {E^2}_1
    {E^2}_3
    {E^3}_2
    \cdot
    \ket{1}
  \)
  \nonumber
  \\
  & \stackrel{\eqref{Botia32}}{=} &
  -
  q
  \beta_4
  \(
    {\[\alpha+1\]}_q
    {E^2}_1
    {E^3}_2
    \cdot
    \ket{1}
    -
    {\[\alpha\]}_q
    {E^3}_2
    {E^2}_1
    \cdot
    \ket{1}
  \)
  \nonumber
  \\
  & \stackrel{\eqref{E21Annihilatesket1}}{=} &
  -
  q
  {\[\alpha+1\]}_q
  \beta_4
  {E^2}_1
  {E^3}_2
  \cdot
  \ket{1}
  \stackrel{\eqref{Botia23}}{=}
  {\[\alpha+1\]}_q
  \beta_4
  {E^3}_1
  \cdot
  \ket{1}
  \nonumber
  \\
  & \stackrel{(\ref{eq:InitialKetDefinition}b)}{=} &
  {\[\alpha+1\]}_q
  \beta_4
  \beta_3^{-1}
  \ket{3}
  \stackrel{(\ref{eq:NormalisationConstants}b)}{=}
  {\[\alpha+1\]}_q^{1/2}
  \ket{3}.
  \eqlabel{E23onket4}
\ene
Combining the results of \eqref{E23onket1} to \eqref{E23onket4}, we
have, in accordance with (\ref{eq:TheRepresentation}$f$):
\be
  \pi \( {E^2}_3 \)
  =
  {\[\alpha\]}_q^{1/2}
  \ketbra{1}{2}
  +
  {\[\alpha+1\]}_q^{1/2}
  \ketbra{3}{4}.
\ee
Again, from the unitarity of the representation,
$
  \pi \( {E^3}_2 \)
  =
  {\pi \( {E^2}_3 \)}^\dagger
$,
hence, in accordance with (\ref{eq:TheRepresentation}$g$):
\be
  \pi \( {E^3}_2 \)
  =
  {\[\alpha\]}_q^{1/2}
  \ketbra{2}{1}
  +
  {\[\alpha+1\]}_q^{1/2}
  \ketbra{4}{3}.
\ee

\pagebreak

We do not really need to  discover the forms of $\pi\({E^1}_3\)$ and
$\pi\({E^3}_1\)$; however we mention them here for completeness.
In accordance with (\ref{eq:TheRepresentationExtended}$d,e$), we have:
\be
  \pi\({E^1}_3\)
  & \stackrel{(\ref{eq:NonSimpleGenerators}a)}{=} &
  \pi\( {E^1}_2 {E^2}_3 - q^{-1} {E^2}_3 {E^1}_2 \)
  \\
  \eq
  \pi \( {E^1}_2 \)
  \pi \( {E^2}_3 \)
  -
  q^{-1}
  \pi \( {E^2}_3 \)
  \pi \( {E^1}_2 \)
  \\
  & \stackrel{(\ref{eq:TheRepresentation}d,f)}{=} &
  -
  \ketbra{2}{3}
  \cdot
  \(
    {\[\alpha\]}_q^{1/2} \ketbra{1}{2}
    +
    {\[\alpha+1\]}_q^{1/2} \ketbra{3}{4}
  \)
  \\
  & &
  \qquad \qquad
  +
  q^{-1}
  \(
    {\[\alpha\]}_q^{1/2} \ketbra{1}{2}
    +
    {\[\alpha+1\]}_q^{1/2} \ketbra{3}{4}
  \)
  \cdot
  \ketbra{2}{3}
  \\
  \eq
  -
  {\[\alpha+1\]}_q^{1/2}
  \ketbra{2}{4}
  +
  q^{-1}
  {\[\alpha\]}_q^{1/2} \ketbra{1}{3}.
  \\
  \pi \( {E^3}_1 \)
  & \stackrel{(\ref{eq:NonSimpleGenerators}b)}{=} &
  \pi \( {E^3}_2 \)
  \pi \( {E^2}_1 \)
  -
  q
  \pi \( {E^2}_1 \)
  \pi \( {E^3}_2 \)
  \\
  & \stackrel{(\ref{eq:TheRepresentation}e,g)}{=} &
  -
  \(
    {\[\alpha\]}_q^{1/2} \ketbra{2}{1}
    +
    {\[\alpha+1\]}_q^{1/2} \ketbra{4}{3}
  \)
  \cdot
  \ketbra{3}{2}
  \\
  & &
  \qquad \qquad
  +
  q
  \ketbra{3}{2}
  \cdot
  \(
    {\[\alpha\]}_q^{1/2} \ketbra{2}{1}
    +
    {\[\alpha+1\]}_q^{1/2} \ketbra{4}{3}
  \)
  \\
  \eq
  -
  {\[\alpha+1\]}_q^{1/2} \ketbra{4}{2}
  +
  q
  {\[\alpha\]}_q^{1/2}
  \ketbra{3}{1}.
\ee

\pagebreak


\subsection{$U_q[gl(2|1)]$ as a Hopf Superalgebra}
\addtocontents{toc}{\protect\vspace{-2.5ex}}
\seclabel{Uqgl21asaHopfSuperalgebra}

$U_q\[gl(2|1)\]$ may be regarded as a Hopf superalgebra, when equipped
with suitable coproduct $\Delta$, counit $\varepsilon$ and
antipode $S$ structures \cite[p1971]{Zhang:92}. Here, we only describe
the coproduct, as that is all that we shall require below. Details for
the counit and antipode may be found in \secref{Uqglmn}.

\subsubsection{Coproduct $\Delta$}

We define a coproduct (a.k.a. comultiplication) structure, which is a
$\mathbb{Z}_2$ graded algebra homomorphism
$
  \Delta : U_q\[gl(2|1)\] \to U_q\[gl(2|1)\] \otimes U_q\[gl(2|1)\]
$,
by:
\bse
  \lefteqn{
    \hspace{-10mm}
    \left.
    \begin{array}{lrcl}
      (a) &
      \Delta ( {E^a}_a)
      \eq
      I \otimes {E^a}_a + {E^a}_a \otimes I,
      \qquad \qquad \qquad \qquad \qquad
      a = 1, 2, 3
      \\
      (b) &
      \Delta ( {E^1}_2 )
      \eq
      {E^1}_2 \otimes q^{-\frac{1}{2} ( {E^1}_1 - {E^2}_2 )}
      +
      q^{\frac{1}{2} ( {E^1}_1 - {E^2}_2 )} \otimes {E^1}_2
      \\
      (c) &
      \Delta ( {E^2}_1 )
      \eq
      {E^2}_1 \otimes q^{-\frac{1}{2} ( {E^1}_1 - {E^2}_2)}
      +
       q^{\frac{1}{2} ( {E^1}_1 - {E^2}_2 )} \otimes {E^2}_1
      \\
      (d) &
      \Delta ( {E^2}_3 )
      \eq
      {E^2}_3 \otimes q^{-\frac{1}{2} ( {E^2}_2 + {E^3}_3 )}
      +
      q^{\frac{1}{2}( {E^2}_2 + {E^3}_3 )} \otimes {E^2}_3
      \\
      (e) &
      \Delta ( {E^3}_2 )
      \eq
      {E^3}_2 \otimes q^{-\frac{1}{2}( {E^2}_2 + {E^3}_3 )}
      +
      q^{\frac{1}{2}( {E^2}_2 + {E^3}_3 )} \otimes {E^3}_2,
    \end{array}
    \qquad
    \right\}
  }
  \eqlabel{Uqgl21Coproduct}
\ese
and extended to an algebra homomorphism on all of $U_q\[gl(2|1)\]$.  As
$\Delta$ is a homomorphism, we have of course that
$\Delta(1)=1\otimes1$, and that for
$x,y \in U_q\[gl(2|1)\]$:
\bse
  \Delta(x y)
  =
  \Delta(x)
  \Delta(y).
  \eqlabel{CoproductGradedHomo}
\ese
$\Delta$ being graded means that it preserves grading, i.e.
$\[\Delta(x)\] = \[x\]$ for homogeneous $x \in U_q\[gl(2|1)\]$.
Of course, (\ref{eq:Uqgl21Coproduct}$a$) implies that:
\be
  \Delta ( q^{{E^a}_a} )
  \eq
  q^{{E^a}_a} \otimes q^{{E^a}_a}.
\ee


\subsubsection{The Twist Map $T$ and
               the Alternative Coproduct $\overline{\Delta}$}

As well as this standard coproduct, there exists another possible
coproduct structure:
$ \overline{\Delta} $, defined by
$ \overline{\Delta} = T \cdot \Delta $, where
$
  T
  :
  U_q[gl(2|1)] \otimes U_q[gl(2|1)]
  \to
  U_q[gl(2|1)] \otimes U_q[gl(2|1)]
$
is the twist map, defined for homogeneous elements
$ x, y \in U_q[gl(2|1)] $:
\be
  T \( x \otimes y \)
  =
  {\( - \)}^{\[ x \] \[ y \]}
  \( y \otimes x \).
\ee
We shall require this operator below, in
\secref{TheQuantumYang--BaxterEquation}. More generally, the notion of
a twist map is applicable to the tensor product of \emph{any} two
vector spaces.

\pagebreak


\subsubsection{$\Delta({E^3}_1)$ and $\Delta({E^1}_3)$}

In \secref{ABasisforV_1} and \secref{ABasisforV_3}, we shall
require the services of $\Delta({E^3}_1)$ and $\Delta({E^1}_3)$
respectively.
\bne
  \hspace{-15pt}
  \Delta({E^3}_1)
  & \stackrel{(\ref{eq:NonSimpleGenerators}b)}{=} &
  \Delta
  (
    {E^3}_2
    {E^2}_1
    -
    q
    {E^2}_1
    {E^3}_2
  )
  \nonumber
  \\
  & \stackrel{\eqref{CoproductGradedHomo}}{=} &
  \Delta({E^3}_2)
  \Delta({E^2}_1)
  -
  q
  \Delta({E^2}_1)
  \Delta({E^3}_2)
  \nonumber
  \\
  & \stackrel{(\ref{eq:Uqgl21Coproduct}c,e)}{=} &
  \(
    {E^3}_2 \otimes q^{-\frac{1}{2}( {E^2}_2 + {E^3}_3 )}
    +
    q^{\frac{1}{2}( {E^2}_2 + {E^3}_3 )} \otimes {E^3}_2
  \)
  \cdot
  \nonumber
  \\
  & & \qquad \qquad
  \(
    {E^2}_1 \otimes q^{-\frac{1}{2} ( {E^1}_1 - {E^2}_2)}
    +
     q^{\frac{1}{2} ( {E^1}_1 - {E^2}_2 )} \otimes {E^2}_1
  \)
  \nonumber
  \\
  & & -
  q
  \(
    {E^2}_1 \otimes q^{-\frac{1}{2} ( {E^1}_1 - {E^2}_2)}
    +
     q^{\frac{1}{2} ( {E^1}_1 - {E^2}_2 )} \otimes {E^2}_1
  \)
  \cdot
  \nonumber
  \\
  & & \qquad \qquad
  \(
    {E^3}_2 \otimes q^{-\frac{1}{2}( {E^2}_2 + {E^3}_3 )}
    +
    q^{\frac{1}{2}( {E^2}_2 + {E^3}_3 )} \otimes {E^3}_2
  \)
  \nonumber
  \\
  & \stackrel{(\ref{eq:ExplicitInterchanges2}d,g,h,k)}{=} &
  \(
    {E^3}_2 {E^2}_1
    -
    q^{1}
    {E^2}_1 {E^3}_2
  \)
  \otimes
  q^{-\frac{1}{2} ( {E^1}_1 + {E^3}_3 )}
  \nonumber
  \\
  & &
  +
  q^{\frac{1}{2} ( {E^1}_1 + {E^3}_3 )}
  \otimes
  \(
    {E^3}_2 {E^2}_1
    -
    q^{1}
    {E^2}_1 {E^3}_2
  \)
  \nonumber
  \\
  & &
  -
  \( q - q^{-1} \)
  q^{\frac{1}{2}( {E^1}_1 - {E^2}_2 )} {E^3}_2
  \otimes
  {E^2}_1 q^{-\frac{1}{2} ( {E^2}_2 + {E^3}_3)}
  \eqlabel{DeltaE31}
  \\
  & \stackrel{(\ref{eq:NonSimpleGenerators}b)}{=} &
  {E^3}_1
  \otimes
  q^{-\frac{1}{2} ( {E^1}_1 + {E^3}_3 )}
  +
  q^{\frac{1}{2} ( {E^1}_1 + {E^3}_3 )}
  \otimes
  {E^3}_1
  \nonumber
  \\
  & &
  -
  \( q - q^{-1} \)
  q^{\frac{1}{2}( {E^1}_1 - {E^2}_2 )} {E^3}_2
  \otimes
  {E^2}_1 q^{-\frac{1}{2} ( {E^2}_2 + {E^3}_3)}.
  \nonumber
\ene
\bne
  \hspace{-29pt}
  \Delta({E^1}_3)
  & \stackrel{(\ref{eq:NonSimpleGenerators}a)}{=} &
  \Delta
  (
    {E^1}_2
    {E^2}_3
    -
    q^{-1}
    {E^2}_3
    {E^1}_2
  )
  \nonumber
  \\
  & \stackrel{\eqref{CoproductGradedHomo}}{=} &
  \Delta({E^1}_2)
  \Delta({E^2}_3)
  -
  q^{-1}
  \Delta({E^2}_3)
  \Delta({E^1}_2)
  \nonumber
  \\
  & \stackrel{(\ref{eq:Uqgl21Coproduct}b,d)}{=} &
  \(
    {E^1}_2 \otimes q^{-\frac{1}{2} ( {E^1}_1 - {E^2}_2 )}
    +
    q^{\frac{1}{2} ( {E^1}_1 - {E^2}_2 )} \otimes {E^1}_2
  \)
  \cdot
  \nonumber
  \\
  & & \qquad \qquad
  \(
    {E^2}_3 \otimes q^{-\frac{1}{2} ( {E^2}_2 + {E^3}_3 )}
    +
    q^{\frac{1}{2}( {E^2}_2 + {E^3}_3 )} \otimes {E^2}_3
  \)
  \nonumber
  \\
  & &
  -
  q^{-1}
  \(
    {E^2}_3 \otimes q^{-\frac{1}{2} ( {E^2}_2 + {E^3}_3 )}
    +
    q^{\frac{1}{2}( {E^2}_2 + {E^3}_3 )} \otimes {E^2}_3
  \)
  \cdot
  \nonumber
  \\
  & & \qquad \qquad
  \(
    {E^1}_2 \otimes q^{-\frac{1}{2} ( {E^1}_1 - {E^2}_2 )}
    +
    q^{\frac{1}{2} ( {E^1}_1 - {E^2}_2 )} \otimes {E^1}_2
  \)
  \nonumber
  \\
  & \stackrel{(\ref{eq:ExplicitInterchanges2}b,e,f,i)}{=} &
  \(
    {E^1}_2 {E^2}_3
    -
    q^{-1}  {E^2}_3 {E^1}_2
  \)
  \otimes
  q^{-\frac{1}{2} ( {E^1}_1 + {E^3}_3)}
  \nonumber
  \\
  & &
  +
  q^{\frac{1}{2}( {E^1}_1 + {E^3}_3 )}
  \otimes
  \(
    {E^1}_2 {E^2}_3
    -
    q^{-1}
    {E^2}_3 {E^1}_2
  \)
  \nonumber
  \\
  & &
  +
  \( q - q^{-1} \)
  q^{\frac{1}{2} ( {E^2}_2 + {E^3}_3 )} {E^1}_2
  \otimes
  {E^2}_3 q^{-\frac{1}{2}( {E^1}_1 - {E^2}_2 )}
  \eqlabel{DeltaE13}
  \\
  & \stackrel{(\ref{eq:NonSimpleGenerators}a)}{=} &
  {E^1}_3
  \otimes
  q^{-\frac{1}{2} ( {E^1}_1 + {E^3}_3)}
  +
  q^{\frac{1}{2}( {E^1}_1 + {E^3}_3 )}
  \otimes
  {E^1}_3
  \nonumber
  \\
  & &
  +
  \( q - q^{-1} \)
  q^{\frac{1}{2} ( {E^2}_2 + {E^3}_3 )} {E^1}_2
  \otimes
  {E^2}_3 q^{-\frac{1}{2}( {E^1}_1 - {E^2}_2 )}.
  \nonumber
\ene

\pagebreak


\subsection{The Tensor Product Module $V \otimes V$}
\addtocontents{toc}{\protect\vspace{-2.5ex}}

From the $4$ dimensional representation on the space $V$, we may
construct a $16$ dimensional representation on the tensor product space
$V\otimes V$. The grading on the basis for $V$ induces a
natural grading on the basis for $V\otimes V$:
\be
  \[\ket{i}\otimes\ket{j}\]
  \defeq
  \[\ket{i}\]
  +
  \[\ket{j}\].
\ee
The $\mathbb{Z}_2$ graded action of
$ U_q\[gl(2|1)\] \otimes U_q[gl(2|1)] $
on $ V \otimes V $
is defined for generators $ x, y \in U_q[gl(2|1)]$ and vectors
$u,v\in V$ in the following manner:
\bse
  \( x \otimes y \)
  \cdot
  \( u \otimes v \)
  \defeq
  {\(-\)}^{\[y\]\[u\]}
  \( x \cdot u \otimes y \cdot v \),
  \qquad
  \mathrm{homogeneous~}
  y, u,
  \eqlabel{Botia105}
\ese
and naturally extended by linearity.


\subsubsection{Decomposition of the Tensor Product Module}

The tensor product module $ V \otimes V $ has the following
decomposition with respect to the coproduct for generic values of
$\alpha$:
\bse
  V \otimes V
  =
  V_1 \oplus V_2 \oplus V_3,
  \eqlabel{TPdecomp}
\ese
where:
\be
  \begin{array}{
         ccl@{}r@{\hspace{1mm}}r@{\hspace{1mm}}c@{\hspace{1mm}}l@{}rcl
               }
    V_1 & \textrm{has~highest~weight} & ( & 0, & 0 & | & 2\alpha   & ) &
      \textrm{and~dimension} & 4 \\
    V_2 &                             & ( & 0,& -1 & | & 2\alpha+1 & ) &
                             & 8 \\
    V_3 &                             & ( & -1,& -1 & | & 2\alpha+2& ) &
                             & 4.
  \end{array}
\ee
This decomposition is found in \cite[p2770]{BrackenGouldLinksZhang:95};
which, strictly speaking, refers to the nonquantised case $gl\(2|1\)$,
but which extends naturally to the quantised case $U_q[gl(2|1)]$,
cf. \cite[p191]{LinksGould:92b}.

\pagebreak

We shall construct symmetry adapted orthonormal bases
$\{\ket{\Psi^1_j}\}_{j=1}^4$ and $\{\ket{\Psi^3_j}\}_{j=1}^4$, for the
spaces $V_1$ and $V_3$ respectively, in terms of the basis elements
$\{\ket{i}\}_{i=1}^4$ of $V$.  In the construction process, we shall
introduce arbitrary phase factors; these will be mentioned in passing.
The bases are:

\bse
  \lefteqn{
    \hspace{-13mm}
    \left.
    \begin{array}{lrcl}
      (a) &
      \ket{\Psi^1_1}
      \eq
      \ket{1} \otimes \ket{1}
      \\
      (b) &
      \ket{\Psi^1_2}
      \eq
      {(q^{\alpha}+q^{-\alpha})}^{-\frac{1}{2}}
      \(
        q^{\frac{\alpha}{2}} \ket{1} \otimes \ket{2}
        +
        q^{-\frac{\alpha}{2}} \ket{2} \otimes \ket{1}
      \)
      \\
      (c) &
      \ket{\Psi^1_3}
      \eq
      (q^{\alpha}+q^{-\alpha})^{-\frac{1}{2}}
      \(
        q^{\frac{\alpha}{2}} \ket{1} \otimes \ket{3}
        +
        q^{-\frac{\alpha}{2}} \ket{3} \otimes \ket{1}
      \)
      \\
      (d) &
      \ket{\Psi^1_4}
      \eq
      {\( q^{\alpha} + q^{-\alpha} \)}^{-\frac{1}{2}}
      {[2\alpha+1]}_q^{-\frac{1}{2}}
      \cdot
      \nonumber
      \\
      & & &
      \hspace{-18mm}
      \[
        {[\alpha+1]}_q^{\frac{1}{2}}
        (
          q^{\alpha} \ket{1} \otimes \ket{4}
          +
          q^{-\alpha} \ket{4} \otimes \ket{1}
        )
        -
        {[\alpha]}_q^{\frac{1}{2}}
        (
          q^{\frac{1}{2}} \ket{3} \otimes \ket{2}
          -
          q^{-\frac{1}{2}} \ket{2} \otimes \ket{3}
        )
      \]
    \end{array}
    \right\}
  }
  \eqlabel{V1KetBasis}
\ese

\bse
  \lefteqn{
    \hspace{-13mm}
    \left.
    \begin{array}{lrcl}
      (a) &
      \ket{\Psi^3_1}
      \eq
      {(q^{\alpha+1}+q^{-\alpha-1})}^{-\frac{1}{2}}
      {[2\alpha+1]}_q^{-\frac{1}{2}}
      \cdot
      \nonumber
      \\
      & & &
      \hspace{-25mm}
      \[
        {[\alpha]}_q^{\frac{1}{2}}
        (
          q^{\alpha+1} \ket{4} \otimes \ket{1}
          +
          q^{-\alpha-1} \ket{1} \otimes \ket{4}
        )
        +
        [\alpha+1]_q^{\frac{1}{2}}
        (
          q^{\frac{1}{2}} \ket{3} \otimes \ket{2}
          -
          q^{-\frac{1}{2}} \ket{2} \otimes \ket{3}
        )
      \]
      \\
      (b) &
      \ket{\Psi^3_2}
      \eq
      {(q^{\alpha+1}+q^{-\alpha-1})}^{-\frac{1}{2}}
      \(
        q^{\frac{\alpha+1}{2}} \ket{4} \otimes \ket{2}
        +
        q^{-\frac{\alpha+1}{2}} \ket{2} \otimes \ket{4}
      \)
      \\
      (c) &
      \ket{\Psi^3_3}
      \eq
      {(q^{\alpha+1}+q^{-\alpha-1})}^{-\frac{1}{2}}
      \(
        q^{\frac{\alpha+1}{2}} \ket{4} \otimes \ket{3}
        +
        q^{-\frac{\alpha+1}{2}} \ket{3} \otimes \ket{4}
      \)
      \\
      (d) &
       \ket{\Psi^3_4}
      \eq
       \ket{4} \otimes \ket{4}.
    \end{array}
    \right\}
  }
  \eqlabel{V3KetBasis}
\ese
Observe that these submodules do not have purely odd or even bases,
but that each basis vector is homogeneous.

Our basis is similar to those presented in
\cite{GouldHibberdLinksZhang:96} and \cite{DeWitKauffmanLinks:98},
except that we have made some different choices in the definitions of
the $\ket{i}$. Of those two papers, the basis presented in the former
contains a slight error as well. In particular, to get from the basis
in \cite{DeWitKauffmanLinks:98} to the above, we must interchange the
definition of $\ket{2}$ with $\ket{3}$ (and also $\ket{\Psi^1_2}$ with
$\ket{\Psi^1_3}$ and $\ket{\Psi^3_2}$ with $\ket{\Psi^3_3}$), and
replace $\ket{4}$ with $-\ket{4}$. These changes naturally lead to
some different phases in the respective normalisation constants.

Bases dual to
$\{\ket{\Psi^k_j}\}_{j=1}^4$ (for $k=1,3$), labelled
$\{\bra{\Psi^k_j}\}_{j=1}^4$ respectively, are found from the
definitions:
\bne
  \bra{\Psi^k_j}
  \eq
  \ket{\Psi^k_j}^{\dagger},
  \qquad \qquad \qquad \qquad
  \;
  k = 1, 3,
  \qquad
  j = 1, \dots, 4,
  \eqlabel{dualrule1}
  \\
  {( \ket{i} \otimes \ket{j} )}^{\dagger}
  \eq
  {\( - \)}^{\[ \ket{i} \] \[ \ket{j} \]}
  \( \bra{i} \otimes \bra{j} \),
  \qquad
  i, j = 1, \dots, 4.
  \eqlabel{dualrule2}
\ene

\pagebreak

Observe that the general form of each of the basis vectors
$\ket{\Psi^k_j}$ is:
\be
  \ket{\Psi^k_j}
  =
  \sum_{m}
    {\theta^{k j}}_{m}
    \( \ket{{x^{k j}}_{m}} \otimes \ket{{y^{k j}}_{m}} \),
\ee
where the coefficients ${\theta^{k j}}_{m}$ are in general complex
scalar functions of $q$ and $\alpha$. To simplify the discussion, we
shall limit $q$ and $\alpha$ to be real and positive, which ensures
that the ${\theta^{k j}}_{m}$ are real, and hence equal to their own
complex conjugates. In that case, from \eqref{dualrule1} and
\eqref{dualrule2}, the duals of these vectors $\ket{\Psi^k_j}$ are
given by:
\be
  \bra{\Psi^k_j}
  =
  \sum_{m}
    {\( - \)}^{ \[ \ket{{x^{k j}}_{m}} \] \[ \ket{{y^{k j}}_{m}} \]}
    {\theta^{k j}}_{m}
    \( \bra{{x^{k j}}_{m}} \otimes \bra{{y^{k j}}_{m}} \).
\ee

Explicitly, we have:
\bse
  \lefteqn{
    \hspace{-13mm}
    \left.
    \begin{array}{lrcl}
      (a) &
      \bra{\Psi^1_1}
      \eq
      \bra{1} \otimes \bra{1}
      \\
      (b) &
      \bra{\Psi^1_2}
      \eq
      (q^{\alpha}+q^{-\alpha})^{-\frac{1}{2}}
      \(
        q^{\frac{\alpha}{2}} \bra{1} \otimes \bra{2}
        +
        q^{-\frac{\alpha}{2}} \bra{2} \otimes \bra{1}
      \)
      \\
      (c) &
      \bra{\Psi^1_3}
      \eq
      (q^{\alpha}+q^{-\alpha})^{-\frac{1}{2}}
      \(
        q^{\frac{\alpha}{2}} \bra{1} \otimes \bra{3}
        +
        q^{-\frac{\alpha}{2}} \bra{3} \otimes \bra{1}
      \)
      \\
      (d) &
      \bra{\Psi^1_4}
      \eq
      {\( q^{\alpha} + q^{-\alpha} \)}^{-\frac{1}{2}}
      [2\alpha+1]_q^{-\frac{1}{2}}
      \cdot
      \nonumber
      \\
      & & &
      \hspace{-18mm}
      \[
        [\alpha+1]_q^{\frac{1}{2}}
        (
          q^{\alpha} \bra{1} \otimes \bra{4}
          +
          q^{-\alpha} \bra{4} \otimes \bra{1}
        )
        +
        [\alpha]_q^{\frac{1}{2}}
        (
          q^{\frac{1}{2}} \bra{3} \otimes \bra{2}
          -
          q^{-\frac{1}{2}} \bra{2} \otimes \bra{3}
        )
      \]
    \end{array}
    \right\}
  }
  \eqlabel{V1BraBasis}
\ese

\bse
  \lefteqn{
    \hspace{-13mm}
    \left.
    \begin{array}{lrcl}
      (a) &
      \bra{\Psi^3_1}
      \eq
      (q^{\alpha+1}+q^{-\alpha-1})^{-\frac{1}{2}}
      [2\alpha+1]_q^{-\frac{1}{2}}
      \cdot
      \nonumber
      \\
      & & &
      \hspace{-25mm}
      \[
        [\alpha]_q^{\frac{1}{2}}
        (
          q^{\alpha+1} \bra{4} \otimes \bra{1}
          +
          q^{-\alpha-1} \bra{1} \otimes \bra{4}
        )
        -
        [\alpha+1]_q^{\frac{1}{2}}
        (
          q^{\frac{1}{2}} \bra{3} \otimes \bra{2}
          -
          q^{-\frac{1}{2}} \bra{2} \otimes \bra{3}
        )
      \]
      \\
      (b) &
      \bra{\Psi^3_2}
      \eq
      {\( q^{\alpha+1}+q^{-\alpha-1} \)}^{-\frac{1}{2}}
      \(
        q^{\frac{\alpha+1}{2}} \bra{4} \otimes \bra{2}
        +
        q^{-\frac{\alpha+1}{2}} \bra{2} \otimes \bra{4}
      \)
      \\
      (c) &
      \bra{\Psi^3_3}
      \eq
      (q^{\alpha+1}+q^{-\alpha-1})^{-\frac{1}{2}}
      \(
        q^{\frac{\alpha+1}{2}} \bra{4} \otimes \bra{3}
        +
        q^{-\frac{\alpha+1}{2}} \bra{3} \otimes \bra{4}
      \)
      \\
      (d) &
      \bra{\Psi^3_4}
      \eq
      \bra{4} \otimes \bra{4}.
    \end{array}
    \right\}
  }
  \eqlabel{V3BraBasis}
\ese

The multiplication operations between the dual bases for $V \otimes V$
are given by:
\bse
  \hspace{-12mm}
  \left.
  \begin{array}{lrcl}
    (a) &
    \( \ket{i} \otimes \ket{j} \) \( \bra{k} \otimes \bra{l} \)
    \eq
    {\( - \)}^{\[ \ket{j} \] \[ \bra{k} \]}
    \( \ketbra{i}{k} \otimes \ketbra{j}{l} \)
    \\
    (b) &
    \( \bra{i} \otimes \bra{j} \) \( \ket{k} \otimes \ket{l} \)
    \eq
    {\( - \)}^{\[ \bra{j} \] \[ \ket{k} \]}
    \delta^i_k \delta^j_l
  \end{array}
  \right\}
  \hspace{2mm}
  i, j, k, l = 1, 2, 3, 4.
  \eqlabel{DualBasesMultRules}
\ese

So, where do the bases $\ket{\Psi^k_j}$ come from, and why aren't there
basis vectors for $V_2$? These things are explained in
\S\ref{sec:ABasisforV_1} to \S\ref{sec:NoBasisforV_2}.


\subsubsection{An Orthonormal Basis for $V_1$}
\seclabel{ABasisforV_1}

The highest weight vector of $V$ is $\ket{1}$, of weight
$\(0,0\,|\,\alpha\)$, hence the highest weight vector of $V \otimes V$
is $\ket{1} \otimes \ket{1}$, of weight $\(0,0\,|\,2\alpha\)$.  We
shall label by $V_1$ the submodule of $V \otimes V$ generated by this
vector $\ket{1} \otimes \ket{1}$, i.e.  $V_1$ has highest weight
$\(0,0\,|\,2\alpha\)$, from which we may determine, using Kac's
dimension formula, that it has dimension $4$. We shall construct an
orthonormal basis $\{\ket{\Psi^1_j}\}_{j=1}^4$ for $V_1$, using the Kac
induced module construction. Initially, as already defined in
(\ref{eq:V1KetBasis}$a$):
\be
  \ket{\Psi^1_1}
  \defeq
  \ket{1}
  \otimes
  \ket{1}.
\ee
The normality of $\ket{1}$ implies the normality of
$\ket{1}\otimes\ket{1}$:
\be
  \hspace{-2pt}
  \(\bra{1}\otimes\bra{1}\)
  \cdot
  \(\ket{1}\otimes\ket{1}\)
  \stackrel{(\ref{eq:DualBasesMultRules}b)}{=}
  {\(-\)}^{\[\bra{1}\]\[\ket{1}\]}
  \( \braket{1}{1}\otimes\braket{1}{1} \)
  \stackrel{\eqref{Braketii=1}}{=}
  {\(-\)}^{0 \cdot 0}
  1
  =
  1,
\ee
For
$V_1$, the equivalent of \eqref{E21Annihilatesket1} for $V\otimes V$ is
$\Delta({E^2}_1)\cdot\ket{\Psi^1_1}=0$. This fact will be used
implicitly.  The other basis vectors for $V_1$ are generated from the
action of (the coproduct of) the odd lowering generators of
$U_q[gl(2|1)]$, in much the same manner as the basis
$\{\ket{i}\}_{i=1}^4$ for $V$ itself was generated from $\ket{1}$, cf.
\eqref{InitialKetDefinition}:
\bse
  \left.
  \begin{array}{lrcl}
    (a) &
    \ket{\Psi^1_2}
    & \defeq &
    \beta^1_2
    \Delta({E^3}_2)
    \cdot
    \(\ket{1} \otimes \ket{1}\)
    \\
    (b) &
    \ket{\Psi^1_3}
    & \defeq &
    \beta^1_3
    \Delta({E^3}_1)
    \cdot
    \(\ket{1} \otimes \ket{1}\)
    \\
    (c) &
    \ket{\Psi^1_4}
    & \defeq &
    \beta^1_4
    \Delta({E^3}_2 {E^3}_1)
    \cdot
    \(\ket{1} \otimes \ket{1}\),
  \end{array}
  \qquad \qquad \qquad \qquad \qquad \; \;
  \right\}
  \eqlabel{InitialKetPsi1Definition}
\ese
where the $ \beta^1_j $, for $ j = 2,3,4 $, are suitable normalisation
constants.  As the coproduct is a homomorphism
\eqref{CoproductGradedHomo}, we have
$\Delta\({E^3}_2{E^3}_1\)=\Delta\({E^3}_2\)\Delta\({E^3}_1\)$.  Again,
the basis vectors are ordered in terms of decreasing weight:
\be
  \begin{array}{
       ccl@{}r@{\hspace{1mm}}r@{\hspace{1mm}}c@{\hspace{1mm}}l@{}l@{}c
               }
    \ket{\Psi^1_1} & \textrm{has~weight} &
      ( &  0, &  0 & | & 2 \alpha   & ) & \\
    \ket{\Psi^1_2} &                     &
      ( &  0, & -1 & | & 2 \alpha+1 & ) & \\
    \ket{\Psi^1_3} &                     &
      ( & -1, &  0 & | & 2 \alpha+1 & ) & \\
    \ket{\Psi^1_4} &                     &
      ( & -1, & -1 & | & 2 \alpha+2 & ) & .
  \end{array}
\ee
Without further ado, let us discover $\ket{\Psi^1_2}$. Firstly, we have:
\be
  \lefteqn{
    \Delta\({E^3}_2\)
    \cdot
    \(\ket{1} \otimes \ket{1}\)
  }
  \\
  & \stackrel{(\ref{eq:Uqgl21Coproduct}e)}{=} &
  {E^3}_2 \cdot \ket{1}
  \otimes
  q^{-\frac{1}{2}( {E^2}_2 + {E^3}_3 )} \cdot \ket{1}
  +
  q^{\frac{1}{2}( {E^2}_2 + {E^3}_3 )} \cdot \ket{1}
  \otimes
  {E^3}_2 \cdot \ket{1}
  \\
  & \stackrel{\eqref{ExplicitExpCartanActiononket1}}{=} &
  {E^3}_2 \cdot \ket{1}
  \otimes
  q^{-\frac{\alpha}{2}} \ket{1}
  +
  q^{\frac{\alpha}{2}} \ket{1}
  \otimes
  {E^3}_2 \cdot \ket{1}
  \\
  & \stackrel{(\ref{eq:TheRepresentation}g)}{=} &
  {[\alpha]}_q^{-1/2}
  \(
    q^{\frac{\alpha}{2}}
    \ket{1} \otimes \ket{2}
    +
    q^{-\frac{\alpha}{2}}
    \ket{2} \otimes \ket{1}
  \).
\ee
To deduce $ \beta^1_2 $, we demand
$\braket{\Psi^1_2}{\Psi^1_2}=1$. Investigating:
\be
  \hspace{-3pt}
  \[
    q^{\frac{\alpha}{2}}
    \bra{1} \otimes \bra{2}
    +
    q^{-\frac{\alpha}{2}}
    \bra{2} \otimes \bra{1}
  \]
  \cdot
  \[
    q^{\frac{\alpha}{2}}
    \ket{1} \otimes \ket{2}
    +
    q^{-\frac{\alpha}{2}}
    \ket{2} \otimes \ket{1}
  \]
  \stackrel{(\ref{eq:DualBasesMultRules}b)}{=}
  \( q^{\alpha} + q^{-\alpha} \),
\ee
hence we require, up to a phase factor
$
  \beta^1_2
  =
  {\( q^{\alpha} + q^{-\alpha} \)}^{-\frac{1}{2}}
$.
Thus, (\ref{eq:InitialKetPsi1Definition}$a$) yields,
in accordance with (\ref{eq:V1KetBasis}$b$):
\be
  \ket{\Psi^1_2}
  =
  {\( q^{\alpha} + q^{-\alpha} \)}^{-\frac{1}{2}}
  \[
    q^{\frac{\alpha}{2}}
    \ket{1} \otimes \ket{2}
    +
    q^{-\frac{\alpha}{2}}
    \ket{2} \otimes \ket{1}
  \].
\ee

To build $\ket{\Psi^1_3}$, we evaluate:
\bne
  \lefteqn{
    \Delta\({E^3}_1\)
    \cdot
    \(\ket{1} \otimes \ket{1}\)
  }
  \nonumber
  \\
  & \stackrel{\eqref{DeltaE31}}{=} &
  \(
    {E^3}_2 {E^2}_1
    -
    q^{1}
    {E^2}_1 {E^3}_2
  \)
  \cdot
  \ket{1}
  \otimes
  q^{-\frac{1}{2} ( {E^1}_1 + {E^3}_3 )}
  \cdot
  \ket{1}
  \nonumber
  \\
  & &
  +
  q^{\frac{1}{2} ( {E^1}_1 + {E^3}_3 )}
  \cdot
  \ket{1}
  \otimes
  \(
    {E^3}_2 {E^2}_1
    -
    q^{1}
    {E^2}_1 {E^3}_2
  \)
  \cdot
  \ket{1}
  \nonumber
  \\
  & &
  -
  \( q - q^{-1} \)
  q^{\frac{1}{2}( {E^1}_1 - {E^2}_2 )} {E^3}_2
  {E^3}_2
  \cdot
  \ket{1}
  \otimes
  {E^2}_1 q^{-\frac{1}{2} ( {E^2}_2 + {E^3}_3)}
  \cdot
  \ket{1}
  \nonumber
  \\
  & \stackrel{(\ref{eq:Botia111},\ref{eq:E21Annihilatesket1})}{=} &
  -
  q^{1}  {E^2}_1 {E^3}_2
  \cdot
  \ket{1}
  \otimes
  q^{-\frac{\alpha}{2}}
  \ket{1}
  -
  q^{1+\frac{\alpha}{2}}
  \ket{1}
  \otimes
  {E^2}_1 {E^3}_2
  \cdot
  \ket{1}
  \nonumber
  \\
  & &
  -
  \( q - q^{-1} \)
  q^{\frac{1}{2}( {E^1}_1 - {E^2}_2 )} {E^3}_2
  {E^3}_2
  \cdot
  \ket{1}
  \otimes
  {E^2}_1
  \cdot
  \ket{1}
  \nonumber
  \\
  & \stackrel{(\ref{eq:E21Annihilatesket1},
               \ref{eq:TheRepresentation}g)}{=} &
  -
  q^{1}
  {\[\alpha\]}_q^{1/2}
  \(
    q^{-\frac{\alpha}{2}}
    {E^2}_1
    \cdot
    \ket{2}
    \otimes
    \ket{1}
    +
    q^{\frac{\alpha}{2}}
    \ket{1}
    \otimes
    {E^2}_1
    \cdot
    \ket{2}
  \)
  \nonumber
  \\
  & \stackrel{(\ref{eq:TheRepresentation}e)}{=} &
  q
  {\[\alpha\]}_q^{\frac{1}{2}}
  \[
    q^{\frac{\alpha}{2}}
    \ket{1} \otimes \ket{3}
    +
    q^{-\frac{\alpha}{2}}
    \ket{3} \otimes \ket{1}
  \].
  \eqlabel{Botia55}
\ene

To deduce $ \beta^1_3 $, we demand
$\braket{\Psi^1_3}{\Psi^1_3}=1$. Investigating:
\be
  \hspace{-3pt}
  \[
    q^{\frac{\alpha}{2}}
    \bra{1} \otimes \bra{3}
    +
    q^{-\frac{\alpha}{2}}
    \bra{3} \otimes \bra{1}
  \]
  \cdot
  \[
    q^{\frac{\alpha}{2}}
    \ket{1} \otimes \ket{3}
    +
    q^{-\frac{\alpha}{2}}
    \ket{3} \otimes \ket{1}
  \]
  \stackrel{(\ref{eq:DualBasesMultRules}b)}{=}
  \( q^{\alpha} + q^{-\alpha} \),
\ee
hence we require, up to a phase factor
$
  \beta^1_3
  =
  {\( q^{\alpha} + q^{-\alpha} \)}^{-\frac{1}{2}}
$.
Thus, (\ref{eq:InitialKetPsi1Definition}$b$) yields,
in accordance with (\ref{eq:V1KetBasis}$c$):
\be
  \ket{\Psi^1_3}
  =
  {\( q^{\alpha} + q^{-\alpha} \)}^{-\frac{1}{2}}
  \[
    q^{\frac{\alpha}{2}}
    \ket{1} \otimes \ket{3}
    +
    q^{-\frac{\alpha}{2}}
    \ket{3} \otimes \ket{1}
  \].
\ee

\pagebreak

In building $\ket{\Psi^1_4}$, we require:
\be
  \lefteqn{
    \Delta({E^3}_2)
    \Delta({E^3}_1)
    \cdot
    \(\ket{1} \otimes \ket{1}\)
  }
  \\
  & \stackrel{(\ref{eq:Uqgl21Coproduct}e,\ref{eq:Botia55})}{=} &
  q
  {\[\alpha\]}_q^{\frac{1}{2}}
  \(
    {E^3}_2 \otimes q^{-\frac{1}{2}( {E^2}_2 + {E^3}_3 )}
    +
    q^{\frac{1}{2}( {E^2}_2 + {E^3}_3 )} \otimes {E^3}_2
  \)
  \cdot
  \\
  & &
  \qquad \qquad \qquad \qquad
  \(
    q^{-\frac{\alpha}{2}}
    \ket{3} \otimes \ket{1}
    +
    q^{\frac{\alpha}{2}}
    \ket{1} \otimes \ket{3}
  \)
  \\
  & \stackrel{\eqref{Botia105}}{=} &
  q
  {\[\alpha\]}_q^{\frac{1}{2}}
  \[
    \begin{array}{l}
      q^{-\frac{\alpha}{2}}
      {E^3}_2 \cdot \ket{3}
      \otimes
      q^{-\frac{1}{2}( {E^2}_2 + {E^3}_3 )} \cdot \ket{1}
      \\
      \quad -
      q^{-\frac{\alpha}{2}}
      q^{\frac{1}{2}( {E^2}_2 + {E^3}_3 )} \cdot \ket{3}
      \otimes
      {E^3}_2 \cdot \ket{1}
      \\
      \quad +
      q^{\frac{\alpha}{2}}
      {E^3}_2 \cdot \ket{1}
      \otimes
      q^{-\frac{1}{2}( {E^2}_2 + {E^3}_3 )} \cdot \ket{3}
      \\
      \quad +
      q^{\frac{\alpha}{2}}
      q^{\frac{1}{2}( {E^2}_2 + {E^3}_3 )} \cdot \ket{1}
      \otimes
      {E^3}_2 \cdot \ket{3}
    \end{array}
  \]
  \\
  & \stackrel{(\ref{eq:Botia111},\ref{eq:TheRepresentation}g)}{=} &
  q
  {\[\alpha\]}_q^{\frac{1}{2}}
  \[
    \begin{array}{l}
      q^{\frac{\alpha}{2}}
      {\[\alpha+1\]}_q^{\frac{1}{2}}
      \ket{1}
      \otimes
      \ket{4}
      +
      q^{-\frac{\alpha}{2}}
      {\[\alpha+1\]}_q^{\frac{1}{2}}
      \ket{4}
      \otimes
      \ket{1}
      \\
      \quad -
      q^{-\frac{\alpha}{2}}
      q^{\frac{1}{2}( {E^2}_2 + {E^3}_3 )} \cdot \ket{3}
      \otimes
      {\[\alpha\]}_q^{\frac{1}{2}}
      \ket{2}
      \\
      \quad +
      q^{\frac{\alpha}{2}}
      {\[\alpha\]}_q^{\frac{1}{2}}
      \ket{2}
      \otimes
      q^{-\frac{1}{2}( {E^2}_2 + {E^3}_3 )} \cdot \ket{3}
    \end{array}
  \]
  \\
  & \stackrel{(\ref{eq:TheRepresentationExtended}b,c)}{=} &
  q
  {\[\alpha\]}_q^{\frac{1}{2}}
  \[
    \begin{array}{l}
      {\[\alpha+1\]}_q^{\frac{1}{2}}
      \(
        q^{\frac{\alpha}{2}}
        \ket{1} \otimes \ket{4}
        +
        q^{-\frac{\alpha}{2}}
        \ket{4} \otimes \ket{1}
      \)
      \\
      \quad -
      {\[\alpha\]}_q^{\frac{1}{2}}
      \(
        q^{-\frac{\alpha}{2}}
        q^{\frac{1}{2}+\frac{\alpha}{2}}
        \ket{3} \otimes \ket{2}
        -
        q^{\frac{\alpha}{2}}
        q^{-\frac{1}{2}-\frac{\alpha}{2}}
        \ket{2} \otimes \ket{3}
      \)
    \end{array}
  \]
  \\
  \eq
  q
  {\[\alpha\]}_q^{\frac{1}{2}}
  \[
    \begin{array}{l}
      {\[\alpha+1\]}_q^{\frac{1}{2}}
      (
        q^{\alpha}
        \ket{1} \otimes \ket{4}
        +
        q^{-\alpha}
        \ket{4} \otimes \ket{1}
      )
      \\
      \qquad -
      {\[\alpha\]}_q^{\frac{1}{2}}
      (
        q^{\frac{1}{2}}
        \ket{3} \otimes \ket{2}
        -
        q^{-\frac{1}{2}}
        \ket{2} \otimes \ket{3}
      )
    \end{array}
  \].
\ee

To deduce $ \beta^1_4 $, we demand
$\braket{\Psi^1_4}{\Psi^1_4}=1$. Investigating:
\be
  \lefteqn{
    \m{[}{l}
      {\[\alpha+1\]}_q^{\frac{1}{2}}
      (
        q^{\alpha}
        \bra{1} \otimes \bra{4}
        +
        q^{-\alpha}
        \bra{4} \otimes \bra{1}
      )
      \\
      \qquad+
      {\[\alpha\]}_q^{\frac{1}{2}}
      (
        q^{\frac{1}{2}}
        \bra{3} \otimes \bra{2}
        -
        q^{-\frac{1}{2}}
        \bra{2} \otimes \bra{3}
      )
    \me
    \cdot
  }
  \\
  & &
  \qquad
  \m{[}{l}
    {\[\alpha+1\]}_q^{\frac{1}{2}}
    (
      q^{\alpha}
      \ket{1} \otimes \ket{4}
      +
      q^{-\alpha}
      \ket{4} \otimes \ket{1}
    )
    \\
    \qquad -
    {\[\alpha\]}_q^{\frac{1}{2}}
    (
      q^{\frac{1}{2}}
      \ket{3} \otimes \ket{2}
      -
      q^{-\frac{1}{2}}
      \ket{2} \otimes \ket{3}
    )
  \me
  \\
  & \stackrel{(\ref{eq:DualBasesMultRules}b)}{=} &
  \[
    {\[\alpha+1\]}_q
    \( q^{2\alpha} + q^{-2 \alpha} \)
    +
    {\[\alpha\]}_q
    \( q^{1} + q^{-1} \)
  \]
  \\
  \eq
  {\( q - q^{-1} \)}^{-1}
  \[
    q^{3\alpha+1}
    -
    q^{-3\alpha-1}
    +
    q^{\alpha+1}
    -
    q^{-\alpha-1}
  \]
  \\
  \eq
  {\( q - q^{-1} \)}^{-1}
  \(
    q^{\alpha}
    +
    q^{-\alpha}
  \)
  \(
    q^{2\alpha+1}
    -
    q^{-2\alpha-1}
  \)
  \\
  \eq
  {\[2\alpha+1\]}_q
  \( q^{\alpha} + q^{-\alpha} \),
\ee
hence we require, up to a phase factor
$
  \beta^1_4
  =
  q^{-1}
  {\[\alpha\]}_q^{-\frac{1}{2}}
  {\( q^{\alpha} + q^{-\alpha} \)}^{-\frac{1}{2}}
  {\[2\alpha+1\]}_q^{-\frac{1}{2}}
$.
Thus, (\ref{eq:InitialKetPsi1Definition}$c$) yields,
in accordance with (\ref{eq:V1KetBasis}$d$):
\be
  \! \! \! \! \! \! \! \! \! \!
  \! \! \! \! \! \! \! \! \! \!
  \ket{\Psi^1_4}
  \eq
  {\( q^{\alpha} + q^{-\alpha} \)}^{-\frac{1}{2}}
  {[2\alpha+1]}_q^{-\frac{1}{2}}
  \cdot
  \\
  & &
  \! \! \! \! \! \! \! \!
  \! \! \! \! \! \! \! \!
  \[
    {[\alpha+1]}_q^{\frac{1}{2}}
    \(
      q^{\alpha} \ket{1} \otimes \ket{4}
      +
      q^{-\alpha} \ket{4} \otimes \ket{1}
    \)
    -
    {[\alpha]}_q^{\frac{1}{2}}
    \(
      q^{\frac{1}{2}} \ket{3} \otimes \ket{2}
      -
      q^{-\frac{1}{2}} \ket{2} \otimes \ket{3}
    \)
  \].
\ee


\subsubsection{An Orthonormal Basis for $V_3$}
\seclabel{ABasisforV_3}

The construction of a normalised basis for $V_3$ is, in a sense, dual
to that of $V_1 $.
The \emph{lowest} weight vector of $V$ is $\ket{4}$, of weight
$\(-1,-1\,|\,\alpha+2\)$ (see \eqref{WeightsofKets}), hence the lowest
weight vector of $V \otimes V$ is $\ket{4} \otimes \ket{4}$, of weight
$\(-2,-2\,|\,2\alpha+4\)$.

We shall regard the submodule of $V \otimes V$ generated by the action
of (the coproducts of) the odd \emph{raising} generators on this vector
$\ket{4} \otimes \ket{4}$ as $V_3$, i.e.  $V_3$ has lowest weight
$\(-2,-2\,|\,-2\alpha+4\)$, from which we may determine, again using
Kac's dimension formula that $V_3$ has dimension $4$ (as did $V_1$).

We shall construct an orthonormal basis $\{\ket{\Psi^3_j}\}_{j=1}^4$
for $V_3$, using the Kac induced module construction.
Initially, as already defined in (\ref{eq:V1KetBasis}$a$):
\be
  \ket{\Psi^3_4}
  \defeq
  \ket{4}
  \otimes
  \ket{4}.
\ee
As for $\ket{1}\otimes\ket{1}$, the normality of $\ket{4}$ implies the
normality of $\ket{4}\otimes\ket{4}$:
\be
  \hspace{-2pt}
  \(\bra{4}\otimes\bra{4}\)
  \cdot
  \(\ket{4}\otimes\ket{4}\)
  =
  {\(-\)}^{\[\bra{4}\]\[\ket{4}\]}
  \( \braket{4}{4}\otimes\braket{4}{4} \)
  =
  {\(-\)}^{0 \cdot 0}
  1
  =
  1.
\ee

The other vectors in $V_3$ are generated from the action of (the
coproduct of) the odd raising generators of $U_q[gl(2|1)]$ on this
lowest weight state $ \ket{\Psi^3_4} = \ket{4} \otimes \ket{4} $:
\bse
  \left.
  \begin{array}{lrcl}
    (a) &
    \ket{\Psi^3_3}
    & \defeq &
    \beta^3_3
    \Delta({E^2}_3)
    \cdot
    \(\ket{4} \otimes \ket{4}\)
    \\
    (b) &
    \ket{\Psi^3_2}
    & \defeq &
    \beta^3_2
    \Delta({E^1}_3)
    \cdot
    \(\ket{4} \otimes \ket{4}\)
    \\
    (c) &
    \ket{\Psi^3_1}
    & \defeq &
    \beta^3_1
    \Delta({E^2}_3 {E^1}_3)
    \cdot
    \(\ket{4} \otimes \ket{4}\),
  \end{array}
  \qquad \qquad \qquad \qquad \qquad \; \;
  \right\}
  \eqlabel{InitialKetPsi3Definition}
\ese
where the $ \beta^3_j $, for $ j = 3,2,1 $ are suitable normalisation
constants.
Again, the vectors are ordered in terms of decreasing weight:
\be
  \begin{array}{
        ccl@{}r@{\hspace{1mm}}r@{\hspace{1mm}}c@{\hspace{1mm}}l@{}l@{}c
               }
    \ket{\Psi^3_1} & \textrm{has~weight} &
      ( & -1, & -1 & | & 2 \alpha+2 & ) & \\
    \ket{\Psi^3_2} &                     &
      ( & -1, & -2 & | & 2 \alpha+3 & ) & \\
    \ket{\Psi^3_3} &                     &
      ( & -2, & -1 & | & 2 \alpha+3 & ) & \\
    \ket{\Psi^3_4} &                     &
      ( & -2, & -2 & | & 2 \alpha+4 & ) & .
  \end{array}
\ee
Recall that $\ket{4}$ is a \emph{lowest weight} vector of $V_1$. This
means that it is annihilated by all lowering generators.  Without
further ado, let us discover $\ket{\Psi^3_3}$. Firstly, we have:
\be
  \lefteqn{
    \hspace{-5mm}
    \Delta ( {E^2}_3 )
    \cdot
    \( \ket{4} \otimes \ket{4} \)
  }
  \\
  & \stackrel{(\ref{eq:Uqgl21Coproduct}d)}{=} &
  {E^2}_3 \cdot \ket{4}
  \otimes
  q^{-\frac{1}{2} ( {E^2}_2 + {E^3}_3 )}
  \cdot
  \ket{4}
  +
  q^{\frac{1}{2}( {E^2}_2 + {E^3}_3 )}
  \cdot
  \ket{4}
  \otimes
  {E^2}_3
  \cdot
  \ket{4}
  \\
  & \stackrel{(\ref{eq:TheRepresentation}f)}{=} &
  {\[\alpha+1\]}_q^{1/2}
  \ket{3}
  \otimes
  q^{-\frac{1}{2} ( {E^2}_2 + {E^3}_3 )}
  \cdot
  \ket{4}
  +
  q^{\frac{1}{2}( {E^2}_2 + {E^3}_3 )}
  \cdot
  \ket{4}
  \otimes
  {\[\alpha+1\]}_q^{1/2}
  \ket{3}
  \\
  & \stackrel{(\ref{eq:TheRepresentationExtended}b,c)}{=} &
  {\[\alpha+1\]}_q^{1/2}
  \[
    q^{\frac{\alpha+1}{2}}
    \ket{4} \otimes \ket{3}
    +
    q^{-\frac{\alpha+1}{2}}
    \ket{3} \otimes \ket{4}
  \].
\ee
To deduce $ \beta^3_3 $, we demand
$\braket{\Psi^3_3}{\Psi^3_3}=1$. Investigating:
\be
  \lefteqn{
    \hspace{-7mm}
    \[
      q^{\frac{\alpha+1}{2}}
      \bra{4} \otimes \bra{3}
      +
      q^{-\frac{\alpha+1}{2}}
      \bra{3} \otimes \bra{4}
    \]
    \cdot
    \[
      q^{\frac{\alpha+1}{2}}
      \ket{4} \otimes \ket{3}
      +
      q^{-\frac{\alpha+1}{2}}
      \ket{3} \otimes \ket{4}
    \]
  }
  \\
  & \stackrel{(\ref{eq:DualBasesMultRules}b)}{=} &
  \(
    q^{\alpha+1}
    +
    q^{-\alpha-1}
  \),
\ee
hence we require, up to a phase factor
$
  \beta^3_3
  =
  {\[\alpha+1\]}_q^{1/2}
  {\(
    q^{\alpha+1}
    +
    q^{-\alpha-1}
  \)}^{-1/2}
$.
Thus, (\ref{eq:InitialKetPsi3Definition}$a$) yields,
in accordance with (\ref{eq:V3KetBasis}$c$):
\be
  \ket{\Psi^3_3}
  =
  {\(
    q^{\alpha+1}
    +
    q^{-\alpha-1}
  \)}^{-1/2}
  \[
    q^{\frac{\alpha+1}{2}}
    \ket{4} \otimes \ket{3}
    +
    q^{-\frac{\alpha+1}{2}}
    \ket{3} \otimes \ket{4}
  \].
\ee

To build $\ket{\Psi^3_2}$, we evaluate:
\bne
  \lefteqn{
    \Delta\({E^1}_3\)
    \cdot
    \(\ket{4} \otimes \ket{4}\)
  }
  \nonumber
  \\
  & \stackrel{\eqref{DeltaE13}}{=} &
  \(
    {E^1}_2 {E^2}_3
    -
    q^{-1}  {E^2}_3 {E^1}_2
  \)
  \cdot
  \ket{4}
  \otimes
  q^{-\frac{1}{2} ( {E^1}_1 + {E^3}_3)}
  \cdot
  \ket{4}
  \nonumber
  \\
  & &
  +
  q^{\frac{1}{2}( {E^1}_1 + {E^3}_3 )}
  \cdot
  \ket{4}
  \otimes
  \(
    {E^1}_2 {E^2}_3
    -
    q^{-1}
    {E^2}_3 {E^1}_2
  \)
  \cdot
  \ket{4}
  \nonumber
  \\
  & &
  +
  \( q - q^{-1} \)
  q^{\frac{1}{2} ( {E^2}_2 + {E^3}_3 )} {E^1}_2
  \cdot
  \ket{4}
  \otimes
  {E^2}_3 q^{-\frac{1}{2}( {E^1}_1 - {E^2}_2 )}
  \cdot
  \ket{4}
  \nonumber
  \\
  & \stackrel{(\ref{eq:TheRepresentation}d,f,
               \ref{eq:TheRepresentationExtended}a,c)}{=} &
  {\[\alpha+1\]}_q^{1/2}
  \(
    q^{-\frac{\alpha}{2}}
    {E^1}_2
    \cdot
    \ket{3}
    \otimes
    \ket{4}
    +
    q^{\frac{\alpha}{2}}
    \ket{4}
    \otimes
    {E^1}_2
    \cdot
    \ket{3}
  \)
  \nonumber
  \\
  & \stackrel{(\ref{eq:TheRepresentation}d)}{=} &
  -
  {\[\alpha+1\]}_q^{1/2}
  \[
    q^{\frac{\alpha+1}{2}}
    \ket{4} \otimes \ket{2}
    +
    q^{-\frac{\alpha+1}{2}}
    \ket{2} \otimes \ket{4}
  \].
  \eqlabel{Botia103}
\ene
To deduce $ \beta^3_2 $, we demand
$\braket{\Psi^3_2}{\Psi^3_2}=1$. Investigating:
\be
  \lefteqn{
    \hspace{-7mm}
    \[
      q^{\frac{\alpha+1}{2}}
      \bra{4} \otimes \bra{2}
      +
      q^{-\frac{\alpha+1}{2}}
      \bra{2} \otimes \bra{4}
    \]
    \cdot
    \[
      q^{\frac{\alpha+1}{2}}
      \ket{4} \otimes \ket{2}
      +
      q^{-\frac{\alpha+1}{2}}
      \ket{2} \otimes \ket{4}
    \]
  }
  \\
  & \stackrel{(\ref{eq:DualBasesMultRules}b)}{=} &
  \(
    q^{\alpha+1}
    +
    q^{-\alpha-1}
  \),
\ee
hence we require, up to a phase factor
$
  \beta^3_2
  =
  {\[\alpha+1\]}_q^{1/2}
  {\(
    q^{\alpha+1}
    +
    q^{-\alpha-1}
  \)}^{-1/2}
$.
Thus, (\ref{eq:InitialKetPsi3Definition}$b$) yields,
in accordance with (\ref{eq:V3KetBasis}$b$):
\be
  \ket{\Psi^3_2}
  =
  {\(
    q^{\alpha+1}
    +
    q^{-\alpha-1}
  \)}^{-1/2}
  \[
    q^{\frac{\alpha+1}{2}}
    \ket{4} \otimes \ket{2}
    +
    q^{-\frac{\alpha+1}{2}}
    \ket{2} \otimes \ket{4}
  \].
\ee

\pagebreak

In building $\ket{\Psi^3_1}$, we require:
\be
  \lefteqn{
    \hspace{-5mm}
    \Delta({E^2}_3 {E^1}_3)
    \cdot
    \(\ket{4} \otimes \ket{4}\)
  }
  \\
  & \stackrel{(\ref{eq:CoproductGradedHomo})}{=} &
  \Delta({E^2}_3)
  \Delta({E^1}_3)
  \cdot
  \(\ket{4} \otimes \ket{4}\)
  \\
  & \stackrel{(\ref{eq:Botia103},
               \ref{eq:Uqgl21Coproduct}d)}{=} &
  -
  {\[\alpha+1\]}_q^{1/2}
  \[
    {E^2}_3 \otimes q^{-\frac{1}{2} ( {E^2}_2 + {E^3}_3 )}
    +
    q^{\frac{1}{2}( {E^2}_2 + {E^3}_3 )} \otimes {E^2}_3
  \]
  \cdot
  \\
  & &
  \qquad
  \[
    q^{\frac{1}{2}(\alpha+1)}
    \ket{4} \otimes \ket{2}
    +
    q^{-\frac{1}{2}(\alpha+1)}
    \ket{2} \otimes \ket{4}
  \]
  \\
  & \stackrel{\eqref{Botia105}}{=} &
  -
  {\[\alpha+1\]}_q^{1/2}
  \m{[}{l}
    q^{\frac{1}{2}(\alpha+1)}
    {E^2}_3
    \cdot
    \ket{4}
    \otimes
    q^{-\frac{1}{2} ( {E^2}_2 + {E^3}_3 )}
    \cdot
    \ket{2}
    \\
    \quad +
    q^{\frac{1}{2}(\alpha+1)}
    q^{\frac{1}{2}( {E^2}_2 + {E^3}_3 )}
    \cdot
    \ket{4}
    \otimes
    {E^2}_3
    \cdot
    \ket{2}
    \\
    \quad +
    q^{-\frac{1}{2}(\alpha+1)}
    {E^2}_3
    \cdot
    \ket{2}
    \otimes
    q^{-\frac{1}{2} ( {E^2}_2 + {E^3}_3 )}
    \cdot
    \ket{4}
    \\
    \quad -
    q^{-\frac{1}{2}(\alpha+1)}
    q^{\frac{1}{2}( {E^2}_2 + {E^3}_3 )}
    \cdot
    \ket{2}
    \otimes
    {E^2}_3
    \cdot
    \ket{4}
  \me
  \\
  & \stackrel{(\ref{eq:TheRepresentation}f,
               \ref{eq:TheRepresentationExtended}b,c)}{=} &
  -
  {\[\alpha+1\]}_q^{1/2}
  \[
    \begin{array}{l}
      {\[\alpha\]}_q^{1/2}
      \[
        q^{\alpha+1}
        \ket{4} \otimes \ket{1}
        +
        q^{-\alpha-1}
        \ket{1} \otimes \ket{4}
      \]
      \\
      \hspace{5mm}
      +
      {\[\alpha+1\]}_q^{1/2}
      \[
        q^{\frac{1}{2}}
        \ket{3} \otimes \ket{2}
        -
        q^{-\frac{1}{2}}
        \ket{2} \otimes \ket{3}
      \]
    \end{array}
  \].
\ee
To deduce $ \beta^3_1 $, we demand
$\braket{\Psi^3_1}{\Psi^3_1}=1$. Investigating:
\be
  \lefteqn{
    \hspace{-5mm}
    \m{[}{l}
      {\[\alpha\]}_q^{1/2}
      \[
        q^{\alpha+1}
        \bra{4} \otimes \bra{1}
        +
        q^{-\alpha-1)}
        \bra{1} \otimes \bra{4}
      \]
      \\
      \hspace{5mm}
      -
      {\[\alpha+1\]}_q^{1/2}
      \[
        q^{\frac{1}{2}}
        \bra{3} \otimes \bra{2}
        -
        q^{-\frac{1}{2}}
        \bra{2} \otimes \bra{3}
      \]
    \me
    \cdot
  }
  \\
  & &
  \m{[}{l}
    {\[\alpha\]}_q^{1/2}
    \[
      q^{\alpha+1}
      \ket{4} \otimes \ket{1}
      +
      q^{-\alpha-1}
      \ket{1} \otimes \ket{4}
    \]
    \\
    \hspace{5mm}
    +
    {\[\alpha+1\]}_q^{1/2}
    \[
      q^{\frac{1}{2}}
      \ket{3} \otimes \ket{2}
      -
      q^{-\frac{1}{2}}
      \ket{2} \otimes \ket{3}
    \]
  \me
  \\
  & \stackrel{(\ref{eq:DualBasesMultRules}b)}{=} &
  {\[\alpha\]}_q
  \(
    q^{2 \alpha + 2}
    +
    q^{- 2 \alpha - 2}
  \)
  +
  {\[\alpha+1\]}_q
  \(
    q^{1}
    +
    q^{-1}
  \)
  \\
  \eq
  {\( q - q^{-1} \)}^{-1}
  \m{[}{l}
    \( q^\alpha - q^{-\alpha} \)
    \(
      q^{2 \alpha + 2}
      +
      q^{- 2 \alpha - 2}
    \)
    \\
    \hspace{5mm}
    +
    \( q^{\alpha+1} - q^{-\alpha-1} \)
    \(
      q^{1}
      +
      q^{-1}
    \)
  \me
  \\
  \eq
  {\( q - q^{-1} \)}^{-1}
  \[
    q^{3 \alpha + 2}
    -
    q^{- 3 \alpha - 2}
    +
    q^{\alpha}
    -
    q^{-\alpha}
  \]
  =
  \(
    q^{\alpha + 1}
    +
    q^{-\alpha-1}
  \)
  {\[2\alpha+1\]}_q.
\ee
Thus, up to a phase factor
$
  \beta^3_1
  =
  {\(
    q^{\alpha + 1}
    +
    q^{-\alpha-1}
  \)}^{-1/2}
  {\[\alpha+1\]}_q^{-1/2}
  {\[2\alpha+1\]}_q^{-1/2}
$,
so (\ref{eq:InitialKetPsi3Definition}$c$) yields,
in accordance with (\ref{eq:V3KetBasis}$a$):
\be
  \hspace{-10pt}
  \ket{\Psi^3_1}
  \eq
  {(q^{\alpha+1}+q^{-\alpha-1})}^{-\frac{1}{2}}
  {[2\alpha+1]}_q^{-\frac{1}{2}}
  \cdot
  \nonumber
  \\
  & &
  \hspace{-20mm}
  \[
    {[\alpha]}_q^{\frac{1}{2}}
    \(
      q^{\alpha+1} \ket{4} \otimes \ket{1}
      +
      q^{-\alpha-1} \ket{1} \otimes \ket{4}
    \)
    -
    [\alpha+1]_q^{\frac{1}{2}}
    \(
      q^{-\frac{1}{2}} \ket{2} \otimes \ket{3}
      -
      q^{\frac{1}{2}} \ket{3} \otimes \ket{2}
    \)
  \].
\ee

\pagebreak


\subsubsection{No Basis for $V_2$?}
\seclabel{NoBasisforV_2}

Although it is not needed below, it would be more complete to find such
a basis.  $V_2$ is necessarily $8$ dimensional, and its highest weight
is $(0,-1\,|\,2\alpha+1)$, determined as the highest remaining weight
after the weights used in the basis for $V_1$ have been eliminated.  In
order to build a basis for $V_2$, we would have to define a highest
weight vector $\ket{\Psi^2_1}$ as a linear combination of weight
vectors of this weight, and then build the rest of the basis by the
action of (the coproduct of) the lowering generators on
$\ket{\Psi^2_1}$. In fact, we may write down $\ket{\Psi^2_1}$
immediately. The only possible combinations of weights that sum to
$\(0,-1\,|\,2\alpha+1\)$ are those of the vectors
$\ket{1}\otimes\ket{2}$ and $\ket{2}\otimes\ket{1}$. This combination
already appears once, in $\ket{\Psi^1_2}$; examination of that vector
(\ref{eq:V1KetBasis}$b$) shows that the only other possible orthogonal,
normalised combination is (as usual, up to a phase factor): 
\be
  \ket{\Psi^2_1}
  =
  {(q^{\alpha}+q^{-\alpha})}^{-\frac{1}{2}}
  \(
    q^{-\frac{\alpha}{2}} \ket{1} \otimes \ket{2}
    -
    q^{\frac{\alpha}{2}} \ket{2} \otimes \ket{1}
  \).
\ee
In this case, the underlying ${gl\(2|1\)}_{\overline{0}}$ module is
$2$, not $1$ dimensional, i.e.  the action of the even lowering
generator ${E^2}_1$ (viz $\Delta({E^2}_1)$) on $\ket{\Psi^2_1}$ is
\emph{not} to annihilate it, cf. \eqref{E21Annihilatesket1}.  Thus, the
construction of $V_2$ would involve requiring the next lowest weight
vector of the top level to be:
\be
  \ket{\Psi^2_3}
  \defeq
  \beta^2_3
  \Delta({E^2}_1)
  \cdot
  \ket{\Psi^2_1}.
\ee
This vector is actually of a lower weight (viz $(-1,0\,|\,2\alpha+1)$)
than a vector in the next level, hence its numbering is peculiar.
This would complete a basis for the top level of $V_2$.

The next level would have a basis
$
  \{
    \ket{\Psi^2_2},
    \ket{\Psi^2_4},
    \ket{\Psi^2_5},
    \ket{\Psi^2_6}
  \}
$
consisting of the $4$ vectors formed
by the action of ${E^3}_2$ and ${E^3}_1$ on $\ket{\Psi^2_1}$ and
$\ket{\Psi^2_3}$, where our vectors are naturally ordered in terms of
decreasing weight. So far, we have:
\be
  \begin{array}{
       ccccl@{}r@{\hspace{1mm}}r@{\hspace{1mm}}c@{\hspace{1mm}}l@{}l@{}c
               }
    \ket{\Psi^2_1} &  &      &\textrm{of~weight} &
      ( &  0, & -1 & | & 2 \alpha+1 & ) & \\
    \ket{\Psi^2_2} \eq \beta^2_2 \Delta({E^3}_2) \cdot \ket{\Psi^2_1} &&
      ( &  0, & -2 & | & 2 \alpha+2 & ) & \\
    \ket{\Psi^2_3} \eq \beta^2_3 \Delta({E^2}_1) \cdot \ket{\Psi^2_1} &&
      ( & -1, &  0 & | & 2 \alpha+1 & ) & \\
    \ket{\Psi^2_4} \eq \beta^2_4 \Delta({E^3}_1) \cdot \ket{\Psi^2_1} &&
      ( & -1, & -1 & | & 2 \alpha+2 & ) & \\
    \ket{\Psi^2_5} \eq \beta^2_5 \Delta({E^3}_2) \cdot \ket{\Psi^2_3} &&
      ( & -1, & -1 & | & 2 \alpha+2 & ) & \\
    \ket{\Psi^2_6} \eq \beta^2_6 \Delta({E^3}_1) \cdot \ket{\Psi^2_3} &&
      ( & -2, &  0 & | & 2 \alpha+2 & ) & .
  \end{array}
\ee
Observe that the natural definition of $V_2$ in this manner yields a
basis that contains a pair of vectors of the \emph{same} weight, viz
$\ket{\Psi^2_4}$ and $\ket{\Psi^2_5}$. This pair will not naturally be
orthogonal, so we \emph{may} have to use the Gram--Schmidt process to
orthonormalise them (although we might be lucky!).

The construction of a basis for $V_2$ would then continue for one more
level.

\pagebreak


\subsection{Projectors onto the Subspaces $V_k\subset V\otimes V$}
\addtocontents{toc}{\protect\vspace{-2.5ex}}

From the sets of basis vectors $\{\ket{\Psi^k_j}\}$ and their duals
$\{\bra{\Psi^k_j}\}$,
we may construct projectors $P_k:V\otimes V\to V_k$ by:
\bse
  P_k
  =
  \sum_{j=1}^4
    \ket{\Psi^k_j} \bra{\Psi^k_j},
  \qquad \qquad
  k = 1, 3,
  \eqlabel{ProjectorDefinition}
\ese
where the multiplications are defined by \eqref{DualBasesMultRules}.
Of course, we could also define $P_2$, but we haven't defined a basis
for $V_2$ \dots
The gentle reader is reminded that projectors satisfy:
$ P_i P_j = \delta_{i j} P_j $.

Now let $ I $ be the identity operator on the tensor product space
$V\otimes V$, defined as:
\be
  I
  =
  \sum_{i, j=1}^4
    \ketbra{i}{i} \otimes \ketbra{j}{j}.
\ee

As $ P_1 + P_2 + P_3 = I $, we thus do not need to explicitly
construct $ P_2 $ (or even a basis for $ V_2 $); we simply set:
\bse
  P_2
  \defeq
  I - P_1 - P_3.
  \eqlabel{P2def}
\ese
An alternate notation is to let $ {e^i}_j \defeq \ketbra{i}{j}$ be an
elementary rank $2$ tensor, which ensures that we must write
the elementary rank $4$ tensors as:
\bse
  {e^{i k}}_{j l}
  =
  {e^i}_j \otimes {e^k}_l
  =
  \ketbra{i}{j} \otimes \ketbra{k}{l}.
  \eqlabel{KetBrasasTensors}
\ese
We then have:
\bse
  I
  =
  \sum_{i, j=1}^4
    {e^i}_i \otimes {e^j}_j
  =
  \sum_{i, j=1}^4
    {e^{i j}}_{i j}.
  \eqlabel{Idef}
\ese

\pagebreak
We may substitute the definitions of the dual bases for $V\otimes V$
(\ref{eq:V1KetBasis}, \ref{eq:V3KetBasis}, \ref{eq:V1BraBasis} and
\ref{eq:V3BraBasis}) into the definition of the projectors $P_1$ and
$P_3$ \eqref{ProjectorDefinition}, paying heed to the definitions of
the tensor multiplications (\ref{eq:DualBasesMultRules}$a$) to yield the
following explicit expressions (checked with \textsc{Mathematica}!):
\bne
  P_1
  \eq
  {e^{1 1}}_{1 1}
  \nonumber
  \\
  & & +
  (q^{\alpha}+q^{-\alpha})^{-1}
  \cdot
  \nonumber
  \\
  & &
  \quad
  \(
    \begin{array}{l}
      \(
        q^{\alpha}
        {e^{1 2}}_{1 2}
        +
        q^{-\alpha}
        {e^{2 1}}_{2 1}
      \)
      +
      \(
        {e^{2 1}}_{1 2}
        -
        {e^{1 2}}_{2 1}
      \)
      \\
      \quad
      +
      \(
        q^{\alpha}
        {e^{1 3}}_{1 3}
        +
        q^{-\alpha}
        {e^{3 1}}_{3 1}
      \)
      +
      \(
        {e^{3 1}}_{1 3}
        -
        {e^{1 3}}_{3 1}
      \)
    \end{array}
  \)
  \nonumber
  \\
  & & +
  (q^{\alpha}+q^{-\alpha})^{-1}
  {[2\alpha+1]}_q^{-1}
  \cdot
  \nonumber
  \\
  & &
  \left\{
    \begin{array}{l}
      {[\alpha+1]}_q
      \(
        \(
          q^{2\alpha}
          {e^{1 4}}_{1 4}
          +
          q^{-2\alpha}
          {e^{4 1}}_{4 1}
        \)
        +
        \(
          {e^{4 1}}_{1 4}
          +
          {e^{1 4}}_{4 1}
        \)
      \)
      \\
      +
      {[\alpha]}_q^{\frac{1}{2}}
      [\alpha+1]_q^{\frac{1}{2}}
      \cdot
      \\
      \quad
      \(
        \begin{array}{l}
          q^{\alpha+\frac{1}{2}}
          {e^{1 4}}_{3 2}
          -
          q^{-\alpha-\frac{1}{2}}
          {e^{4 1}}_{2 3}
          \\
          \quad
          +
          q^{\alpha-\frac{1}{2}}
          {e^{2 3}}_{1 4}
          -
          q^{-\alpha+\frac{1}{2}}
          {e^{3 2}}_{4 1}
          \\
          \quad
          +
          q^{-\alpha-\frac{1}{2}}
          {e^{2 3}}_{4 1}
          -
          q^{\alpha+\frac{1}{2}}
          {e^{3 2}}_{1 4}
          \\
          \quad
          +
          q^{-\alpha+\frac{1}{2}}
          {e^{4 1}}_{3 2}
          -
          q^{\alpha-\frac{1}{2}}
          {e^{1 4}}_{2 3}
        \end{array}
      \)
      \\
      +
      {[\alpha]}_q
      \(
        \(
          q^{1}
          {e^{3 2}}_{3 2}
          +
          q^{-1}
          {e^{2 3}}_{2 3}
        \)
        -
        \(
          {e^{3 2}}_{2 3}
          +
          {e^{2 3}}_{3 2}
        \)
      \)
    \end{array}
  \right\}.
  \eqlabel{ExplicitP1}
\ene
\bne
  P_3
  \eq
  {(q^{\alpha+1} + q^{-\alpha-1})}^{-1}
  {[2\alpha+1]}_q^{-1}
  \cdot
  \nonumber
  \\
  & &
  \left\{
  \begin{array}{l}
    {[\alpha]}_q
    \(
      \(
        q^{2\alpha+2}
        {e^{4 1}}_{4 1}
        +
        q^{-2\alpha-2}
        {e^{1 4}}_{1 4}
      \)
      +
      \(
        {e^{4 1}}_{1 4}
        +
        {e^{1 4}}_{4 1}
      \)
    \)
    \\
    +
    {[\alpha]}_q^{\frac{1}{2}}
    {[\alpha+1]}_q^{\frac{1}{2}}
    \cdot
    \nonumber
    \\
    \quad
    \(
      \begin{array}{l}
        q^{\alpha+\frac{1}{2}}
        {e^{4 1}}_{2 3}
        -
        q^{-\alpha-\frac{1}{2}}
        {e^{1 4}}_{3 2}
        \\
        \quad
        +
        q^{\alpha+\frac{3}{2}}
        {e^{3 2}}_{4 1}
        -
        q^{-\alpha-\frac{3}{2}}
        {e^{2 3}}_{1 4}
        \\
        \quad
        +
        q^{-\alpha-\frac{1}{2}}
        {e^{3 2}}_{1 4}
        -
        q^{\alpha+\frac{1}{2}}
        {e^{2 3}}_{4 1}
        \\
        \quad
        +
        q^{-\alpha-\frac{3}{2}}
        {e^{1 4}}_{2 3}
        -
        q^{\alpha+\frac{3}{2}}
        {e^{4 1}}_{3 2}
      \end{array}
    \)
    \\
    +
    {[\alpha+1]}_q
    \(
      \(
        q^{1}
        {e^{3 2}}_{3 2}
        +
        q^{-1}
        {e^{2 3}}_{2 3}
      \)
      -
      \(
        {e^{3 2}}_{2 3}
        +
        {e^{2 3}}_{3 2}
      \)
    \)
  \end{array}
  \right\}
  \nonumber
  \\
  & &
  +
  {(q^{\alpha+1} + q^{-\alpha-1})}^{-1}
  \cdot
  \nonumber
  \\
  & &
  \quad
  \(
    \begin{array}{l}
      \(
        q^{\alpha+1}
        {e^{4 2}}_{4 2}
        +
        q^{-\alpha-1}
        {e^{2 4}}_{2 4}
      \)
      +
      \(
        {e^{2 4}}_{4 2}
        -
        {e^{4 2}}_{2 4}
      \)
      \\
      \quad
      +
      \(
        q^{\alpha+1}
        {e^{4 3}}_{4 3}
        +
        q^{-\alpha-1}
        {e^{3 4}}_{3 4}
      \)
      +
      \(
        {e^{3 4}}_{4 3}
        -
        {e^{4 3}}_{3 4}
      \)
    \end{array}
  \)
  \nonumber
  \\
  & &
  + {e^{4 4}}_{4 4}.
  \eqlabel{ExplicitP3}
\ene

\pagebreak


\subsection{The Quantum R Matrix and the Braid Generator}
\addtocontents{toc}{\protect\vspace{-2.5ex}}

\subsubsection{The Quantum Yang--Baxter Equation}
\seclabel{TheQuantumYang--BaxterEquation}

For any classical Lie superalgebra $\mathfrak{g}$, the associated
quantum superalgebra $U_q\[\mathfrak{g}\]$ is an example of a
quasitriangular Hopf superalgebra, which necessarily admits a
\emph{universal R matrix}
$R\in U_q\[\mathfrak{g}\]\otimes U_q\[\mathfrak{g}\]$ satisfying (among
other relations) the \emph{quantum Yang--Baxter equation}:
\bse
  R_{12} R_{13} R_{23}
  =
  R_{23} R_{13} R_{12},
  \eqlabel{QYB}
\ese
where the subscripts refer to the (standard) embedding of $R$ acting on
the triple tensor product space
$
  U_q[\mathfrak{g}] \otimes U_q[\mathfrak{g}] \otimes U_q[\mathfrak{g}]
$.

From any representation of $U_q[\mathfrak{g}]$, we may obtain a tensor
solution of \eqref{QYB} by replacing the superalgebra elements with
their matrix representatives.  Observe that the $R$ `matrix' is in fact
a rank $4$ tensor.  Similarly to \eqref{DualBasesMultRules},
multiplication of the tensor product of \emph{matrices} $ a, b, c, d $
is governed by:
\be
  \( a \otimes b \)
  \( c \otimes d \)
  =
  {\( - \)}^{[b][c]}
  \( a c \otimes b d \),
  \qquad
  \textrm{homogeneous} \;
  b, c.
\ee
Recall that we were careful to only define our basis for $V\otimes V$
for $q$ and $\alpha$ real and positive.  As it is known that the R
matrix is unique,%
\footnote{
 The R matrix is unique up to several considerations, see
 \cite[p236]{LinksScheunertGould:94} or \cite{KhoroshkinTolstoy:91}.
}
analytic continuation makes our R matrix valid for any complex $q$ and
$\alpha$.


\subsubsection{The Graded Permutation Operator $P$}

We introduce the \emph{graded permutation operator}
$P:V\otimes V\to V\otimes V$, defined for graded basis vectors
$v^k,v^l\in V$ by:
\bse
  P ( v^k \otimes v^l )
  =
  {\( - \)}^{[k][l]}
  ( v^l \otimes v^k ),
  \eqlabel{Pdef}
\ese
and extended by linearity. (We use the shorthand $[v^k]\equiv\[k\]$.)
Specialising to the case of the $U_q[gl(2|1)]$ representation
$\pi\equiv\pi_{\lambda}$, in terms of the decomposition
\eqref{TPdecomp} of $V\otimes V$, we find%
\footnote{%
  This is taken from \cite{LinksGould:92b}, but with a
  slight change of notation and a convenient choice of normalisation.
}
that as $q\to 1$, $P$ is the identity operator $I$ on $V_1$ and $V_3$,
and $-I$ on $V_2$:
\bse
  P
  =
  \lim_{q\to1}
  \[
    P_1 - P_2 + P_3.
  \]
  \eqlabel{Pproj}
\ese
We will require $P$ in \secref{ObtainingtheRMatrixRandsigma}.

\vfill


\subsubsection{Obtaining the R Matrix $R$ and the
               Braid Generator $\sigma$}
\seclabel{ObtainingtheRMatrixRandsigma}

The \emph{quantum} R matrix $R$ mentioned in
\secref{TheQuantumYang--BaxterEquation} is a specialisation of a
\emph{trigonometric R matrix} $R(u)$, which corresponds to the
associated \emph{affine} quantum superalgebra $U_q[\mathfrak{g}^{(1)}]$
not the basic $U_q[\mathfrak{g}]$, where the parameter $u$ in
$R(u)$ is introduced in the process of the extension of
$U_q[\mathfrak{g}]$ to $U_q[\mathfrak{g}^{(1)}]$. As any representation
of the affine extension is a representation of the underlying quantum
superalgebra, $R(u)$ may also be regarded as acting on
$V\otimes V$, indeed $R(u)$ satisfies the
\emph{spectral parameter-dependent Yang--Baxter equation}
\cite{GouldHibberdLinksZhang:96} (cf. \eqref{QYB}):
\bse
  {R(u)}_{12}
  {R(u+v)}_{13}
  {R(v)}_{23}
  =
  {R(v)}_{23}
  {R(u+v)}_{13}
  {R(u)}_{12}.
  \eqlabel{SpectralQYB}
\ese
Setting $\check{R}(u)=PR(u)$, in
\cite{BrackenGouldLinksZhang:95,BrackenGouldZhangDelius:94b} it was
shown that in our case (i.e. $\mathfrak{g}=gl(2|1)$), $\check{R}(u)$
has the form:
\be
  \check{R}(u)
  \propto
  \frac {q^u-q^{2\alpha}}{1-q^{u+2\alpha}}P_1
  +
  P_2
  +
  \frac{1-q^{u+2\alpha +2}}{q^u-q^{2\alpha+2}} P_3.
\ee
In \cite{GouldHibberdLinksZhang:96}, the constant of proportionality
was chosen as $ 1 $. Here we choose $ -1 $; which will ensure (see
below) that \eqref{Pproj} is satisfied. Thus:
\bse
  \check{R}(u)
  =
  -
  \frac{q^u-q^{2\alpha}}{1-q^{u+2\alpha}}
  P_1
  -
  P_2
  -
  \frac{1-q^{u+2\alpha +2}}{q^u-q^{2\alpha+2}}
  P_3.
  \eqlabel{checkRu}
\ese
In passing from $U_q[\mathfrak{g}^{(1)}]$ to $U_q[\mathfrak{g}]$,
we take the limit as $u\to\infty$, and our R matrix does likewise, viz:
\bse
  \check{R}
  \defeq
  \lim_{u\to\infty}
    \check{R}(u)
  =
  q^{-2 \alpha} P_1
  -
  P_2
  +
  q^{2 \alpha + 2} P_3.
  \eqlabel{PRdef}
\ese
Now, we define \cite{LinksGould:92b,LinksGouldZhang:93}:
\be
  \sigma
  \defeq
  \check{R}.
\ee
Taking $\lim_{u\to\infty}$ in the spectral parameter-dependent
Yang--Baxter equation \eqref{SpectralQYB}:
\be
  \( \sigma \otimes I \)
  \( I \otimes \sigma \)
  \( \sigma \otimes I \)
  =
  \( I \otimes \sigma \)
  \( \sigma \otimes I \)
  \( I \otimes \sigma \).
\ee
It follows that $\sigma$ is a representation of the \emph{braid
generator}, which is a sufficient requirement that $\sigma$ may be
chosen as the (positive) crossing matrix when computing our link
invariant.

Using \eqref{P2def} and \eqref{Idef}, \eqref{PRdef} becomes:
\bne
  \sigma
  \eq
  \( 1 + q^{-2 \alpha} \) P_1
  +
  \( 1 + q^{2 \alpha + 2} \) P_3
  -
  \sum_{ij}
    {e^{i j}}_{i j}.
  \nonumber
  \\
  \eq
  q^{-\alpha} \( q^{\alpha} + q^{-\alpha} \) P_1
  +
  q^{\alpha+1} \( q^{\alpha+1} + q^{-\alpha-1} \) P_3
  -
  \sum_{ij}
    {e^{i j}}_{i j}.
  \eqlabel{sigmadef}
\ene
From the explicit expressions for the projectors
(\ref{eq:ExplicitP1},\ref{eq:ExplicitP3}), we may now construct an
explicit expression for $\sigma$.  To simplify its presentation, we use
the substitution $p\defeq q^\alpha$, which systematically replaces
$\alpha$, so that we may regard $\sigma$ as a function of $q$ and $p$
rather than $q$ and $\alpha$. We also set:
\bne
  Y
  & \defeq &
  {\( p^{-2} - q^2 + p^2 q^2 - 1 \)}^{1/2}
  =
  {\[ q \( p - p^{-1} \) \( p q - p^{-1} q^{-1} \) \]}^{1/2}
  \nonumber
  \\
  \eq
  \( q^{\frac{1}{2}} - q^{-\frac{1}{2}} \)
  {[\alpha]}_q^{\frac{1}{2}} {[\alpha+1]}_q^{\frac{1}{2}}.
  \eqlabel{Ydefn}
\ene
With the help of \textsc{Mathematica}, we find:
\bne
  \sigma
  \eq
    p^{-2}          {e^{1 1}}_{1 1}
  + (p^{-2} - 1)    {e^{2 1}}_{2 1}
  -                 {e^{2 2}}_{2 2}
  \nonumber
  \\
  & &
  + (p^{-2} - 1)    {e^{3 1}}_{3 1}
  + (q^2 - 1)       {e^{3 2}}_{3 2}
  -                 {e^{3 3}}_{3 3}
  \nonumber
  \\
  & &
  + Y^2             {e^{4 1}}_{4 1}
  + (p^2 q^2 - 1)   {e^{4 2}}_{4 2}
  + (p^2 q^2 - 1)   {e^{4 3}}_{4 3}
  + p^2 q^2         {e^{4 4}}_{4 4}
  \nonumber
  \\
  & &
  + p^{-1} ({e^{2 1}}_{1 2} - {e^{1 2}}_{2 1})
  + p^{-1} ({e^{3 1}}_{1 3} - {e^{1 3}}_{3 1})
  \nonumber
  \\
  & &
  +        ({e^{4 1}}_{1 4} + {e^{1 4}}_{4 1})
  - q      ({e^{3 2}}_{2 3} + {e^{2 3}}_{3 2})
  \nonumber
  \\
  & &
  + Y      ({e^{4 1}}_{2 3} - {e^{2 3}}_{4 1})
  + q Y    ({e^{3 2}}_{4 1} - {e^{4 1}}_{3 2})
  \nonumber
  \\
  & &
  + p q      ({e^{2 4}}_{4 2} - {e^{4 2}}_{2 4})
  + p q      ({e^{3 4}}_{4 3} - {e^{4 3}}_{3 4}).
  \eqlabel{Explicitsigma}
\ene
Inspection of the basis vectors $\ket{\Psi^k_j}$ for $V\otimes V$ shows
that we may expect a maximum of $36$ nonzero entries in $\sigma$. In
fact, $ \sigma $ has only $ 26 $ nonzero entries; ten entries
(fortuitously) cancel to zero.

\pagebreak


\subsubsection{Aside: Evaluation of $\check{R}(u)$}

Having constructed projectors, it is a simple matter to apply
\eqref{checkRu} to them to yield the trigonometric R matrix
$\check{R}(u)$, of significance in the study of exactly solvable models
\cite{GouldHibberdLinksZhang:96}. Again, we let \textsc{Mathematica} do
the work:
\bne
  \check{R}(u)
  \eq
  -
  {[\alpha-\textstyle{\frac{u}{2}}]}_q
  {[\alpha+\textstyle{\frac{u}{2}}]}_q^{-1}
  {e^{1 1}}_{1 1}
  -
  ({e^{2 2}}_{2 2} + {e^{3 3}}_{3 3})
  \nonumber
  \\
  & &
  -
  {[\alpha+1+\textstyle{\frac{u}{2}}]}_q
  {[\alpha+1-\textstyle{\frac{u}{2}}]}_q^{-1}
  {e^{4 4}}_{4 4}
  \nonumber
  \\
  & &
  -
  {[\alpha]}_q
  {[\alpha+\textstyle{\frac{u}{2}}]}_q^{-1}
  \[
    q^{\frac{u}{2}}
    ({e^{2 1}}_{2 1} + {e^{3 1}}_{3 1})
    +
    q^{-\frac{u}{2}}
    ({e^{1 2}}_{1 2} + {e^{1 3}}_{1 3})
  \]
  \nonumber
  \\
  & &
  -
  {[\alpha]}_q
  {[\alpha+1]}_q
  {[\alpha+\textstyle{\frac{u}{2}}]}_q^{-1}
  {[\alpha+1-\textstyle{\frac{u}{2}}]}_q^{-1}
  (
    q^{u}
    {e^{4 1}}_{4 1}
    +
    q^{-u}
    {e^{1 4}}_{1 4}
  )
  \nonumber
  \\
  & &
  -
  {[\textstyle{\frac{u}{2}}]}_q
  {[\alpha+\textstyle{\frac{u}{2}}]}_q^{-1}
  (
    {e^{1 2}}_{2 1} - {e^{2 1}}_{1 2}
    +
    {e^{1 3}}_{3 1} - {e^{3 1}}_{1 3}
  )
  \nonumber
  \\
  & &
  +
  {[\alpha]}_q^{\frac{1}{2}}
  {[\alpha+1]}_q^{\frac{1}{2}}
  {[\textstyle{\frac{u}{2}}]}_q
  {[\alpha+\textstyle{\frac{u}{2}}]}_q^{-1}
  {[1+\alpha-\textstyle{\frac{u}{2}}]}_q^{-1}
  \times
  \nonumber
  \\
  & &
  \qquad
  \[
    \begin{array}{l}
      -
      q^{-\frac{1}{2}-\frac{u}{2}}
      ({e^{1 4}}_{2 3} - {e^{2 3}}_{1 4})
      -
      q^{\frac{1}{2}-\frac{u}{2}}
      ({e^{1 4}}_{3 2} - {e^{3 2}}_{1 4})
      \\
      +
      q^{\frac{u}{2}-\frac{1}{2}}
      ({e^{2 3}}_{4 1} - {e^{4 1}}_{2 3})
      +
      q^{\frac{u}{2}+\frac{1}{2}}
      ({e^{3 2}}_{4 1} - {e^{4 1}}_{3 2})
    \end{array}
  \]
  \nonumber
  \\
  & &
  +
  {[\textstyle{\frac{u}{2}}]}_q
  {[\textstyle{1-\frac{u}{2}}]}_q
  {[\alpha+\textstyle{\frac{u}{2}}]}_q^{-1}
  {[\alpha+1-\textstyle{\frac{u}{2}}]}_q^{-1}
  ({e^{1 4}}_{4 1} + {e^{4 1}}_{1 4})
  \nonumber
  \\
  & &
  +
  {(q-q^{-1})}^{-2}
  {[\alpha+\textstyle{\frac{u}{2}}]}_q^{-1}
  {[1+\alpha-\textstyle{\frac{u}{2}}]}_q^{-1}
  \times
  \nonumber
  \\
  & &
  \qquad
  \[
    \begin{array}{l}
      \(2q^{-1} -q^{2\alpha+1} - q^{-2\alpha-1} + q^{-u+1} - q^{-u-1} \)
      {e^{2 3}}_{2 3}
      \\
      +
      \(2 q^1   -q^{-2\alpha-1} - q^{2\alpha+1}  + q^{u-1}  - q^{u+1} \)
      {e^{3 2}}_{3 2}
    \end{array}
  \]
  \nonumber
  \\
  & &
  -
  {[\alpha+1]}_q
  {[\alpha+1-\textstyle{\frac{u}{2}}]}_q^{-1}
  \[
    q^{\frac{u}{2}}
    ({e^{4 2}}_{4 2} + {e^{4 3}}_{4 3})
    +
    q^{-\frac{u}{2}}
    ({e^{2 4}}_{2 4} - {e^{3 4}}_{3 4})
  \]
  \nonumber
  \\
  & &
  +
  {[\textstyle{\frac{u}{2}}]}_q^2
  {[\alpha+\textstyle{\frac{u}{2}}]}_q^{-1}
  {[\alpha+1-\textstyle{\frac{u}{2}}]}_q^{-1}
  ({e^{2 3}}_{3 2} + {e^{3 2}}_{2 3})
  \nonumber
  \\
  & &
  -
  {[\textstyle{\frac{u}{2}}]}_q
  {[\alpha+1-\textstyle{\frac{u}{2}}]}_q^{-1}
  (
    {e^{2 4}}_{4 2}
    -
    {e^{4 2}}_{2 4}
    +
    {e^{3 4}}_{4 3}
    -
    {e^{4 3}}_{3 4}
  ).
  \eqlabel{ExplicitCheckRu}
\ene
This tensor has $36$ nonzero entries: the fortuitous cancellations that
caused $\sigma$ to have only $26$ nonzero entries do not occur here.
Observe that $\lim_{u\to\infty}\check{R}(u)=\check{R}=\sigma$.  From
this expression, one may obtain the Hamiltonian $H^Q(U)$ defined in
\cite{GouldHibberdLinksZhang:96}.

\pagebreak


\subsubsection{A Skein Relation Satisfied by $\sigma$}

From \eqref{sigmadef}, it is
straightforward to deduce that $\sigma$ satisfies the polynomial
identity:
\bse
  q^{-1} \sigma^3
  +
  \( q^{-1} - q^{- 2 \alpha - 1} - q^{2 \alpha + 1} \) \sigma^2
  +
  \( q - q^{- 2 \alpha - 1} - q^{2 \alpha + 1} \) \sigma
  +
  q I
  =
  0,
  \eqlabel{Uqgl21SkeinRelation}
\ese
known as a \emph{skein relation}. To justify
\eqref{Uqgl21SkeinRelation}, begin by observing that from
\eqref{PRdef}, the eigenvalues of $\sigma$
are $q^{-2\alpha}$, $-1$ and $q^{2\alpha+2}$, hence:
\be
  \(\sigma-q^{-2\alpha}I\)
  \(\sigma+I\)
  \(\sigma-q^{2\alpha+2}I\)
  =
  0,
\ee
which may be expanded to yield:
\be
  \sigma^3
  +
  \( 1 - q^{-2\alpha} - q^{2\alpha+2} \)
  \sigma^2
  +
  \( q^2 - q^{-2\alpha} - q^{2\alpha+2} \)
  \sigma
  +
  q^2 I
  =
  0,
\ee
which may be immediately seen to be equivalent to
\eqref{Uqgl21SkeinRelation}.

We shall evaluate our link invariant using the state model based on the
components of $\sigma$. The invariant would also be able to be
evaluated in some cases using \eqref{Uqgl21SkeinRelation}, although not
in all cases since it is of third order.  (The invariant may also be
directly evaluated for a class of links using quantum superalgebra
representation theoretic results \cite{GouldLinksZhang:96b}.)


\subsubsection{Construction of $R$ Itself?}

We have defined $\sigma=PR$, shortcutting the actual finding of $R$.
For good measure, we might obtain $R$ itself from $\sigma$ in the
following manner.  Observe that where $\ket{\Psi^k_j (q)}$ is a basis
for $V\otimes V$, as per (\ref{eq:V1KetBasis},\ref{eq:V3KetBasis}) and
$P$ is as per \eqref{Pdef}, we have:
\be
  P \ket{\Psi^k_j (q)}
  =
  \ket{\Psi^k_j (q^{-1})},
\ee
thus where:
\be
  \sigma
  =
  \sum_{k=1}^3
    \rho_k
    \sum_j
      \ket{\Psi^k_j (q)}
      \bra{\Psi^k_j (q)},
\ee
for the eigenvalues $\rho_1 = q^{-2\alpha}, \rho_2 = -1,
\rho_3 = q^{2\alpha+2}$, we have:
\be
  R
  \eq
  P \check{R}
  =
  P
  \sum_k
    \rho_k
    \sum_j
      \ket{\Psi^k_j (q)}
      \bra{\Psi^k_j (q)}
  =
  \sum_k
    \rho_k
    \sum_j
      P
      \ket{\Psi^k_j (q)}
      \bra{\Psi^k_j (q)}
  \\
  \eq
  \sum_k
    \rho_k
    \sum_j
      \ket{\Psi^k_j (q^{-1})}
      \bra{\Psi^k_j (q)}.
\ee

\cleardoublepage


\section{Explicit Construction of $\sigma$, $\sigma^{-1}$, $\Omega^\pm$
         and $\mho^\pm$}
\seclabel{RepresentationofSigma}
\addtocontents{toc}{\protect\vspace{-2.5ex}}

\subsection{The Origin of $\sigma$}
\addtocontents{toc}{\protect\vspace{-2.5ex}}

The explicit form of the representation of the braid generator $\sigma$
on the graded module $V\otimes V$ was presented in
\eqref{Explicitsigma}. The grading on $V$ means that this tensor
actually satisfies a graded version of the quantum Yang--Baxter
equation \eqref{QYB}, which looks similar, except that it includes
strings of parity factors. Our abstract tensor formalism for the
construction of link invariants also requires tensor multiplications
that do not include such parity factors. Thus, for our purposes, we
shall require a version of the above $\sigma$, as it appears after the
grading on $V$ is removed. We shall refer to that tensor here as
$\overline{\sigma}$, however it is still labelled ``$\sigma$'' as far
as the link invariant material is concerned \dots


\subsubsection{Components of Rank $4$ Graded Tensors}

Let $A$ represent an arbitrary graded rank $4$ tensor acting on
$V\otimes V$, then in terms of scalar coefficients $A^{i~k}_{j~l}$:
\bse
  A
  =
  \sum_{i j k l}
    A^{i~k}_{j~l}
    {e^{i k}}_{j l}
  =
  \sum_{i j k l}
    A^{i~k}_{j~l}
    \( {e^i}_j \otimes {e^k}_l \).
  \eqlabel{StarStar}
\ese
Note that we are representing the scalar coefficients as
$A^{i~k}_{j~l}$ rather than the more standard ${A^{i~k}}_{j~l}$; this
is for consistency with the way that the abstract tensors for the link
diagrams are presented.  This does not apply to the elementary
tensors.

Let us write $\{v^k\}$ as the (standard) basis for $V$, where of course
the $v^k$ satisfy:
\bse
  {e^i}_j v^k
  =
  \delta_{k j} v^i,
  \eqlabel{Star}
\ese
hence the action of $A$ on $V\otimes V$ is:
\bne
  A ( v^j \otimes v^l )
  & \stackrel{\eqref{StarStar}}{=} &
  \sum_{i k m n}
    A^{i~k}_{m~n}
    \( {e^i}_m \otimes {e^k}_n \)
    ( v^j \otimes v^l )
  \nonumber
  \\
  & \stackrel{\eqref{Botia105}}{=} &
  \sum_{i k m n}
    A^{i~k}_{m~n}
    {(-)}^{[j]([k]+[n])}
    \( {e^i}_m v^j \otimes {e^k}_n v^l \)
  \nonumber
  \\
  & \stackrel{\eqref{Star}}{=} &
  \sum_{i k m n}
    A^{i~k}_{m~n}
    {(-)}^{[j]([k]+[n])}
    \( \delta_{j m} v^i \otimes \delta_{l n} v^k \)
  \nonumber
  \\
  \eq
  \sum_{i k}
    A^{i~k}_{j~l}
    {(-)}^{[j]([k]+[l])}
    \( v^i \otimes v^k \).
  \eqlabel{Aaction}
\ene


\subsubsection{Aside: Components of $P$}

Specialising \eqref{Aaction} to the case where the operator $A$ is $P$,
we have:
\be
  \sum_{i k}
    P^{i~k}_{j~l}
    {(-)}^{[j]([k]+[l])}
    \( v^i \otimes v^k \)
  & \stackrel{\eqref{Aaction}}{=} &
  P
  \( v^j \otimes v^l \)
  \stackrel{\eqref{Pdef}}{=}
  {(-)}^{[j][l]}
  \( v^l \otimes v^j \)
  \\
  & \stackrel{\eqref{Star}}{=} &
  \sum_{i k}
    {(-)}^{[j][l]}
    \delta_{i l} \delta_{k j}
    \( v^i \otimes v^k \).
\ee
Equating components:
$
  P^{i~k}_{j~l}
  {(-)}^{[j]([k]+[l])}
  =
  {(-)}^{[j][l]}
  \delta_{i l} \delta_{k j}
$,
so that:
\bse
  P^{i~k}_{j~l}
  =
  {(-)}^{[j]}
  \delta_{i l} \delta_{k j}.
  \eqlabel{Pcomponents}
\ese
This expression simplifies; the only nonzero entries are:
\be
  P^{i~j}_{j~i}
  =
  {(-)}^{[j]},
  \qquad \; \;
  i, j = 1, \dots, 4.
\ee
From this we may reconstruct $P$:
\be
  P
  & \stackrel{\eqref{StarStar}}{=} &
  \sum_{i j k l}
    P^{i~k}_{j~l}
    \( {e^i}_j \otimes {e^k}_l \)
  \stackrel{\eqref{Pcomponents}}{=}
  \sum_{i j k l}
    {(-)}^{[j]}
    \delta_{i l} \delta_{k j}
    \( {e^i}_j \otimes {e^k}_l \)
  =
  \sum_{i j}
    {(-)}^{[j]}
    {e^{i j}}_{j i}
  \\
  \eq
  {e^{1 1}}_{1 1}
  -
  {e^{1 2}}_{2 1}
  -
  {e^{1 3}}_{3 1}
  +
  {e^{1 4}}_{4 1}
  +
  {e^{2 1}}_{1 2}
  -
  {e^{2 2}}_{2 2}
  -
  {e^{2 3}}_{3 2}
  +
  {e^{2 4}}_{4 2}
  \\
  & &
  +
  {e^{3 1}}_{1 3}
  -
  {e^{3 2}}_{2 3}
  -
  {e^{3 3}}_{3 3}
  +
  {e^{3 4}}_{4 3}
  +
  {e^{4 1}}_{1 4}
  -
  {e^{4 2}}_{2 4}
  -
  {e^{4 3}}_{3 4}
  +
  {e^{4 4}}_{4 4}.
\ee
A check on the consistency of our results is that from
\eqref{sigmadef}, we expect $\lim_{q\to1}\sigma=P$. Observing that
$\lim_{q\to1}Y=0$, we see that indeed this is the case.

\pagebreak


\subsubsection{Components After Removal of the Grading on $V$}

In \eqref{Aaction}, the parity factor is constructed from the degrees
of vectors; in the ungraded case, there would be no such factor, indeed
we would have:
\be
  \overline{A} (v^j \otimes v^l)
  \eq
  \sum_{i k}
    \overline{A}^{i~k}_{j~l}
    (v^i \otimes v^k).
\ee
Note that we are expressly \emph{not} dropping the grading on $A$,
although we shall refer to the new operator as $\overline{A}$.  This
motivates us to set:
\bse
  \overline{A}^{i~k}_{j~l}
  =
  {(-)}^{[j]([k]+[l])}
  A^{i~k}_{j~l}.
  \eqlabel{overlineA}
\ese
Removal of the grading on $V$ has the net effect of changing the sign
of one in four of the entries in the respective tensors.


\subsubsection{Explicit Components of $\overline{\sigma}$}

Applying \eqref{overlineA} to the explicit expression for $\sigma$
(i.e.  \eqref{Explicitsigma}) yields the ungraded edition of $\sigma$:
\be
  \overline{\sigma}
  \eq
    p^{-2}          {e^{1 1}}_{1 1}
  + (p^{-2} - 1)    {e^{2 1}}_{2 1}
  -                 {e^{2 2}}_{2 2}
  \\
  & &
  + (p^{-2} - 1)    {e^{3 1}}_{3 1}
  + (q^2 - 1)       {e^{3 2}}_{3 2}
  -                 {e^{3 3}}_{3 3}
  \\
  & &
  + Y^2             {e^{4 1}}_{4 1}
  + (p^2 q^2 - 1)   {e^{4 2}}_{4 2}
  + (p^2 q^2 - 1)   {e^{4 3}}_{4 3}
  + p^2 q^2         {e^{4 4}}_{4 4}
  \\
  & &
  + p^{-1} ({e^{2 1}}_{1 2} + {e^{1 2}}_{2 1})
  + p^{-1} ({e^{3 1}}_{1 3} + {e^{1 3}}_{3 1})
  \\
  & &
  +        ({e^{4 1}}_{1 4} + {e^{1 4}}_{4 1})
  - q      ({e^{3 2}}_{2 3} + {e^{2 3}}_{3 2})
  \\
  & &
  - Y      ({e^{4 1}}_{2 3} + {e^{2 3}}_{4 1})
  + q Y    ({e^{3 2}}_{4 1} + {e^{4 1}}_{3 2})
  \\
  & &
  + p q      ({e^{2 4}}_{4 2} + {e^{4 2}}_{2 4})
  + p q      ({e^{3 4}}_{4 3} + {e^{4 3}}_{3 4}),
\ee
which is the tensor that we shall use to represent the positive
crossings in the abstract tensor formalism.

\pagebreak


\subsubsection{Explicit Components of $\overline{\sigma}^{-1}$}

Similarly to \eqref{PRdef}, we may find ${\sigma}^{-1}$ from the
knowledge of the eigenvalues of $\sigma$. Where the eigenvalues of
$\sigma$ are $q^{-2 \alpha}$, $-1$ and $q^{2\alpha+2}$, the eigenvalues
of $\sigma^{-1}$ are $q^{2 \alpha}$, $-1$ and $q^{-2\alpha-2}$
respectively, thus:
\be
  \sigma^{-1}
  =
  q^{2 \alpha} P_1
  -
  P_2
  +
  q^{-2 \alpha - 2} P_3.
\ee
To get to $\overline{\sigma}^{-1}$, the version of $\sigma^{-1}$ after
the grading on the indices has been removed, we again invoke
\eqref{overlineA}. This yields:
\be
  {\overline{\sigma}}^{-1}
  \eq
    p^2                {e^{1 1}}_{1 1}
  + (p^2 - 1)          {e^{1 2}}_{1 2}
  + (p^2 - 1)          {e^{1 3}}_{1 3}
  + q^{-2}Y^2          {e^{1 4}}_{1 4}
  \\
  & &
  -                    {e^{2 2}}_{2 2}
  + (q^{-2}-1)         {e^{2 3}}_{2 3}
  + (p^{-2}q^{-2}-1)   {e^{2 4}}_{2 4}
  \\
  & &
  -                    {e^{3 3}}_{3 3}
  + (p^{-2}q^{-2}-1)   {e^{3 4}}_{3 4}
  + p^{-2}q^{-2}       {e^{4 4}}_{4 4}
  \\
  & &
  + p                 ({e^{2 1}}_{1 2} + {e^{1 2}}_{2 1})
  + p                 ({e^{3 1}}_{1 3} + {e^{1 3}}_{3 1})
  \\
  & &
  +                   ({e^{4 1}}_{1 4} + {e^{1 4}}_{4 1})
  - q^{-1}            ({e^{3 2}}_{2 3} + {e^{2 3}}_{3 2})
  \\
  & &
  - q^{-1}Y           ({e^{1 4}}_{3 2} + {e^{3 2}}_{1 4})
  + q^{-2}Y           ({e^{1 4}}_{2 3} + {e^{2 3}}_{1 4})
  \\
  & &
  + p^{-1}q^{-1}      ({e^{4 2}}_{2 4} + {e^{2 4}}_{4 2})
  + p^{-1}q^{-1}      ({e^{4 3}}_{3 4} + {e^{3 4}}_{4 3}),
\ee
which is the tensor that we shall use to represent the negative
crossings in the abstract tensor formalism.  Perhaps not surprisingly,
$\overline{\sigma}^{-1}$ also has $26$ nonzero components.


\subsubsection{Writing $\overline{\sigma}$ and $\overline{\sigma}^{-1}$
               as Matrices}

It may be easier on the eyes to view $\overline{\sigma}$ and
$\overline{\sigma}^{-1}$ in the more `graphical' form of matrices
rather than via lists of components.  Our rank $4$ tensors have $4$
indices, hence we cannot present them elegantly on the page in matrix
form as we can do with rank $2$ tensors. We can however, present a rank
$4$ tensor on the two dimensional medium of a page as an ordinary
matrix using certain standard conventions. (This really amounts to an
isomorphism.)

To begin with, we shall represent rank $2$ tensors as matrices, that
is, the elementary rank $2$ tensor $ {e^i}_j $ will be represented by
the elementary matrix $ e_{i,j} $.  Recall that our underlying carrier
space $V$ is $4$ dimensional. This means that all indices will run from
$1, \dots, 4$, i.e. we intend $ e_{i,j} $  to be a $ 4 \times 4 $
matrix.

We adopt the (usual) convention that the elementary rank
$4$ tensor $ {e^{i k}}_{j l}\defeq{e^i}_j \otimes {e^k}_l$ is
constructed by insertion of a copy of the elementary rank $2$ tensor
${e^k}_l$ at each location of ${e^i}_j$ (i.e.  each element of
${e^i}_j$ is multiplied by the whole of ${e^k}_l$).  This means that
${e^{i k}}_{j l}$ is represented by the elementary ($16 \times 16$)
matrix $e_{4(i-1)+k,4(j-1)+l}$.

Our convention then tells us that $ A^{i~k}_{j~l} $ appears as
the $(4(i-1)+k,4(j-1)+l)$ entry of
$\psi(A)$ (the matrix representing $A$), explicitly:
\bse
  {\psi(A)}_{4(i-1)+k,4(j-1)+l}
  =
  A^{i~k}_{j~l}.
  \eqlabel{matrixentryofA}
\ese

Using \eqref{matrixentryofA}, we evaluate the expressions for the
matrices $\psi(\overline{\sigma})$ and $\psi(\overline{\sigma}^{-1})$:
\def\scs{\scriptstyle}
\be
  \m{[}{@{}*{3}{*{3}{c@{}@{}}c|}*{3}{c@{}@{}}c@{}}
    {\scs p^{-2}} &.&.&.&.&.&.&.&.&.&.&.&.&.&.&.\\
    .&.&.&.&{\scs p^{-1}} &.&.&.&.&.&.&.&.&.&.&.\\
    .&.&.&.&.&.&.&.&{\scs p^{-1}} &.&.&.&.&.&.&.\\
    .&.&.&.&.&.&.&.&.&.&.&.&{\scs 1} &.&.&.\\
    \hline
    .&{\scs p^{-1}} &.&.&{\scs p^{-2}-1} &.&.&.&.&.&.&.&.&.&.&.\\
    .&.&.&.&.&{\scs -1} &.&.&.&.&.&.&.&.&.&.\\
    .&.&.&.&.&.&.&.&.&{\scs -q} &.&.&{\scs -Y} &.&.&.\\
    .&.&.&.&.&.&.&.&.&.&.&.&.&{\scs p q} &.&.\\
    \hline
    .&.&{\scs p^{-1}}&.&.&.&.&.&{\scs p^{-2}-1}&.&.&.&.&.&.&.\\
    .&.&.&.&.&.&{\scs -q}&.&.&{\scs q^2-1}&.&.&{\scs qY}&.&.&.\\
    .&.&.&.&.&.&.&.&.&.&{\scs -1}&.&.&.&.&.\\
    .&.&.&.&.&.&.&.&.&.&.&.&.&.&{\scs p q}&.\\
    \hline
    .&.&.&{\scs 1}&.&.&{\scs -Y}&.&.&{\scs qY}&.&.&{\scs Y^2}&.&.&.\\
    .&.&.&.&.&.&.&{\scs p q} &.&.&.&.&.&{\scs p^2 q^2 -1} &.&.\\
    .&.&.&.&.&.&.&.&.&.&.&{\scs p q} &.&.&{\scs p^2 q^2 -1} &.\\
    .&.&.&.&.&.&.&.&.&.&.&.&.&.&.&{\scs p^2 q^2}
  \me
\ee

\be
  \m{[}{@{}*{3}{*{3}{c@{}@{}}c|}*{3}{c@{}@{}}c@{}}
    {\scs p^2}&.&.&.&.&.&.&.&.&.&.&.&.&.&.&.\\
    .&{\scs p^2-1}&.&.&{\scs p}&.&.&.&.&.&.&.&.&.&.&.\\
    .&.&{\scs p^2-1}&.&.&.&.&.&{\scs p}&.&.&.&.&.&.&.\\
    .&.&.&{\scs q^{-2}Y^2}&.&.&{\scs q^{-2}Y}&.&.&
                {\scs -q^{-1}Y}&.&.&{\scs 1}&.&.&.\\
    \hline
    .&{\scs p}&.&.&.&.&.&.&.&.&.&.&.&.&.&.\\
    .&.&.&.&.&{\scs -1}&.&.&.&.&.&.&.&.&.&.\\
    .&.&.&{\scs q^{-2}Y}&.&.&{\scs q^{-2}-1}&.&.&{\scs -q^{-1}}&.
                &.&.&.&.&.\\
    .&.&.&.&.&.&.&{\scs p^{-2}q^{-2}-1}&.&.&.&.&.
                &{\scs p^{-1}q^{-1}}&.&.\\
    \hline
    .&.&{\scs p}&.&.&.&.&.&.&.&.&.&.&.&.&.\\
    .&.&.&{\scs -q^{-1}Y}&.&.&{\scs -q^{-1}}&.&.&.&.&.&.&.&.&.\\
    .&.&.&.&.&.&.&.&.&.&{\scs -1}&.&.&.&.&.\\
    .&.&.&.&.&.&.&.&.&.&.&{\scs p^{-2}q^{-2}-1}&.&.
                &{\scs p^{-1}q^{-1}}&.\\
    \hline
    .&.&.&{\scs 1}&.&.&.&.&.&.&.&.&.&.&.&.\\
    .&.&.&.&.&.&.&{\scs p^{-1}q^{-1}}&.&.&.&.&.&.&.&.\\
    .&.&.&.&.&.&.&.&.&.&.&{\scs p^{-1}q^{-1}}&.&.&.&.\\
    .&.&.&.&.&.&.&.&.&.&.&.&.&.&.&{\scs p^{-2}q^{-2}}
  \me.
\ee

\pagebreak


\subsection{Caps and Cups $ \Omega^\pm $ and $ \mho^\pm $}
\addtocontents{toc}{\protect\vspace{-2.5ex}}
\seclabel{CapsandCups}

We will use:
\be
  \Omega^+
  \eq
  \mho^+
  =
  I_4
  \\
  \Omega^-
  \eq
  \m{[}{c@{}c@{}c@{}c}
    {\scs q^{- 2 \alpha}} & . & . & . \\
    . & {\scs - q^{- 2 \alpha}} & . & . \\
    . & . & {\scs - q^{- 2 \( \alpha + 1 \)}} & . \\
    . & . & . & {\scs q^{- 2 \( \alpha + 1 \)}}
  \me
  =
  \m{[}{cccc}
    {\scs p^{-2}} & . & . & . \\
    . & {\scs - p^{-2}} & . & . \\
    . & . & {\scs - p^{-2} q^{-2}} & . \\
    . & . & . & {\scs p^{-2} q^{-2}}
  \me
  \\
  \qquad
  \mho^-
  \eq
  \m{[}{c@{}c@{}c@{}c}
    {\scs q^{2 \alpha}} & . & . & . \\
    . & {\scs - q^{2 \alpha}} & . & . \\
    . & . & {\scs - q^{2 \( \alpha + 1 \)}} & . \\
    . & . & . & {\scs q^{2 \( \alpha + 1 \)}}
  \me
  =
  \m{[}{cccc}
    {\scs p^2} &            . &                . & . \\
    .          & {\scs - p^2} &                . & . \\
    .          &            . & {\scs - p^2 q^2} & . \\
    .          &            . &                . & {\scs p^2 q^2}
  \me.
\ee
These choices for $ \Omega^\pm $ and $ \mho^\pm $ are sufficient
although not unique. To justify them, we begin with the comment that
the graphical consistency considerations of \figref{ReidemeisterI}
(invariance under the first Reidemeister move) require:
\be
  {(\Omega^\pm)}_{x~a} \cdot {(\mho^\pm)}^{a~y}
  =
  \delta^y_x
  =
  {(\mho^\pm)}^{y~a} \cdot {(\Omega^\pm)}_{a~x}.
\ee
This means that $\Omega^\pm={(\mho^\pm)}^{-1}$, so we need
only find $\mho^\pm$.

\begin{figure}[htbp]
  \begin{center}
    \input{graphics/ReidemeisterI.pstex_t}
    \caption{%
      Invariance under the first Reidemeister move.
    }
    \addtocontents{lof}{\protect\vspace{-2.5ex}}
    \figlabel{ReidemeisterI}
  \end{center}
\end{figure}

\pagebreak


\subsubsection{Choosing $\mho^+$ and $\Omega^+$}

$\mho^+$ may be chosen from graphical consistency considerations
involving the quarter-turn rotations in \figref{XLandXR}. The choice:
\be
  {\( \mho^+ \)}^{a~b}
  =
  \delta^a_b
\ee
i.e. $\mho^+ = I_4$,
ensures that the definition \eqref{XRdef}, i.e.
\be
  {\( X_r \)}^{a~c}_{b~d}
  \defeq
  {(X)}^{c~g}_{f~b}
  \cdot
  {\( \mho^+ \)}^{a~f}
  \cdot
  {\( \Omega^+ \)}_{g~d}
\ee
(where $ X $ is either $ \sigma $ or $ \sigma^{-1} $),
simplifies to the elegant form:
$
  {\( X_r \)}^{a~c}_{b~d}
  =
  {(X)}^{c~d}_{a~b}
$,
and indeed we also have
$
  {\( X_d \)}^{a~c}_{b~d}
  =
  {(X)}^{d~b}_{c~a}
$.
Clearly, we then have $\Omega^+=I_4$, viz
\bse
  {(\Omega^+)}_{a~b}
  =
  \delta^a_b.
  \eqlabel{ComponentsofOmegaplus}
\ese


\subsubsection{Choosing $\mho^-$ and $\Omega^-$}

For the choice of $\mho^-$, we begin by imposing the
graphical consistency requirement of \figref{Loop} (invariance under
ambient isotopy), i.e. \eqref{AmbientIsotopy}:
\be
  \sum_{abc}
    {(\sigma)}^{y~a}_{x~b}
    \cdot
    {(\Omega^+)}_{a~c}
    \cdot
    {(\mho^-)}^{b~c}
  =
  \delta^y_x.
\ee
The choice of $\Omega^+$ in \eqref{ComponentsofOmegaplus}
simplifies this to:
\bse
  \sum_{ab}
    {(\sigma)}^{y~a}_{x~b}
    \cdot
    {(\mho^-)}^{b~a}
  =
  \delta^y_x.
  \eqlabel{LoopRemovalConsistencySimplified}
\ese
We also invoke the following result
\cite[Lemma~2,~p354]{LinksGouldZhang:93}
(and see also \cite{LinksGould:92b}):
\bse
  \( I \otimes \mathrm{str} \)
  [ ( I \otimes q^{2 h_\rho} ) \sigma ]
  =
  K I,
  \eqlabel{ResultFromLGZ93}
\ese
for some constant $K$ depending on the normalisation of $\sigma$, where
$ \mathrm{str} $ denotes the supertrace. By expanding the LHS of
\eqref{ResultFromLGZ93}, we shall find an expression for $\mho^-$
sufficient to satisfy \eqref{LoopRemovalConsistencySimplified}. In a
sense, the result \eqref{ResultFromLGZ93} amounts to satisfying the
requirement that our invariant is an invariant of ambient isotopy.

To summarise, we shall find:
\bse
  {\( \mho^- \)}^{b a}
  =
  {\( - \)}^{\[ b \]}
  \pi {( q^{-2 h_\rho} )}^b_a,
  \eqlabel{mho-defn}
\ese
which, as $\Omega^- = {(\mho^-)}^{-1}$, yields:
\be
  {\( \Omega^- \)}_{b a}
  =
  {\( - \)}^{\[ b \]}
  \pi {( q^{2 h_\rho} )}^b_a,
\ee
so that we shall require the evaluation of the components
${\pi(q^{\pm 2 h_\rho})}^b_a$.

To begin the expansion of the LHS of \eqref{ResultFromLGZ93}, we have:
\be
  \( I \otimes q^{-2 h_{\rho}} \) \sigma
  \eq
  \sum_{ijkl}
    {\[ \( I \otimes q^{-2 h_{\rho}} \) \sigma \]}^{i~j}_{k~l}
    \( {e^i}_k \otimes {e^j}_l \),
\ee
so we must deduce the form of
$
  {\[ \( I \otimes q^{-2 h_{\rho}} \) \sigma \]}^{i~j}_{k~l}
$.
To this end, we have:
\be
  I \otimes q^{-2 h_{\rho}}
  \eq
  \(
    \sum_{im}
      {\delta}^i_m
      {e^i}_m
  \)
  \otimes
  \(
    \sum_{jn}
      {\(q^{-2 h_{\rho}}\)}^j_n
      {e^j}_n
  \)
  \\
  \eq
  \sum_{ijmn}
    {\delta}^i_m
    {\(q^{-2 h_{\rho}}\)}^j_n
    \( {e^i}_m \otimes {e^j}_n \),
\ee
i.e. in components:
\bse
  {(I \otimes q^{-2 h_{\rho}})}^{i~j}_{m~n}
  =
  {\delta}^i_m
  {\(q^{-2 h_{\rho}}\)}^j_n.
  \eqlabel{Iotimesq-2hrhocomponents}
\ese
Also:
\bse
  \sigma
  =
  \sum_{mnkl}
    {(\sigma)}^{m~n}_{k~l}
    \( {e^m}_k \otimes {e^n}_l \),
  \eqlabel{sigma}
\ese

\pagebreak

Let $A$ and $B$ be general (ungraded) rank $4$ tensors, viz
$
  A = \sum_{ijmn} A^{i~j}_{m~n} \( {e^i}_m \otimes {e^j}_n \)
$
and
$
  B = \sum_{mnkl} B^{m~n}_{k~l} \( {e^m}_k \otimes {e^n}_l \)
$,
we have:
\be
  A B
  \eq
  \sum_{ijklmn}
    A^{i~j}_{m~n}
    B^{m~n}_{k~l}
    \( {e^i}_m \otimes {e^j}_n \)
    \( {e^m}_k \otimes {e^n}_l \)
  \\
  \eq
  \sum_{ijklmn}
    A^{i~j}_{m~n}
    B^{m~n}_{k~l}
    \( {e^i}_m {e^m}_k \otimes {e^j}_n {e^n}_l \)
  \\
  \eq
  \sum_{ijklmn}
    A^{i~j}_{m~n}
    B^{m~n}_{k~l}
    \( {e^i}_k \otimes {e^j}_l \),
\ee
as, for \emph{each} $m$ (i.e. without summation over $m$), we have
${e^i}_m{e^m}_k={e^i}_k$. Thus:
\bse
  {\( A B \)}^{i~j}_{k~l}
  =
  \sum_{mn}
    A^{i~j}_{m~n}
    B^{m~n}_{k~l}.
  \eqlabel{ABcomponents}
\ese
Where $A$ is
$I \otimes q^{-2 h_{\rho}}$ and $B$ is $\sigma$, we have the components:
\bse
  {\[
    \( I \otimes q^{-2 h_{\rho}} \)
    \sigma
  \]}^{i~j}_{k~l}
  \stackrel{(%
    \ref{eq:Iotimesq-2hrhocomponents},%
    \ref{eq:sigma},%
    \ref{eq:ABcomponents}%
  )}{=}
  \sum_{mn}
    \delta^i_m
    {(q^{-2 h_{\rho}})}^j_n
    {(\sigma)}^{m~n}_{k~l}.
  \eqlabel{ItensorqSigmacomponents}
\ese

If $M$ is a rank $2$ tensor (i.e. a matrix), viz:
$M=\sum_{ij}M^i_j {e^i}_j$, then we have:
\be
  \str (M)
  \defeq
  \sum_j
    {(-)}^{[j]}
    M^j_j,
\ee
hence if $A$ is a rank $4$ tensor, viz
$A=\sum_{ijkl}A^{i~j}_{k~l} ({e^i}_k \otimes {e^j}_l)$:
\be
  (I \otimes \str) (A)
  =
  \sum_{ijk}
    {(-)}^{[j]}
    A^{i~j}_{k~j}
    {e^i}_k,
\ee
i.e. componentwise:
\bse
  {\[ (I \otimes \str) (A) \]}^i_k
  =
  \sum_{j}
    {(-)}^{[j]}
    A^{i~j}_{k~j}.
  \eqlabel{ITensorSupertraceonAcomponents}
\ese

\pagebreak

Applying the definition of the action of $I\otimes\str$ to the
expression $\(I\otimes q^{-2 h_{\rho}} \) \sigma$:
\be
  {\left\{
    \(I \otimes \str\)
    \[ \( I \otimes q^{-2 h_{\rho}} \) \sigma \]
  \right\}}^y_x
  & \stackrel{\eqref{ITensorSupertraceonAcomponents}}{=} &
  \sum_b
    {(-)}^{[b]}
    {\[ \( I \otimes q^{-2 h_{\rho}} \) \sigma \]}^{y~b}_{x~b}
  \\
  & \stackrel{\eqref{ItensorqSigmacomponents}}{=} &
  \sum_{abm}
    {(-)}^{[b]}
    \delta^y_m
    {(q^{-2 h_{\rho}})}^b_a
    {(\sigma)}^{m~a}_{x~b}
  \\
  \eq
  \sum_{ab}
    {(-)}^{[b]}
    {(q^{-2 h_{\rho}})}^b_a
    {(\sigma)}^{y~a}_{x~b}%
    .
\ee
Now, by \eqref{ResultFromLGZ93}, we have:
$
  {\left\{
    \(I \otimes \str\)
    \[ \( I \otimes q^{-2 h_{\rho}} \) \sigma \]
  \right\}}^y_x
  =
  K
  \delta^y_x
$, for some constant $K$, that is we have:
\bse
  \sum_{ab}
    {(-)}^{[b]}
    {(q^{-2 h_{\rho}})}^b_a
    {(\sigma)}^{y~a}_{x~b}
  =
  K
  \delta^y_x.
  \eqlabel{Drosophilia}
\ese
Thus, we have $K=1$. Choosing:
\bse
  {(\mho^-)}^{b~a}
  =
  {(-)}^{[b]}
  {(q^{-2 h_{\rho}})}^b_a
  \eqlabel{mho-componentsintermsofq-2hrho}
\ese
converts \eqref{Drosophilia} to
\eqref{LoopRemovalConsistencySimplified}, and we have thus justified
\eqref{mho-defn}.


\subsubsection{Evaluation of the Components ${(q^{-2 h_{\rho}})}^b_a$}

Now we evaluate the components ${(q^{-2 h_{\rho}})}^b_a$. To
this end, let $H$ be the Cartan subalgebra of $ gl(2|1)$, with dual
the $gl(2|1)$ root space $H^*$. A basis for $H^*$ is:
\be
  \{
    \varepsilon_1 = \(1,0\,|\,0\),
    \varepsilon_2 = \(0,1\,|\,0\),
    \varepsilon_3 = \(0,0\,|\,1\)
  \},
\ee
on which we have an invariant bilinear form
$ ( \cdot, \cdot ) : H^* \times H^* \to \BC $ defined by:
\bse
  (\varepsilon_i, \varepsilon_j)
  \defeq
  {(-)}^{[i]}
  \delta_{i j}.
  \eqlabel{InvariantBilinearForm}
\ese

Now, $gl(2|1)$ has the two simple positive roots
$\alpha_1=(1,-1\,|\,0)=\varepsilon_1-\varepsilon_2$ (even) and
$\alpha_0=(0,1\,|\,-1)=\varepsilon_2-\varepsilon_3$ (odd), and another,
nonsimple, odd positive root
$\alpha_1+\alpha_0=(1,0\,|\,-1)=\varepsilon_1-\varepsilon_3$.  From
these we have the half-sums of all even and all odd roots
$\rho_0=(\frac{1}{2},-\frac{1}{2}\,|\,0)$ and
$\rho_1=(\frac{1}{2},\frac{1}{2}\,|\,-1)$ respectively, hence their
\emph{graded} half-sum $\rho\defeq\rho_0-\rho_1$ is:
\be
  \rho
  =
  (0,-1\,|\,1)
  =
  -\varepsilon_2 + \varepsilon_3,
\ee
from which we may proceed to determine the associated Cartan element
$h_{\rho}$.

\pagebreak

Firstly, as $H$ is a dual vector space to $H^*$, we may implicitly
define generators $\{{E^i}_i\}_{i=1}^3$ of $H$ (i.e. $gl(2|1)$ Cartan
generators) by the action:
\bse
  {E^i}_i ( \varepsilon_j )
  \defeq
  \delta_{i j},
  \qquad
  i, j = 1, 2, 3.
  \eqlabel{HactiononHstar}
\ese
Next, we define $ h_{\rho} \in H $ to be such that:
\bse
  h_{\rho}
  ( \alpha_i)
  =
  (\rho, \alpha_i),
  \qquad
  i = 0, 1.
  \eqlabel{hrhodefn}
\ese
Examination of the definition of $\rho$ shows that for the two simple
roots $ \alpha_0 = \varepsilon_2 - \varepsilon_3 $ and
$ \alpha_1 = \varepsilon_1 - \varepsilon_2 $, we have:
\be
  (\rho, \alpha_i)
  =
  {\textstyle \frac{1}{2}}
  (\alpha_i, \alpha_i),
  \qquad
  i = 0, 1,
\ee
thus \eqref{hrhodefn} may be more usefully written as:
\bse
  h_{\rho}
  ( \alpha_i)
  =
  {\textstyle \frac{1}{2}}
  (\alpha_i, \alpha_i),
  \qquad
  i = 0, 1.
  \eqlabel{hrhousefuldefn}
\ese
For $gl(2|1)$ (as opposed to $sl(2|1)$, for which it is unique),
$h_{\rho}$ is in fact only determined up to addition of multiples of
the first order Casimir $C={E^1}_1+{E^2}_2+{E^3}_3$, which satisfies
$C(\alpha_i)=0$ for $i = 0,1$.  Thus, in setting
$
  h_{\rho}
  =
  \beta_1 {E^1}_1
  +
  \beta_2 {E^2}_2
  +
  \beta_3 {E^3}_3
$,
for complex coefficients $ \beta_i $ to be determined by
\eqref{hrhousefuldefn}; we may choose $\beta_1 = 0$. Thus, using:
\be
  h_{\rho}
  (\varepsilon_j)
  =
  \sum_i
    \beta_i
    {E^i}_i
  (\varepsilon_j)
  \stackrel{\eqref{HactiononHstar}}{=}
  \sum_i
    \beta_i
    \delta_{i j}
  =
  \beta_j,
  \qquad
  i = 1, 2, 3,
\ee
we may write:
\be
  h_{\rho}
  (\alpha_0)
  \eq
  h_{\rho}
  (\varepsilon_2)
  -
  h_{\rho}
  (\varepsilon_3)
  =
  \beta_2 - \beta_3
  \\
  h_{\rho}
  (\alpha_1)
  \eq
  h_{\rho}
  (\varepsilon_1)
  -
  h_{\rho}
  (\varepsilon_2)
  =
  \beta_1 - \beta_2.
\ee
This information may be combined with data obtained by expansion of
\eqref{hrhousefuldefn}:
\be
  h_{\rho}
  (\alpha_0)
  \eq
  {\textstyle \frac{1}{2}}
  (\alpha_0, \alpha_0)
  =
  {\textstyle \frac{1}{2}}
  (\varepsilon_2 - \varepsilon_3, \varepsilon_2 - \varepsilon_3)
  \stackrel{\eqref{InvariantBilinearForm}}{=}
  {\textstyle \frac{1}{2}}
  (1 - 1)
  =
  0
  \\
  h_{\rho}
  (\alpha_1)
  \eq
  {\textstyle \frac{1}{2}}
  (\alpha_1, \alpha_1)
  =
  {\textstyle \frac{1}{2}}
  (\varepsilon_1 - \varepsilon_2, \varepsilon_1 - \varepsilon_2)
  \stackrel{\eqref{InvariantBilinearForm}}{=}
  {\textstyle \frac{1}{2}}
  (1 + 1)
  =
  1
\ee
to yield a linear system for the $\beta_i$:
\be
  \beta_2 - \beta_3
  \eq
  0
  \\
  \beta_1 - \beta_2
  \eq
  1.
\ee
Thus, we have:
$
  \beta_2
  =
  \beta_3
  =
  \beta_1 - 1 = -1
$,
and so, as $ \beta_1 = 0 $:
\be
  h_{\rho}
  =
  - {E^2}_2 - {E^3}_3.
\ee
Note that this is expressly \emph{not} $-{E^2}_2+{E^3}_3$ as might be
na\"{\i}vely assumed from inspection of $\rho$.

The same $h_{\rho}$ holds for $U_q[gl(2|1)]$, hence we have:
\be
  q^{-2 h_{\rho}}
  =
  q^{2 {E^2}_2 + 2 {E^3}_3},
\ee
and:
\be
  \pi \( q^{-2 h_{\rho}} \)
  \eq
  \pi ( q^{2 {E^2}_2} )
  \pi ( q^{2 {E^3}_3} )
  \\
  & \stackrel{(\ref{eq:TheRepresentationExtended}b,c)}{=} &
  \[
    \ketbra{1}{1}
    +
    q^{-2}
    \ketbra{2}{2}
    +
    \ketbra{3}{3}
    +
    q^{-2}
    \ketbra{4}{4}
  \]
  \cdot
  \\
  & & \quad
  \[
    q^{2 \alpha}
    \ketbra{1}{1}
    +
    q^{2 (\alpha+1)}
    \ketbra{2}{2}
    +
    q^{2 (\alpha+1)}
    \ketbra{3}{3}
    +
    q^{2 (\alpha+2)}
    \ketbra{4}{4}
  \]
  \\
  \eq
    q^{2 \alpha}
    \ketbra{1}{1}
    +
    q^{2 \alpha}
    \ketbra{2}{2}
    +
    q^{2 (\alpha+1)}
    \ketbra{3}{3}
    +
    q^{2 (\alpha+1)}
    \ketbra{4}{4},
\ee
i.e., expressed as a matrix:
\bse
  \pi (-q^{2 h_\rho} )
  =
  \m{[}{cccc}
    q^{2\alpha} &           . &               . &              . \\
              . & q^{2\alpha} &               . &              . \\
              . &           . & q^{2(\alpha+1)} &              . \\
              . &           . &               . & q^{2(\alpha+1)}
  \me
  =
  \m{[}{cccc}
    p^{2} &      . &           . &            . \\
         . & p^{2} &           . &            . \\
         . &     . & p^{2} q^{2} &           . \\
         . &     . &           . & p^{2} q^{2}
  \me.
  \eqlabel{matrixofq-2hrho}
\ese

Note that this result disagrees slightly with that presented in the
published work \cite{DeWitKauffmanLinks:98}, where the definitions of
$\ket{2}$ and $\ket{3}$ are interchanged with respect to the above.
This disagreement will not affect the evaluation of the invariant.

Lastly, \eqref{matrixofq-2hrho} allows us to use
\eqref{mho-componentsintermsofq-2hrho} to state:
\be
  \mho^-
  =
  \m{[}{cccc}
    p^{2} &       . &             . &          . \\
        . & - p^{2} &             . &          . \\
        . &       . & - p^{2} q^{2} &          . \\
        . &       . &             . & p^{2} q^{2}
  \me,
\ee
and hence $\Omega^- = (\mho^-)^{-1}$ is:
\be
  \Omega^-
  =
  \m{[}{cccc}
    p^{-2} &        . &               . &            . \\
         . & - p^{-2} &               . &            . \\
         . &        . & - p^{-2} q^{-2} &            . \\
         . &        . &               . & p^{-2} q^{-2}
  \me.
\ee

\cleardoublepage


\section{Implementation in \protect\textsc{Mathematica}}
\addtocontents{toc}{\protect\vspace{-2.5ex}}

\subsection{Construction of the Braid Generator}
\addtocontents{toc}{\protect\vspace{-2.5ex}}

The program \texttt{Uqgl21.m} (see \secref{Uqgl21.m}) consists of a
suite of functions that evaluate the explicit representation of
$\sigma$, starting from the basis vectors $\ket{\Psi^k_j}$ of $V\otimes
V$. It also contains auxiliary functions required to implement
noncommutativity and a few other functions to construct objects related
to $\sigma$.

The code is self-explanatory. To use it, the program must first be
imported into a \textsc{Mathematica} session via a command line such as
\texttt{Get["Uqgl21.m"]}. After the functions are loaded, typing
\texttt{sigma} will evaluate $\sigma$, and typing
\texttt{ungradedsigma} will evaluate $\overline{\sigma}$.  Other
functions behave similarly, and are named in an obvious fashion.
Perusal of the code listing in \secref{Uqgl21.m} will confirm what
functions are available.  Note that once $\sigma$ has been evaluated
within a given session, \textsc{Mathematica} remembers its value, and
when called upon again, it will not have to be recomputed.

Of some interest is the code \texttt{CheckRu} that computes
$\check{R}(u)$. Whilst the output from that function is correct,
considerable effort went into converting screen-based output to the
human readable form of \eqref{ExplicitCheckRu}.

\pagebreak


\subsection{Evaluation of the Invariant}
\addtocontents{toc}{\protect\vspace{-2.5ex}}

The program \texttt{LinksGouldInvariant.m} (see
\secref{LinksGouldInvariant.m}) consists of a suite of functions that
evaluate the Links--Gould invariant for various links, starting from
$\overline{\sigma}$ (computed by a function within \texttt{Uqgl21.m})
and a set of abstract tensors. It has been designed to compute the
Links--Gould invariant and the (unnormalised) bracket polynomial (and
thence the Jones polynomial). As it depends on \texttt{Uqgl21.m}, both
that file and \texttt{LinksGouldInvariant.m} must be imported.

The command \texttt{LinkInvariant[LinksGould, Link]} returns the
polynomial invariant corresponding to the link \texttt{Link}, which may
be a pretzel \texttt{Pretzel[pp, qq, rr]} or any or all of the links in
\texttt{KnotList}, which is:
\be
  \begin{array}{l}
    \mathtt{Unknot}, \mathtt{HopfLink}, \mathtt{Trefoil},
    \mathtt{FigureEight},
    \\
    \mathtt{Cinquefoil}, \mathtt{FiveTwo}, \mathtt{WhiteheadLink},
    \\
    \mathtt{SixOne}, \mathtt{SixTwo}, \mathtt{SixThree},
    \\
    \mathtt{Septfoil}, \mathtt{SevenTwo},
    \\
    \mathtt{EightSeventeen}, \mathtt{NineFortyTwo},
    \mathtt{TenFortyEight}, \mathtt{KT}, \mathtt{KTI}.
  \end{array}
\ee
The program also contains several auxiliary functions and a few other
functions to evaluate the properties of the invariant. In particular,
we have:
\begin{itemize}
\item
  \texttt{YangBaxterChecker[Invariant]} asks if the $\sigma$ being used
  for the invariant satisfies the quantum Yang--Baxter equation.

\item
  \texttt{LGCantDetectNoninvertibilityQ[Link]} answers the question:
  ``Is is that $LG$ cannot detect the noninvertibility of \texttt{Link}
  (if it is noninvertible)?'' by checking whether $LG$ for
  \texttt{Link} displays the symmetry that denies it the ability to
  distinguish \texttt{Link} from its inverse.
  \texttt{Link} must be an element of
  \texttt{KnotList}, or be a \texttt{Pretzel[pp, qq, rr]}.

\item
  \texttt{LGDetectsChiralityofLinkQ} answers the question:  ``Does $LG$
  detect the chirality of \texttt{Link} (assuming it is chiral)?'' by
  checking whether $LG$ for \texttt{Link} is not palindromic.  Again,
  \texttt{Link} must be an element of \texttt{KnotList}, or be a
  \texttt{Pretzel[pp, qq, rr]}.

\item
  \texttt{ListofPretzels[N]} returns a list of the noninvertible pretzel
  knots $(p,q,r)$ for all odd, distinct combinations of
  $3\leqslant p<q<r\leqslant N$ (odd). 

\item
  \texttt{LGCantDetectNoninvertibilityofPretzelsQ[N]} returns a list of
  answers to \texttt{LGCantDetectNoninvertibilityQ[Pretzel[pp,qq,rr]]}
  for all the pretzels \texttt{Pretzel[pp,qq,rr]} in
  \texttt{ListofPretzels[N]}.

\item
  \texttt{LGDetectsChiralityofPretelsQ[N]} returns a list of answers to
  the question
  \texttt{LGDetectsChiralityofLinkQ[Pretzel[pp,qq,rr]]}, for all the
  pretzels \texttt{Pretzel[pp,qq,rr]} in \texttt{ListofPretzels[N]}.

\item
  \texttt{JonesPolynomial[Link]} evaluates the Jones polynomial
  for the link in question.
\end{itemize}

\cleardoublepage


\section{Experimental Results and Analysis}
\addtocontents{toc}{\protect\vspace{-2.5ex}}

\subsection{Analysis of the Behaviour of the Invariant}
\addtocontents{toc}{\protect\vspace{-2.5ex}}

Before proceeding to evaluate the invariant for our test examples,
we digress to consider its properties so as to be able to
make predictions about the evaluations. We are interested in
whether $LG$ distinguishes the chirality of links, the noninvertibility
of (single-component) knots and mutations of knots.


\subsubsection{Chirality: The Effect of Reflection on $LG_K$}

Say that we have a link $K$ with associated abstract tensor $T_{K}$ and
invariant $LG_{K}$. If we replace $K$ with its reflection $K^*$, every
positive (respectively negative) crossing in $K$ will have been
replaced by the equivalent negative (respectively positive) crossing in
$K^*$, thus, $T_{K^*}$ will be the same as $T_{K}$ except for the
interchange $\sigma \leftrightarrow \sigma^{-1}$.  The caps
$\Omega^\pm$ and cups $\mho^\pm$ will remain unchanged, as will the
directions of the arrows (recall that chirality is independent of
orientation).

From the uniqueness of the universal R matrix for any quantum
(super)algebra \cite{KhoroshkinTolstoy:91}, the following relation
holds (for appropriate normalisation):
\be
  R^{-1} (q)
  =
  R (q^{-1}),
\ee
hence, as $ \sigma = P R $, we have:
\be
  \sigma^{-1} (q)
  =
  P \sigma (q^{-1}) P.
\ee
(In the case of the Links--Gould invariant, $R$ of course has the extra
variable $\alpha$, but this doesn't affect the above.) Thus, up
to a basis transformation, $\sigma$ and $\sigma^{-1}$ are
interchangeable by the change of variable $q \mapsto q^{-1}$. It then
follows that the invariant for $K^*$ is obtainable from that of $K$ by
the same change of variable, thus we have proven:

\begin{proposition}
  If $K$ is amphichiral then the polynomial ${LG}_K(q)$ is palindromic,
  i.e. invariant under the mapping $q \mapsto q^{-1}$ (and hence
  $p\mapsto p^{-1}$).
  \proplabel{AmphichiralMeansPalindromic}
\end{proposition}

\pagebreak

From the polynomial for a link $K$, we may immediately write down the
polynomial for $K^*$:
\bse
  {LG}_{K^*} (q,p)
  =
  {LG}_K \( q^{-1}, p^{-1} \)
  \eqlabel{Utricularia1}
\ese

As we have:
\be
  K = K^*
  \quad \to \quad
  {LG}_K (q,p)
  =
  {LG}_{K^*} (q,p),
\ee
then we have, conversely, that:
\bse
  {LG}_K (q,p)
  \neq
  {LG}_{K^*} (q,p)
  \quad \to \quad
  K \neq K^*,
  \eqlabel{Utricularia3}
\ese
i.e. if the polynomials corresponding to $K$ and $ K^* $ are
distinct, then $K$ must be chiral. Using the identity
\eqref{Utricularia1}, the test of \eqref{Utricularia3} becomes:
\be
  {LG}_K (q,p)
  \neq
  {LG}_K \( q^{-1}, p^{-1} \)
  \quad \to \quad
  K \neq K^*,
\ee
hence we have demonstrated:

\begin{proposition}
  If ${LG}_K(q,p)$ is \emph{not} palindromic, then $LG$ detects the
  chirality of $K$.
  \proplabel{NotPalindromicImpliesChiral}
\end{proposition}

\pagebreak


\subsubsection{The Effect of Inversion on $LG_K$}

Replacing $K$ with $K^{-1}$ amounts to reversing every arrow in $K$.
$T_{K^{-1}}$ will thus be the same as $T_{K}$ except for the following
changes, best expressed in terms of the auxiliary tensors.  For the
crossings, where $X$ is either $\sigma$ or $\sigma^{-1}$, interchange
$X$ with $X_d$ and $X_l$ with $X_r$; and for the caps and cups,
interchange only the signs, i.e.  interchange $\Omega^\pm$ with
$\Omega^\mp$ and $\mho^\pm$ with $\mho^\mp$.

This has the effect that $T_{K}$ is replaced by $T_{K^{-1}}$, the
dual tensor acting on the dual space \cite{ReshetikhinTuraev:90}.
Recalling that the tensors representing $(1,1)$ tangles act as scalar
multiples of the identity on $V$, the dual tensor has exactly the same
form, thus:

\begin{proposition}
  A knot invariant derived from an irreducible representation of a
  quantum (super)algebra is unable to detect inversion.
  \proplabel{CantDetectInvertibility}
\end{proposition}

Now observe that the representation of $U_q[gl(2|1)]$ acting on the
dual module $V^*$ is given by the replacement
$\alpha \mapsto -(\alpha+1)$ (with an appropriate redefinition of the
Cartan elements).  Thus for a given $(1,1)$ tangle $K$, with invariant
${LG}_K (q,p)$, the invariant ${LG}_{K^{-1}}$ of its inverse $K^{-1}$
is obtained as:
\be
  {LG}_{K^{-1}} (q,p)
  =
  {LG}_K (q,q^{-1} p^{-1}).
\ee
However, in view of \propref{CantDetectInvertibility}, such an
invariant is unable to detect inversion, hence we must have
${LG}_{K^{-1}}={LG}_{K}$, from which it follows that:

\begin{proposition}
  $ {LG}_K (q,p) $ enjoys the symmetry property:
  \bse
    {LG}_K (q,p) = {LG}_K \( q, q^{-1} p^{-1} \).
    \eqlabel{Melanotaenia}
  \ese
  \proplabel{ASymmetryofLG}
\end{proposition}


\subsubsection{The Effect of Mutation on $LG_K$}

The question of whether our invariant is able to distinguish mutants is
answerable in the negative. Theorem 5 of \cite{MortonCromwell:96}
states that (for quantum algebras) if the modules occurring in the
decomposition of $V \otimes V$ each have unit multiplicity, as indeed
\eqref{TPdecomp} shows in our case, then the invariant is unable to
detect mutations.  Whilst this was proved in \cite{MortonCromwell:96}
for the case of quantum algebras, the extension to the case of quantum
superalgebras is quite straightforward. We make a proposition out of
it:

\begin{proposition}
  The Links--Gould invariant is unable to distinguish between mutants.
  \proplabel{CantDistinguishMutants}
\end{proposition}

\pagebreak


\subsection{Invariant Evaluations for the Ordinary Links}
\addtocontents{toc}{\protect\vspace{-2.5ex}}

Evaluations of $LG_K(q,p)$ for our collection of ordinary links are
presented in \tabref{gl21oriented}. We have, as intended, that
$LG_{0_1}=1$.  Observe that $LG_{4_1}$ and $LG_{6_3}$ are palindromic
as $4_1$ and $6_3$ are amphichiral, whilst $LG_K$ is not palindromic
for any other links as the rest are chiral.  Observe also that the
invariant is actually a polynomial in $q^2$ (and hence $p^2$) rather
than just $q$ (and $p$), thus we might write
$LG_K(Q \defeq q^2, P \defeq p^2)$ rather than $LG_K(q,p)$.

\begin{table}[ht]
  \centering
  \begin{tabular}{||c||l||}
    \hline\hline
    & \\[-3mm]
    $K$ &
    \multicolumn{1}{c||}{$ {LG}_K (q,p) $}
    \\[1mm]
    \hline\hline
    & \\[-2mm]
    $ 2^2_1 $ &
    $ -1 + p^{-2} - q^{2} + p^{2} q^{2} $
    \\[2mm]
    \hline
    & \\[-2mm]
    $ 3_1 $ &
    $
      1 + p^{-4} - p^{-2} + 2 q^{2} - p^{-2} q^{2} - p^{2} q^{2}
      - p^{2} q^{4} + p^{4} q^{4}
    $
    \\[2mm]
    \hline
    & \\[-2mm]
    $ 4_1 $ &
    $
        7
      +   \( p^{-4} q^{-2} + p^{4} q^{2} \)
      - 3 \( p^{-2}        + p^{2}     \)
      - 3 \( p^{-2} q^{-2} + p^{2} q^{2} \)
      + 2 \(        q^{-2} +       q^{2} \)
    $
    \\[2mm]
    \hline
    & \\[-2mm]
    $ 5_1 $ &
    $
      \begin{array}{l}
          1
        +    p^{-8}
        -    p^{-6}
        +    p^{-4}
        -    p^{-2}
        + 2         q^{2}
        -    p^{-6} q^{2}
        + 2  p^{-4} q^{2}
        - 2  p^{-2} q^{2}
        -    p^{2}  q^{2}
        + 2         q^{4}
        \\
        \quad
        -    p^{-2} q^{4}
        - 2  p^{2}  q^{4}
        +    p^{4}  q^{4}
        -    p^{2}  q^{6}
        + 2  p^{4}  q^{6}
        -    p^{6}  q^{6}
        -    p^{6}  q^{8}
        +    p^{8}  q^{8}
      \end{array}
    $
    \\[3mm]
    \hline
    & \\[-2mm]
    $ 5_2 $ &
    $
      \begin{array}{l}
          3
        + 3  p^{-4}
        - 5  p^{2}
        + 10        q^{2}
        +    p^{-4} q^{2}
        - 6  p^{-2} q^{2}
        - 5  p^{2}  q^{2}
        + 4         q^{4}
        -    p^{-2} q^{4}
        - 6  p^{2}  q^{4}
        \\
        \quad
        + 3  p^{4}  q^{4}
        -    p^{2}  q^{6}
        +    p^{4}  q^{6}
      \end{array}
    $
    \\[3mm]
    \hline
    & \\[-2mm]
    $ 5^2_1 $ &
    $
      \begin{array}{l}
        - 10
        +   p^{-6}  q^{-2}
        - 3 p^{-4}
        - 3 p^{-4}  q^{-2}
        + 4 p^{-2}  q^{-2}
        + 9 p^{-2}
        - 2         q^{-2}
        \\
        \quad
        - 8         q^{2}
        + 2 p^{-2}  q^{2}
        + 9 p^{2}   q^{2}
        + 4 p^{2}
        + 2 p^{2}   q^{4}
        - 3 p^{4}   q^{2}
        - 3 p^{4}   q^{4}
        +   p^{6}   q^{4}
      \end{array}
    $
    \\[2mm]
    \hline
    & \\[-2mm]
    $ 6_1 $ &
    $
      \begin{array}{l}
          17
        +    p^{-4} q^{-4}
        + 3  p^{-4} q^{-2}
        - 3  p^{-2} q^{-4}
        - 10 p^{-2} q^{-2}
        - 7  p^{-2}
        + 2         q^{-4}
        + 10        q^{-2}
        \\
        \quad
        + 4         q^{2}
        - 3  p^{2}  q^{-2}
        - 10 p^{2}
        - 7  p^{2}  q^{2}
        +    p^{4}
        + 3  p^{4}  q^{2}
      \end{array}
    $
    \\[2mm]
    \hline
    & \\[-2mm]
    $ 6_2 $ &
    $
      \begin{array}{l}
          11
        - 3  p^{-6}
        + 9  p^{-4}
        - 12 p^{-2}
        - 4  p^{2}
        + 2         q^{-2}
        +    p^{-8} q^{-2}
        - 3  p^{-6} q^{-2}
        + 4  p^{-4} q^{-2}
        \\
        \quad
        - 4  p^{-2} q^{-2}
        + 14        q^{2}
        + 2  p^{-4} q^{2}
        - 8  p^{-2} q^{2}
        - 12 p^{2}  q^{2}
        + 4  p^{4}  q^{2}
        + 2         q^{4}
        - 8  p^{2}  q^{4}
        \\
        \quad
        + 9  p^{4}  q^{4}
        - 3  p^{6}  q^{4}
        + 2  p^{4}  q^{6}
        - 3  p^{6}  q^{6}
        +    p^{8}  q^{6}
      \end{array}
    $
    \\[3mm]
    \hline
    & \\[-2mm]
    $ 6_3 $ &
    $
      \begin{array}{l}
          25
        +    \( p^{-8} q^{-4} + p^{8}  q^{4}  \)
        + 10 \(        q^{-2} +        q^{2}  \)
        + 11 \( p^{-4} q^{-2} + p^{4}  q^{2}  \)
        + 4  \( p^{-4}        + p^{4}         \)
        \\
        \quad
        + 4  \( p^{-4} q^{-4} + p^{4}  q^{4}  \)
        - 2  \( p^{-2} q^{2}  + p^{2}  q^{-2} \)
        - 16 \( p^{-2}        + p^{2}         \)
        \\
        \quad
        - 16 \( p^{-2} q^{-2} + p^{2}  q^{2}  \)
        - 2  \( p^{-2} q^{-4} + p^{2}  q^{4}  \)
        - 3  \( p^{-6} q^{-2} + p^{6}  q^{2}  \)
        \\
        \quad
        - 3  \( p^{-6} q^{-4} + p^{6}  q^{4}  \)
      \end{array}
    $
    \\[3mm]
    \hline
    & \\[-2mm]
    $ 7_1 $ &
    $
      \begin{array}{l}
        1
        +    p^{-12}
        -    p^{-10}
        +    p^{-8}
        -    p^{-6}
        +    p^{-4}
        -    p^{-2}
        + 2         q^{2}
        -    p^{10} q^{2}
        + 2  p^{-8} q^{2}
        - 2  p^{-6} q^{2}
        \\
        \quad
        + 2  p^{-4} q^{2}
        - 2  p^{-2} q^{2}
        -    p^{2}  q^{2}
        + 2         q^{4}
        -    p^{-6} q^{4}
        + 2  p^{-4} q^{4}
        - 2  p^{-2} q^{4}
        - 2  p^{2}  q^{4}
        \\
        \quad
        +    p^{4}  q^{4}
        + 2         q^{6}
        -    p^{-2} q^{6}
        - 2  p^{2}  q^{6}
        + 2  p^{4}  q^{6}
        -    p^{6}  q^{6}
        -    p^{2}  q^{8}
        + 2  p^{4}  q^{8}
        - 2  p^{6}  q^{8}
        \\
        \quad
        +    p^{8}  q^{8}
        -    p^{6}  q^{10}
        + 2  p^{8}  q^{10}
        -    p^{10} q^{10}
        -    p^{10} q^{12}
        +    p^{12} q^{12}
      \end{array}
    $
    \\[3mm]
    \hline
    & \\[-2mm]
    $ 7_2 $ &
    $
      \begin{array}{l}
          5
        + 5  p^{-4}
        - 9  p^{-2}
        + 20        q^{2}
        + 3  p^{-4} q^{2}
        - 14 p^{-2} q^{2}
        - 9  p^{2}  q^{2}
        + 14        q^{4}
        +    p^{-4} q^{4}
        \\
        \quad
        - 6  p^{-2} q^{4}
        - 14 p^{2}  q^{4}
        + 5  p^{4}  q^{4}
        + 4         q^{6}
        -    p^{-2} q^{6}
        - 6  p^{2}  q^{6}
        + 3  p^{4}  q^{6}
        -    p^{2}  q^{8}
        +    p^{4}  q^{8}
      \end{array}
    $
    \\[5mm]
    \hline\hline
  \end{tabular}
  \caption[%
    ${LG}_K (q,p)$, evaluated for various oriented links $K$.
  ]{%
    ${LG}_K (q,p)$, evaluated for various oriented links $K$.  Data for
    $8_{17}$ are presented in \secref{EightSeventeen} and for $9_{42}$
    and $10_{48}$ in \secref{NineFortyTwoandTenFortyEight}.
  }
  \addtocontents{lot}{\protect\vspace{-2.5ex}}
  \tablabel{gl21oriented}
\end{table}

Evaluations of the Jones polynomials $V_K(t)$ for our example links are
presented in \tabref{JonesPolynomials}.  These polynomials agree with
published data \cite{Adams:94,Jones:85,Kauffman:93,Kawauchi:96}.  Their
correctness confirms that the evaluations of the Links--Gould invariant
are correct.

\begin{table}[ht]
  \centering
  \begin{tabular}{||c||l||}
    \hline\hline
    & \\[-3mm]
    $K$ &
    \multicolumn{1}{c||}{$ {V}_K (t) $}
    \\[1mm]
    \hline\hline
    & \\[-2mm]
    $ 2^2_1 $ &
    $
      -    t^{\frac{1}{2}}
      -    t^{\frac{5}{2}}
    $
    \\[2mm]
    \hline
    & \\[-2mm]
    $ 3_1 $ &
    $
           t
      +    t^{3}
      -    t^{4}
    $
    \\[2mm]
    \hline
    & \\[-2mm]
    $ 4_1 $ &
    $
           \( t^{-2} +    t^{2} \)
      -    \( t^{-1} +    t     \)
      + 1
    $
    \\[2mm]
    \hline
    & \\[-2mm]
    $ 5_1 $ &
    $
           t^{2}
      +    t^{4}
      -    t^{5}
      +    t^{6}
      -    t^{7}
    $
    \\[2mm]
    \hline
    & \\[-2mm]
    $ 5_2 $ &
    $
           t
      -    t^{2}
      + 2  t^{3}
      -    t^{4}
      +    t^{5}
      -    t^{6}
    $
    \\[2mm]
    \hline
    & \\[-2mm]
    $ 5^2_1 $ &
    $
      -    t^{-\frac{3}{2}}
      +    t^{-\frac{1}{2}}
      - 2  t^{\frac{1}{2}}
      +    t^{\frac{3}{2}}
      - 2  t^{\frac{5}{2}}
      +    t^{\frac{7}{2}}
    $
    \\[2mm]
    \hline
    & \\[-2mm]
    $ 6_1 $ &
    $
           t^{-4} 
      -    t^{-3}
      +    t^{-2}
      - 2  t^{-1}
      + 2
      -    t
      +    t^{2}
    $
    \\[2mm]
    \hline
    & \\[-2mm]
    $ 6_2 $ &
    $
           t^{-1}
      - 1
      + 2  t
      - 2  t^{2}
      + 2  t^{3}
      - 2  t^{4}
      +    t^{5}
    $
    \\[2mm]
    \hline
    & \\[-2mm]
    $ 6_3 $ &
    $
      -    \( t^{-3} + t^{3} \)
      + 2  \( t^{-2} + t^{2} \)
      - 2  \( t^{-1} + t     \)
      + 3
    $
    \\[2mm]
    \hline
    & \\[-2mm]
    $ 7_1 $ &
    $
           t^{3}
      +    t^{5}
      -    t^{6}
      +    t^{7}
      -    t^{8}
      +    t^{9}
      -    t^{10}
    $
    \\[2mm]
    \hline
    & \\[-2mm]
    $ 7_2 $ &
    $
           t
      -    t^{2}
      + 2  t^{3}
      - 2  t^{4}
      + 2  t^{5}
      -    t^{6}
      +    t^{7}
      -    t^{8}
    $
    \\[2mm]
    \hline\hline
    & \\[-2mm]
    $ 8_{17} $ &
    $
           \( t^{-4} + t^{4} \)
      - 3  \( t^{-3} + t^{3} \)
      + 5  \( t^{-2} + t^{2} \)
      - 6  \( t^{-1} + t     \)
      + 7
    $
    \\[2mm]
    \hline\hline
    & \\[-2mm]
    $ 9_{42} $ &
    $
           \( t^{-3} +    t^{3} \)
      -    \( t^{-2} +    t^{2} \)
      +    \( t^{-1} +    t     \)
      - 1
    $
    \\[2mm]
    \hline
    & \\[-2mm]
    $ 10_{48} $ &
    $
      \begin{array}{l}
        -   \( t^{-5} + t^{5} \)
        + 2 \( t^{-4} + t^{4} \)
        - 4 \( t^{-3} + t^{3} \)
        + 6 \( t^{-2} + t^{2} \)
        \\
        \quad
        - 7 \( t^{-1} + t \)
        + 9
      \end{array}
    $
    \\[2mm]
    \hline\hline
    & \\[-2mm]
    $ KT, KTI $ &
    $
           t^{-6}
      - 2  t^{-5}
      + 2  t^{-4}
      - 2  t^{-3}
      +    t^{-2}
      + 2  t
      - 2  t^{2}
      + 2  t^{3}
      -    t^{4}
    $
    \\[5mm]
    \hline\hline
  \end{tabular}
  \caption{
    ${V}_K (t)$, evaluated for various oriented links $K$.
  }
  \addtocontents{lot}{\protect\vspace{-2.5ex}}
  \tablabel{JonesPolynomials}
\end{table}

\pagebreak


\subsection{The Noninvertibility of $ 8_{17} $ is not Detected}
\addtocontents{toc}{\protect\vspace{-2.5ex}}
\seclabel{EightSeventeen}

Recall that $ 8_{17} $ is the smallest noninvertible knot.
We find:
\begin{eqnarray*}
  \lefteqn{
    {LG}_{8_{17}} (q,p)
    =
  }
  \\
  & &
    139
  +     \( p^{-12} q^{-6} + p^{12} q^6    \)
  -   4 \( p^{-10} q^{-6} + p^{10} q^6    \)
  -   4 \( p^{-10} q^{-4} + p^{10} q^4    \)
  \\
  & &
  +   7 \( p^{-8}  q^{-6} + p^8    q^6    \)
  +  18 \( p^{-8}  q^{-4} + p^8    q^4    \)
  +   7 \( p^{-8}  q^{-2} + p^8    q^2    \)
  \\
  & &
  -   7 \( p^{-6}  q^{-6} + p^6    q^6    \)
  -  36 \( p^{-6}  q^{-4} + p^6    q^4    \)
  -  36 \( p^{-6}  q^{-2} + p^6    q^2    \)
  \\
  & &
  -   7 \( p^{-6}         + p^6           \)
  +   3 \( p^{-4}  q^{-6} + p^4    q^6    \)
  +  40 \( p^{-4}  q^{-4} + p^4    q^4    \)
  \\
  & &
  +  82 \( p^{-4}  q^{-2} + p^4    q^2    \)
  +  40 \( p^{-4}         + p^4           \)
  +   3 \( p^{-4}  q^2    + p^4    q^{-2} \)
  \\
  & &
  -  22 \( p^{-2}  q^{-4} + p^2    q^4    \)
  - 102 \( p^{-2}  q^{-2} + p^2    q^2    \)
  - 102 \( p^{-2}         + p^2           \)
  \\
  & &
  -  22 \( p^{-2}  q^2    + p^2    q^{-2} \)
  +   4 \( q^{-4}         + q^4           \)
  +  68 \( q^{-2}         + q^2           \).
\end{eqnarray*}
We observe the invariance of \propref{ASymmetryofLG}, viz
$
  {LG}_{8_{17}} (q,p)
  =
  {LG}_{8_{17}} (q, q^{-1} p^{-1})
$,
illustrating that $LG$ \emph{doesn't} detect the noninvertibility of
$8_{17}$. Also note that as $ 8_{17} $ is amphichiral, the polynomial
is palindromic, as predicted by \propref{AmphichiralMeansPalindromic}.


\subsection{The Chirality of $9_{42}$ and $10_{48}$ is Detected}
\addtocontents{toc}{\protect\vspace{-2.5ex}}
\seclabel{NineFortyTwoandTenFortyEight}

We have:
\begin{eqnarray*}
  \lefteqn{
    {LG}_{9_{42}} (q,p)
    =
  }
  \\
  & &
      3
    +   p^{-8} q^{-6}
    - 2 p^{-6} q^{-6}
    - 2 p^{-6} q^{-4}
    +   p^{-4} q^{-6}
    + 3 p^{-4} q^{-4}
    +   p^{-4} q^{-2}
    +   p^{-4}
  \\
  & &
    -   p^{-2} q^{-4}
    -   p^{-2} q^{-2}
    - 3 p^{-2}
    - 3 p^{-2} q^2
    + 6        q^2
    + 2        q^4
    -   p^2    q^{-2}
    -   p^2
    - 3 p^2    q^2
  \\
  & &
    - 3 p^2    q^4
    +   p^4    q^{-2}
    + 3 p^4
    +   p^4    q^2
    +   p^4    q^4
    - 2 p^6
    - 2 p^6    q^2
    +   p^8    q^2
  \\
  \lefteqn{
    {LG}_{10_{48}} (q,p)
    =
  }
  \\
  & &
      165
    + 5   p^{-8}
    - 25  p^{-6}
    + 68  p^{-4}
    - 129 p^{-2}
    - 132 p^{2}
    + 67  p^{4}
    - 22  p^{6}
    + 4   p^{8}
  \\
  & &
    +     p^{-16} q^{-8}
    - 3   p^{-14} q^{-8}
    + 4   p^{-12} q^{-8}
    - 4   p^{-10} q^{-8}
    + 4   p^{-8}  q^{-8}
    - 2   p^{-6}  q^{-8}
  \\
  & &
    - 3   p^{-14} q^{-6}
    + 12  p^{-12} q^{-6}
    - 21  p^{-10} q^{-6}
    + 24  p^{-8}  q^{-6}
    - 22  p^{-6}  q^{-6}
    + 13  p^{-4}  q^{-6}
  \\
  & &
    - 3   p^{-2}  q^{-6}
    + 16          q^{-4}
    + 5   p^{-12} q^{-4}
    - 23  p^{-10} q^{-4}
    + 50  p^{-8}  q^{-4}
    - 69  p^{-6}  q^{-4}
  \\
  & &
    + 67  p^{-4}  q^{-4}
    - 43  p^{-2}  q^{-4}
    - 3   p^{2}   q^{-4}
    + 94          q^{-2}
    - 6   p^{-10} q^{-2}
    + 29  p^{-8}  q^{-2}
  \\
  & &
    - 72  p^{-6}  q^{-2}
    + 119 p^{-4}  q^{-2}
    - 132 p^{-2}  q^{-2}
    - 43  p^{2}   q^{-2}
    + 13  p^{4}   q^{-2}
    - 2   p^{6}   q^{-2}
  \\
  & &
    + 88          q^{2}
    - 2   p^{-6}  q^{2}
    + 12  p^{-4}  q^{2}
    - 39  p^{-2}  q^{2}
    - 129 p^{2}   q^{2}
    + 119 p^{4}   q^{2}
    - 69  p^{6}   q^{2}
  \\
  & &
    + 24  p^{8}   q^{2}
    - 4   p^{10}  q^{2}
    + 12          q^{4}
    - 2   p^{-2}  q^{4}
    - 39  p^{2}   q^{4}
    + 68  p^{4}   q^{4}
    - 72  p^{6}   q^{4}
    + 50  p^{8}   q^{4}
  \\
  & &
    - 21  p^{10}  q^{4}
    + 4   p^{12}  q^{4}
    - 2   p^{2}   q^{6}
    + 12  p^{4}   q^{6}
    - 25  p^{6}   q^{6}
    + 29  p^{8}   q^{6}
    - 23  p^{10}  q^{6}
  \\
  & &
    + 12  p^{12}  q^{6}
    - 3   p^{14}  q^{6}
    - 2   p^{6}   q^{8}
    + 5   p^{8}   q^{8}
    - 6   p^{10}  q^{8}
    + 5   p^{12}  q^{8}
    - 3   p^{14}  q^{8}
    +     p^{16}  q^{8}.
  \\
  & &
\end{eqnarray*}
As neither of these polynomials are palindromic, $LG$ distinguishes the
chirality of both these knots.


\subsection{The Noninvertible Pretzels are not Distinguished}
\addtocontents{toc}{\protect\vspace{-2.5ex}}
\seclabel{NoninvertiblePretzelsarenotDistinguished}

Experiments show that the Links--Gould invariant for this class of
noninvertible knots always displays the symmetry of
\eqref{Melanotaenia}, for all $p,q,r\leqslant 67$. This amounts to
$5456$ knots%
\footnote{
  More generally, if we sample for all $p,q,r\leqslant N$, for some odd
  $N\geqslant7$, we will have tested $(N-1)(N-3)(N-5)/48$ pretzels, the
  largest of which is a knot of $3N-6$ crossings.
}
(sheer bloody-mindedness!),
the smallest being the $(3,5,7)$ pretzel, a knot of $3+5+7=15$
crossings, and the largest being the $(63,65,67)$ pretzel, a knot of
$63+65+67=195$ crossings.  This illustrates the assertion of
\propref{CantDetectInvertibility} that our invariant cannot detect the
noninvertibility of \emph{any} knot.
The smallest of these polynomials is:
\begin{eqnarray*}
  \lefteqn{
    {LG}_{(3,5,7)} (q,p)
    =
  }
  \\
  & &
    57
  + 57  p^{-4}
  - 113 p^{-2}
  + 308 q^2
  + 84  q^2 p^{-4}
  - 279 q^2 p^{-2}
  - 113 p^2 q^2
  + 468 q^4
  \\
  & &
  + 83  q^4  p^{-4}
  - 329 q^4 p^{-2}
  - 279 p^2 q^4
  + 57  p^4 q^4
  + 464 q^6
  + 60  q^6 p^{-4}
  - 279 q^6 p^{-2}
  \\
  & &
  - 329 p^2 q^6
  + 84  p^4 q^6
  + 338 q^8
  + 30  q^8 p^{-4}
  - 172 q^8 p^{-2}
  - 279 p^2 q^8
  + 83  p^4 q^8
  \\
  & &
  + 174 q^{10}
  + 9   q^{10} p^{-4}
  - 71  q^{10} p^{-2}
  - 172 p^2 q^{10}
  + 60  p^4 q^{10}
  + 56  q^{12}
  +     q^{12} p^{-4}
  \\
  & &
  - 16  q^{12} p^{-2}
  - 71  p^2 q^{12}
  + 30  p^4 q^{12}
  + 8   q^{14}
  -     q^{14} p^{-2}
  - 16  p^2 q^{14}
  + 9   p^4 q^{14}
  \\
  & &
  -     p^2 q^{16}
  +     p^4 q^{16}.
\end{eqnarray*}
We also discover that all these $5456$ pretzels are chiral as none are
palindromic, lending weight to Trotter's assertion that all such
pretzels are chiral.


\subsection{The Kinoshita--Terasaka Mutants are not Distinguished}
\addtocontents{toc}{\protect\vspace{-2.5ex}}

The polynomials for \emph{both} mutants are:
\begin{eqnarray*}
  \lefteqn{
    {LG}_{KT} (q,p)
    =
  }
  \\
  & &
    - 23
    -    p^{-6} q^{-8}
    -    p^{-6} q^{-6}
    +  2 p^{-6} q^{-4}
    +    p^{-6} q^{-2}
    -    p^{-6}
    +    p^{-4} q^{-8}
    +  6 p^{-4} q^{-6}
  \\
  & &
    -  3 p^{-4} q^{-4}
    -  9 p^{-4} q^{-2}
    +  2 p^{-4}
    +  3 p^{-4} q^2
    -  7 p^{-2} q^{-6}
    -  7 p^{-2} q^{-4}
    + 18 p^{-2} q^{-2}
  \\
  & &
    +  9 p^{-2}
    - 11 p^{-2} q^2
    -  2 p^{-2} q^4
    +  2        q^{-6}
    + 14        q^{-4}
    -  8        q^{-2}
    +  6        q^2
    + 10        q^4
    -  7 p^2    q^{-4}
  \\
  & &
    -  7 p^2    q^{-2}
    + 18 p^2
    +  9 p^2    q^2
    - 11 p^2    q^4
    -  2 p^2    q^6
    +    p^4    q^{-4}
    +  6 p^4    q^{-2}
    -  3 p^4
    -  9 p^4    q^2
  \\
  & &
    +  2 p^4    q^4
    +  3 p^4    q^6
    -    p^6    q^{-2}
    -    p^6
    +  2 p^6    q^2
    +    p^6    q^4
    -    p^6    q^6,
\end{eqnarray*}
hence $LG$ \emph{does not} distinguish between these mutants.  As
predicted by the theorem of \cite{MortonCromwell:96} (i.e.
\propref{CantDistinguishMutants}), the tensors $KTA$ and $KTA'$ are in
fact \emph{identical}, which explains why the pair of mutants yield the
same invariant.

\vfill

\cleardoublepage


\section{Conclusions and Extensions}
\addtocontents{toc}{\protect\vspace{-2.5ex}}
\seclabel{Conclusions}

\subsection{Conclusions}
\addtocontents{toc}{\protect\vspace{-2.5ex}}

The Links--Gould invariant is at least as powerful as the HOMFLY
and Kauffman (two-variable) invariants in the following sense.  Neither
of those invariants can detect the chirality of $9_{42}$, although they
can do so for \emph{all} preceding knots. Beyond $9_{42}$, there are
other knots for which they fail to detect
chirality.  $10_{48}$ is the next example for which the HOMFLY
invariant fails to detect chirality (although the Kauffman does detect
it). However, although our invariant \emph{does} detect the chirality
of $9_{42}$ and $10_{48}$, there is as yet no evidence to tell us
whether it detects the chirality of \emph{all preceding} knots, that
is, for all we know, it may fail on $9_{41}$ (which is chiral,
according to \cite[p385]{Jones:87}, \cite[p265]{Kawauchi:96}).  This
information is summarised in \tabref{chiralitydetectcomparison}.

\begin{table}[ht]
  \centering
  \begin{tabular}{||c||c|c|c||}
    \hline\hline
    Knot & HOMFLY & Kauffman & Links-Gould \\
    \hline\hline
    & & & \\[-4mm]
    $ 3_1$ to $7_2$     & Yes    & Yes   & Yes       \\
    & & & \\[-4mm]
    \hline
    & & & \\[-4mm]
    $ 7_3$ to $9_{41}$     & Yes    & Yes   & ?       \\
    & & & \\[-4mm]
    \hline
    & & & \\[-4mm]
    $ 9_{42}$              & No   & No  & Yes \\
    & & & \\[-4mm]
    \hline
    & & & \\[-4mm]
    $ 9_{43}$ to $10_{47}$ & Yes    & Yes   & ?       \\
    & & & \\[-4mm]
    \hline
    & & & \\[-4mm]
    $10_{48}$              & No   & Yes   & Yes \\
    & & & \\[-4mm]
    \hline
    & & & \\[-4mm]
    $10_{49}$ and above    & sometimes  & sometimes & ?       \\[2mm]
    \hline\hline
  \end{tabular}
  \caption{%
    Chirality detection properties of two-variable link invariants.
  }
  \addtocontents{lot}{\protect\vspace{-2.5ex}}
  \tablabel{chiralitydetectcomparison}
\end{table}

We have as yet found no chiral knot for which $LG$ fails to detect its
chirality, so whether our invariant is a complete invariant for
chirality is an open question.

Recall that the method used to evaluate our invariant for a particular
link involves the deduction of an appropriate abstract tensor that is
suitable for use with efficient shortcuts using auxiliary tensors, and
this requires the labour-intensive drawing of a suitable presentation
(i.e. finding an appropriate quasi-Morse function).  It is thus
infeasible to use this method to systematically evaluate the invariant
for all links up to any particular point. This situation is in marked
contrast to that for the HOMFLY and Kauffman polynomials, for which
extensive tables have been (automatically) computed. The question of
whether the Links--Gould invariant is superior to the other
two-variable invariants is therefore currently unanswerable by
experiment. To this end, it would be useful to conduct research to
develop a method to automate the deduction of the abstract tensor
directly from a braid presentation.  (The abstract tensor would then be
applicable to the evaluation of \emph{any} other invariants based on
state models.) Work on this program is in preparation \cite{DeWit:99}.

\pagebreak

As explained in \propref{CantDetectInvertibility}, no invariant of
this type%
\footnote{
  That is, an invariant based on a state space model where the crossing
  tensor $\sigma$ is based on the R matrix associated with a quantum
  (super)algebra representation.
}
can be expected to detect inversion. Our early hopes that our invariant
might do so were shattered by experimental results; indeed this led us
to the literature on invariants constructed from quantum algebras,
which led us to propose that the same considerations naturally carry
over to quantum superalgebras.

Also, as explained in \propref{CantDistinguishMutants}, no state space
invariant of this type%
\footnote{
  That is, an invariant based on a state space model where the crossing
  tensor $\sigma$ is based on the R matrix associated with a quantum
  (super)algebra representation in which the decomposition of the
  tensor product contains only unit multiplicities.
}
can be expected to distinguish mutants; again negative experimental
results caused us to make the considerations leading to the theoretical
proof of this statement.


\subsection{Extension -- A $U_q[gl(3|1)]$ Two-Variable Link Invariant?}
\addtocontents{toc}{\protect\vspace{-2.5ex}}

Having investigated the Links--Gould invariant based on the
$U_q[gl(2|1)]$ family of representations labelled $\(0,0\,|\,\alpha\)$,
it is immediately attractive to ask if other two-variable invariants
may be constructed from the families $\(\dot{0}_m\,|\,\alpha\)$ of
$U_q[gl(m|1)]$ representations, for general $m\geqslant 2$. The
immediate answer is that yes, of course they exist.  Their actual
evaluation would require a lot more computation as the R matrices would
be larger, but the process would in principle be the same. In practice,
explicit construction of the R matrices is hampered by the algebraic
complexities involved in finding an explicit orthogonal basis for the
tensor product module.


\subsubsection{An Orthonormal Basis for $V\otimes V$}

We illustrate the situation with the example of the $U_q[gl(3|1)]$
family of representations $\(0,0,0\,|\,\alpha\)$. This family is $8$
dimensional (cf. the family for $U_q[gl(2|1)]$ is $4$ dimensional),
which means that the R matrix has $8^4=4096$ rather than $4^4=256$
(albeit mostly zero) entries, so the computational effort in evaluating
abstract tensors would be orders of magnitude larger.

A basis for the decomposition of $V\otimes V$ is provided in
\cite{GeGouldZhangZhou:98a}, but this is not orthogonal.  The
decomposition involves two modules of dimension $8$, for which a fully
orthogonal basis is given, and two modules of dimension $24$, of which
a basis for only one (called $V_2$) is provided, and that basis is not
orthogonal. Using the Gram--Schmidt process to (manually!)
orthogonalise that basis for $V_2$ yields reasonably tractable
expressions for all but the last vector (originally called
$\ket{\Psi^2_{24}}$).  Attempting to normalise this basis would further
complicate the expressions involved, and it appears that the resulting
R matrix would contain expressions similar to $Y$ (see \eqref{Ydefn})
in our R matrix, although somewhat more complicated. The complexity of
$Y$ doesn't propagate through to evaluations of $LG$, and we anticipate
that evaluations of the $U_q[gl(3|1)$ invariant would behave
similarly.

\vfill


\subsubsection{Sour Grapes}

Perhaps there is nothing lost by the lack of success in constructing a
state model for an invariant based on $U_q[gl(3|1)]$ family of
representations $\(0,0,0\,|\,\alpha\)$. Such an invariant would have
the same inability to detect inversion and mutation as does $LG$, and
there is no particular reason to expect that it would always detect
chirality.


\subsection{Automatic Construction of $U_q[gl(m|n)]$ Representations}
\addtocontents{toc}{\protect\vspace{-2.5ex}}

The process used to build the explicit $U_q[gl(2|1)]$ representation
was presented in rigorous detail so that it might be analysed in order
to consider the automatic construction of $U_q[gl(m|n)]$ (generally
non-parametric) representations. (A description of $U_q[gl(m|n)]$ is
provided in \secref{Uqglmn}.)

A general element of $U_q[gl(m|n)]$ is a sum of (complex) multiples of
products of $U_q[gl(m|n)]$ generators.  As $U_q[gl(m|n)]$ \emph{is} a
(graded) vector space, we may meaningfully call these general elements
``vectors''.

So, given a vector involving weighted strings of $U_q[gl(m|n)]$
generators, where each string may contain nonsimple generators, how do
we convert this to some kind of normal ordering?  Well, firstly we must
(recursively) replace every nonsimple generator in the vector with a
combination of simple generators using
\eqref{UqglmnNonSimpleGenerators}, until our vector has only simple
generators. Secondly, we must (recursively) examine every product of
two generators ``$AB$'' in each summand in the vector.  If the pair of
generators are disordered (i.e. $A > B$), we use the relations of
\secref{UqglmnRelations} to replace the product with another, ordered
string.

This process provides us with all the necessary machinery for the
manipulation of vectors encountered when constructing concrete
representations of $U_q[gl(m|n)]$.

\cleardoublepage


\appendix

\section{The Quantum Superalgebra $U_q[gl(m|n)]$}
\addtocontents{toc}{\protect\vspace{-2.5ex}}
\seclabel{Uqglmn}

The origin of the quantum superalgebra $U_q[gl(m|n)]$
(we intend $m,n \in \BZ^+$) as a \emph{quantum deformation} of
$gl(m|n)$ is described in
\cite{%
  BrackenGouldZhang:90,%
  DeliusGouldLinksZhang:95a,%
  ScheunertNahmRittenberg:77a,%
  Yamane:91,%
  Yamane:94%
},
and in the book by Chari and Pressley
\cite[see~\S6.5]{ChariPressley:94}.  Here, we set the context for
\secref{Uqgl21} (with obvious extensions of the definitions of the
graded commutator and $q$ bracket).  The material is taken from
\cite[pp1237-1238]{Zhang:93}, except that we have substituted the
notation $q^{{\(-\)}^{\[a\]}{E^a}_a}$ for $K_a$ and modified the
definition of the coproduct and antipode.


\subsection{$ \protect \mathbb{Z}_2 $ Grading on $U_q[gl(m|n)]$
  Indices}
\addtocontents{toc}{\protect\vspace{-2.5ex}}

Firstly, we define a $\mathbb{Z}_2$ grading
$ \[\cdot\] : \{ 1, \dots, m+n \} \to \mathbb{Z}_2 $
on the $m+n$ $gl\(m|n\)$ \emph{indices}:
\be
  \[ a \]
  =
  \left\{
  \begin{array}{lll}
    0 \quad & \textrm{if} \quad a = 1, \dots, m \qquad & \textrm{(even)}
    \\
    1       & \textrm{if} \quad a = m+1, \dots, m+n    & \textrm{(odd)}.
  \end{array}
  \right.
\ee
Throughout, we shall use dummy indices $ a, b = 1, \dots, m+n $ where
meaningful.


\subsection{$U_q[gl(m|n)]$ Generators}
\addtocontents{toc}{\protect\vspace{-2.5ex}}

A set of generators for $U_q[gl(m|n)]$ is:
\be
  \hspace{-25pt}
  \left\{
    \begin{array}{lll}
      q^{ {(-)}^{[a]} {E^a}_a}, \quad & a = 1, \dots, m+n \qquad &
        \textrm{($m+n$ Cartan)} \\
      {E^{a}}_{b},       & a<b = 1, \dots, m+n &
         \textrm{($(m+n)(m+n-1)/2$ raising)} \\
      {E^{b}}_{a},       & a<b = 1, \dots, m+n &
         \textrm{($(m+n)(m+n-1)/2$ lowering)}.
    \end{array}
  \right\}
\ee
We write the Cartan generators as $q^{ {(-)}^{[a]} {E^a}_a}$; they might
also be written as $q^{{E^a}_a}$; the notation including the parity
factor ensures consistency with the literature and what follows below.
More significantly, we are \emph{not} writing the Cartan generators as
${E^a}_a$ (as per $gl(m|n)$).

On all the ${(m+n)}^2$ \emph{generators} we define
a natural $\mathbb{Z}_2$ grading in terms of the grading on the
indices:
\be
  [{E^a}_b]
  \defeq
  [a] + [b]
  \qquad
  \( \mod~2 \),
\ee
where the definition of the exponential yields
$
  [ q^{{E^a}_a} ] = \[ {E^a}_a \] = 0
$.
The product of homogeneous $x,y\in U_q[gl(m|n)]$ has degree:
\be
  \[ x y \]
  \defeq
  \[x\] + \[y\]
  \qquad
  \( \mod~2 \).
\ee


\subsubsection{$U_q[gl(m|n)]$ Simple Generators}

The $U_q[gl(m|n)]$ \emph{simple} generators are a subset of $3(m+n)-2$
elements of the above generators:
\be
  \left\{
    \begin{array}{lll}
      q^{ {(-)}^{[a]} {E^a}_a}, \qquad & a = 1, \dots, m+n \qquad &
        \textrm{($m+n$ Cartan)} \\
      {E^{a}}_{a+1},     & a,a+1 = 1, \dots, m+n &
        \textrm{($m+n-1$ raising)} \\
      {E^{a+1}}_{a},     & a,a+1 = 1, \dots, m+n &
        \textrm{($m+n-1$ lowering)},
    \end{array}
  \right\}
\ee
where by the expression ${E^{a}}_{a+1}$, $ a,a+1 = 1, \dots, m+n $, we
of course intend to exclude the nonmeaningful ${E^{0}}_{1}$ and
${E^{m+n}}_{m+n+1}$. Note that there are only two \emph{odd} simple
generators:  ${E^{m}}_{m+1}$ (raising) and ${E^{m+1}}_{m}$ (lowering).
The fact that there are $m+n-1$ simple raising generators indicates
that $U_q[gl(m|n)]$ has rank $m+n-1$.


\subsubsection{$U_q[gl(m|n)]$ Nonsimple Generators}

In the $gl(m|n)$ case, the remaining nonsimple (non Cartan) generators
satisfy the same commutation relations as do the simple generators.
Writing ${e^a}_b$ as the ${(m+n)}^2$ $ gl(m|n)$ generators, to parallel
the ${E^a}_b$ of $U_q[gl(m|n)]$, we have:
\bse
  \[ {e^a}_b, {e^c}_d \]
  =
  {\delta^c}_b {e^a}_d
  -
  {(-)}^{([a]+[b])([c]+[d])}
  {\delta^a}_d {e^c}_b,
  \eqlabel{glmnCommutationRelations}
\ese
and this is true for \emph{all} generators ${e^a}_b, {e^c}_d$, not just
the simple ones. To be sure, Zhang \cite[p1238]{Zhang:93} obtains
these generators by taking the limit as $q \to 1$ of $ U_q[gl(m|n)]$
to yield $U[gl(m|n)]$, with simple generators:
\be
  \hspace{-11mm}
  \left\{
    \begin{array}{lll}
      {e^a}_{a+1}
      \triangleq
      \lim_{q\to1}
        {E^a}_{a+1},
      \quad & a = 1, \dots, m+n-1 \quad &
        \textrm{($m+n-1$ raising)}
      \\
      {e^a}_{a-1}
      \triangleq
      \lim_{q\to1}
        {E^a}_{a-1},
      \quad & a = 2, \dots, m+n \qquad &
        \textrm{($m+n-1$ lowering)}
      \\
      {e^a}_a
      \triangleq
      \lim_{q\to1}
      {\[ {E^a}_a \]}_q,
      & a = 1, \dots, m+n \qquad &
        \textrm{($m+n$ Cartan)}.
    \end{array}
  \right\}
\ee
(Of course, we have $\lim_{q\to1} q^{ \pm {E^a}_{a} } = 1$.)
Having defined these, the $gl(m|n)$ commutation relations
\eqref{glmnCommutationRelations} serve to define the nonsimple
generators.

The situation is different in $U_q[gl(m|n)]$.  The remaining
(nonsimple) generators do not satisfy the same commutation relations,
in fact they are recursively defined in terms of the simple generators
\cite[p1971,~(3)]{Zhang:92} and
\cite[p1238,~(2)]{Zhang:93}:

\bse
  \left.
  \begin{array}{lrcll}
    (a) &
    {E^a}_b
    & \defeq &
      {E^a}_c
      {E^c}_b
      -
      q^{-{(-)}^{[c]}}
      {E^c}_b
      {E^a}_c
    \qquad
    & \textrm{(raising)}
    \\
    (b) &
    {E^b}_a
    & \defeq &
      {E^b}_c
      {E^c}_a
      -
      q^{{(-)}^{[c]}}
      {E^c}_a
      {E^b}_c
    &
    \textrm{(lowering)},
  \end{array}
  \quad
  \right\}
  \eqlabel{UqglmnNonSimpleGenerators}
\ese
where $a<c<b$.
To evaluate a graded commutator involving a nonsimple generator, that
generator must first be expanded using
\eqref{UqglmnNonSimpleGenerators}, and then the graded commutator
expanded by linearity.


\subsection{$U_q[gl(m|n)]$ Relations}
\addtocontents{toc}{\protect\vspace{-2.5ex}}
\seclabel{UqglmnRelations}

For the \emph{simple} $U_q[gl(m|n)]$ generators, we have the
following $U_q[gl(m|n)]$ \emph{relations}:

\begin{itemize}
\item
  The Cartan generators all commute:
  \be
    q^{\pm{E^a}_a}
    q^{\pm{E^b}_b}
    =
    q^{\pm{E^b}_b}
    q^{\pm{E^a}_a},
    \qquad \qquad
    q^{{E^a}_a}
    q^{-{E^a}_a}
    =
    1,
  \ee

\item
  The Cartan generators commute with the simple raising and lowering
  generators in the following manner:
  \bse
    q^{{\(-\)}^{\[a\]}{E^a}_a} {E^b}_{b\pm1} q^{-{\(-\)}^{\[a\]}{E^a}_a}
    =
    q^{{\(-\)}^{\[a\]} \( \delta^a_b - \delta^a_{b\pm1} \)}
    {E^b}_{b\pm1}.
    \eqlabel{UqglmnCartanRaisingCommutation}
  \ese
  From \eqref{UqglmnCartanRaisingCommutation}, we have the following
  useful interchange:
  \bse
    q^{{\(-\)}^{\[a\]}{E^a}_a}
    {E^b}_{b\pm1}
    =
    q^{{\(-\)}^{\[a\]} \( \delta^a_b - \delta^a_{b\pm1} \)}
    {E^b}_{b\pm1}
    q^{{\(-\)}^{\[a\]}{E^a}_a}.
    \eqlabel{UqglmnUsefulInterchangePositive}
  \ese
  Replacing $q$ with $q^{-1}$ in \eqref{UqglmnUsefulInterchangePositive}
  yields the equivalent:
  \be
    q^{-{\(-\)}^{\[a\]}{E^a}_a}
    {E^b}_{b\pm1}
    =
    q^{- {\(-\)}^{\[a\]} \( \delta^a_b - \delta^a_{b\pm1} \)}
    {E^b}_{b\pm1}
    q^{-{\(-\)}^{\[a\]}{E^a}_a}.
  \ee

\item
  The squares of the odd simple generators are zero:
  \be
    {\( {E^{m}}_{m+1} \)}^2
    =
    {\( {E^{m+1}}_{m} \)}^2
    =
    0.
  \ee
  It may be shown that this implies that the squares of the nonsimple
  odd generators are also zero.

\item
  The non-Cartan generators satisfy the following commutation relations
  (this is the really interesting part!):
  \be
    \[ {E^a}_{a+1}, {E^{b+1}}_b \]
    =
    {\delta^a}_b
    {\[
      {E^{a}}_{a}
      -
      {(-)}^{[a]+[a+1]}
      {E^{a+1}}_{a+1}
    \]}_q.
  \ee

  We also have, for $|a-b| > 1$, the commutations:
  \be
    {E^{a}}_{a+1}
    {E^{b}}_{b+1}
    =
    {E^{b}}_{b+1}
    {E^{a}}_{a+1}
    \qquad
    \textrm{and}
    \qquad
    {E^{a+1}}_{a}
    {E^{b+1}}_{b}
    =
    {E^{b+1}}_{b}
    {E^{a+1}}_{a}.
  \ee

\pagebreak

\item
  There are two additional relations, known as the
  \emph{Serre relations}. Their inclusion ensures that the algebra is
  reduced enough to be \emph{simple}. For $a\neq m$:
  \be
    \begin{array}{l}
      {\( {E^{a}}_{a+1} \)}^2 {E^{a\pm1}}_{a\pm1+1}
      -
      \( q + q^{-1} \) {E^{a}}_{a+1} {E^{a\pm1}}_{a\pm1+1} {E^{a}}_{a+1}
      \\
      \hspace{70mm}
      +
      {E^{a\pm1}}_{a\pm1+1} {\( {E^{a}}_{a+1} \)}^2
      =
      0
      \\
      {\( {E^{a+1}}_{a} \)}^2 {E^{a\pm1+1}}_{a\pm1}
      -
      \( q + q^{-1} \) {E^{a+1}}_{a} {E^{a\pm1+1}}_{a\pm1} {E^{a+1}}_{a}
      \\
      \hspace{70mm}
      +
      {E^{a\pm1+1}}_{a\pm1} {\( {E^{a+1}}_{a} \)}^2
      =
      0.
    \end{array}
  \ee
  To complement these, for a limited set of odd generators, we have:
  \bse
     \[ {E^{m-1}}_{m+2}, {E^{m}}_{m+1} \]
     =
     \[ {E^{m+2}}_{m-1}, {E^{m+1}}_{m} \]
     =
     0,
     \eqlabel{Uqglmnspecialoddcommutation}
  \ese
  where, as all the elements are odd, the graded commutator devolves
  to a simple anticommutator, viz \eqref{Uqglmnspecialoddcommutation}
  may be written as the exchanges:
  \be
    {E^{m-1}}_{m+2} {E^{m}}_{m+1}
    \eq
    -
    {E^{m}}_{m+1} {E^{m-1}}_{m+2},
    \\
    {E^{m+2}}_{m-1} {E^{m+1}}_{m}
    \eq
    -
    {E^{m+1}}_{m} {E^{m+2}}_{m-1}.
  \ee
  Note the presence of the nonsimple
  generators $ {E^{m-1}}_{m+2} $ and $ {E^{m+2}}_{m-1} $ here.
  The index range ensures that these relations are irrelevant to
  $U_q[gl(2|1)]$, which is why we didn't encounter them in
  \secref{Uqgl21Relations}.
\end{itemize}

It must be emphasised that expressions involving nonsimple generators
may in general only be manipulated by prior expansion using
\eqref{UqglmnNonSimpleGenerators}.

\pagebreak


\subsection{$U_q[gl(m|n)]$ as a Hopf Superalgebra}
\addtocontents{toc}{\protect\vspace{-2.5ex}}

$U_q\[gl(m|n)\]$ may be regarded as a Hopf superalgebra when equipped
with the following coproduct $\Delta$, counit $\varepsilon$ and
antipode $S$ structures \cite[p1238]{Zhang:93}. It is in fact
quasitriangular (i.e. it possesses an R matrix).

\subsubsection{Coproduct $\Delta$}

We define a coproduct (a.k.a. comultiplication) structure,
which is a
$\mathbb{Z}_2$ graded algebra homomorphism
$
  \Delta : U_q\[gl(m|n)\] \to U_q\[gl(m|n)\] \otimes U_q\[gl(m|n)\]
$
by:
\be
  \lefteqn{
    \begin{array}{lrcl}
      (a) &
      \Delta ({E^{a}}_{a+1})
      \eq
      {E^{a}}_{a+1}
      \otimes
      q^{
        -\frac{1}{2}
        \(
          {(-)}^{[a]} {E^{a}}_{a}
          -
          {(-)}^{[a+1]} {E^{a+1}}_{a+1}
        \)
      }
      \\
      & & &
      \qquad
      \qquad
      +
      q^{
        \frac{1}{2}
        \(
          {(-)}^{[a]} {E^{a}}_{a}
          -
          {(-)}^{[a+1]} {E^{a+1}}_{a+1}
        \)
      }
      \otimes
      {E^{a}}_{a+1}
      \\
      (b) &
      \Delta ({E^{a+1}}_{a})
      \eq
      {E^{a+1}}_{a}
      \otimes
      q^{
        -
        \frac{1}{2}
        \(
          {(-)}^{[a]} {E^{a}}_{a}
          -
          {(-)}^{[a+1]} {E^{a+1}}_{a+1}
        \)
      }
      \\
      & & &
      \qquad
      \qquad
      +
      q^{
        \frac{1}{2}
        \(
          {(-)}^{[a]} {E^{a}}_{a}
          -
          {(-)}^{[a+1]} {E^{a+1}}_{a+1}
        \)
      }
      \otimes
      {E^{a+1}}_{a}
      \\
      (c) &
      \Delta (q^{\pm {(-)}^{[a]} {E^{a}}_{a}})
      \eq
      q^{\pm {(-)}^{[a]} {E^{a}}_{a}}
      \otimes
      q^{\pm {(-)}^{[a]} {E^{a}}_{a}},
    \end{array}
  }
\ee
and extended to an algebra homomorphism on $U_q\[gl(m|n)\]$.
$\Delta$ being graded means that it preserves grading, viz for
homogeneous $x \in U_q\[gl(m|n)\]$, we have $\[\Delta(x)\] = \[x\]$.
$\Delta$ being a homomorphism means that $\Delta(1)=1\otimes1$, and
for $x,y \in U_q\[gl(m|n)\]$:
\be
  \Delta(x y)
  =
  \Delta(x)
  \Delta(y).
\ee

As well as this standard coproduct, there exists another possible
coproduct structure:  $\overline{\Delta}$, defined by
$\overline{\Delta}=T\cdot\Delta$, where the \emph{twist map} $T$, an
operator on the tensor product $U_q[gl(m|n)]\otimes U_q[gl(m|n)]$, is
defined for homogeneous $x,y\in U_q[gl(m|n)]$ by:
\be
  T \( x \otimes y \)
  =
  {\( - \)}^{\[ x \] \[ y \]}
  \( y \otimes x \).
\ee


\subsubsection{Counit $\varepsilon$}

We define a counit $\varepsilon:U_q\[gl(m|n)\]\to\BC$, also a
$\mathbb{Z}_2$ graded algebra homomorphism:
\be
  \varepsilon ( {E^{a}}_{a\pm1} )
  \eq
  0,
  \\
  \varepsilon ( q^{ {(-)}^{[a]} {E^a}_a} )
  \eq
  \varepsilon ( q^{ - {(-)}^{[a]} {E^a}_a} )
  =
  \varepsilon ( 1 )
  =
  1,
\ee
and extended to an algebra homomorphism on all of $U_q\[gl(m|n)\]$.


\subsubsection{Antipode $S$}

We define an antipode $S:U_q\[gl(m|n)\]\to U_q[gl(m|n)]$, which is a
$\mathbb{Z}_2$ graded algebra \emph{anti}homomorphism (actually an
anti\emph{auto}morphism):
\be
  S ( {E^{a}}_{{a+1}} )
  \eq
  -
  q^{
    -
    \frac{1}{2}
    \[
      {(-)}^{[a]}
      +
      {(-)}^{[a+1]}
    \]
  }
  {E^{a}}_{{a+1}}
  \\
  S ( {E^{a+1}}_{a} )
  \eq
  -
  q^{
    \frac{1}{2}
    \[
      {(-)}^{[a]}
      +
      {(-)}^{[a+1]}
    \]
  }
  {E^{a+1}}_{{a}}
  \\
  S (q^{ {(-)}^{[a]} {E^{a}}_{a}})
  \eq
  q^{- {(-)}^{[a]} {E^a}_a},
\ee
and
extended as an algebra homomorphism on all of $U_q\[gl(m|n)\]$. Because
$S$ is an \emph{anti}automorphism, i.e. for homogeneous $x,y\in
U_q[gl(m|n)]$:
\be
  S \( x y \)
  =
  {(-)}^{[x] [y]}
  S \( y \) S \( x \),
\ee
we may deduce that $S\(1\)=1$.

\cleardoublepage


\section{\textsc{Mathematica} Code}
\addtocontents{toc}{\protect\vspace{-2.5ex}}

\subsection{\texttt{Uqgl21.m}}
\addtocontents{toc}{\protect\vspace{-2.5ex}}
\seclabel{Uqgl21.m}

\footnotesize
\begin{verbatim}

(* ------------------------------------------------------------------

  Uqgl21.m

  Mathematica code to build projectors and the braid generator.
  We are interested in the 4 dimensional representation LA = (0,0|a)
  of the Quantum Superalgebra U_q[gl(2|1)].

  David De Wit

  20       May  1996  --   3      June  1996
  20     March  1998  --  27     March  1998
  08   October  1998  --  20  November  1998

   ------------------------------------------------------------------ *)

(*  Basis vectors of V \otimes V *)

ket[1, 1] = Kron[ket[1], ket[1]];

ket[1, 2] =
  (q^a + q^(-a))^(-1/2) *
  (q^(a/2) * Kron[ket[1], ket[2]] + q^(-a/2) * Kron[ket[2], ket[1]]);

ket[1, 3] =
  (q^a + q^(-a))^(-1/2) *
  (q^(a/2) * Kron[ket[1], ket[3]] + q^(-a/2) * Kron[ket[3], ket[1]]);

ket[1, 4] =
  (q^a + q^(-a))^(-1/2) *
  Qbracket[q, 2 a + 1]^(-1/2) *
  (
    Qbracket[q, a + 1]^(1/2) *
    (q^a * Kron[ket[1], ket[4]] + q^(-a) * Kron[ket[4], ket[1]])
    -
    Qbracket[q, a]^(1/2) *
    (q^(1/2) * Kron[ket[3], ket[2]] - q^(-1/2) * Kron[ket[2], ket[3]])
  );












ket[3, 1] =
  (q^(a+1) + q^(-a-1))^(-1/2) *
  Qbracket[q, 2 a + 1]^(-1/2) *
  (
    Qbracket[q, a]^(1/2) *
    (q^(a+1) * Kron[ket[4], ket[1]] + q^(-a-1) * Kron[ket[1], ket[4]])
    +
    Qbracket[q, a + 1]^(1/2) *
    (q^(1/2) * Kron[ket[3], ket[2]] - q^(-1/2) * Kron[ket[2], ket[3]])
  );

ket[3, 2] =
  (q^(a+1) + q^(-a-1))^(-1/2) *
  (
    q^((a+1)/2) * Kron[ket[4], ket[2]]
    +
    q^(-(a+1)/2) * Kron[ket[2], ket[4]]
  );

ket[3, 3] =
  (q^(a+1) + q^(-a-1))^(-1/2) *
  (
    q^((a+1)/2) * Kron[ket[4], ket[3]]
    +
    q^(-(a+1)/2) * Kron[ket[3], ket[4]]
  );

ket[3, 4] = Kron[ket[4], ket[4]];

Grading[ket[k_]] := (Grading[ket[k]] = If[MemberQ[{1,4},k],0,1]);

Grading[bra[k_]] := (Grading[bra[k]] = Grading[ket[k]]);

DimTPSubModule[1] := 4;

DimTPSubModule[3] := 4;

bra[k_, j_] :=
  (bra[k, j] =
    Expand[
      ket[k, j] //.
      Kron[ket[h_], ket[i_]] :>
        (-1)^(Grading[ket[h]]*Grading[ket[i]]) * Kron[bra[h], bra[i]]
    ]
  ) /;
  (MemberQ[{1,3},k] && (j <= DimTPSubModule[k]));

(* ------------------------------------------------------------------ *)










(* Utility functions *)

KroneckerDelta[i_, j_] := If[i == j, 1, 0];

OperatorQ[x_NonCommutativeMultiply] := Map[OperatorQ, Apply[Or, x]];
OperatorQ[x_Times] := Map[OperatorQ, Apply[Or, x]];
OperatorQ[x_Plus] := Map[OperatorQ, Apply[Or, x]];

OperatorQ[ket[_]] := True;       OperatorQ[bra[_]] := True;
OperatorQ[Kron[_, _]] := True;
OperatorQ[_] := False;

SetAttributes[NonCommutativeMultiply, {Listable, Flat, OneIdentity}];

Unprotect[NonCommutativeMultiply];

  NonCommutativeMultiply[Kron[A_, B_], Kron[C_, D_]] :=
    (-1)^(Grading[B] * Grading[C]) * Kron[A ** C, B ** D];

  (X_) ** ((alpha_) (Y_)) := alpha X ** Y /; NumberQ[alpha];
  (X_) ** (alpha_) := alpha X /; NumberQ[alpha];
  0 ** X_ := 0 /; OperatorQ[X];
  ((a_)*(b_)) ** (c_) := a*b ** c /; !OperatorQ[a];
  (c_) ** ((a_)*(b_)) := a*c ** b /; !OperatorQ[a];
  ((a_) + (b_)) ** (c_) := a ** c + b ** c;
  (a_) ** ((b_) + (c_)) := a ** b + a ** c;

Protect[NonCommutativeMultiply];

(* ------------------------------------------------------------------ *)

(* Build projectors: *)

M = 4;

Iden := Outer[Times, IdentityMatrix[M], IdentityMatrix[M]];

Projector[K_] :=
  (Projector[K] =
    Table[
      Simplify[
        Coefficient[
          Sum[ket[K, j] ** bra[K, j], {j,DimTPSubModule[K]}],
          Kron[ket[i] ** bra[j], ket[k] ** bra[l]]
        ]
      ],
      {i,M},{j,M},{k,M},{l,M}
    ]) /; MemberQ[{1,3},K]

Projector[2] := Simplify[Iden - (Projector[1] + Projector[3])];

(* ------------------------------------------------------------------ *)






(*  Build the braid generator sigma and associated objects *)

sigma := sigma =
  Module[
    {
      temp =
        q^(-2a) Projector[1] - Projector[2] + q^(2a+2) Projector[3],
      ruleqX = {-1 - q^2 + q^(-2a) + q^(2+2a) :> X},
      ruleqatop =
        {
          q^(-a) :> p^(-1),       q^(-2a) :> p^(-2),
          q^(1+a) :> q p,         q^(2+2a) :> q^2 p^2
        }
    },
    temp = temp //.  Qbracket[q_, x_] :> (q^x - q^(-x))/(q - q^(-1));
    temp = Simplify[Expand[Simplify[temp]] //. ruleqX];
    temp = Simplify[PowerExpand[temp]];
    temp[[4,3,1,2]] = -q Sqrt[Expand[(temp[[4,3,1,2]]/q)^2] //. ruleqX];
    temp[[3,4,2,1]] = q Sqrt[Expand[(temp[[3,4,2,1]]/q)^2] //. ruleqX];
    temp[[2,4,3,1]] = -Sqrt[Expand[temp[[2,4,3,1]]^2] //. ruleqX];
    temp[[4,2,1,3]] = Sqrt[Expand[temp[[4,2,1,3]]^2] //. ruleqX];
    Expand[temp //. Sqrt[X] -> Y //. X -> Y^2] //. ruleqatop
  ];

ungradedsigma :=
    Table[
      (-1)^(Grading[ket[j]]*(Grading[ket[k]]+Grading[ket[l]])) *
        sigma[[i,j,k,l]],
      {i,M},{j,M},{k,M},{l,M}
    ];

matrixofungradedsigma :=
  Table[
    ungradedsigma[[Ceiling[m/M], Ceiling[n/M],
                   m+M(1-Ceiling[m/M]), n+M(1-Ceiling[n/M])]],
    {m,M^2},{n,M^2}
  ];

matrixofungradedsigmainv :=
  Simplify[Inverse[matrixofungradedsigma]];

ungradedsigmainv := ungradedsigmainv =
  Table[
    matrixofungradedsigmainv[[M(i-1)+k,M(j-1)+l]],
    {i,M},{j,M},{k,M},{l,M}
  ];

P := P =
  Table[
    (-1)^Grading[ket[j]]*KroneckerDelta[i,l]*KroneckerDelta[k,j],
    {i,M},{j,M},{k,M},{l,M}
  ];

ls := Limit[sigma //. Y->Sqrt[-1-q^2+p^(-2)+q^2 p^2] //. p->q^a, q->1];

(* ------------------------------------------------------------------ *)


CheckRu := CheckRu = 
  Module[
    {
      temp =
      (
        - (q^u-q^(2a))/(1-q^(u+2a)) Projector[1]
        - Projector[2]
        - (1-q^(u+2a+2))/(q^u-q^(2a+2)) Projector[3]
      )
    },
    temp = temp //.  Qbracket[q_, x_] :> (q^x - q^(-x))/(q - q^(-1));
    temp = Factor[temp];

    temp[[2,1,3,4]] = Sqrt[Simplify[Expand[temp[[2,1,3,4]]^2]]];
    temp[[1,2,4,3]] = - temp[[2,1,3,4]];
    temp[[3,1,2,4]] = Sqrt[Simplify[Expand[temp[[3,1,2,4]]^2]]];
    temp[[1,3,4,2]] = - temp[[3,1,2,4]];
    temp[[2,4,3,1]] = Sqrt[Simplify[Expand[temp[[2,4,3,1]]^2]]];
    temp[[4,2,1,3]] = - temp[[2,4,3,1]];
    temp[[3,4,2,1]] = Sqrt[Simplify[Expand[temp[[3,4,2,1]]^2]]];
    temp[[4,3,1,2]] = - temp[[3,4,2,1]];

    temp = Simplify[temp];
    temp = temp //. ((q^(x_) - 1) :> q^(x/2) Q[x/2] (q-q^(-1)));
    temp = temp //. ((q^u-q^x_) :> q^(u/2+x/2) Q[u/2-x/2] (q-q^(-1)));
    temp = temp //. ((q^x_-q^u) :> q^(x/2+u/2) Q[x/2-u/2] (q-q^(-1)));
    temp = temp //. (Q[-a+x_] :> -Q[a-x]);
    Simplify[temp]
  ];

Format[Q[x_]] := SequenceForm["[", x, "]_q"];

(* ------------------------------------------------------------------ *)

\end{verbatim}
\normalsize

\pagebreak

\subsection{\texttt{LinksGouldInvariant.m}}
\addtocontents{toc}{\protect\vspace{-2.5ex}}
\seclabel{LinksGouldInvariant.m}

\footnotesize
\begin{verbatim}

(* ------------------------------------------------------------------

  LinksGouldInvariant.m

  Mathematica code to compute various link polynomial invariants for
  various (oriented) links.  We define the positive crossing tensor
  sigma (here called RMatrix), its inverse (here called SMatrix),
  several matrices derived from these, and the `cup-cap' matrices
  OMMatrix, OPMatrix, UMMatrix, UPMatrix. Requires "Uqgl21.m".

  David De Wit, Louis H Kauffman, Jon R Links

  04  August     1997  --  21  August    1997
  24  September  1997  --  04  December  1997
  06  October    1998  --  24  November  1998

  ------------------------------------------------------------------- *)

(* DimV is the range of the indices in the tensors, i.e. the
   dimension of the underlying representation.  *)

DimV[RawBracket] := 2;

DimV[LinksGould] := 4;

(* ------------------------------------------------------------------ *)

UMMatrix[LinksGould] := DiagonalMatrix[{p^2, -p^2, -q^2 p^2, q^2 p^2}];

UMMatrix[RawBracket] := {{0, I A}, {- I A^(-1), 0}};

OMMatrix[Invariant_] := OMMatrix[Invariant] =
  Inverse[UMMatrix[Invariant]];

OPMatrix[RawBracket] := OMMatrix[RawBracket];

OPMatrix[LinksGould] := IdentityMatrix[DimV[LinksGould]];

UPMatrix[Invariant_] := UPMatrix[Invariant] =
  Inverse[OPMatrix[Invariant]];

(* ------------------------------------------------------------------ *)












YangBaxterChecker[Invariant_] :=
  Module[
    {
      M = DimV[Invariant],   r = RMatrix[Invariant],   temp
    },
    temp =
      Table[
        Sum[
          r[[a,i,b,j]] r[[j,k,c,f]] r[[i,d,k,e]]
          -
          r[[b,j,c,i]] r[[a,d,j,k]] r[[k,e,i,f]],
          {i,M},{j,M},{k,M}
        ],
        {a,M},{b,M},{c,M},{d,M},{e,M},{f,M}
      ];
    Apply[And, Map[(# == 0)&, Flatten[Expand[temp]]]]
  ];

(* ------------------------------------------------------------------ *)

RMatrix[RawBracket] := RMatrix[RawBracket] =
  With[
    {
      M = DimV[RawBracket], OM = OMMatrix[RawBracket],
      UM = OMMatrix[RawBracket]
    },
    Table[
      A^(-1)      OM[[i,k]] UM[[j,l]]
      +
      A      KroneckerDelta[i,j] KroneckerDelta[k,l],
      {i,M},{j,M},{k,M},{l,M}
    ]
  ];

RMatrix[LinksGould] := RMatrix[LinksGould] =
  Expand[ungradedsigma /. Y -> Sqrt[p^(-2) - q^2 + p^2 q^2 - 1]];

SMatrix[Invariant_] := SMatrix[Invariant] =
  Module[
    {
      M = DimV[Invariant],     matrixofr, matrixofs,
      r = RMatrix[Invariant]
    },
    matrixofr =
      Table[
        r[[Ceiling[m/M], Ceiling[n/M],
           m+M(1-Ceiling[m/M]), n+M(1-Ceiling[n/M])]],
        {m,M^2},{n,M^2}
      ];
    matrixofs = Inverse[matrixofr];
    Table[
      Expand[matrixofs[[M(i-1)+k,M(j-1)+l]]],
      {i,M},{j,M},{k,M},{l,M}
    ]
  ];

(* ------------------------------------------------------------------ *)

XTensor[Invariant_, X_] :=
  (XTensor[Invariant, X] =
    Switch[
      X,
      R, RMatrix[Invariant],         S, SMatrix[Invariant],
      RD, XDTensor[Invariant, R],    SD, XDTensor[Invariant, S]
    ]
  ) /; MemberQ[{R, S, RD, SD}, X];

XLTensor[Invariant_, X_] :=
  (XLTensor[Invariant, X] =
    With[
      {
        M = DimV[Invariant],      Xint = XTensor[Invariant, X],
        OM = OMMatrix[Invariant], UM = UMMatrix[Invariant]
      },
      Table[
        Expand[Sum[Xint[[e,d,a,h]] OM[[b,e]] UM[[h,c]], {e,M},{h,M}]],
        {a,M},{b,M},{c,M},{d,M}
      ]
    ]
  ) /; MemberQ[{R, S}, X];

XRTensor[Invariant_, X_] :=
  (XRTensor[Invariant, X] =
    With[
      {
        M = DimV[Invariant],      Xint = XTensor[Invariant, X],
        OP = OPMatrix[Invariant], UP = UPMatrix[Invariant]
      },
      Table[
        Expand[Sum[Xint[[c,f,g,b]] UP[[a,f]] OP[[g,d]], {f,M},{g,M}]],
      {a,M},{b,M},{c,M},{d,M}
      ]
    ]
  ) /; MemberQ[{R, S}, X];

XDTensor[Invariant_, X_] :=
  (XDTensor[Invariant, X] =
    With[
      {
        M = DimV[Invariant],      Xint = XTensor[Invariant, X],
        OP = OPMatrix[Invariant], UP = UPMatrix[Invariant]
      },
      Table[
        Expand[Sum[
          Xint[[e,f,g,h]] UP[[a,h]] OP[[g,b]] UP[[c,f]] OP[[e,d]],
          {e,M},{f,M},{g,M},{h,M}
        ]],
        {a,M},{b,M},{c,M},{d,M}
      ]
    ]
  ) /; MemberQ[{R, S}, X];

(* ------------------------------------------------------------------ *)



PowerofXTensor[Invariant_, power_Integer?Positive, X_] :=
  (PowerofXTensor[Invariant, power, X] =
    Module[
      {
        M = DimV[Invariant],  Xpminus, Xint = XTensor[Invariant, X]
      },
      If[
        power == 1,
        Xint,
        Xpminus = PowerofXTensor[Invariant, power - 1, X];
        Table[
          Expand[Sum[Xint[[a,e,c,f]] Xpminus[[e,b,f,d]], {e,M},{f,M}]],
          {a,M},{b,M},{c,M},{d,M}
        ]
      ]
    ]
  ) /; MemberQ[{R, RD, S, SD}, X];

(* ------------------------------------------------------------------ *)

XUXDTensor[Invariant_, X_] :=
  (XUXDTensor[Invariant, X] =
    With[
      {
        M = DimV[Invariant],        Xint = XTensor[Invariant, X],
        UM = UMMatrix[Invariant],   OP = OPMatrix[Invariant],
        XD = XDTensor[Invariant, X]
      },
      Table[
        Expand[Sum[
          Xint[[a,b,e,f]] XD[[g,h,c,d]] OP[[e,g]] UM[[f,h]],
          {e,M},{f,M},{g,M},{h,M}
        ]],
        {a,M},{b,M},{c,M},{d,M}
      ]
    ]
  ) /; MemberQ[{R, S}, X];

XDXUTensor[Invariant_, X_] :=
  (XDXUTensor[Invariant, X] =
    With[
      {
        M = DimV[Invariant],        Xint = XTensor[Invariant, X],
        OM = OMMatrix[Invariant],   UP = UPMatrix[Invariant],
        XD = XDTensor[Invariant, X]
      },
      Table[
        Expand[Sum[
          XD[[a,b,e,f]] Xint[[g,h,c,d]] OM[[e,g]] UP[[f,h]],
          {e,M},{f,M},{g,M},{h,M}
        ]],
        {a,M},{b,M},{c,M},{d,M}
      ]
    ]
  ) /; MemberQ[{R, S}, X];

(* ------------------------------------------------------------------ *)

XLXRTensor[Invariant_, X_] :=
  (XLXRTensor[Invariant, X] =
    With[
      {
        M = DimV[Invariant],
        XL = XLTensor[Invariant, X], XR = XRTensor[Invariant, X]
      },
      Table[
        Expand[Sum[XL[[a,e,c,f]] XR[[e,b,f,d]], {e,M},{f,M}]],
        {a,M},{b,M},{c,M},{d,M}
      ]
    ]
  ) /; MemberQ[{R, S}, X];

XRXLTensor[Invariant_, X_] :=
  (XRXLTensor[Invariant, X] =
    With[
      {
        M = DimV[Invariant],
        XR = XRTensor[Invariant, X], XL = XLTensor[Invariant, X]
      },
      Table[
        Expand[Sum[XR[[a,e,c,f]] XL[[e,b,f,d]], {e,M},{f,M}]],
        {a,M},{b,M},{c,M},{d,M}
      ]
    ]
  ) /; MemberQ[{R, S}, X];

(* ------------------------------------------------------------------ *)

XUDUTensor[Invariant_, 1, X_] :=
  (XUDUTensor[Invariant, 1, X] =
    XTensor[Invariant, X]
  ) /; MemberQ[{R, S}, X];

XUDUTensor[Invariant_, N_, X_] :=
  (XUDUTensor[Invariant, N, X] =
    With[
      {
        M = DimV[Invariant],      XDXU = XDXUTensor[Invariant, X],
        UM = UMMatrix[Invariant], OP = OPMatrix[Invariant],
        XUDUTensorLess = XUDUTensor[Invariant, N-2, X]
      },
      Table[
        Expand[Sum[
          XUDUTensorLess[[a,b,e,f]] XDXU[[g,h]] OP[[e,g]] UM[[f,h]],
          {e,M},{f,M},{g,M},{h,M}
        ]],
        {a,M},{b,M}
      ]
    ]
  ) /; (OddQ[N] && (N >= 3) && MemberQ[{R, S}, X]);

(* ------------------------------------------------------------------ *)




AbstractTensor[Invariant_, Unknot, x_, y_] := KroneckerDelta[x,y];

AbstractTensor[Invariant_, HopfLink, x_, y_] :=
  With[
    {
      M = DimV[Invariant],
      RSquared = PowerofXTensor[Invariant, 2, R],
      UM = UMMatrix[Invariant], OP = OPMatrix[Invariant]
    },
    Sum[RSquared[[y,x,a,b]] OP[[a,c]] UM[[b,c]], {a,M},{b,M},{c,M}]
  ];

AbstractTensor[Invariant_, Trefoil, x_, y_] :=
  With[
    {
      M = DimV[Invariant],
      RCubed = PowerofXTensor[Invariant, 3, R],
      UM = UMMatrix[Invariant], OP = OPMatrix[Invariant]
    },
    Sum[RCubed[[y,x,c,d]] OP[[c,f]] UM[[d,f]], {c,M},{d,M},{f,M}]
  ];

AbstractTensor[Invariant_, FigureEight, x_, y_] :=
  With[
    {
      M = DimV[Invariant],
      RLRR = XLXRTensor[Invariant, R], SUSD = XUXDTensor[Invariant, S],
      UP = UPMatrix[Invariant],        OM = OMMatrix[Invariant]
    },
    Sum[
      RLRR[[y,a,b,c]] SUSD[[a,x,c,d]] OM[[b,e]] UP[[d,e]],
      {a,M},{b,M},{c,M},{d,M},{e,M}
    ]
  ];

AbstractTensor[Invariant_, Cinquefoil, x_, y_] :=
  With[
    {
      M = DimV[Invariant],
      RFifth = PowerofXTensor[Invariant, 5, R],
      UM = UMMatrix[Invariant], OP = OPMatrix[Invariant]
    },
    Sum[RFifth[[y,x,c,d]] OP[[c,f]] UM[[d,f]], {c,M},{d,M},{f,M}]
  ];

AbstractTensor[Invariant_, Septfoil, x_, y_] :=
  With[
    {
      M = DimV[Invariant],
      RSeventh = PowerofXTensor[Invariant, 7, R],
      UM = UMMatrix[Invariant], OP = OPMatrix[Invariant]
    },
    Sum[RSeventh[[y,x,c,d]] OP[[c,f]] UM[[d,f]], {c,M},{d,M},{f,M}]
  ];




AbstractTensor[Invariant_, FiveTwo, x_, y_] :=
  With[
    {
      M = DimV[Invariant],      ZThree = XUDUTensor[Invariant, 3, R],
      RSquared = PowerofXTensor[Invariant, 2, R],
      UP = UPMatrix[Invariant], OM = OMMatrix[Invariant]
    },
    Sum[
      ZThree[[b,c,d,x]] RSquared[[a,b,y,d]] OM[[e,a]] UP[[e,c]],
      {a,M},{b,M},{c,M},{d,M},{e,M}
    ]
  ];

AbstractTensor[Invariant_, WhiteheadLink, x_, y_] :=
  Module[
    {
      M = DimV[Invariant],      RRRL = XRXLTensor[Invariant, R], W,
      RDSquared = PowerofXTensor[Invariant, 2, RD],
      SSquared = PowerofXTensor[Invariant, 2, S],
      OP = OPMatrix[Invariant], UP = UPMatrix[Invariant],
      UM = UMMatrix[Invariant]
    },
    W =
      Table[
        Expand[Sum[
          SSquared[[c,j,e,f]] RDSquared[[g,h,i,d]] OP[[e,g]] UM[[f,h]],
          {e,M},{f,M},{g,M},{h,M}
        ]],
        {c,M},{j,M},{i,M},{d,M}
      ];
    Sum[
      RRRL[[a,i,y,b]] W[[c,x,i,d]] OP[[c,a]] UP[[d,b]],
      {a,M},{b,M},{c,M},{d,M},{i,M}
    ]
  ];

AbstractTensor[Invariant_, SixOne, x_, y_] :=
  Module[
    {
      M = DimV[Invariant],      ZThree = XUDUTensor[Invariant, 3, S],
      SD = XDTensor[Invariant, S], RRRL = XRXLTensor[Invariant, R],
      UM = UMMatrix[Invariant], OP = OPMatrix[Invariant],
      UP = UPMatrix[Invariant], OM = OMMatrix[Invariant], SOA
    },
    SOA =
      Table[
        Sum[
          SD[[b,c,f,h]] ZThree[[g,i,d,x]] OM[[f,g]] UP[[h,i]],
          {f,M},{g,M},{h,M},{i,M}
        ],
        {b,M},{c,M},{d,M}
      ];
    Sum[
      SOA[[b,c,d]] RRRL[[a,b,y,d]] OP[[e,a]] UM[[e,c]],
      {a,M},{b,M},{c,M},{d,M},{e,M}
    ]
  ];

AbstractTensor[Invariant_, SixTwo, x_, y_] :=
  Module[
    {
      M = DimV[Invariant],      SRSL = XRXLTensor[Invariant, S],
      RCubed = PowerofXTensor[Invariant, 3, R],
      RD = XDTensor[Invariant, R],
      UM = UMMatrix[Invariant], OP = OPMatrix[Invariant],
      UP = UPMatrix[Invariant], OM = OMMatrix[Invariant], STA
    },
    STA =
      Table[
        Sum[
          RCubed[[e,f,y,g]] RD[[a,b,c,d]] OM[[c,e]] UP[[d,f]],
          {c,M},{d,M},{e,M},{f,M}
        ],
        {a,M},{b,M},{g,M}
      ];
    Sum[
      STA[[a,b,g]] SRSL[[b,h,g,x]] OP[[i,a]] UM[[i,h]],
      {a,M},{b,M},{g,M},{h,M},{i,M}
    ]
  ];

AbstractTensor[Invariant_, SixThree, x_, y_] :=
  Module[
    {
      M = DimV[Invariant],
      Rint = XTensor[Invariant, R], Sint = XTensor[Invariant, S],
      RSquared = PowerofXTensor[Invariant, 2, R],
      SSquared = PowerofXTensor[Invariant, 2, S],
      UP = UPMatrix[Invariant], OM = OMMatrix[Invariant], STA, STB, ST
    },
    STA =
      Table[
        Sum[SSquared[[a,b,e,f]] Rint[[d,e,y,i]], {e,M}],
        {a,M},{b,M},{d,M},{f,M},{i,M}
      ];
    STB =
      Table[
        Sum[Sint[[b,c,g,h]] RSquared[[f,g,i,x]], {g,M}],
        {b,M},{c,M},{f,M},{h,M},{i,M}
      ];
    ST =
      Table[
        Sum[STA[[a,b,d,f,i]] STB[[b,c,f,h,i]], {b,M},{f,M},{i,M}],
        {a,M},{c,M},{d,M},{h,M}
      ];
    Sum[
      ST[[a,c,d,h]] OM[[j,d]] UP[[j,h]] OM[[k,a]] UP[[k,c]],
      {a,M},{c,M},{d,M},{h,M},{j,M},{k,M}
    ]
  ];






AbstractTensor[Invariant_, SevenTwo, x_, y_] :=
  With[
    {
      M = DimV[Invariant],      ZFive = XUDUTensor[Invariant, 5, R],
      RSquared = PowerofXTensor[Invariant, 2, R],
      UP = UPMatrix[Invariant], OM = OMMatrix[Invariant]
    },
    Sum[
      ZFive[[b,c,d,x]] RSquared[[a,b,y,d]] OM[[e,a]] UP[[e,c]],
      {a,M},{b,M},{c,M},{d,M},{e,M}
    ]
  ];

AbstractTensor[Invariant_, EightSeventeen, x_, y_] :=
  Module[
    {
      M = DimV[Invariant],          EA, EB, EC, ED, ES,
      Rint = XTensor[Invariant, R], Sint = XTensor[Invariant, S],
      RSquared = PowerofXTensor[Invariant, 2, R],
      SSquared = PowerofXTensor[Invariant, 2, S],
      UM = UMMatrix[Invariant], OP = OPMatrix[Invariant]
    },
    EA =
      Table[
        Expand[Sum[RSquared[[a,b,g,d]] SSquared[[c,g,e,f]], {g,M}]],
        {a,M},{b,M},{c,M},{d,M},{e,M},{f,M}
      ];
    EC =
      Table[
        Expand[Sum[Sint[[d,k,f,l]] Rint[[b,m,k,n]], {k,M}]],
        {b,M},{m,M},{d,M},{n,M},{f,M},{l,M}
      ];
    ED =
      Table[
         Expand[Sum[Sint[[n,o,l,j]] Rint[[m,h,o,i]], {o,M}]],
        {m,M},{h,M},{n,M},{i,M},{l,M},{j,M}
      ];
    EB =
      Table[
        Expand[Sum[
          EC[[b,m,d,n,f,l]] ED[[m,h,n,i,l,j]],
          {l,M},{m,M},{n,M}
        ]],
        {b,M},{h,M},{d,M},{i,M},{f,M},{j,M}
      ];
    ES =
      Table[
        Expand[Sum[
          EA[[a,b,c,d,e,f]] EB[[b,h,d,i,f,j]],
          {b,M},{d,M},{f,M}
        ]],
        {a,M},{h,M},{c,M},{i,M},{e,M},{j,M}
      ];
    Sum[
      ES[[y,x,c,i,e,j]] OP[[c,r]] UM[[i,r]] OP[[e,s]] UM[[j,s]],
      {c,M},{e,M},{i,M},{j,M},{s,M},{r,M}
    ]
  ];
AbstractTensor[Invariant_, NineFortyTwo, x_, y_] :=
  Module[
    {
      M = DimV[Invariant], SCubed = PowerofXTensor[Invariant, 3, S],
      RDSquared = PowerofXTensor[Invariant, 2, RD], N,
      RURD = XUXDTensor[Invariant, R], SDSU = XDXUTensor[Invariant, S],
      OM = OMMatrix[Invariant],        UP = UPMatrix[Invariant],
      OP = OPMatrix[Invariant]
    },
    N =
      Table[
        Expand[Sum[
          RDSquared[[a,b,c,d]] SCubed[[e,f,g,h]] OM[[c,e]] UP[[d,f]],
          {c,M},{d,M},{e,M},{f,M}
        ]],
        {a,M},{b,M},{g,M},{h,M}
      ];
    Sum[
      N[[a,b,y,h]] SDSU[[b,i,h,j]] RURD[[k,x,i,m]] UP[[m,j]] OP[[k,a]],
      {a,M},{b,M},{h,M},{i,M},{j,M},{k,M},{m,M}
    ]
  ];

AbstractTensor[Invariant_, TenFortyEight, x_, y_] :=
  Module[
    {
      M = DimV[Invariant], TA, TB, TT,
      Rint = XTensor[Invariant, R],
      SSquared = PowerofXTensor[Invariant, 2, S],
      SCubed = PowerofXTensor[Invariant, 3, S],
      RFourth = PowerofXTensor[Invariant, 4, R],
      UP = UPMatrix[Invariant], OP = OPMatrix[Invariant],
      OM = OMMatrix[Invariant], UM = UMMatrix[Invariant]
    },
    TA =
      Table[
        Expand[Sum[SSquared[[a,b,y,f]] RFourth[[f,d,g,h]], {f,M}]],
        {a,M},{b,M},{d,M},{g,M},{h,M}
      ];
    TB =
      Table[
        Expand[Sum[SCubed[[b,c,d,e]] Rint[[e,x,h,i]], {e,M}]],
        {b,M},{c,M},{d,M},{h,M},{i,M}
      ];
    TT =
      Table[
        Expand[Sum[
          TA[[a,b,d,g,h]] TB[[b,c,d,h,i]],
          {b,M},{d,M},{h,M}
        ]],
        {a,M},{c,M},{g,M},{i,M}
      ];
    Sum[
      TT[[a,c,g,i]] OM[[j,a]] UP[[j,c]] OP[[g,k]] UM[[i,k]],
      {a,M},{c,M},{g,M},{i,M},{j,M},{k,M}
    ]
  ];

AbstractTensor[Invariant_, Pretzel[pp_, qq_, rr_], x_, y_] :=
  Module[
    {
      M = DimV[Invariant],              P1, P2, P3,
      Zp = XUDUTensor[Invariant,pp,R],  Zq = XUDUTensor[Invariant,qq,R],
      Zr = XUDUTensor[Invariant,rr,R],
      OM = OMMatrix[Invariant],         UP = UPMatrix[Invariant]
    },
    P1 =
      Expand[Table[
        Sum[Zp[[a,b,y,e]] * Zq[[b,c,e,f]], {b,M},{e,M}],
        {a,M},{c,M},{f,M}
      ]];
    P2 =
      Expand[Table[
        Sum[P1[[a,c,f]] * Zr[[c,d,f,x]], {c,M},{f,M}],
        {a,M},{d,M}
      ]];
    P3 =
      Expand[Table[
        Sum[P2[[a,d]] * OM[[g,a]], {a,M}],
        {d,M},{g,M}
      ]];
    Expand[Sum[P3[[d,g]] * UP[[g,d]], {d,M},{g,M}]]
  ] /; ((OddQ[pp] && OddQ[qq] && OddQ[rr]) && (3 <= pp < qq < rr));

KTATensor[Invariant_] :=
  Module[
    {
      M = DimV[Invariant],             SD = XDTensor[Invariant, S],
      SSquared = PowerofXTensor[Invariant, 2, S],
      OP = OPMatrix[Invariant],        UM = UMMatrix[Invariant],
      RURD = XUXDTensor[Invariant, R], temp
    },
    temp =
      Table[
        Expand[Sum[
          SSquared[[d,s,f,g]] SD[[h,i,e,c]] OP[[f,h]] UM[[g,i]],
          {f,M},{g,M},{h,M},{i,M}
        ]],
        {d,M},{s,M},{e,M},{c,M}
      ];
    Table[
      Expand[Sum[RURD[[a,d,b,e]] temp[[d,s,e,c]], {d,M},{e,M}]],
     {a,M},{s,M},{b,M},{c,M}
    ]
  ];











KTAITensor[Invariant_] :=
  Module[
    {
      M = DimV[Invariant],             SD = XDTensor[Invariant, S],
      SSquared = PowerofXTensor[Invariant, 2, S],
      OP = OPMatrix[Invariant],        UM = UMMatrix[Invariant],
      RURD = XUXDTensor[Invariant, R], temp
    },
    temp =
      Table[
        Expand[Sum[
          SSquared[[a,d,f,g]] SD[[h,i,b,e]] OP[[f,h]] UM[[g,i]],
          {f,M},{g,M},{h,M},{i,M}
        ]],
      {a,M},{d,M},{b,M},{e,M}
    ];
    Table[
      Expand[Sum[temp[[a,d,b,e]] RURD[[d,s,e,c]], {d,M},{e,M}]],
      {a,M},{s,M},{b,M},{c,M}
    ]
  ];

KTCTensor[Invariant_] :=
  Module[
    {
      M = DimV[Invariant],               KTB,
      OP = OPMatrix[Invariant],          OM = OMMatrix[Invariant],
      UP = UPMatrix[Invariant],          UM = UMMatrix[Invariant],
      SLSR = XLXRTensor[Invariant, S],   Rint = XTensor[Invariant, R],
      RDSquared = PowerofXTensor[Invariant, 2, RD],
      SUSD = XUXDTensor[Invariant, S]
    },
    KTB =
      Table[
        Expand[Sum[
          Rint[[d,a,b,c]] RDSquared[[l,m,f,n]] SUSD[[a,e,n,g]] *
          OP[[b,l]] UM[[c,m]],
          {a,M},{b,M},{c,M},{l,M},{m,M},{n,M}
        ]],
        {d,M},{e,M},{f,M},{g,M}
      ];
    Table[
      Expand[Sum[
        KTB[[d,e,f,g]] SLSR[[h,i,j,k]] OM[[f,h]] UP[[g,i]],
        {f,M},{g,M},{h,M},{i,M}
      ]],
      {d,M},{e,M},{j,M},{k,M}
    ]
  ];









AbstractTensor[Invariant_, KT, x_, y_] :=
  With[
    {
      M = DimV[Invariant],        KTC = KTCTensor[Invariant],
      KTA = KTATensor[Invariant], UP = UPMatrix[Invariant],
      OP = OPMatrix[Invariant],   OM = OMMatrix[Invariant]
    },
    Sum[
      KTA[[a,x,b,c]] KTC[[d,e,j,k]] *
      OM[[b,d]] UP[[c,e]] OP[[a,j]] UP[[k,y]],
      {a,M},{b,M},{c,M},{d,M},{e,M},{j,M},{k,M}
    ]
  ];

AbstractTensor[Invariant_, KTI, x_, y_] :=
  With[
    {
      M = DimV[Invariant],          KTC = KTCTensor[Invariant],
      KTAI = KTAITensor[Invariant], UP = UPMatrix[Invariant],
      OP = OPMatrix[Invariant],     OM = OMMatrix[Invariant]
    },
    Sum[
      KTAI[[a,x,b,c]] KTC[[d,e,j,k]] *
      OM[[b,d]] UP[[c,e]] OP[[a,j]] UP[[k,y]],
      {a,M},{b,M},{c,M},{d,M},{e,M},{j,M},{k,M}
    ]
  ];

(* ------------------------------------------------------------------ *)





























KnotList =
  {
    Unknot, HopfLink, Trefoil, FigureEight,
    Cinquefoil, FiveTwo, WhiteheadLink,
    SixOne, SixTwo, SixThree, Septfoil, SevenTwo,
    EightSeventeen, NineFortyTwo, TenFortyEight, KT, KTI
  };

LinkPolynomial::usage = "LinkPolynomial[Invariant, Link] evaluates the
link polynomial invariant of name Invariant for the knot described by
Link.  Invariant must be one of `RawBracket', or `LinksGould', and Link
must be one of the elements of `KnotList' or `Pretzel[pp, qq, rr]'."

Attributes[LinkPolynomial] = {Listable};

LinkPolynomial[RawBracket, Link_] :=
  LinkPolynomial[RawBracket, Link] =
    With[
      {
        M = DimV[RawBracket],
        OM = OMMatrix[RawBracket], UP = UPMatrix[RawBracket]
      },
      Print["Building LinkPolynomial[RawBracket, ", Link, "]."];
      Expand[Sum[
        AbstractTensor[RawBracket, Link, x, y] OM[[z,y]] UP[[z,x]],
        {x,M},{y,M},{z,M}
      ]]
    ];

LinkPolynomial[LinksGould, Link_] :=
  LinkPolynomial[LinksGould, Link] =
  (
    Print["Building LinkPolynomial[LinksGould, ", Link, "]."];
    Expand[AbstractTensor[LinksGould, Link, 1, 1]]
  );

(* ------------------------------------------------------------------ *)

Attributes[LGCantDetectNoninvertibilityQ] = {Listable};

LGCantDetectNoninvertibilityQ[Link_] :=
  With[
    {poly = LinkPolynomial[LinksGould, Link] /. p -> q^a},
    Expand[poly - (poly /. a -> -(1+a))] === 0
  ];

Attributes[LGDetectsChiralityQ] = {Listable};

LGDetectsChiralityQ[Link_] :=
  With[
    {poly = LinkPolynomial[LinksGould, Link] /. p -> q^a},
    PowerExpand[poly - (poly /. q-> q^(-1))] =!= 0
  ];

(* ------------------------------------------------------------------ *)



ListofPretzels[N_] :=
  Flatten[
    Table[Pretzel[pp,qq,rr], {pp,3,N,2}, {qq,pp+2,N,2}, {rr,qq+2,N,2}]
  ];

LGCantDetectNoninvertibilityofPretzelsQ[N_] :=
  LGCantDetectNoninvertibilityQ[ListofPretzels[N]];

LGDetectsChiralityofPretzelsQ[N_] :=
  LGDetectsChiralityQ[ListofPretzels[N]];

(* ------------------------------------------------------------------ *)

Writhe[Pretzel[pp_, qq_, rr_]] := pp + qq + rr;

Writhe[Unknot] =          0;
Writhe[HopfLink] =        2;    Writhe[Trefoil] =         3;
Writhe[FigureEight] =     0;    Writhe[Cinquefoil] =      5;
Writhe[FiveTwo] =         5;    Writhe[WhiteheadLink] =   2;
Writhe[SixOne] =         -2;    Writhe[SixTwo] =          2;
Writhe[SixThree] =        0;    Writhe[Septfoil] =        7;
Writhe[SevenTwo] =        7;    Writhe[EightSeventeen] =  0;
Writhe[NineFortyTwo] =   -1;    Writhe[TenFortyEight] =   0;
Writhe[KT] =             -2;    Writhe[KTI] =            -2;

Attributes[JonesPolynomial] = {Listable};

JonesPolynomial[Link_] :=
  With[
    {
      bu = LinkPolynomial[RawBracket, Unknot],
      b = LinkPolynomial[RawBracket, Link]
    },
    Expand[(-A)^(-3 Writhe[Link]) Expand[Factor[b/bu]]] /. A -> t^(-1/4)
  ];

TableofJonesCoeffs[ListofLinks_] :=
  Module[
    {
      J = JonesPolynomial[ListofLinks], powers
    },
    powers = -Exponent[J, t^(-1)];
    Transpose[
     {
       ListofLinks,
       powers,
       CoefficientList[Expand[t^(-powers)] J, t]
     }
    ]
  ];

(* ------------------------------------------------------------------ *)

\end{verbatim}
\normalsize

\cleardoublepage


\def\JMP{Journal of Mathematical Physics}
\def\JPA{Journal of Physics A. Mathematical and General}
\def\LMP{Letters in Mathematical Physics}
\def\NPB{Nuclear Physics B}
\def\RMP{Reports on Mathematical Physics}
\def\CMP{Communications in Mathematical Physics}
\def\JKTR{Journal of Knot Theory and its Ramifications}



\addcontentsline{toc}{section}{References}

\bibliographystyle{plain}
\bibliography{thesis}

\end{document}